\documentclass{gtart_h}  

\def\ifplaintex{\expandafter\ifx\csname documentclass\endcsname\relax}

\def\ifplaintex{\expandafter\ifx\csname documentclass\endcsname\relax}

\ifplaintex 
\else
\fi

\expandafter\ifx\csname epsfbox\endcsname\relax\input epsf\fi

\def\gt{{\mathsurround=0pt\it $\cal G\mskip-2mu$eometry \&\ 
$\cal T\!\!$opology}}        %  journal title in recommended style

\def\gtp{{\mathsurround=0pt\it $\cal G\mskip-2mu$eometry \&\ 
$\cal T\!\!$opology $\cal P\!$ublications}}  % GT publications

\def\lognumber#1{\def\thelognumber{#1}}
\def\volumenumber#1{\def\thevolumenumber{#1}}
\def\papernumber#1{\def\thepapernumber{#1}}
\def\volumeyear#1{\def\thevolumeyear{#1}}

\def\pagenumbers#1#2{\def\startpage{#1}\def\finishpage{#2}}
\def\published#1{\def\publishdate{#1}}
\def\proposed#1{\def\theproposer{#1}}
\def\seconded#1{\def\theseconders{#1}}
\def\received#1{\def\receiveddate{#1}}
\def\revised#1{\def\reviseddate{#1}}
\def\accepted#1{\def\accepteddate{#1}}

\def\asciiurl#1{\def\theasciiurl{#1}}

\let\\\par\let\thelognumber\relax
\let\thevolumenumber\relax\let\thepapernumber\relax
\let\thevolumeyear\relax\let\thesamplenumber\relax\let\startpage\relax
\let\finishpage\relax\let\publishdate\relax\let\receiveddate\relax
\let\reviseddate\relax\let\accepteddate\relax\let\theasciititle\relax
\let\theasciiauthors\relax
\let\theasciiabstract\relax
\let\theasciiemail\relax\let\theshortauthors\relax\let\theshorttitle\relax
\let\theasciiurl\relax

\long\def\maketitlep{   % start of definition of \maketitlep

\count0=\startpage

\gt\hfill      %   Journal title (top left) 
\hbox to 77pt{\vbox to 0pt{\vglue -15pt\epsfbox{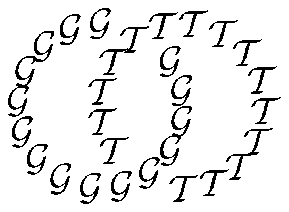}\vss}\hss}
\break
{\small\ifx\thesamplenumber\relax % sample?  
Volume \else Sample
\fi\thevolumenumber\ (\thevolumeyear)
\startpage--\finishpage\nl
Published: \publishdate}
\vglue 0.5truein plus 0.4fil minus 0.1truein

{\parskip=0pt\leftskip 0pt plus 1fil\def\\{\par\smallskip}{\ifplaintex\large
\else\Large\fi\bf\thetitle}\par\medskip}   

\vglue 0pt plus 0.1fil 

{\parskip=0pt\leftskip 0pt plus 1fil\def\\{\par}{\sc\theauthors}
\par\medskip}

\vglue 0pt plus 0.1fil 

{\small\parskip=0pt\let\newline\\
{\leftskip 0pt plus 1fil\def\\{\par}{\sl\theaddress}\par}
\expandafter\ifx\theemail\relax    % email address?
\relax\else\vglue 5pt plus 0.02fil minus 2pt\def\\{\stdspace{\rm 
and}\stdspace} 
\cl{Email:\stdspace\tt\theemail}\fi
\ifx\theurl\relax                  % URL given?
\relax\else\vglue 5pt plus 0.02fil minus 2pt\def\\{\stdspace{\rm 
and}\stdspace}
\cl{URL:\stdspace\tt\theurl}\fi\par}

\vglue 7pt plus 0.3fil minus 3pt

{\bf Abstract}
\vglue 5pt plus 0.1fil minus 2pt

\theabstract

\vglue 7pt plus 0.3fil minus 3pt

{\bf AMS Classification numbers}\quad Primary:\quad \theprimaryclass

Secondary:\quad \thesecondaryclass

\vglue 5pt plus 0.3fil minus 2pt

{\bf Keywords:}\quad \thekeywords

\vglue 10pt plus 0.5fil minus 5pt

{\small  Proposed: \theproposer\hfill Received: \receiveddate\nl
Seconded: \theseconders\hfill 
\ifx\reviseddate\relax                         % paper revised?
Accepted: \accepteddate                        % no
\else
Revised: \reviseddate                          % yes
\fi}
\eject
}       %  end of definition of \maketitlep

\font\phead=cmsl9 scaled 950
\font=cmsl9 scaled 1050
\font\pnum=cmbx10 scaled 913
\font\lnum=cmbx10 
\font\pfoot=cmsl9 scaled 950
\font=cmsl9 scaled 1050
\ifplaintex
\headline{\vbox to 0pt{\vskip -4.5mm\line{\small\phead\ifnum
\count0=\startpage ISSN 1364-0380 (on line)
1465-3060 (printed) \hfill {\pnum\folio}\else\ifodd\count0\def\\{ }% 
\ifx\theshorttitle\relax\thetitle\else\theshorttitle\fi\hfill{\pnum\folio}
\else\def\\{ and }{\pnum\folio}\hfill\ifx\theshortauthors\relax\theauthors
\else\theshortauthors\fi\fi\fi}\vss}}
\footline{\vbox to 0pt{\vglue 0mm\line{\small\pfoot\ifnum\count0=\startpage
\copyright\ \gtp\hfill\else
\gt, Volume \thevolumenumber\ (\thevolumeyear)\hfill\fi}\vss
}}
\else
\makeatletter
\def\@oddhead{{\small\lhead\ifnum\count0=\startpage ISSN 1364-0380 (on line)
1465-3060 (printed) \hfill {\lnum\number\count0}\else\ifodd\count0
\def\\{ }\ifx\theshorttitle\relax \thetitle \else\theshorttitle\fi\hfill
{\lnum\number\count0}\else\def\\{ and }{\lnum\number\count0}
\hfill\ifx\theshortauthors\relax 
\theauthors\else\theshortauthors\fi\fi\fi}}\def\@evenhead{\@oddhead}
\def\@oddfoot{\small\lfoot\ifnum\count0=\startpage\copyright\ \gtp\hfill\else
\gt, Volume \thevolumenumber\ (\thevolumeyear)\hfill\fi}
\def\@evenfoot{\@oddfoot}
\makeatother
\fi

\newwrite\gtoutfile
\long\gdef\makeheadfile{  %%% start of definition of \makeheadfile
{\def\\{, }\def\s{ }
\immediate\openout\gtoutfile head.xxx
\immediate\write\gtoutfile{Proxy-for: \ifx\theasciiauthors\relax
\theauthors\else\theasciiauthors\fi\s<\ifx\theasciiemail\relax\theemail\else\theasciiemail\fi>}
\immediate\write\gtoutfile{\noexpand\\}
\immediate\write\gtoutfile{Authors: \ifx\theasciiauthors\relax
\theauthors\else\theasciiauthors\fi}
{\def\\{ }\immediate\write\gtoutfile{Title: \ifx\theasciititle\relax
\thetitle\else\theasciititle\fi}}
\immediate\write\gtoutfile{Subj-class: GT or SG or MG etc}
\immediate\write\gtoutfile{MSC-class: \theprimaryclass\ifx\thesecondaryclass\relax\else, \thesecondaryclass\fi}
\immediate\write\gtoutfile{Journal-ref: Geom. Topol. \thevolumenumber
(\thevolumeyear) \startpage-\finishpage}
\immediate\write\gtoutfile{Comments: Published by Geometry and Topology at}
\immediate\write\gtoutfile{\s\s http://www.maths.warwick.ac.uk/gt/GTVol\thevolumenumber/paper\thepapernumber.abs.html}
\immediate\write\gtoutfile{\noexpand\\}
\immediate\write\gtoutfile{}
\ifx\theasciiabstract\relax
\immediate\write\gtoutfile{\theabstract}\else
\immediate\write\gtoutfile{\theasciiabstract}\fi
\immediate\write\gtoutfile{}
\immediate\write\gtoutfile{\noexpand\\}
\immediate\write\gtoutfile{}
\immediate\closeout\gtoutfile}}  %%% end of definition of \makeheadfile

\def\maketitlepage{\maketitlep\makeheadfile}
\let\maketitle\maketitlepage

\lognumber{465}
\received{2 August 2003}
\volumenumber{9}\papernumber{16}\volumeyear{2005}
\pagenumbers{571}{697}   
\revised{26 February 2005}
\published{19 April 2005}
\accepted{29 March 2005}
\proposed{Frances Kirwan}
\seconded{Ralph Cohen, Gang Tian}

\usepackage{amsmath,amssymb,pstricks,pst-coil,pst-node}

\numberwithin{equation}{section}
\newtheorem{thm}{Theorem}[section] 
\newtheorem{prp}[thm]{Proposition}
\newtheorem{lmm}[thm]{Lemma}   
\newtheorem{crl}[thm]{Corollary} 
\newtheorem{dfn}[thm]{Definition} 

\newtheorem*{remark}{Remark}
\newtheorem*{remarks}{Remarks}

\newtheorem{alprp}{Proposition}

\newtheorem{ques}{Question}

\renewcommand{\frak}{\mathfrak}
\renewcommand{\Bbb}{\mathbb}

\def\I{\mathfrak i}
\def\M{\mathfrak M}

\def\under#1{\underline{#1}}
\def\ov#1{\overline{#1}}
\def\ti#1{\widetilde{#1}}
\def\tilde#1{\widetilde{#1}}
\def\hat#1{\widehat{#1}}
\def\sm#1{\begin{small}#1\end{small}}

\def\lra{\longrightarrow}
\def\Lra{\Longrightarrow}
\def\Llra{\Longleftrightarrow}

\def\lan{\langle}
\def\ran{\rangle}
\def\lr#1{\lan{#1}\ran}
\def\blr#1{\big\lan{#1}\big\ran}
\def\llan{\langle\langle}
\def\rran{\rangle\rangle}
\def\llan{\lan\!\lan}
\def\rran{\ran\!\ran}
\def\llrr#1{\llan{#1}\rran}
\def\Llan{\big\lan\!\big\lan}
\def\Rran{\big\ran\!\big\ran}
\def\LlRr#1{\Llan{#1}\Rran}

\def\al{\alpha}
\def\be{\beta}
\def\de{\delta}
\def\ep{\epsilon}
\def\ga{\gamma}
\def\io{\iota}
\def\ka{\kappa}
\def\la{\lambda}

\def\om{\omega}
\def\si{\sigma}
\def\th{\theta}
\def\ve{\varepsilon}
\def\ups{\upsilon}
\def\vph{\varphi}

\def\vr{\varrho}
\def\vt{\vartheta}

\def\De{\Delta}
\def\Ga{\Gamma}

\def\Om{\Omega}
\def\Si{\Sigma}

\def\and{\textnormal{and}}
\renewcommand{\rk}{\operatorname{rk}}

\newcommand{\codim}{\operatorname{codim}}
\renewcommand{\Im}{\operatorname{Im}}

\newcommand{\ev}{\operatorname{ev}}
\newcommand{\const}{\operatorname{const}}
\newcommand{\ord}{\operatorname{ord}}

\def\P{\Bbb{P}^n}
\def\PP{\Bbb{P}^2}
\def\PPP{\Bbb{P}^3}
\def\i{\infty}
\def\eset{\emptyset}
\def\w{\wedge}

\def\overK{\overline{{\mathcal{K}}}}
\def\overM{{}\mskip4mu\overline{\mskip-4mu{\mathcal{M}}}}
\def\overR{{}\overline{\mathbb{R}\mskip-2mu}\mskip2mu}
\def\overS{{}\mskip2mu\overline{\mskip-2mu{\mathcal{S}}}}
\def\overT{\overline{{\mathcal{T}}}}
\def\overU{\overline{{\mathcal{U}}}}
\def\overV{\overline{{\mathcal{V}}}}
\def\overZ{\overline{{\mathcal{Z}}}}

\begin{document}

\title[Counting rational curves of arbitrary shape]
{Counting rational curves of arbitrary shape\\in projective spaces}
\author{Aleksey Zinger}

\address{Department of Mathematics, Stanford University\\Stanford, CA 94305-2125, USA}
\email{azinger@math.stanford.edu}
\urladdr{http://math.stanford.edu/~azinger/}
\asciiurl{http://math.stanford.edu/ azinger/}

\begin{abstract}
We present an approach to a large class of enumerative problems
concerning rational curves in projective spaces.  This approach uses
analysis to obtain topological information about moduli spaces of
stable maps.  We demonstrate it by enumerating one-component rational
curves with a triple point or a tacnodal point in the
three-dimensional projective space and with a cusp in any projective
space.
\end{abstract}


\maketitlepage

\section{Introduction}
\label{intro_sec}

\subsection{Background}
\label{back_subs}

Enumerative geometry of algebraic varieties is a field of mathematics that
has been a subject of ongoing research since at least the nineteenth
century.  Problems in this field are geometric in nature.  It has
connections to other fields of mathematics as well as to theoretical
physics.  The general goal of enumerative algebraic geometry is to
determine the number of geometric objects that satisfy pre-specified
geometric conditions.  The objects are often, but not always, complex
curves in a smooth algebraic manifold.  Such curves would be required
to represent a given homology class, to have certain singularities,
and to satisfy various contact conditions with respect to a collection
of subvarieties.  One of the most well-known examples of an enumerative
problem is:

\begin{ques}
\label{g0n2_ques}
If $d$ is a positive integer, what is the number $n_d$ of degree--$d$ rational 
curves that pass through $3d{-}1$ points in general position
in the complex projective plane $\PP$?
\end{ques}

Since the number of lines through any two distinct points is one,
$n_1{=}1$.  A little bit of algebraic geometry and topology
gives $n_2{=}1$ and $n_3{=}12$.  It is far harder to find that
$n_4 = 620$, but this number was computed as early as the middle of
the nineteenth century; \hbox{see~\cite[page 378]{Ze}}.  The higher-degree
numbers remained unknown until the early~1990s, when a recursive formula
for the numbers~$n_d$ was announced; see~\cite{KM} and~\cite{RT}.

For more than a hundred years, tools of algebraic geometry had been
the dominant force behind progress in enumerative algebraic geometry.
However, in~\cite{G}, Gromov initiated the study of pseudoholomorphic
curves in symplectic manifolds and demonstrated their usefulness
by obtaining a number of important results in symplectic topology.
Since then moduli spaces of stable maps, ie of the parameterizations of
pseudoholomorphic curves, have evolved into a powerful tool in enumerative
geometry and have become a central object in algebraic geometry.
In particular, these moduli spaces lie behind the derivation of the
recursive formula for the numbers $n_d$ in \cite{KM} and \cite{RT}.
The latter work, in fact, gives a recursive-formula solution to the
natural generalization of Question~\ref{g0n2_ques} to projective spaces
of arbitrary dimension:

\begin{ques}
\label{g0_ques}
Suppose $n\ge2$, $d$, and $N$ are positive integers, 
and 
$$\mu\equiv(\mu_1,\ldots,\mu_N)$$
is an $N$--tuple of proper subvarieties  of~$\P$ in general position such that 
\begin{equation}\label{g0_ques_e}
\sum_{l=1}^{N}\codim_{\Bbb{C}}\mu_l =d(n + 1)+n-3+N.
\end{equation}
What is the number $n_d(\mu)$ of degree--$d$ rational curves 
that pass through the subvarieties $\mu_1,\ldots,\mu_N$?
\end{ques}

Condition~\eqref{g0_ques_e} is necessary to ensure that 
the expected answer is finite and not clearly zero.
For straightforward geometric reasons 
(the same ones as for the linearity property and the divisor equation
in \cite[Section~1]{RT}), it is sufficient to solve 
Question~\ref{g0_ques}, as well as other similar questions, 
for tuples $\mu$ of linear subspaces of~$\P$ of codimension at least~two.
Thus, the enumerative formulas given below are stated and proved only for 
constraints~$\mu$ that are points in $\PP$ or points and lines in~$\PPP$.
However, analogous formulas hold for arbitrary constraints~$\mu$.

Following \cite{G} and \cite{K}, moduli spaces of stable maps 
into algebraic manifolds became subjects of much research in algebraic geometry.
Algebraic geometers usually denote by $\ov\M_{0,N}(\P,d)$
the stable-map compactification of the space $\M_{0,N}(\P,d)$ of  equivalence classes of 
degree--$d$ holomorphic maps from $\Bbb{P}^1$ with $N$~marked points into~$\P$.
These spaces are described as algebraic stacks in~\cite{FP}.
While their cohomology is not entirely understood, it is shown in~\cite{P} that 
the intersections of \textit{tautological} cohomology classes in 
$\ov\M_{0,N}(\P,d)$ can be computed via explicit recursive formulas.
These cohomology classes include all cohomology classes that arise 
through natural geometric constructions.
As an application to enumerative geometry, \cite{P} expresses
the number $|{\mathcal{S}}_1(\mu)|$  of Question~\ref{g0n2cusp_ques}
in terms of intersections of tautological classes in
$\ov\M_{0,N}(\P,d)$ and then in terms of the numbers~$n_d$.

\begin{ques}
\label{g0n2cusp_ques}
If $d$ is a positive integer, what is the number $\big|{\mathcal{S}}_1(\mu)\big|$ 
of degree--$d$ rational curves that have a cusp and pass through
a tuple~$\mu$ of $3d{-}2$ points in general position in~$\PP$?
\end{ques}

A similar approach to enumerative geometry of plane curves  is taken
in~\cite{V}.  Using relationships derived in~\cite{DH}, \cite{V} expresses
the ``codimension-one'' enumerative numbers of rational plane curves,
such as those of Questions~\ref{g0n2cusp_ques}--\ref{g0n2tac_ques},
in terms of intersection numbers of tautological classes in $\ov{\frak
M}_{0,N}(\P,d)$ and the latter in terms of the numbers~$n_d$.

\begin{ques}
\label{g0n2p3_ques}
If $d$ is a positive integer, what is the number 
$\frac{1}{6}\big|{\mathcal{V}}_1^{(2)}(\mu)\big|$ 
of degree--$d$ rational curves that have a triple point and pass through
a tuple~$\mu$ of $3d{-}2$ points in general position in~$\PP$?
\end{ques}

\begin{ques}
\label{g0n2tac_ques}
If $d$ is a positive integer, what is the number 
$\frac{1}{2}\big|{\mathcal{S}}_1^{(1)}(\mu)\big|$ 
of degree--$d$ rational curves that have a tacnode and pass through
a tuple~$\mu$ of $3d{-}2$ points in general position in~$\PP$?
\end{ques}

Questions~\ref{g0n2_ques} and \ref{g0n2cusp_ques}--\ref{g0n2tac_ques}
can actually be solved using more classical methods of algebraic geometry,
as is done in~\cite{R1} and~\cite{R2}.
However, the derivations in~\cite{R1} and~\cite{R2} involve 
fairly complicated algebraic geometry.
In contrast, the computations in~\cite{P} and~\cite{V} involve 
much less algebraic geometry and rely on known results,
obtained via fairly complicated algebraic geometry elsewhere,
including~\cite{DH}, \cite{FP}, and~\cite{K}.

The method given in this paper can be used in a straightforward, if
somewhat laborious, manner to express the number of rational curves
in a complex projective space, that have a $k$--fold point, for example,
and pass through a set of constraints in general position, in terms of
intersections of tautological classes in the moduli spaces of stable
rational maps.  In Subsection~\ref{general_subs}, we describe in more
detail the scope of the applicability of this method.  Its application
makes practically no use of algebraic geometry.  The method itself relies
on a number of technical results, only some of which are contained in
this paper, and the rest elsewhere, including~\cite{LT,MS,RT,Z1,Z2}.

The author would like to thank Tomasz Mrowka and Jason Starr for
helpful conversations during the preparation of this manuscript and
Izzet Coskun, Joachim Kock, Ravi Vakil, and the referee for comments
on early versions of this paper.  The author was partially supported
by the Clay Mathematics Institute and an NSF Postdoctoral Fellowship.
Most of this work was completed at~MIT.

\subsection{Outline of the method}
\label{method_subs}

The first step in our approach is to describe a subset ${\mathcal{Z}}$ of
a moduli space of stable rational maps, or of a closely related space,
such that the cardinality of ${\mathcal{Z}}$ is a known multiple of the number
we are looking~for.  We would also like the subset ${\mathcal{Z}}$ to be the
zero set of a reasonably well-behaved section~$s$ of a bundle $V$ over
a reasonably nice submanifold~${\mathcal{S}}$ of the ambient space~$\overM$.  For example, in the case of Question~\ref{g0n2cusp_ques}, we might
take ${\mathcal{S}}$ to be the subset of $\M_{0,1}(\PP,d)$ consisting of the
equivalence classes of maps whose images pass through the $3d{-}2$ points
in~$\PP$ and take ${\mathcal{Z}}$ to be the subset of ${\mathcal{S}}$ consisting
of the equivalence classes of maps whose differential vanishes at the
marked point.  Alternatively, we can also allow ${\mathcal{Z}}$ to be the
preimage under a reasonably well-behaved map $h\co{\mathcal{S}} \lra \mathcal{
X}$ of a submanifold $\De$ of~${\mathcal{X}}$.  For example, in the case of
Question~\ref{g0n2p3_ques}, we might take ${\mathcal{S}}$ to be the subset
of $\M_{0,3}(\PP,d)$ consisting of the equivalence classes of maps~$b$
whose images pass through the $3d{-}2$ points such that $\ev_1(b) =
\ev_2(b)$, where $\ev_1$ and $\ev_2$ are the evaluation maps at the
first and second marked points of~$\M_{0,3}(\PP,d)$. We could then take
$${\mathcal{Z}}=\{\ev_1 \times \ev_3\}^{-1}(\De_{\PP\times\PP})\cap{\mathcal{S}},$$
where $\De_{\PP\times\PP}$ denotes the diagonal in $\PP \times \PP$.
In the case of Question~\ref{g0n2tac_ques}, we might take the
ambient space to be the projectivization of a natural rank-two bundle
over~$\ov\M_{0,2}(\PP,d)$.  However, in practice, we will keep track of
the points on $\Bbb{P}^1$ that get mapped to the constraints, ie there
will be marked points labeled by the positive integers $1,\ldots,N$,
where $N$ is the number of constraints.  The marked points of the domain
of a stable map that describe the singularities of the image curve will
be labeled by $\hat{1}$, $\hat{2}$, etc.

If ${\overS}$ is a smooth compact oriented manifold,
$V \lra {\overS}$  is a smooth oriented vector bundle
of the same rank as the dimension of ${\mathcal{S}}$,
and $\ti{s}\co{\overS} \lra V$ is a smooth section, 
which is transverse to the zero set in $V$, then
\begin{equation}\label{method_subs_e1}
^{\pm}\big|\ti{s}^{-1}(0)\big|=
\big\lan e(V),{\overS}\big\ran,
\end{equation}
where $^{\pm}\big|\ti{s}^{-1}(0)\big|$ is the signed cardinality
of the set $\ti{s}^{-1}(0)$.
Equation~\eqref{method_subs_e1} is valid under more general circumstances.
In the cases of interest to us, 
the ambient ${\overM}$ is an oriented stratified topological orbifold
and ${\mathcal{S}}$ is a smooth submanifold of the main stratum~${\mathcal{M}}$
such that ${\overS} - {\mathcal{S}}$ is contained 
in a finite union of smooth manifolds
of dimension less than the dimension of~${\mathcal{S}}$.
Under these assumption, ${\mathcal{S}}$ determines a homology class
in~${\overM}$.
Furthermore, if $\ti{s}$ is a continuous section of $V$ over 
${\overS}$ and $e(V)$ is the restriction of 
a cohomology class on ${\overM}$, then 
equality~\eqref{method_subs_e1} still holds.
By~\eqref{method_subs_e1},
if $s$ is any continuous section of $V$ 
over ${\overS}$ such that $s|{\mathcal{S}}$ is transverse to the zero set
and \hbox{${\mathcal{Z}} \equiv s^{-1}(0)\cap{\mathcal{S}}$} is a finite set, then
\begin{equation}\label{method_subs_e2}
^{\pm}\big|{\mathcal{Z}}\big|=\big\lan e(V),{\overS}\big\ran-
{\mathcal{C}}_{\partial{\overS}}(s),
\end{equation}
where ${\mathcal{C}}_{\partial{\overS}}(s)$ is
the \textit{$s$--contribution of $\partial{\overS}$ to the euler class of~$V$}.
In other words,  ${\mathcal{C}}_{\partial{\overS}}(s)$ is the signed number
of zeros of a small generic perturbation~$\tilde{s}$ of~$s$
that lie near~$\partial{\overS}$.
If the behavior of $s$ near $\partial{\overS}$ can be understood, 
it is reasonable to hope that the number ${\mathcal{C}}_{\partial{\overS}}(s)$
can be computed, at least in terms of evaluations of some cohomology classes.
On the other hand, in the case of Question~\ref{g0n2cusp_ques},
${\overS}$ is a tautological class in the appropriate moduli
space of stable maps.
Thus, if $e(V)$ is also a tautological class, 
$\lr{e(V),{\mathcal{S}}}$ is computable, and we are done.
Most of the time, however, 
we will have to describe $\lr{e(V),{\overS}}$ as 
the signed cardinality of a subset ${\mathcal{Z}}'$ 
of a space which is a step closer to a tautological class than ${\mathcal{S}}$
and apply equation~\eqref{method_subs_e2} with~${\mathcal{Z}}'$.
Eventually, we will end up with intersections of tautological classes
in moduli spaces of stable rational maps.

The topological setup of the previous paragraph is only slightly more general 
than that of \cite[Section~3]{Z1}.
However, it is not sufficient for our purposes.
We now present two significant generalizations of this setup.
The first is that equation~\eqref{method_subs_e2} makes sense 
even if the section $s$ is defined only over ${\mathcal{S}}$
and does not extend over ${\overS} - {\mathcal{S}}$.
In such a case, we can use a cutoff function to define a
new section $s'$ that vanishes on a neighborhood of 
${\overS} - {\mathcal{S}}$ and thus extends to a continuous
section over~${\overS}$.
The term ${\mathcal{C}}_{\partial{\overS}}(s)$ is then the signed number of zeros 
of a small generic perturbation~$\ti{s}$ of~$s'$
that lie near~$\partial{\overS}$.
If we can understand the behavior of $s$ near $\partial{\overS}$ 
and choose the cutoff function carefully, it is again reasonable to hope that 
we can determine the number ${\mathcal{C}}_{\partial{\overS}}(s)$.

The second generalization has a very different flavor.
Suppose ${\mathcal{S}}$ and ${\overM}$ are as above and
${\mathcal{X}}$ is a smooth compact oriented manifold.
If $h\co{\overM} \lra {\mathcal{X}}$ is continuous map such that
the restriction of $h$ to every stratum of ${\overM}$ is smooth,
then $h|{\mathcal{S}}$ is a pseudocycle in the sense 
of~\cite{MS} and~\cite{RT}, ie it determines an element 
of~$H_*({\mathcal{X}};\Bbb{Z})$.
In particular, if $\De$ is an immersed compact oriented submanifold
of~${\mathcal{X}}$ such that $\dim{\mathcal{S}} + \dim\De = \dim{\mathcal{X}}$, 
there is a well-defined homology-intersection number
$$\LlRr{\{h|{\mathcal{S}}\}^{-1}(\De)}\equiv
\LlRr{h^{-1}(\De),{\overS}}.$$
If ${\mathcal{Y}}$ is an immersed compact oriented submanifold
of ${\mathcal{X}}$ such~that 
$$[{\mathcal{Y}}]=[\De]\in H_*({\mathcal{X}};\Bbb{Z}),\qquad
h\big(\partial{\overS}\big)\cap{\mathcal{Y}}=\eset,$$
and $h$ is transversal to ${\mathcal{Y}}$ on ${\mathcal{S}}$, then
\begin{equation}\label{method_subs_e3}
\LlRr{\{h|{\mathcal{S}}\}^{-1}(\De)}=
\LlRr{\{h|{\mathcal{S}}\}^{-1}({\mathcal{Y}})}  \equiv\,
^{\pm} \big|\{h|{\mathcal{S}}\}^{-1}({\mathcal{Y}})\big|.
\end{equation}
Alternatively, if $\th$ is a small perturbation of $h$
on a neighborhood of ${\overS}$ in~${\overM}$,
$$\LlRr{\{h|{\mathcal{S}}\}^{-1}(\De)} =\,
^{\pm} \big|\{\th|{\mathcal{S}}\}^{-1}(\De)\big|.$$
Thus, if $h\co{\overM} \lra {\mathcal{X}}$ is a continuous map
as above such that $h|{\mathcal{S}}$ is transversal to~$\De$
and \hbox{${\mathcal{Z}} \equiv \{h|{\mathcal{S}}\}^{-1}(\De)$} is a finite~set,
\begin{equation}\label{method_subs_e4}
^{\pm} \big|{\mathcal{Z}}\big|=
\LlRr{\{h|{\mathcal{S}}\}^{-1}(\De)}-
{\mathcal{C}}_{\partial{\overS}}(h,\De),
\end{equation}
where ${\mathcal{C}}_{\partial{\overS}}(h,\De)$ denotes
the \textit{$(h,\De)$--contribution to the intersection number
$\llrr{\{h|{\mathcal{S}}\}^{-1}(\De)}$}, ie the signed cardinality
of the subset of $\th^{-1}(\De)\cap{\mathcal{S}}$ consisting of
the points that lie near $\partial{\overS}$ 
for a small generic perturbation~$\th$ of $h$ near~$\partial{\overS}$.
If the image of  a stratum ${\mathcal{Z}}_i$ of $\partial{\overS}$
under~$h$ is disjoint from $\De$, then clearly ${\mathcal{Z}}_i$
does not contribute to~${\mathcal{C}}_{\partial{\overS}}(h,\De)$.
If $h$~maps ${\mathcal{Z}}_i$  into~$\De$,
on a neighborhood of~${\mathcal{Z}}_i$ we can view $h$ and $\th$
as vector-bundle sections.
Thus, if we can understand the behavior of $h$ near~${\overS}$,
computing ${\mathcal{C}}_{\partial{\overS}}(h,\De)$ is no different
than computing ${\mathcal{C}}_{\partial{\overS}}(s)$ in 
the topological setup presented first.
On the other hand, in the cases of interest to us, 
we will be able to find a submanifold~${\mathcal{Y}}$ as in~\eqref{method_subs_e3}
such that $^{\pm} |\{h|{\mathcal{S}}\}^{-1}({\mathcal{Y}})|$
can be expressed as evaluation of tautological classes
on~${\overS}$; see Subsection~\ref{n3p3_sum_subs}, for example.
If ${\overS}$ itself is not a tautological class in
a moduli space of rational stable maps, 
we will have to describe $^{\pm} |\{h|{\mathcal{S}}\}^{-1}({\mathcal{Y}})|$
as the signed cardinality of a subset ${\mathcal{Z}}'$ 
of a space which is a step closer to a tautological class than ${\mathcal{S}}$
and apply  equation~\eqref{method_subs_e2} or~\eqref{method_subs_e4}
with~${\mathcal{Z}}'$.
Eventually, we will end up with intersections of tautological classes
on moduli spaces of stable rational maps.

In Subsection~\ref{top_dfn_subs}, we describe our topological assumptions
on ${\mathcal{S}}$, ${\overM}$, and the behavior of $s$ or $h$ near~$\partial{\overS}$.
These assumptions imply that the sets  $s^{-1}(0)\cap{\mathcal{S}}$
and $\{h|{\mathcal{S}}\}^{-1}(\De)$ are finite.
Roughly speaking, we require that $\partial{\overS}$ be contained
in a finite union of smooth manifolds ${\mathcal{Z}}_i$ such that near
each ${\mathcal{Z}}_i$ the section~$s$ or the map~$h$ can be approximated 
by a polynomial map between vector bundles over~${\mathcal{Z}}_i$. 
The polynomial map may contain terms of negative degree.
Propositions~\ref{zeros_prp} and~\ref{euler_prp} of 
Subsection~\ref{top_comp_subs}  give an inductive procedure for computing the contribution 
from each space ${\mathcal{Z}}_i$ to ${\mathcal{C}}_{\partial{\overS}}(s)$ or
to ${\mathcal{C}}_{\partial{\overS}}(h,\De)$ in good cases.
The two propositions describe how to set up a finite tree with 
topological intersection numbers assigned to the nodes and  
with integer weights assigned to the edges.
The root of the tree is assigned the first term on the right-hand side 
of~\eqref{method_subs_e2} or~\eqref{method_subs_e4}.
The number on the left-hand side of~\eqref{method_subs_e2} or~\eqref{method_subs_e4}
is a weighted sum of the numbers at the nodes.
The weight of the number assigned to a node is the product of 
the weights assigned to the edges between the node and the root.

\begin{remark}
{\rm The method presented in Subsection~\ref{top_dfn_subs}
is an improvement over that of Section~3 in~\cite{Z1} even for the
basic topological setup of the second paragraph of this subsection.
In particular, its use does not require applications of the Implicit
Function Theorem (IFT) to describe a neighborhood of $\partial{\overS}$
in~${\overS}$.  The complexity of applying the IFT increases rapidly with the
dimension of the boundary strata, as a comparison between
\cite[Subsection~5.4]{Z1} and \cite[Subsection~2.3]{Z3} suggests.}
\end{remark}

In order to apply the topological method of this paper to enumerative
problems, we use Lemma~\ref{str_lmm} and Proposition~\ref{str_prp}.
The former is a rather elementary result in complex geometry
and implies that various bundle sections over smooth strata
of moduli spaces of stable rational maps are transverse to the zero~set.
The latter depends on the explicit construction of the gluing map in~\cite{Z2}
and describes the behavior of these bundle sections near the boundary 
of each stratum.
In many cases, Proposition~\ref{str_prp}, combined with Lemma~\ref{str_lmm},
implies that natural submanifolds ${\mathcal{S}}$ of moduli spaces of stable rational maps, 
or of closely related spaces, that are needed for counting singular rational curves are
well-behaved near $\partial{\overS}$ and that the behavior near $\partial{\overS}$
of various natural vector-bundle sections over ${\mathcal{S}}$ can be approximated by polynomials.

\subsection{Computed examples}
\label{results_subs}

We now describe the main enumerative results derived in this paper
using the computational method outlined above.
These are the enumerations of triple-pointed and of tacnodal rational one-component 
curves in $\PPP$  and of rational one-component cuspidal curves in~$\P$
that pass through a collection of constraints in general position.
The reason we choose these examples is that they illustrate all aspects of our method 
and lead to new results.
The numerical values of some low-degree numbers can be found at the end of the paper.
Note that our low-degree numbers pass the standard classical checks;
see Section~\ref{tables_sec}.

We start by giving a formula describing the number of cuspidal curves
in~$\P$.
This is actually the least interesting example of the three mentioned,
as it should have really been done in~\cite{Z4}.
However, the solution to this example is easier to state and explain
than the answers to the two other primary examples.

\begin{thm}
\label{cusps_thm}
Suppose $n \ge 2$, $d \ge 1$, $N \ge 0$, and
$\mu = (\mu_1,\ldots,\mu_N)$ is 
an $N$--tuple of proper subvarieties of~$\P$ in general position
such~that
$$\sum_{l=1}^{N}\codim_{\Bbb{C}}\mu_l=d(n + 1)-2+N.$$
The number of rational cuspidal degree--$d$ curves 
that pass through the constraints~$\mu$ is given~by
$$\big|{\mathcal{S}}_1(\mu)\big|
=\sum_{k=1}^{2k\le n+2}   (-1)^{k-1}(k - 1)!
 \sum_{l=0}^{n+2-2k} \binom{n + 1}{l}
\blr{a_{\hat{0}}^l\,\eta_{\hat{0},n+2-2k-l},{\overV}_k(\mu)}.$$
\end{thm}

We now explain the notation involved in the statement of
Theorem~\ref{cusps_thm}.  The compact oriented topological
manifold~${\overV}_k(\mu)$, which in general may be an orbifold, consists
of unordered $k$--tuples of stable rational maps of total degree~$d$.
Each map comes with a special marked point~$(i,\i)$.  All these marked
points are mapped to the same point in~$\P$.  In particular, there is
a well-defined evaluation map
$$\ev_{\hat{0}}\co {\overV}_k(\mu)\lra\P$$
which sends each tuple of stable maps to the value at (any) one of
the special marked points.  We also require that the union of the
images of the maps in each tuple intersect each of the constraints
\hbox{$\mu_1,\ldots,\mu_N$}.  In fact, the elements in the tuple carry
a total of~$N$ marked points, $y_1,\ldots,y_N$,  in addition to the $k$
special marked points.  These marked points are mapped to the constraints
$\mu_1,\ldots,\mu_N$, respectively.  Roughly speaking, each element of
${\overV}_k(\mu)$ corresponds to a degree--$d$ rational curve in~$\P$,
which has at least $k$ irreducible components, and $k$ of the components
meet at the same point in~$\P$.  The precise definition of the spaces
${\overV}_k(\mu)$ can be found in Subsection~\ref{notation_subs}.

The cohomology classes~$a_{\hat{0}}$ and~$\eta_{\hat{0},l}$ 
are tautological classes in~${\overV}_k(\mu)$.
In~fact,
$$a_{\hat{0}}=\ev_{\hat{0}}^*c_1\big({\mathcal{O}}_{\P}(1)\big).$$
Let ${\overV}_k'(\mu)$ be the oriented topological orbifold
defined as~${\overV}_k(\mu)$, except without specifying
the marked points \hbox{$y_1,\ldots,y_N$} mapped to the constraints
\hbox{$\mu_1,\ldots,\mu_N$}.
Then, there is well-defined forgetful~map,
$$\pi_k\co{\overV}_k(\mu)\lra{\overV}_k'(\mu),$$
which drops the marked points $y_1,\ldots,y_N$
and contracts the unstable components.
Let 
$$\eta_{\hat{0},l}' \in H^{2l}({\overV}_k'(\mu))$$
be the sum of all degree--$l$ monomials in  
$$\psi_{(1,\i)},\ldots,\psi_{(k,\i)},$$
where $\psi_{(i,\i)}$ is the first chern class of the universal
cotangent line bundle for the marked point $(i,\i) \in \Bbb{P}^1$.
Since ${\overV}_k'(\mu)$ is a collection of \textit{un}ordered
$k$--tuples, a priori $\psi_{(i,\i)}$ may not be well defined as
an element of ${\overV}_k'(\mu)$.  However, it is easy to see that
every symmetric polynomial in $\psi_{(1,\i)},\ldots,\psi_{(k,\i)}$ is
well defined.  We put
$$\eta_{\hat{0},l}=\pi_k^*\eta_{\hat{0},l}' \in
H^{2l}({\overV}_k(\mu)).$$
In Subsection~\ref{notation_subs}, we give a definition of~$\eta_l$ that
does not involve the projection map~$\pi_k$.  The algorithm of~\cite{P}
for computing intersections of tautological classes in $\ov\M_{0,N}(d,\P)$
applies, with no change, to computing the intersection numbers involved
in the statement of Theorem~\ref{cusps_thm}.

We will call top intersections of tautological classes on
$\ov\M_{0,k}(\P,d)$ and on closely related spaces, such as
${\overV}_k(\mu)$ and projectivizations of natural vector bundles over
${\overV}_k(\mu)$, \textit{level~0 numbers}.  All such numbers can be
computed using the algorithm of~\cite{P}.  Counts of rational curves
with $s$ basic singularity conditions will be called \textit{level $s$
numbers}.  Every level $s$ number, with $s > 0$, can be written in
the form~\eqref{method_subs_e2} or~\eqref{method_subs_e4} such that
the middle term is a level~$(s{-}1)$ number.  For example, the number
$|{\mathcal{S}}_1(\mu)|$ of Theorem~\ref{cusps_thm} is a level~1 number.

Counts of rational curves with a triple point or a tacnode are level~2 numbers.
Indeed, the sets of such curves are subsets of the space of one-component rational curves
with a node.
Counts of such curves are level~1 numbers,
since the next level down are the rational curves that pass through the given constraints;
see the first paragraph of Subsection~\ref{method_subs}.
Thus, the first two theorems below express level~2 numbers in terms of level~1 numbers.
After stating them, we give some clarification on the notation involved and  
then state several lemmas that express the relevant level~1 numbers 
in terms of level~0 numbers.
 
\begin{thm}
\label{n3p3_thm}
Let $d$, $p$, and $q$ be nonnegative integers such that $2p + q = 4d - 3$.
The number of rational one-component degree--$d$ curves 
that have a triple point and 
pass through a tuple $\mu$ of $p$~points and $q$~lines 
in general position in~$\PPP$ is $\frac{1}{6}|{\mathcal{V}}_1^{(2)}(\mu)|$, where
\begin{equation*}\begin{split}
\big|{\mathcal{V}}_1^{(2)}(\mu)\big|&=
\big|{\mathcal{V}}_1^{(1)}(\mu + H^0)\big|+
\big\lan a_{\hat{0}},{\overV}_1^{(1)}(\mu + H^1)\big\ran+
\big\lan 16a_{\hat{0}} + 8\eta_{\hat{0},1},
{\overS}_1(\mu)\big\ran\\
&\qquad+2\big|{\mathcal{V}}_2^{(1)}(\mu)\big|
-\big\lan (12 - d)a_{\hat{0}}^2 + 8a_{\hat{0}}\eta_{\hat{0},1}
 + 2\eta_{\hat{0},1}^2,{\overV}_1^{(1)}(\mu)\big\ran
-2\big|{\mathcal{S}}_2(\mu)\big|.
\end{split}\end{equation*}
\end{thm}

\begin{thm}
\label{n3tac_thm}
Let $d$, $p$, and $q$ be nonnegative integers such that $2p{+}q
= 4d{-}3$.  The number of rational one-component degree--$d$
curves that have a tacnodal point and pass through a tuple $\mu$
of $p$~points and $q$~lines in general position in~$\PPP$ is
$\frac{1}{2}|{\mathcal{S}}_1^{(1)}(\mu)|$, where
\begin{equation*}\begin{split}
\big|{\mathcal{S}}_1^{(1)}(\mu)\big|
= \big\lan 6a_{\hat{0}}^2 + \eta_{\hat{0},1}^2,
{\overV}_1^{(1)}(\mu)\big\ran+
\big\lan 4a_{\hat{0}} +\frac{1}{2}\eta_{\hat{0},1},
{\overV}_2^{(1,1)}(\mu)\big\ran+7\big|{\mathcal{S}}_2(\mu)\big| \qquad& \\
-\big\lan 20a_{\hat{0}} + 19\eta_{\hat{0},1},
{\overS}_1(\mu)\big\ran
-2\big|{\mathcal{V}}_2^{(1)}(\mu)\big|.& 
\end{split}\end{equation*}
\end{thm}

\begin{figure}[ht!]\small
\begin{pspicture}(-1.6,-2)(10,0)
\psset{unit=.36cm}
% 1st diagram
\psarc(-5,-3){5}{0}{35}\psarc(5,-3){5}{145}{180}\pscircle*(0,-3){.17}
\rput(0,-4.3){${\mathcal{S}}_1$}
% 2nd diagram
\psarc(5,-3.6){1.9}{15}{165}\psarc(5,.2){1.9}{195}{345}
\pscircle*(5,-1.7){.17}
\rput(5,-4.3){${\mathcal{S}}_2$}
% 3rd diagram
\psarc(11.91,-.79){2}{180}{260}\psarc(9.09,-.79){2}{280}{0}
\psarc(10.5,-.79){.59}{0}{180}\pscircle*(10.5,-2.2){.17}
\rput(10.7,-4){${\mathcal{V}}_1^{(1)}$}
% 4th diagram
\psarc(16.91,-.79){2}{180}{260}\psarc(14.09,-.79){2}{280}{0}
\psarc(15.5,-.79){.59}{0}{180}\psline(13.5,-2.2)(17.5,-2.2)\pscircle*(15.5,-2.2){.17}
\rput(15.7,-4){${\mathcal{V}}_2^{(1)}$}
% 5th diagram
\psarc(21.42,-1.35){.95}{116.53}{243.47}\psarc(20.58,-1.35){.95}{-63.47}{63.47}
\psarc(20.58,.35){.95}{-63.47}{-25}\psarc(21.42,.35){.95}{205}{243.47}
\psarc(20.58,-3.05){.95}{25}{63.47}\psarc(21.42,-3.05){.95}{116.53}{155}
\pscircle*(21,-2.2){.17}\pscircle*(21,-.5){.17}
\rput(21.3,-4){${\mathcal{V}}_2^{(1,1)}$}
% 6th diagram
\psline(25.1,-.14)(26.9,-2.66)\psline(27.9,-.14)(26.1,-2.66)
\psline(24.9,-.7)(28.1,-.7)
\pscircle*(26.5,-2.1){.17}\pscircle*(25.5,-.7){.17}\pscircle*(27.5,-.7){.17}
\rput(26.95,-4){${\mathcal{V}}_{2,(0,1)}^{(1;0,1)}$}
\end{pspicture}
\caption{Images of typical elements of ${\mathcal{S}}_*$ and ${\mathcal{V}}_*^*$}
\label{images_fig}
\end{figure}
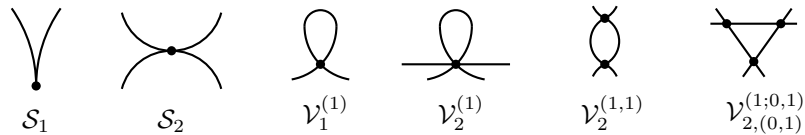

We define the spaces ${\overV}_k^{(1)}(\mu)$ as follows.  Let
${\mathcal{V}}_k^{(1)}(\mu)$ be the space of $k$--tuples of stable maps
as in the construction of the space ${\overV}_k(\mu)$, but with the
following exceptions.  Every element of each $k$--tuple lies in the main
stratum of the appropriate moduli space of stable maps, ie the domain
of the map is~$\Bbb{P}^1$.  Furthermore, one of the elements of each
$k$--tuple~$b$ carries a special marked point, labeled by~$\hat{1}$,
and the value of the map at this point is $\ev_{\hat{0}}(b)$.  In the
space \smash{${\mathcal{V}}_2^{(1,1)}(\mu)$} each of the two components carries
a special marked point, one of which is labeled  by $\hat{1}$ and the
other by~$\hat{2}$.  Furthermore, $\ev_{\hat{1}}(b) = \ev_{\hat{2}}(b)$
for all $2$--tuples $b$ in~${\mathcal{V}}_2^{(1,1)}(\mu)$.  The spaces
${\overV}_k^{(1)}(\mu)$ and ${\overV}_2^{(1,1)}(\mu)$  are the closures
of the spaces \smash{${\mathcal{V}}_k^{(1)}(\mu)$} and
\smash{${\mathcal{V}}_2^{(1,1)}(\mu)$}
in the unions of the appropriate products of moduli spaces of stable
rational maps.  We denote by ${\mathcal{S}}^*_*(\mu)$ the subspace of
${\mathcal{V}}^*_*(\mu)$ consisting of tuples of maps with the simplest
possible additional natural singularity.  For example, the differential
of every element of ${\mathcal{S}}_1(\mu)$ vanishes at~$(1,\i)$.  The set
${\mathcal{S}}_2(\mu)$ is described in detail by Lemma~\ref{n3c2tac_lmm}.
Figure~\ref{images_fig} depicts the images of typical elements of these
spaces as well as of \smash{${\mathcal{V}}_{2,(0,1)}^{(1;0,1)}(\mu)$},
which appears in a relationship between level~1 numbers; see the remark
following the proof of Lemma~\ref{n3tac_contr_lmm4b}.  We give formal
definitions of all these spaces in Subsections~\ref{n3p3_sum_subs},
\ref{n3p3_corr_subs}, and~\ref{n3tac_sum_subs}.  Finally, $\mu + H^r$
denotes the $(N + 1)$--tuple of constrains $(\mu_1,\ldots,\mu_N,H^r)$,
where $H^r$ is a generic linear subspace of $\P$ of complex dimension~$r$.

\begin{lmm}
\label{n3circ3_lmm}
Suppose $d$, $p$, and $q$ are nonnegative integers such that $2p{+}q = 4d{-}3$
and $\mu$ is a tuple of $p$~points and $q$~lines in general position in~$\PPP$.
The number of rational connected two-component degree--$d$ curves that pass through
the constraints $\mu$ and such that one of the components of each curve has a node and
the other component is attached at the node of the first component as depicted in Figure~\ref{images_fig} is $\frac{1}{2}|{\mathcal{V}}_2^{(1)}(\mu)|$, where
\begin{equation*}\begin{split}
\big|{\mathcal{V}}_2^{(1)}(\mu)\big|=
\big|{\mathcal{V}}_2(\mu + H^0)\big|+
\big\lan a_{\hat{0}},{\overV}_2(\mu + H^1)\big\ran
+3\big|{\mathcal{V}}_3(\mu)\big| \qquad&\\
-\big\lan (12 - d)a_{\hat{0}}^2 + 4a_{\hat{0}}\eta_{\hat{0},1}
 + 2\eta_{\hat{0},2} - \eta_{\hat{0},1}^2,{\overV}_2(\mu)\big\ran.&
\end{split}\end{equation*}
\end{lmm}

\begin{lmm}
\label{n3c2tac_lmm}
Suppose $d$, $p$, and $q$ are nonnegative integers such that $2p{+}q = 4d{-}3$
and $\mu$ is a tuple of $p$~points and $q$~lines in general position in~$\PPP$.
The number of rational connected two-component degree--$d$ curves 
that pass through the constraints $\mu$ and have a tacnodal point is given~by
$$|{\mathcal{S}}_2(\mu)|=
\big\lan 6a_{\hat{0}}^2 + 4a_{\hat{0}}\eta_{\hat{0},1} + \eta_{\hat{0},2},
{\overV}_2(\mu)\big\ran-3\big|{\mathcal{V}}_3(\mu)\big|.$$
\end{lmm}

\begin{lmm}
\label{n3circ2_lmm}
Suppose $d$, $p$, and $q$ are nonnegative integer such that $2p + q = 4d - 3$
and $\mu$ is a tuple of $p$~points and $q$~lines in general position in~$\PPP$.

{\rm (1)}\qua The number of rational connected two-component degree--$d$
curves, with the components arranged in a circle as in Figure~\ref{images_fig},
that pass through the constraints~$\mu$ and have one of the nodes
on a generic hyperplane 
is given~by
$$\frac{1}{2}\big\lan a_{\hat{0}},{\overV}_2^{(1,1)}(\mu)\big\ran=
\big\lan a_{\hat{0}},{\overV}_2(\mu + \{H^1 : H^2\})\big\ran
-\big\lan 4a_{\hat{0}}^2 + a_{\hat{0}}\eta_{\hat{0},1},
{\overV}_2(\mu)\big\ran.$$
{\rm (2)}\qua Furthermore,
\begin{equation*}\begin{split}
\frac{1}{2}
\big\lan \eta_{\hat{0},1},{\overV}_2^{(1,1)}(\mu)\big\ran =&
\big\lan \eta_{\hat{0},1},{\overV}_2(\mu + \{H^1  : H^2\})\big\ran
+ \big|{\mathcal{V}}_2(\mu + H^0)\big|\\
&\qquad+\big\lan a_{\hat{0}},{\overV}_2(\mu + H^1)\big\ran
+ d\big\lan a_{\hat{0}}^2,{\overV}_2(\mu)\big\ran
-3\big|{\mathcal{V}}_3(\mu)\big|.
\end{split}\end{equation*}
\end{lmm}

\begin{lmm}
\label{n3cusps_lmm}
Suppose $d$, $p$, and $q$ are nonnegative integers such that $2p + q = 4d - 3$ 
and $\mu$ is a tuple of $p$ points and $q$ lines in general position in~$\PPP$.

{\rm (1)}\qua The number of rational degree--$d$ curves that pass through
the constraints $\mu$
and have a cusp on a generic hyperplane is given~by
$$\lan a_{\hat{0}},{\overS}_1(\mu)\ran=
\big\lan 6a_{\hat{0}}^3\eta_{\hat{0},1} + 
4a_{\hat{0}}^2\eta_{\hat{0},1}^2 + 
a_{\hat{0}}\eta_{\hat{0},1}^3,{\overV}_1(\mu)\big\ran-
\big\lan 4a_{\hat{0}}^2 + a_{\hat{0}}\eta_{\hat{0},1},
{\overV}_2(\mu)\big\ran.$$
{\rm (2)}\qua Furthermore,
$$\lan \eta_{\hat{0},1},{\overS}_1(\mu)\ran=
\big\lan 4a_{\hat{0}}^3\eta_{\hat{0},1} + 
6a_{\hat{0}}^2\eta_{\hat{0},1}^2 + 4a_{\hat{0}}\eta_{\hat{0},1}^3 + 
\eta_{\hat{0},1}^4,{\overV}_1(\mu)\big\ran
-\big|{\mathcal{V}}_3(\mu)\big|.$$
\end{lmm}

\begin{lmm}
\label{n3p2_lmm}
Suppose $d$, $p$, and $q$ are nonnegative integers and $\mu$ is a 
tuple of $p$~points and $q$~lines in general position in~$\PPP$.

{\rm (1)}\qua If $2p + q = 4d - 1$, the number of rational one-component
degree--$d$ curves that  pass through the constraints~$\mu$ and have a node
is $\frac{1}{2}|{\mathcal{V}}_1^{(1)}(\mu)|$, where
$$|{\mathcal{V}}_1^{(1)}(\mu)|=
\big\lan (2d - 6)a_{\hat{0}}^2 - 4a_{\hat{0}}\eta_{\hat{0},1}
 - \eta_{\hat{0},1}^2,{\overV}_1(\mu)\big\ran
+\big|{\mathcal{V}}_2(\mu)\big|.$$
{\rm (2)}\qua If $2p + q = 4d - 2$, the number of rational one-component
degree--$d$ curves that pass through the constraints~$\mu$ and have a
node on a generic hyperplane is $\frac{1}{2}\lan
a_{\hat{0}},{\overV}_1^{(1)}(\mu)\ran$, where
\begin{multline*}
\lan a_{\hat{0}},{\overV}_1^{(1)}(\mu)\ran=
\big\lan(2d - 6)a_{\hat{0}}^3 - 
4a_{\hat{0}}^2\eta_{\hat{0},1} - 
a_{\hat{0}}\eta_{\hat{0},1}^2,{\overV}_1(\mu)\big\ran \\
+\big\lan a_{\hat{0}}^2,{\overV}_1(\mu + H^1)\big\ran
+\big\lan a_{\hat{0}},{\overV}_2(\mu)\big\ran.
\end{multline*}
{\rm (3a)}\qua If $2p + q = 4d - 3$, the number of rational one-component
degree--$d$ curves that pass through the constraints~$\mu$
and have a node on a generic line  is $\frac{1}{2}\lan
a_{\hat{0}}^2,{\overV}_1^{(1)}(\mu)\ran$, where
$$\lan a_{\hat{0}}^2,{\overV}_1^{(1)}(\mu)\ran=
2\big\lan a_{\hat{0}}^3,{\overV}_1(\mu + H^1)\big\ran
-\big\lan 4a_{\hat{0}}^3\eta_{\hat{0},1} + 
a_{\hat{0}}^2\eta_{\hat{0},1}^2,{\overV}_1(\mu)\big\ran
+\big\lan a_{\hat{0}}^2,{\overV}_2(\mu)\big\ran.$$
{\rm (3b)}\qua Furthermore,
\begin{equation*}\begin{split}
\big\lan a_{\hat{0}}\eta_{\hat{0},1},
                      {\overV}_1^{(1)}(\mu)\big\ran
=&~ \big\lan a_{\hat{0}}\eta_{\hat{0},1},
                         {\overV}_1(\mu + H^0)\big\ran
+\big\lan a_{\hat{0}}^2\eta_{\hat{0},1},
                         {\overV}_1(\mu + H^1)\big\ran\\
&\qquad\qquad +d\big\lan a_{\hat{0}}^3\eta_{\hat{0},1},
           {\overV}_1(\mu)\big\ran
-\big\lan 4a_{\hat{0}}^2 + 
a_{\hat{0}}\eta_{\hat{0},1},{\overV}_2(\mu)\big\ran\\
\big\lan \eta_{\hat{0},1}^2,{\overV}_1^{(1)}(\mu)\big\ran
=&~ \big\lan \eta_{\hat{0},1}^2,{\overV}_1(\mu + H^0)\big\ran
+ \big\lan a_{\hat{0}}\eta_{\hat{0},1}^2,
                                  {\overV}_1(\mu + H^1)\big\ran\\
&\qquad\qquad +\big\lan 4a_{\hat{0}}^3\eta_{\hat{0},1}
+d\cdot a_{\hat{0}}^2\eta_{\hat{0},1}^2,
	                           {\overV}_1(\mu)\big\ran
-\big|{\mathcal{V}}_3(\mu)\big|.
\end{split}\end{equation*}
\end{lmm}

Every term on the right-hand side of each expression in 
Lemmas~\ref{n3circ3_lmm}--\ref{n3p2_lmm} is a level~0 number,
ie it is a top intersection of tautological classes
in a product of moduli spaces of stable maps and thus
is computable via the explicit formulas of~\cite{P}.
We define the space ${\mathcal{V}}_2(\mu + \{H^1 : H^2\})$
in the same way as ${\mathcal{V}}_2(\mu{+}H^1{+}H^2)$,
with the only exception that we require $H_1$ and $H_2$ 
to lie on different elements of the~tuple~$b$.

\begin{remark}
{\rm If one were to derive a completely general recursive formula for counting 
rational curves with singularities, no separate formula would be necessary for 
generalizations of Lemma~\ref{n3circ2_lmm} part (2), 
Lemma~\ref{n3cusps_lmm} part (2), and Lemma~\ref{n3p2_lmm} part (3b).
Using \cite{P}, one can express all classes $\eta_{*,*}^*$
on products of moduli spaces of stable rational maps in terms of
subspaces of possibly other products of 
moduli spaces of stable rational maps that consist of stable maps 
sending their marked points to various constraints in~$\P$.
In many cases, using Proposition~\ref{str_prp},
one  can thus express evaluations of the classes $\eta_{*,*}^*$
on a space ${\overS}$ of maps that represent curves
with certain singularities in terms of the numbers of singular
curves that pass through various constraints.
Furthermore, the level of the latter numbers will be no higher than
that~of~${\mathcal{S}}$.}
\end{remark}

We prove Theorems~\ref{n3p3_thm} and~\ref{n3tac_thm} 
in Sections~\ref{n3p3_sec} and~\ref{n3tac_sec}.
In particular, we describe the structure of the spaces
\smash{${\overV}_1^{(1)}(\mu)$} and  \smash{${\overV}_2^{(1,1)}(\mu)$}
in Subsection~\ref{n3p3_str_subs} and conclude
that they define a homology class in 
a compact oriented stratified topological orbifold.
Since ${\overS}_1(\mu)$ is shown to be an oriented topological manifold
in \cite[Subsection~5.4]{Z1}, 
it follows that all the terms on the right-hand side of
the formulas in the two propositions are well-defined.
In Subsections~\ref{n3p3_sum_subs} and~\ref{n3tac_sum_subs},
we write ${\mathcal{V}}_1^{(2)}(\mu)$ and ${\mathcal{S}}_1^{(1)}(\mu)$
in the form~\eqref{method_subs_e4} and~\eqref{method_subs_e2},
respectively, and express the first term on the right-hand side
in terms of evaluations of tautological classes
on ${\overV}_1^{(1)}(\mu)$ and on related spaces.
Lemmas~\ref{n3c2tac_lmm} and~\ref{n3cusps_lmm} are proved in~\cite{Z1}.
The first statement of Lemma~\ref{n3p2_lmm} is a special case
of \cite[Theorem~1.1]{Z4}.
The remaining statements of Lemma~\ref{n3p2_lmm}
and Lemmas~\ref{n3circ3_lmm} and~\ref{n3circ2_lmm} 
are proved in Section~\ref{level1nums_sec}.
In Subsections~\ref{n2p3_subs} and~\ref{n2tac_subs},
we show that our method recovers the formulas of~\cite{KQR} and~\cite{V} 
solving Questions~\ref{g0n2p3_ques} and~\ref{g0n2tac_ques},
which are the $\PP$ analogues of the problems addressed by
Theorems~\ref{n3p3_thm} and~\ref{n3tac_thm}.
We conclude by proving Theorem~\ref{cusps_thm} in Subsection~\ref{cusps_subs},
where we construct a tree of contributions and thus
illustrate a point made at the end of Subsection~\ref{method_subs}.

\subsection{General remarks}
\label{general_subs}

Given the claims made in the abstract and at the end of Subsection~\ref{back_subs},
the reader may wonder why this paper is so long,
why the notation is so involved, and why a more general case is not done.
Doing a more general case, instead of the examples we work out,
may in fact shorten this paper.
However, the additional notation needed to describe a general case
is likely to completely obscure the computational method presented~here.

The topological part of our method consists of
Propositions~\ref{zeros_prp} and~\ref{euler_prp}.
The somewhat involved notation of Subsection~\ref{top_dfn_subs}
formally states what it means to take the leading term(s)
of a section along the normal direction to a submanifold.
Proposition~\ref{str_prp} gives power-series expansions for all relevant vector-bundle 
sections near all boundary strata of moduli spaces of stable rational maps.
Describing the terms involved in the power-series expansions 
requires quite a bit of notation.
However, as we will see in later sections, very few boundary strata
actually matter in our computations,
and the expansions of Proposition~\ref{str_prp} corresponding
to such strata are rather simple.
In practice, it is best to draw a tree of these simple strata along with 
all the relevant topological data;
then deriving formulas such as those of Theorems~\ref{n3p3_thm}
and~\ref{n3tac_thm} becomes a nearly-mechanical~task.

The sections and linear maps between vector bundles that we introduce in Subsection~\ref{analysis_subs} are described in an analytic way.
Nevertheless, it is likely that the numbers of zeros of these and related linear maps 
have an algebraic interpretation, and that the same is true of our entire computational approach.
Furthermore, there seem to be some general properties that remain to be explored.
For example, finding any difference between our formulas for $|{\mathcal{S}}_1(\mu)|$
in Theorem~\ref{cusps_thm} and for the genus-one correction term~$CR_1(\mu)$
in \cite[Theorem~1.1]{Z4} requires a rather careful comparison of the~two.
One may also notice some similarity between the expressions
for various level~1 numbers in Lemmas~\ref{n3circ3_lmm}--\ref{n3p2_lmm}.
A topological property of the number of zeros of an affine map
in a simple case is described in \cite[Corollary~4.7]{Z3}.

We conclude this introductory section by describing some classes
of enumerative problems definitely and likely solvable 
by the method described in this paper.
Suppose ${\mathcal{C}}$ is a one-component curve in~$\P$, 
or in another algebraic manifold~$M$.
Let $u\co\ti{\mathcal{C}} \lra {\mathcal{C}}$ be a normalization 
of~${\mathcal{C}}$, ie $\ti{\mathcal{C}}$ is a smooth connected complex curve and
$u\co\ti{\mathcal{C}} \lra M$ is a holomorphic map such 
that the image of $u$ is $\ti{\mathcal{C}}$ and $u$ is one-to-one
outside of a finite set of points of~$\tilde{\mathcal{C}}$.
If $\tilde{p}$ is a point in~$\ti{\mathcal{C}}$, 
let $\si_u(\ti{p})$ be the nonnegative integer such that 
the first $\si_u(\ti{p})$ derivatives of $u$ at $\tilde{p}$ vanish, 
but the derivative of order $\si_u(\ti{p}) + 1$ does not vanish at~$p$.
For example, if $\si_u(\ti{p}) = 0$, $du|_{\ti{p}} \neq 0$,
ie $p$ is a smooth point of the branch of the curve~${\mathcal{C}}$
corresponding to~$\ti{p}$.
If $\si_u(\ti{p}) = 1$, 
$$du|_{\ti{p}} = 0,\qquad\hbox{but}\qquad 
{\mathcal{D}}_{\ti{p}}^{(2)}u\equiv {\mathcal{D}}du|_{\tilde{p}}\neq0,$$ 
where ${\mathcal{D}}$ denotes the covariant differentiation 
with respect to some connection in~$TM$.
In other words, the branch of the curve ${\mathcal{C}}$ at~$p$
corresponding to~$\tilde{p}$ has a cusp at~$p$.
If $p$ is a point of~${\mathcal{C}}$, we denote by $\si_0(p)$ 
the set of branches of ${\mathcal{C}}$ at~$p$;
in particular, $|\si_0(p)| = |u^{-1}(p)|$.
For each $i \in \si_0(p)$, let $\si(p;i) = \si_u(\tilde{p}_i)$
if $\tilde{p}_i$ is the point in $\tilde{\mathcal{C}}$ such that
$u$ maps a small neighborhood of~$\tilde{p}_i$ into the branch~$i$ at~$p$.
Let $\si_0(p) = \bigsqcup\si_{0,k}(p)$ be the partition such~that
$${\mathcal{D}}_{\tilde{p}_i}^{(\si_u(\tilde{p}_i)+1)}u \parallel
{\mathcal{D}}_{\tilde{p}_j}^{(\si_u(\tilde{p}_j)+1)}u 
\qquad\Llra\qquad
i,j \in \si_{0,k}(p)\quad\hbox{for some}~~k.$$
For example, if $\si_0(p) = \si_{0,1}(p) = \{\tilde{p}_1,\tilde{p}_2\}$
and $\si(p;1) = \si(p;2) = 0$, the curve ${\mathcal{C}}$ has a tacnode at~$p$.
We take 
$$\under{\si}(p)=(\si_0(p),\{\si_{0,k}(p)\},\{\si(p;i)\}).$$
The infinite set $\{\under{\si}(p) :p \in {\mathcal{C}}\}$ describes
the singularities of the curve~${\mathcal{C}}$.
However, for all but finitely many points $p$ in~${\mathcal{C}}$,
$\si_0(p) = \{i\}$ is a single-element set and $\si(p;i) = 0$.
Thus, we say that the curve ${\mathcal{C}}$ \textit{has the set of singularities}
$$\{\under{\si}(\al) : \al = 1,\ldots,N\},
\qquad\hbox{where}\qquad
\under{\si}(\al)=\big(\si_0(\al),\{\si_{0,k}(\al)\},\{\si(\al;i)\}\big),$$
if there are distinct points $p_1,\ldots,p_N$ of ${\mathcal{C}}$
such that for all $\al$ and some identification of $\si_0(\al)$
with the set of branches of ${\mathcal{C}}$ at $p_{\al}$,
$\under{\si}({\al}) = \under{\si}(p_{\al})$.
The method of this paper can be used to determine the number of
one-component rational curves in a projective space that have 
any one-element set of singularities of the form 
$\big\{\bigl(\si,\{\si\},\{\si(i)\}\bigr)\big\}$.
In other words, the curves are to have one $|\si|$--fold singular point
and the branches of the curves are to have cusps of
the orders~$\{\si(i)\}$.
In particular, we are not imposing any tacnodal kind of condition.
The singular point of the curves may be required to fall on a subvariety.

Some types of singularities that cause problems for this method are
the flex, two nodes, two cusps, and the tacnode.
The reason is that the expansions of bundle sections
given in Propositions~\ref{str_prp} are not sufficiently fine 
in the cases when the boundary stratum involves curves of 
very low, but positive, degree.
For example, one of the strata of the space of 
one-component degree--$d$ rational curves with one marked point~$y_1$
consists of two-component curves, one of which has degree one
and carries the marked point~$y_1$.
Every element of this stratum has a flex at~$y_1$,
ie is a zero of a certain bundle section~$s$.
However, Proposition~\ref{str_prp} does not give a sufficiently fine 
description of the behavior of~$s$ near this stratum.
On the other hand, a simple dimension-counting argument shows that 
such a boundary stratum cannot occur if the dimension~$n$ 
of the projective space is~two.
If we would like to count one-component curves that have 
two cusps, or two nodes, or a tacnode, 
the problem stratum is the one consisting
of two-component curves, one of which is a double line 
and carries two marked points.
Again, a dimension-counting argument shows that this boundary 
stratum does not occur unless $n$ is at least $6$, $4$, or $5$,
respectively.
Indeed, among the examples worked out in this paper
are the enumerations of tacnodal rational curves in $\PP$ and~$\PPP$.

Many types of curves with multiple components and with first-order tangency conditions 
to disjoint subvarieties can be counted as well.
In fact, as the results described in Subsection~\ref{results_subs} indicate, 
counting one-component singular curves involves
counting multiple-components curves with simpler singularities.
We plan to elaborate more on what types of curves can be counted
and why in a later paper.
Finally, due to the explicit nature of the gluing maps used,
it should be possible to sharpen the expansions of Propositions~\ref{str_prp}
along the few problem strata that appear in more general cases.
If so, every enumerative problem, in the sense described above, will be~solvable.

\begin{remark}
{\rm This paper concerns counting curves in $\P$, but our method may apply to 
counting curves with singularities in some other Kahler manifolds as well.
Our aim is to express curve counts in terms of top intersections of
tautological classes on moduli spaces of stable rational maps.
Thus, in order to obtain actual numbers we make use~of:

{\bf Fact 0}\qua Top intersections of tautological classes on
$\ov\M_{0,k}(\P,d)$ are computable.
It is essential for the method itself that the moduli spaces $\ov\M_{0,k}(\P,d)$
have the expected structure. This is due~to

{\bf Fact 1}\qua If $u\co\Bbb{P}^1 \lra \P$ is a holomorphic map, then
$$H^1(\Bbb{P}^1;u^*T\P \otimes {\mathcal{O}}_{\Bbb{P}^1}(-1))=0.$$
In order for our method to apply to a Kahler manifold $M$,
Fact~1 needs to hold  with $\P$ replaced by~$M$ for maps up to the relevant ``degree", 
ie the homology class of the curves to be counted.
Additional positivity conditions, dependent on the type of singularities involved,
need to be satisfied as well.
For example, in order to count curves with a cusp in $\P$, we rely~on

{\bf Fact 2}\qua If $u\co\Bbb{P}^1 \lra \P$ is a nonconstant holomorphic
map, then
$$H^1(\Bbb{P}^1;u^*T\P \otimes {\mathcal{O}}_{\Bbb{P}^1}(-2))=0.$$
Theorem~\ref{cusps_thm} is valid with $\P$ replaced by~$M$ 
as long as Facts~1 and~2 hold  with $\P$ replaced by~$M$ for maps of
degree up~$d$.}
\end{remark}

\section{Topology}
\label{topology_sec}

\subsection{Monomials maps}
\label{prelim_subs}

This section contains details of the topological aspects of the
computational approach of this paper.  We describe the setting in the
next subsection and state and justify the reductive method for computing
boundary contributions in Subsection~\ref{top_comp_subs}.  The present
subsection collects a few basic facts that are used elsewhere in this
section.  The key statements here are Definition~\ref{monomials_dfn1a}
and Propositions~\ref{monomials_prp1} and~\ref{monomials_prp2}.

We denote by $\overR^+$ the set of nonnegative reals.
Let $\be\co{\overR}^+ \lra [0,1]$ be a smooth cutoff function such that
$$\be(t)=\begin{cases}
0&\hbox{if~}t \le 1,\\
1&\hbox{if~}t \ge 2;
\end{cases}
\qquad\hbox{and}\qquad
\be'(t) > 0~~\hbox{if~}t \in (1,2).$$
If $\de > 0$, let $\be_{\de} \in C^{\i}\bigl({\overR}^+;\Bbb{R}\bigr)$ be given by 
$\be_{\de}(t) = \be\big(\de^{-\frac{1}{2}}t\big)$.
We also denote by $\be_{\de}$ the natural extension of $\be_{\de}$
to~$\Bbb{C}^n$:
$$\be_{\de}\big(\under{z}\big)= \be_{\de}\big(|\under{z}|\big),$$
where $|\under{z}|=\sqrt{|z_1|^2 + \ldots + |z_n|^2}$
if $\under{z} = \big(z_1,\ldots,z_n\big)\in\Bbb{C}^n$.
We write $B_{\de}(0,\Bbb{C}^n)$ for the open ball of radius $\de$
about $0$ in~$\Bbb{C}^n$.
Let $Y_n$ be the union of the $n$ codimension-one coordinate subspaces 
$\Bbb{C}^k \times \{0\} \times \Bbb{C}^{n-1-k}$ in~$\Bbb{C}^n$.

\begin{dfn}
\label{monomials_dfn1a}
{\rm
Suppose $m$ and $n$ are positive integers and 
${\mathcal{A}} = (a)_{ij}$ is an $m \times n$ integer matrix.  Then:

(1)\qua A function $\rho\co\Bbb{C}^n - Y_n \lra \Bbb{C}^m - Y_m$
is a \emph{degree-${\mathcal{A}}$ monomials map} if
$$\rho\big(z_1,\ldots,z_n\big)=
\big(z_1^{a_{1,1}}\ldots z_n^{a_{1,n}},\ldots,
z_1^{a_{m,1}}\ldots z_n^{a_{m,n}}\big)$$
for all $\big(z_1,\ldots,z_n\big) \in \Bbb{C}^n - Y_n$.

(2)\qua A degree--${\mathcal{A}}$ monomials map $\rho$ is
\emph{nondegenerate} if $\rk\rho \equiv \rk{\mathcal{A}}=m$.

(3)\qua If $m = n$, a nondegenerate degree--${\mathcal{A}}$ monomials
map $\rho$ is
\emph{positive} (or \emph{negative}) if all components of
the~vector  ${\mathcal{A}}^{-1}\under{1}$
are positive (or negative).

(4)\qua If $m = n$, a nondegenerate degree--${\mathcal{A}}$ monomials
map $\rho$ is \emph{neutral} if it is neither positive nor negative.}
\end{dfn}

In (3) above, $\under{1}$ denotes the column vector of length $n$ 
consisting of all ones.
If $m = n$, let $\det\rho$ denote the determinant of the square matrix~${\mathcal{A}}$.

\begin{dfn}
\label{monomials_dfn1b}
{\rm
Suppose $n$ is a positive integer,
${\mathcal{A}}$ is an $n \times n$ nondegenerate integer matrix,
and ${\mathcal{A}}_1$ and ${\mathcal{A}}_2$ are row vectors of length~$n$.
Let $\rho$, $\rho_1$, and $\rho_2$ be monomials maps
of degrees ${\mathcal{A}}$, ${\mathcal{A}}_1$, and ${\mathcal{A}}_2$, respectively.

(1)\qua If $\rho$ is a positive monomials map, $\rho_1>_{\rho}\rho_2$ 
if ${\mathcal{A}}_1{\mathcal{A}}^{-1}\under{1}<
{\mathcal{A}}_2{\mathcal{A}}^{-1}\under{1}$.

(2)\qua If $\rho$ is a negative monomials map, $\rho_1>_{\rho}\rho_2$ 
if ${\mathcal{A}}_1{\mathcal{A}}^{-1}\under{1}>
{\mathcal{A}}_2{\mathcal{A}}^{-1}\under{1}$.}
\end{dfn}

\begin{remark}
{\rm The notions of nondegenerate, positive, negative,
and neutral of Definition~\ref{monomials_dfn1a}
are invariant under every reordering of coordinates 
on the domain and/or the target space.
The same is true of the partial-order relation introduced by 
Definition~\ref{monomials_dfn1b}.
We describe geometric consequences of these properties~below.}
\end{remark}

\stepcounter{thm}

\begin{alprp}
\label{monomials_prp1}
If $\rho\co\Bbb{C}^n - Y_n \lra \Bbb{C}^n - Y_n$ is 
a degree--${\mathcal{A}}$ neutral monomials map and ${\mathcal{K}}$ is a compact subset of
$\Bbb{C}^n - Y_n$, 
there exists \hbox{$\de^* = \de^*({\mathcal{A}},{\mathcal{K}}) \in \Bbb{R}^+$} 
such~that
$$\big(B_{\de^*}(0,\Bbb{C}^n) - Y_n\big)\cap
\rho^{-1}(\Bbb{R}^+\cdot{\mathcal{K}})=\eset.$$
\end{alprp}

\begin{alprp}
\label{monomials_prp2}
Suppose $\rho\co\Bbb{C}^n - Y_n \lra \Bbb{C}^n - Y_n$ is 
a degree--${\mathcal{A}}$ positive (or negative) monomials map
and ${\mathcal{K}}$ is a precompact open subset of $\Bbb{C}^n - Y_n$.
Then:

{\rm(1)}\qua The set $\rho^{-1}(\Bbb{R}^+\cdot{\overK})$ is closed in
$\Bbb{C}^n - \{0\}$.

{\rm(2)}\qua For every $\de^* \in \Bbb{R}^+$, $\de \in (0,\de^*)$,
and $\de_+ \in (\de^{1/2},\i)$,
there exists 
$$\ep = \ep({\mathcal{A}},\de^*,\de,\de_+,{\mathcal{K}}) \in\Bbb{R}^+$$
such that for all $t \in (0,\ep)$, the map
$$t^{-1}\be_{\de}\rho\co
\big(B_{2\de^{*1/2}}(0,\Bbb{C}^n) - Y_n\big)\cap
\big\{t^{-1}\be_{\de}\rho\big\}^{-1}({\mathcal{K}})\lra {\mathcal{K}}$$
is a smooth covering projection of 
oriented order \hbox{$|\det\rho|$ ($-|\det\rho|$)}.
Furthermore,
$$\big(B_{2\de^{*1/2}}(0,\Bbb{C}^n) - Y_n\big)\cap
\big\{t^{-1}\be_{\de}\rho\big\}^{-1}({\mathcal{K}})\subset
B_{\de_+}(0,\Bbb{C}^n).$$
{\rm(3)}\qua If $\rho_1$ and $\rho_2$ are monomials maps of degrees
${\mathcal{A}}_1$ and ${\mathcal{A}}_2$ such that $\rho_1 > _{\rho}\rho_2$, 
for every $\ep \in \Bbb{R}^+$, there exists 
$\de = \de({\mathcal{A}},{\mathcal{A}}_1,{\mathcal{A}}_2,\ep)$
such~that
$\big|\rho_2(\under{z})\big|\le\ep\big|\rho_1(\under{z})\big|$
for all $\under{z}\in 
\big(B_{\de}(0,\Bbb{C}^n) - Y_n\big)\cap 
\rho^{-1}(\Bbb{R}^+\cdot{\mathcal{K}})$.
\end{alprp}

The rest of this subsection is devoted to proving these propositions.
Note~that
\begin{equation}\label{monomials_e1}
\det D\rho\big|_{(z_1,\ldots,z_n)}
=(\det\rho)z_1^{A_1-1}\ldots z_n^{A_n-1},
\end{equation}
where $A_j = \sum_{i=1}^{i=n} a_{i,j}$.
Thus, $\Im\rho$ contains an open subset of $\Bbb{C}^n$ 
if (and only if) $\rho$ is nondegenerate in 
the sense of (2) of Definition~\ref{monomials_dfn1a}.
Since $\rho$ is a rational function in complex variables,
it follows that $\Im\rho$ is a dense open subset of~$\Bbb{C}^n$
if $\rho$ is nondegenerate.
Since $\rho$ is given by monomials, $\Im\rho$ is in fact all 
of~$\Bbb{C}^n - Y_n$. 
Thus
$$\rho\co\Bbb{C}^n - Y_n\lra\Bbb{C}^n - Y_n$$
is a local diffeomorphism.
By Lemma~\ref{covering_order1}, this map is in fact 
a covering projection of order~$|\det\rho|$.

\begin{remark}
{\rm The proof of Lemma~\ref{covering_order1} does not rely on 
Lemmas~\ref{monomials_lmm4}--\ref{monomials_lmm2} or
Corollary~\ref{monomials_crl3}.  We postpone its proof until the very
end of this subsection in order to focus on the main aspects of the
proof of Propositions~\ref{monomials_prp1} and~\ref{monomials_prp2}.}
\end{remark}

Denote by $\tilde{Y}_n$ the union of all
the $n$ codimension-one coordinate subspaces 
$\Bbb{R}^k \times \{0\} \times \Bbb{R}^{n-1-k}$ in~$\Bbb{R}^n$.
Let $\tilde{Y}_n^* = \tilde{Y}_n\cap S^{n-1}$.
We identify $\Bbb{C}^n - Y_n$ with 
$$\Bbb{R}^+ \times (S^{n-1} - \tilde{Y}_n^*) \times (S^1)^n$$
by the~map
\begin{eqnarray*}
\Bbb{C}^n - Y_n  & \lra &
\Bbb{R}^+ \times (S^{n-1} - \tilde{Y}_n^*) \times (S^1)^n,\\
\under{z} = (z_1,\ldots,z_n)&\longmapsto &
\bigg(|\under{z}|,\frac{(|z_1|,\ldots,|z_n|)}{|\under{z}|},
\Big(\frac{z_1}{|z_1|},\ldots,\frac{z_n}{|z_n|}\Big)\bigg).
\end{eqnarray*}
With respect to this decomposition,
$$\rho(r,\vph,\th)=\big(f(r,\vph),g(r,\vph),h(\th)\big),$$
where $\big(f(\under{r}),g(\under{r})\big) = \rho(\under{r}) \in
\Bbb{R}^n$ and 
$$h\big(e^{i\th_1},\ldots,e^{i\th_n}\big)
=\rho\big(e^{i\th_1},\ldots,e^{i\th_n}\big)\in(S^1)^n \subset \Bbb{C}^n.$$
Proposition~\ref{monomials_prp1} follows immediately from:

\begin{lmm}
\label{monomials_lmm4}
If $\rho$ is a neutral monomials map,
for every compact subset ${\mathcal{K}}$ of $S^{n-1} - \tilde{Y}_n^*$, 
there exists $\de\in \Bbb{R}^+$ such that for all $r \in (0,\de)$,
$\{g(r,\cdot)\}^{-1}({\mathcal{K}}) = \eset$.
\end{lmm}

\begin{proof}
We use the variables $s_1,\ldots,s_n$ to denote
the standard Euclidean coordinates on the target space $\Bbb{R}^n$,
as well as the corresponding component functions of $(f,g)$.
Let
$$\tilde{\mathcal{A}}=\left(
\begin{array}{ccc}
a_{2,1} - a_{1,1}&\ldots&a_{2,n} - a_{1,n}\\
\vdots&&\vdots\\
a_{n,1} - a_{1,1}&\ldots&a_{n,n} - a_{1,n}
\end{array}\right).$$
It is sufficient to show that, for some $i = 2,\ldots,n$,
$$\limsup_{r\to0}\Big\{\frac{s_i(\under{r})}{s_1(\under{r})} :
|\under{r}| <r \Big\}=0
\quad\hbox{or}\quad
\liminf_{r\to0}\Big\{\frac{s_i(\under{r})}{s_1(\under{r})} :
|\under{r}| <r \Big\}=\i.$$
This condition is equivalent to
\begin{align}
\limsup_{t\to-\i}\Bigl\{\pm&\sum_{j=1}^{n}(a_{i,j}{-}a_{1,j})t_j :
t_j \in (-\i,t)\Bigr\}=-\i\notag\\
\label{monomials_lmm4e3}
\Llra\quad&\left\{\begin{array}{ll}
& a_{i,j} - a_{1,j} \ge 0\ \mbox{for all}\ j = 1,\ldots,n \\
\mbox{or} & a_{i,j} - a_{1,j} \le 0\ \mbox{for all}\ j = 1,\ldots,n.
\end{array}\right.
\end{align}
Note that the two lines above are equivalent 
because ${\mathcal{A}}$ is assumed to be nondegenerate.
If \eqref{monomials_lmm4e3} is not satisfied by any $i = 2,\ldots,n$,
for every nonzero vector $\under{x} \in \Bbb{R}^{n-1}$
there exists a vector $\under{c} \in \Bbb{R}^n$
such~that
$$\under{x}^t\tilde{\mathcal{A}}\under{c}>0
\qquad\hbox{and}\qquad
c_1,\ldots,c_n \ge 0.$$
This means that the image of $\tilde{\mathcal{A}}^t$ contains
no nonzero vector with all components of the same sign.
Thus
\begin{equation}\label{monomials_lmm4e4}
\mbox{there exists}\ (x_1,\ldots,x_n)\ \mbox{in}\ \ker\tilde{\mathcal{A}} - \{0\}
\ \mbox{such that}\ x_1,\ldots,x_n > 0.
\end{equation}
Let $\tilde{\mathcal{A}}_{\hat{j}}$ be the matrix obtained from
$\tilde{\mathcal{A}}$ by removing the $j$th column.
Since $\rho$ is nondegenerate, $\det\tilde{\mathcal{A}}_{\hat{j}} \neq 0$
for some~$j$. Then, by Cramer's Rule,
\begin{equation}\label{monomials_lmm4e5}
x_{j'}=
\frac{(-1)^{j+j'-1}\det\tilde{\mathcal{A}}_{\hat{j}'}}
{\det\tilde{\mathcal{A}}_{\hat{j}}} \, (-x_j)
=\frac{(-1)^{j'}\det\tilde{\mathcal{A}}_{\hat{j}'}}
{(-1)^j\det\tilde{\mathcal{A}}_{\hat{j}}}\, x_j \ \mbox{for all}\ j'.
\end{equation}
Since $\rho$ is neutral, 
$$\frac{(-1)^{j'}\det\tilde{\mathcal{A}}_{\hat{j}'}}
{(-1)^j\det\tilde{\mathcal{A}}_{\hat{j}}}\le0$$
for some $j'$.
Thus \eqref{monomials_lmm4e5} contradicts \eqref{monomials_lmm4e4}.
\end{proof}

\begin{lmm}
\label{monomials_lmm1}
If $\rho$ is a positive or negative monomials map,
then for every $r \in \Bbb{R}^+$, the map
$$g(r,\cdot)\co S^{n-1} - \tilde{Y}_n^* \lra S^{n-1} - \tilde{Y}_n^*,
\ \mbox{given by}\ \th\longmapsto g(r,\th),$$
is a local diffeomorphism.
\end{lmm}

\begin{proof}
We assume that $n \ge 2$; otherwise, there is nothing to prove. 
Suppose $g(r,\cdot)$ is not a local diffeomorphism at
$\under{r} = (r_1,\ldots,r_n) \in \Bbb{R}^n - \tilde{Y}_n$.
Then, there exists $\under{c} \in \Bbb{R}^n - \{0\}$ such that
\begin{equation}
\label{monomials_lmm1e1}
\sum_{j=1}^n r_jc_j=0
\quad \hbox{and}\quad \sum_{j=1}^n \big(a_{i,j}-a_{1,j}\big)r_j^{-1}c_j=0
\ \mbox{for all}\ i = 2,\ldots,n.
\end{equation}
The first equation above is equivalent to the condition
$\under{c} \in T_{\under{r}}S^{n-1}$.
The second equation means that the ratio of the $i$th and the first
Euclidean components of the function $(f,g)$ does not change in 
the direction of $\under{c}$ at~$\under{r}$.
The $n$ conditions~\eqref{monomials_lmm1e1} are equivalent to
${\mathcal{A}}^{(1)}\big(\under{r}^2\big)\under{c}'=0\in\Bbb{R}^n$,
where:
\begin{gather*}
{\mathcal{A}}^{(1)}\big(\under{r}^2\big) \equiv 
{\mathcal{A}}^{(1)}\big(r_1^2,\ldots,r_n^2\big)=
\left(
\begin{array}{ccc}
r_1^2&\ldots& r_n^2\\
a_{2,1} - a_{1,1}&\ldots& a_{2,n} - a_{1,n}\\
\vdots&&\vdots\\
a_{n,1} - a_{1,1}&\ldots& a_{n,n} - a_{1,n}
\end{array}\right)
\end{gather*}
This equation has a nonzero solution only if 
$\det{\mathcal{A}}^{(1)}\big(\under{r}^2\big) = 0$.
However,
$$\det{\mathcal{A}}^{(1)}\big(\under{r}^2\big)
=\sum_{j=1}^n(-1)^{j-1}\big(\det\tilde{\mathcal{A}}_{\hat{j}}\big)r_j^2.$$
Since $\rho$ is positive or negative,
all the elements of the set 
$\big\{(-1)^{j-1}\det\tilde{\mathcal{A}}_{\hat{j}}\big\}$
have the same sign.
Thus, $\det{\mathcal{A}}^{(1)}\big(\under{r}^2\big)$ does not vanish on 
$\Bbb{R}^n - \{0\}$.
It follows that the differential of $g(r,\cdot)$ is an isomorphism 
everywhere on $S^{n-1} - \tilde{Y}_n^*$.
\end{proof}

Let $Y_n^* = Y_n\cap S^{2n-1} \subset \Bbb{C}^n$.
We identify $\Bbb{C}^n - Y_n$ with 
$\Bbb{R}^+ \times (S^{2n-1} - Y_n^*)$ and
denote by 
$$\big(\tilde{f},\tilde{g}\big)\co  \Bbb{R}^+ \times (S^{2n-1} - Y_n^*)
\lra \Bbb{R}^+ \times (S^{2n-1} - Y_n^*)$$
the pair of maps corresponding to $\rho$.

\begin{lmm}
\label{monomials_lmm2}
If $\rho$ is a positive (or negative) monomials map, then
for every $r \in \Bbb{R}^+$, the map
$$\tilde{g}(r,\cdot)\co S^{2n-1} - Y_n^*\lra S^{2n-1} - Y_n^*$$
is a local orientation-preserving (or orientation-reversing)
diffeomorphism, and is a covering projection of oriented order 
$|\det\rho|$ (or $-|\det\rho|$).
\end{lmm}

\begin{proof} To prove the first claim, we assume that $\rho$ is a
positive monomials map; the negative case is proved similarly.
By Lemma~\ref{monomials_lmm1} and 
the decomposition $\tilde{g} = (g,h)$, it follows that
$\tilde{g}(r,\cdot)$ is a local diffeomorphism.
Since $\rho$ is orientation-preserving everywhere,
$\tilde{g}(r,\cdot)$ is orientation-preserving at~$(r,\th)$ if
\begin{equation}\label{monomials_lmm2e1}
\Big\lan   \big\{D\rho\big|_{(r,\th)}\big\}^{-1}
\Big(\frac{\partial}{\partial r}\Big),
\Big(\frac{\partial}{\partial r}\Big)\Big\ran >0.
\end{equation}
Let $\under{c} \in \Bbb{R}^n$ be given by
${\mathcal{A}}\under{c} = \under{1} \in \Bbb{R}^n$. Then
\begin{equation}\label{monomials_lmm2e2}
\rho\big(t^{c_1}z_1,\ldots,t^{c_n}z_n\big)=t\rho\big(z_1,\ldots,z_n\big)
\quad \mbox{and}\quad
\frac{d}{dt}t\rho(z)\Big|_{t=1}\!\!=\frac{\partial}{\partial r}.
\end{equation}
Since $\rho$ is positive, $c_1,\ldots,c_n \in \Bbb{R}^+$ and thus
\begin{equation}\label{monomials_lmm2e3}
\frac{d}{dt}\big|\big(t^{c_1}z_1,\ldots,t^{c_n}z_n\big)\big|>0.
\end{equation}
The desired inequality~\eqref{monomials_lmm2e1} is immediate 
from~\eqref{monomials_lmm2e2} and~\eqref{monomials_lmm2e3}.
The remaining claim follows from Lemma~\ref{covering_order1} 
and~\eqref{monomials_lmm2e2}, since the curves
$$t\longmapsto\big(t^{c_1}z_1,\ldots,t^{c_n}z_n\big)
\quad\hbox{and}\quad
t\longmapsto\big(tw_1,\ldots,tw_n\big),\ \mbox{for}\ t \in (0,\i),$$
with $(z_1,\ldots,z_n)$ and $(w_1,\ldots,w_n) \in S^{2n-1}$,
foliate $\Bbb{C}^n - Y_n$ and intersect
each $(2n - 1)$--sphere $r = \const$ exactly once.
\end{proof}

\begin{remark}
{\rm The first claim of Proposition~\ref{monomials_prp2} now
follows from Lemma~\ref{covering_order1}, the first identity
in~\eqref{monomials_lmm2e2}, and the assumptions that all the exponents
$c_i$ have the same~sign.}
\end{remark}

\begin{crl}
\label{monomials_crl3}
Suppose $\rho$ is a positive (or negative) monomials map,
$\de^* \in \Bbb{R}^+$, $\de \in (0,\de^*)$,
and ${\mathcal{K}}$ is a precompact open subset of $\Bbb{C}^n - Y_n$.
Then, there exists $\ep \in \Bbb{R}^+$ such that
for all $t \in (0,\ep)$, the map
$$t^{-1}\be_{\de}\rho\co 
\big(B_{2\de^{*1/2}}(0,\Bbb{C}^n) - Y_n\big)\cap
\big\{t^{-1}\be_{\de}\rho\big\}^{-1}({\mathcal{K}})\lra {\mathcal{K}}$$
is a smooth covering projection of oriented order
\hbox{$|\det\rho|$ ($-|\det\rho|$)}.
\end{crl}

\begin{proof}
We prove this corollary in the case $\rho$ is negative.
If $\rho$ is positive, a stronger claim can be obtained by 
a similar and somewhat simpler argument.
With notation as above, $\be_{\de}\rho$ corresponds 
to the pair $(\be_{\de}\tilde{f},\tilde{g})$
with respect to the splitting of 
$\Bbb{C}^n - Y_n = \Bbb{R}^+ \times (S^{2n-1} - Y_n^*)$.
Thus
\begin{equation}\label{monomials_crl3e1}
\det D(\be_{\de}\rho)\big|_{(r,\th)}=
\be_{\de}(r)\det D\rho\big|_{(r,\th)}+
\be_{\de}'(r)\tilde{f}(r,\th)
\det\Big(\frac{\partial\tilde{g}}{\partial\th}\Big)\Big|_{(r,\th)}.
\end{equation}
Let $\tilde{\mathcal{K}} \subset S^{2n-1} - Y_n^*$ be the image of
${\mathcal{K}}$
under the projection map onto the second component.
Since the closure of ${\mathcal{K}}$ in $\Bbb{C}^n - Y_n$ is compact,
Lemma~\ref{monomials_lmm2} implies that the set
$$\tilde{U}\equiv\mskip-10mu\bigcup_{\de^{1/2}\le r\le
2\de^{*1/2}}\mskip-20mu
\big\{\tilde{g}(r,\cdot)\big\}^{-1}(\tilde{\mathcal{K}})$$
has compact closure in $S^{2n-1} - Y_n^*$.
Thus, there exists $C > 0$ such that
%\begin{gather}\label{monomials_crl3e2}
$$\det D\rho\big|_{(r,\th)}<C \quad\mbox{and}\quad
\tilde{f}(r,\th)\det\Big(\frac{\partial\tilde{g}}{\partial\th}\Big)
\Big|_{(r,\th)}<-C^{-1}$$
for all $(r,\th)\in\big(\de^{1/2},2\de^{*1/2}\big) \times \tilde{U}$.
%\end{gather}
The second bound is obtained by using Lemma~\ref{monomials_lmm2}.
Choose $\eta > 0$ such that
\begin{equation}\label{monomials_crl3e2b}
\be_{\de}'(r) > 2C^2\be_{\de}(r)
\qquad \mbox{for all}\ r \in \big(\de^{1/2},\de^{1/2} + \eta\big).
\end{equation}
Note that combining \eqref{monomials_crl3e1}--\eqref{monomials_crl3e2b},
we obtain
\begin{equation}\label{monomials_crl3e2c}
\det D(\be_{\de}\rho)\big|_{(r,\th)}<0
\end{equation}
for all $(r,\th) \in \big(\de^{1/2},\de^{1/2} + \eta\big) \times
\tilde{U}$.  Let $\ep > 0$ be such that
\begin{equation}\label{monomials_crl3e3}
\ep\cdot\max\big\{|w| : w \in {\mathcal{K}}\big\}
<{\textstyle\frac{1}{2}}\min\bigl\{\be_{\de}(|\under{z}|)\big|\rho(\under{z})\big| :
\under{z} \in \big(\de^{1/2} + {\textstyle\frac{1}{2}}\eta,2\de^{*1/2}\big)
 \times \tilde{U}\bigr\}.
\end{equation}
We claim that $\ep$ satisfies the required properties.
Suppose that
$t<\ep$, that
$\under{z}=(r,\th) \in B_{2\de^{*1/2}}(0,\Bbb{C}^n)$, and that
$\big\{t^{-1}\be_{\de}\rho\big\}(z) \in {\mathcal{K}}.$
Then
$$r \in \big(\de^{1/2},2\de^{*1/2}\big)\Lra
\th \in \big\{\tilde{g}(r,\cdot)\big\}^{-1}({\mathcal{K}})
\subset\tilde{U}.$$
Assumption \eqref{monomials_crl3e3} on $\ep$ then implies that
$(r,\th) \in \big(\de^{1/2},\de^{1/2} + \frac{1}{2}\eta\big)
 \times \tilde{U}$.
{}From \eqref{monomials_crl3e2c}, we conclude that
the determinant of the derivative of $t^{-1}\be_{\de}\rho$
at $\under{z}$ is negative and the map
$$t^{-1}\be_{\de}\rho\co
\big(B_{2\de^{*1/2}}(0,\Bbb{C}^n) - Y_n\big)\cap
\big\{t^{-1}\be_{\de}\rho\big\}^{-1}({\mathcal{K}})\lra {\mathcal{K}}$$
is a local orientation-reversing diffeomorphism.
It remains to see that for each point $ (s,\vt) \in {\mathcal{K}}$,
$$\big|\big\{(r,\th)
\in \big(\de^{1/2},\de^{1/2} + {\textstyle\frac{1}{2}}\eta\big) \times
\tilde{U} : \be_{\de}(r)\tilde{f}(r,\th) = ts,\tilde{g}(r,\th) = \vt\big\}\big|
=|\det\rho|.$$
By Lemma~\ref{monomials_lmm2}, there are smooth paths
$$\th_i\co \big[\de^{1/2},\de^{1/2} + {\textstyle\frac{1}{2}}\eta\big]
\lra \tilde{U}\ \mbox{for}\ i = 1,\ldots,|\det\rho|$$
such that
$$\big\{\tilde{g}(r,\cdot)\big\}^{-1}(\vt)=
\big\{\th_i(r) :i = 1,\ldots,|\det\rho|\big\}
\quad \mbox{and}\quad \th_i(r) \neq \th_j(r)$$
for all $r \in \big[\de^{1/2},\de^{1/2} + \frac{1}{2}\eta\big]$ and $i \neq j$.
By \eqref{monomials_crl3e2c} and \eqref{monomials_crl3e3},
\begin{gather*}
\{\be_{\de}\tilde{f}\}\big(\de^{1/2},\th_i(\de^{1/2})\big) = 0,\quad
\{\be_{\de}\tilde{f}\}\big(\de^{1/2} + {\textstyle\frac{1}{2}}\eta,
\th_i(\de^{1/2} + {\textstyle\frac{1}{2}}\eta)\big) > ts,\\
\hbox{and}\quad\frac{d}{dr}\{\be_{\de}\tilde{f}\}\big(r,\th_i(r)\big)>0~~
\ \mbox{for all}\ r \in \big(\de^{1/2},\de^{1/2} +
{\textstyle\frac{1}{2}}\eta\big).
\end{gather*}
Thus, for each $i = 1,\ldots,|\det\rho|$, there exists a unique number
$$r_i\in\big(\de^{1/2},\de^{1/2} + {\textstyle\frac{1}{2}}\eta\big)
\ \mbox{such that}\ \{\be_{\de}\tilde{f}\}\big(r_i,\th_i(r_i)\big)=ts,$$
as required.
\end{proof}

Corollary~\ref{monomials_crl3} essentially concludes the proof of 
part (2) of Proposition~\ref{monomials_prp2}.
The claimed inclusion is achieved 
if, in the proof of Corollary~\ref{monomials_crl3}, $\eta$ is chosen so that
$\de^{1/2} + \eta < \de_+$.

We next prove part (3) of Proposition~\ref{monomials_prp2}.
Suppose ${\mathcal{A}}$ is a positive monomials map and
$\tilde{\mathcal{K}} = \rho^{-1}({\overK})$.
Since $\tilde{\mathcal{K}}$ is a compact subset of $\Bbb{C}^n - Y_n$,
there exists $r > 0$ such that
$B_r(0,\Bbb{C}^n)\cap\tilde{\mathcal{K}} = \eset$.
On the other hand, by Lemma~\ref{covering_order1} and 
the first identity in~\eqref{monomials_lmm2e2}, 
if $t \in \Bbb{R}^+$ and $\rho(\under{z}) \in t\cdot{\mathcal{K}}$,
\begin{gather*}
(z_1,\ldots,z_n)=(t^{c_1}w_1,\ldots,t^{c_n}w_n)
\quad\hbox{for some}\quad (w_1,\ldots,w_n)
 \in \rho^{-1}(K)\\ 
\Lra\qquad |\under{z}|\ge t^{\min c_i}r.
\end{gather*}
Thus, if $|\under{z}| \le \de$,
$$\Big|\frac{\rho_2(\under{z})}{\rho_1(\under{z})}\Big|=
t^{{\mathcal{A}}_2\cdot\under{c}-{\mathcal{A}}_1\cdot\under{c}}~
\Big|\frac{\rho_2(\under{w})}{\rho_1(\under{w})}\Big|
\le C\Big(\frac{\de}{r}
\Big)^{({\mathcal{A}}_2\cdot\under{c}-{\mathcal{A}}_1\cdot\under{c})/\min c_i}.$$
Since ${\mathcal{A}}_2\cdot\under{c} > {\mathcal{A}}_1\cdot\under{c}$
and $\min c_i > 0$, the right-hand side above tends to zero 
with~$\de$. If $\rho$ is negative, the proof is similar.

\begin{lmm}
\label{covering_order1}
If $\rho\co \Bbb{C}^n - Y_n \lra \Bbb{C}^n - Y_n$ 
is a nondegenerate monomials map,
$\rho$ is a covering projection of order~$|\det\rho|$.
\end{lmm}

\begin{proof}
By~\eqref{monomials_e1}, we only need to compute the order of the cover.
We can view $\rho$ as a rational map from
$(\Bbb{P}^1)^N$ to~$(\Bbb{P}^1)^N$.
In turn, this rational map induces a holomorphic map
$\tilde{\rho}\co \tilde{M} \lra (\Bbb{P}^1)^N$,
where $\tilde{M}$ is a compact complex manifold obtained
from $(\Bbb{P}^1)^N$ by a sequence of blowups along 
submanifolds disjoint from $\Bbb{C}^n - Y_n$.
Then
\begin{equation}\label{covering_order1e1}\begin{split}
\hbox{ord}~ \rho
=\big\lan\tilde{\rho}^*\big(\om_1\w\ldots\w\om_n\big),[\tilde{M}]\big\ran
&=\int_{\tilde{M}}\tilde{\rho}^*\big(\om_1\w\ldots\w\om_n\big)\\
&=\int_{C^n-Y_n}\rho^*\big(\om_1\w\ldots\w\om_n\big),
\end{split}\end{equation}
where $\om_i$ is the Fubini--Study symplectic form on the $i$th 
$\Bbb{P}^1$--factor of the target space.
Since
$\om_i = \frac{\I}{2\pi}\frac{dz_i\w d\bar{z}_i}{(1+|z_i|^2)^2}$,
see~\cite[p31]{GH}, by \eqref{monomials_e1},
\begin{equation}\label{covering_order1e2}\begin{split}
&\rho^*\big(\om_1\w\ldots\w\om_n\big)\\
&\qquad = |\det\rho|^2\Big(\frac{\I}{2\pi}\Big)^n ~
\frac{r_1^{2A_1-2}\ldots r_1^{2A_n-2}
\, dz_1\w d\bar{z}_1\w\ldots\w dz_n\w d\bar{z}_n}{
(1 + r_1^{2a_{1,1}}\ldots r_n^{2a_{1,n}})^2\ldots
(1 + r_1^{2a_{n,1}}\ldots r_n^{2a_{n,n}})^2}.
\end{split}\end{equation}
Combining \eqref{covering_order1e1} and~\eqref{covering_order1e2} and 
switching to polar coordinates, we obtain
\begin{equation}\label{covering_order1e3}
\ord \rho
=2^n|\det\rho|^2  \int_0^{\i}\!\!\!\!     \ldots   \int_0^{\i}
\frac{r_1^{2A_1-1} \ldots r_1^{2A_n-1} \, dr_1\ldots dr_n}{
(1 + r_1^{2a_{1,1}}  \ldots r_n^{2a_{1,n}})^2\ldots
(1 + r_1^{2a_{n,1}}  \ldots r_n^{2a_{n,n}})^2}.
\end{equation}
The change of variables,
$$(r_1,\ldots,r_n)\lra
\big(r_1^{2a_{1,1}}\ldots r_n^{2a_{1,n}},\ldots,
r_1^{2a_{n,1}}\ldots r_n^{2a_{n,n}}\big),$$
reduces~\eqref{covering_order1e3} to
\begin{equation*}\begin{split}
\ord\rho
&= |\det\rho|\int_0^{\i}\!\!\!\!\ldots\int_0^{\i}
\frac{dr_1\ldots dr_n}{(1 + r_1)^2\ldots(1 + r_n)^2}\\
&=|\det\rho|\bigg(\int_0^{\i}\frac{dr}{(1 + r)^2}\bigg)^n=|\det\rho|,
\end{split}\end{equation*}
as claimed.
\end{proof}

\subsection{Topological setup}
\label{top_dfn_subs}

In this subsection, we give formal definitions of the topological objects
to which the computational method described in the next subsection applies.

We start by extending the concept of monomials maps to vector bundles.
All vector bundles we encounter will be assumed to be complex and normed.
Vector bundles over smooth manifolds will in addition be smooth. 
Given a vector bundle $F \lra X$ and any map $\de\co X \lra \Bbb{R}$, put
$$F_{\de}=\big\{(b;v) \in F : |v|_b<\de(b)\big\}.$$ 
If $F = \bigoplus_{i\in I}F_i$ is the direct sum of nontrivial subbundles
and $I_0 \subset I$, let
$$Y(F;I_0)=\bigcup_{i\in I_0}\Bigl(\bigoplus_{j\in I-\{i\}}  F_j\Bigr)
\subset F.$$

\begin{dfn}
\label{monomials_dfn2}
{\rm Suppose $I_0$, $I$, and $J$ are finite sets, and ${\mathcal{A}} =
(a_{ij})$ is an integer-valued function on $(I_0\sqcup I) \times (I_0\sqcup J)$
such that for all $j \in I_0$,
$a_{ij} = 0$ if $i \neq j$ and $a_{ij} = 1$ if $i = j$.

(1)\qua Suppose $F_j \lra {\mathcal{M}}$ is a vector bundle for each $j \in I_0$
and a line bundle for each $j \in J$,
$$F = \bigoplus_{j\in I_0\sqcup J}   F_j,\quad
\ti{F}_i=\bigotimes_{j\in I_0\sqcup J}   F_j^{\otimes a_{ij}}
\ \mbox{for all}\  i \in I_0 \sqcup I, \quad\mbox{and}\quad
\ti{F} = \bigoplus_{i\in I_0\sqcup I}   \ti{F}_i.$$
A function $\rho\co F - Y(F;J) \lra \tilde{F}$
is a \emph{degree-${\mathcal{A}}$ monomials map on $F$} if
$$\pi_i\rho\big((\ups_j)_{j\in I_0\sqcup J}\big)=
 \bigotimes_{j\in I_0\sqcup J}   \ups_j^{\otimes a_{i,j}}$$
for all $(\ups_j)_{j\in I_0\sqcup J} \in F - Y(F;J)$
and $i \in I_0\sqcup I$,
where $\pi_i\co \ti{F} \lra \ti{F}_i$ is the projection map.

(2)\qua Suppose $\rho$ is as in (1),
$E_i \lra {\mathcal{M}}$ is a vector bundle for each $i \in I$, 
$$E=\bigoplus_{i\in I}E_i,  \qquad\hbox{and}\qquad
\tilde{E}=\bigoplus_{i\in I}E_i \otimes \tilde{F}_i.$$
A function
$$\ti{\rho}\co E \oplus F-Y(E \oplus F;J)\lra\tilde{E}$$
is a \emph{degree-${\mathcal{A}}$ monomials map on $E \oplus F$} if 
$$\pi_i\ti{\rho}\big((w_i)_{i\in I},(\ups)_{j\in I_0\sqcup J}\big)
=w_i\otimes\pi_i\rho(\ups)$$
for all $\big((w_i)_{i\in I},(\ups)_{j\in I_0\sqcup J}\big)
	            \in E \oplus F-Y(E \oplus F;J)$ and $i \in I$.}
\end{dfn}

A monomials map between vector bundles in the sense of
Definition~\ref{monomials_dfn2} part~(1) can be viewed as a pair of
bundle maps
$$\rho_{I,J}\co  F_J \equiv \bigoplus_{j\in J} F_j
\lra \ti{F}_I \equiv \bigoplus_{i\in I}\ti{F}_i
\qquad\hbox{and}\qquad
\rho_{I_0}\co  F_J\oplus \bigoplus_{i\in I_0} F_i
\lra \bigoplus_{i\in I_0} \ti{F}_i.$$
The vector bundles $F_J$ and $\ti{F}_I$ are sums of line bundles and
the restriction of the bundle map $\rho_{I,J}$ to each fiber is 
a monomials map in the sense of Definition~\ref{monomials_dfn1a}.
The degree of this map is~${\mathcal{A}}|I \times J$.
The bundle map~$\rho_{I_0}$ has $I_0$ components, indexed by $i \in I$:
$$\pi_i\rho\co  F_i\oplus F_J \lra \ti{F}_i\equiv
F_i\otimes\bigotimes_{j\in J}F_j^{\otimes a_{ij}}, \qquad
\mbox{given by}\ \big(\ups_i,(\ups_j)_{j\in J}\big)\longmapsto
\ups_i\otimes\bigotimes_{j\in J}\ups_j.$$
A monomials map in the sense of Definition~\ref{monomials_dfn2} part~(2)
is equivalent to a monomials map in the sense of
Definition~\ref{monomials_dfn2} part~(1) with ${\mathcal{A}}$ being a
function on $(I \sqcup \eset) \times (I \sqcup J)$.

If $\rho$ and $\ti{\rho}$ are as in Definition~\ref{monomials_dfn2}
part~(2), and $i \in I$,  we put 
$$\ti{F}_i(\ti{\rho})=\ti{F}_i \qquad\hbox{and}\qquad
\ti{\rho}_i(\ups)=\pi_i\rho(\ups)\in\ti{F}_i(\ti{\rho})
\ \mbox{for all}\ \ups \in F - Y(F;J).$$
If $\rho$ is as in Definition~\ref{monomials_dfn2} part (1),  we call $\rho$ 
\emph{nondegenerate} if the restriction of $\rho$ to a fiber of~$F$ is
\emph{nondegenerate} in the sense of Definition~\ref{monomials_dfn1a},
for some choice of identifications of the sets $I_0\sqcup I$ and $I_0\sqcup J$
with the sets of integers $1,\ldots,|I_0\sqcup I|$ and 
$1,\ldots,|I_0\sqcup J|$.
If $\rho$ is nondegenerate and $|I| = |J|$, we call  $\rho$
\textit{positive} (or \textit{negative}, or \textit{neutral})
if the restriction of $\rho$ to a fiber of~$F$ is
\textit{positive} (or \textit{negative}, or \textit{neutral}).
Similarly, suppose $\rho$ is a positive or negative monomials map,
$I_1$ and $I_2$ are one-element sets, and
${\mathcal{A}}_1$ and ${\mathcal{A}}_2$ are integer-valued functions on
$$I_1 \times (I_0\sqcup J) \qquad\hbox{and}\qquad 
I_2\times(I_0\sqcup J),$$ 
respectively. 
If $\rho_1$ and $\rho_2$ are monomials maps of degrees 
${\mathcal{A}}_1$ and ${\mathcal{A}}_2$, respectively,
we write $\rho_1 >_{\rho} \rho_2$ if
this relation holds for the restriction to a fiber;
see Definition~\ref{monomials_dfn1b}.
Due to the remark following this definition,
the notions of nondegenerate, positive, negative, and neutral
depend only on~${\mathcal{A}}$;
the partial ordering relation depends only on ${\mathcal{A}}$, ${\mathcal{A}}_1$,
and~${\mathcal{A}}_2$.

If $F$ is any (normed) vector bundle and $\de \in \Bbb{R}^+$,
we define the function $\be_{\de}$ on $F$ by:
$$\be_{\de}\co  F\lra[0,1]\subset\Bbb{R}, \qquad
\be_{\de}(\ups)=\be_{\de}\big(|\ups|\big).$$
If  $\rho$ is a positive or negative monomials map between vector bundle 
and $t \in \Bbb{R}^+$, 
we denote by $\deg\rho$ the oriented degree of the map
$t^{-1}\be_{\de}\rho$ given by Proposition~\ref{monomials_prp2}.

The next two definitions characterize the topological spaces with which we
work.  Ms-orbifolds, as described by Definition~\ref{space1_dfn}, include
spaces of stable maps.  Examples of pseudovarieties,  as described by
Definition~\ref{space2_dfn}, that we encounter are subspaces of spaces
of stable maps that consist of elements corresponding to curves with
specified singularities.

\begin{dfn}
\label{space1_dfn}
{\rm A compact topological orbifold ${\overM} = {\mathcal{M}}_n\sqcup
\smash{\bigsqcup\limits_{k=0}\limits^{n-1}}  {\mathcal{M}}_k$
is a \emph{mostly smooth}, or \emph{ms}-orbifold of dimension~$n$ if

(1)\qua ${\mathcal{M}} \equiv {\mathcal{M}}_n$
is an open subset of ${\overM}$, and ${\overM}_k -
{\mathcal{M}}_k\subset\smash{\bigcup\limits_{j<k}}{\mathcal{M}}_j$
for all $k = 0,\ldots,n$;

(2)\qua ${\mathcal{M}}_k$ is a smooth oriented orbifold of dimension $2k$
for all $k = 0,\ldots,n$;

(3)\qua for each $k = 0,\ldots,n - 1$, there exist a smooth
complex vector orbi-bundle ${\mathcal{F}}_k \lra {\mathcal{M}}_k$
and an identification $\phi_k\co U_k \lra V_k$ of neighborhoods of
${\mathcal{M}}_k$ in ${\mathcal{F}}_k$ and in ${\overM}$ such that
\hbox{$\phi_k\co \phi_k^{-1}({\mathcal{M}}) \lra V_k\cap{\mathcal{M}}$}
is an orientation-preserving diffeomorphism.}
\end{dfn}

There is no firm consensus about the correct definition of the orbifold
category.  For our purposes, we put the following, rather strong, 
requirements on the objects involved in Definition~\ref{space1_dfn}.
Each smooth orbifold ${\mathcal{M}}_k$ of Definition~\ref{space1_dfn} 
is the quotient of a smooth manifold $\tilde{\mathcal{M}}_k$
by a smooth action of a compact Lie group~$G_k$.
All points of $\tilde{\mathcal{M}}_k$ have finite stabilizers,
and the set of points with nontrivial stabilizers 
has codimension at least two in~$\tilde{\mathcal{M}}_k$.
In other words, this set is
a finite union of smooth manifolds of dimension at most $2k - 2 + \dim G_k$.
In addition, there exists a vector-bundle splitting 
$$T\tilde{\mathcal{M}}_k=T^v\tilde{\mathcal{M}}_k\oplus
  T^h\tilde{\mathcal{M}}_k,$$
where $T^v\tilde{\mathcal{M}}_k$ is the vertical tangent bundle and
$T^h\tilde{\mathcal{M}}_k$ is a \textit{complex} vector bundle
on which $G_k$ acts by complex-linear automorphisms.  We call
${\mathcal{F}}_k \lra {\mathcal{M}}_k$ a smooth complex vector orbi-bundle
if there exists a smooth complex vector bundle $\tilde{\mathcal{F}}_k
\lra \tilde{\mathcal{M}}_k$ on which $G_k$ acts smoothly.

By a compact topological orbifold ${\overM}$, 
we mean the quotient of a compact Hausdorff topological space
$$\tilde{\mathcal{M}}=
\tilde{\mathcal{M}}_n\sqcup\bigsqcup_{i=0}^{n-1}\tilde{\mathcal{M}}_k'
$$
by a continuous action of a compact Lie group~$G$.
For the purposes of Definition~\ref{space1_dfn},
$\tilde{\mathcal{M}}_k'$ denotes the preimage of ${\mathcal{M}}_k$ under 
the quotient projection map $\tilde{\mathcal{M}} \lra {\overM}$,
$G = G_n$, and the restriction of the continuous $G$--action to
$\tilde{\mathcal{M}}_n$ agrees with the smooth $G_n$--action of
the previous paragraph.
Condition~(3) of Definition~\ref{space1_dfn} means that there exist
\begin{enumerate}
\item[(3a)] a splitting $G_k=G \times G_k'$;
\item[(3b)] a $G_k$--invariant neighborhood $\tilde{U}_k$ 
of $\tilde{\mathcal{M}}_k$ in $\tilde{F}_k$;
\item[(3c)] a $G$--invariant neighborhood $\tilde{V}_k$
of $\tilde{\mathcal{M}}_k'$ in $\tilde{\mathcal{M}}$;
\item[(3d)] a $G$--equivariant topological $G_k'$--fibration 
$\tilde{\phi}_k\co \tilde{U}_k \lra \tilde{V}_k$ such that
\begin{enumerate}
\item[(3d-i)] $\tilde{\phi}_k(\tilde{\mathcal{M}}_k) = \tilde{\mathcal{M}}_k'$
and $\tilde{\phi}_k\co \tilde{\phi}_k^{-1}(\tilde{\mathcal{M}}_n) \lra 
\tilde{V}_k\cap\tilde{\mathcal{M}}_n$ is smooth;
\item[(3d-ii)] for each $x \in
\tilde{\phi}_k^{-1}(\tilde{\mathcal{M}}_n)$, the map
${\pi_{\tilde{\phi}_k(x)}^h} \circ  d\tilde{\phi}_k|_x\co 
T_x^h\tilde{\mathcal{M}}_k \lra
{T_{\tilde{\phi}_k(x)}^h}\tilde{\mathcal{M}}$
is an orientation-preserving isomorphism.
\end{enumerate}
\end{enumerate}

Throughout the rest of the paper by a vector bundle over a smooth orbifold 
we will mean a smooth complex normed vector orbi-bundle.
With notation as above, this means that $\tilde{\mathcal{F}}_k$ carries a
$G_k$--invariant Hermitian inner-product.
Similarly, $V \lra {\overM}$ is a \textit{vector bundle}~if

(1)\qua $V$ is the quotient of a topological Hermitian vector  
bundle $\tilde{V} \lra \tilde{\mathcal{M}}$ by an action of~$G$;

(2)\qua $V|{\mathcal{M}}_k$ is a vector bundle for all $i = 0,\ldots,n$.

In (1), the action of $G$ preserves the Hermitian structure.

If both
${\overM} = {\mathcal{M}}_n\sqcup
\bigsqcup_{k=0}^{n-1}  {\mathcal{M}}_k$
and  ${\overM}' = {\mathcal{M}}_{n'}\sqcup
\bigsqcup_{k=0}^{n'-1}  {\mathcal{M}}_k'$
are ms-orbifolds, a  continuous map 
$\pi\co {\overM} \lra {\overM}'$ will called an 
\textit{ms-map} if
for each $k = 0,\ldots,n$, there exists $k' = 0,\ldots,n'$
such that $\pi\co {\mathcal{M}}_k \lra {\mathcal{M}}_{k'}'$ is a smooth map.

\begin{dfn}
\label{space2_dfn}
{\rm Let ${\overM}$ be an ms-orbifold as in Definition~\ref{space1_dfn}.

(1)\qua A smooth $2m$--dimensional oriented suborbifold ${\mathcal{S}}$
of ${\mathcal{M}}$ is an \emph{$m$--pseu\-do\-cy\-cle} in
${\overM}$ if
${\overS} - {\mathcal{S}}$ is contained in
$\bigsqcup_{k=0}^{n-1}  {\mathcal{M}}_k$ and
${\overS}\cap{\mathcal{M}}_k$ is contained in
a finite union of smooth suborbifolds
of ${\mathcal{M}}_k$ of dimension at \hbox{most $2m - 2$}.

(2)\qua A pseudocycle ${\mathcal{S}}$ is a \emph{pseudovariety} if
${\mathcal{S}}$ is a smooth submanifold of~${\mathcal{M}}$.}
\end{dfn}

This definition of pseudocycle is a variation on that of
\cite[Chapter~7]{MS} and \cite[Section~1]{RT}.
For fairly straightforward topological reasons, every pseudocycle
of~\cite{MS} and~\cite{RT} determines an integral homology class.
For nearly the same reasons, every pseudocycle ${\mathcal{S}}$
of Definition~\ref{space2_dfn} determines an element of
$H_{2m}({\overM};\Bbb{Q})$.  If $\al \in H^{2m}({\overM};\Bbb{Q})$,
we denote the evaluation  of~$\al$ on this homology class by
$\lr{\al,{\overS}}$.  We write $\partial{\overS}$ for the boundary
of~${\mathcal{S}}$, ie the set ${\overS} - {\mathcal{S}}$.

If ${\mathcal{S}}$ is a pseudovariety in ${\overM}$, by a vector bundle 
$V \lra {\overS}$ we will mean a topological Hermitian vector bundle over 
a neighborhood $U_V$ of ${\overS}$ in ${\overM}$ such that for each 
$k = 0,\ldots,n$ the restriction of $V$ to $U_V\cap{\mathcal{M}}_k$ is a smooth Hermitian bundle.
Similarly, we denote by $\Ga({\overS};V)$
the space of continuous sections of $V$ over $U_V$ that restrict
to smooth sections on~$U_V\cap{\mathcal{M}}_k$ for all  $k = 0,\ldots,n$.
For the sake of simplicity, 
we restrict the presentation of our main topological tools 
to vector bundles over pseudovarieties as 
this is sufficient for the purposes of counting rational curves
in projective spaces.

Throughout the paper, we assume that every smooth oriented manifold ${\mathcal{Z}}$
comes with a \textit{system of trivializations}, ie a smooth~map
$$\vt_{\mathcal{Z}}\co T{\mathcal{Z}}_{\de}\lra{\mathcal{Z}},$$
where $\de \in C({\mathcal{Z}};\Bbb{R}^+)$, such that $\vt_{\mathcal{Z}}|T_b{\mathcal{Z}}_{\de(b)}$ is 
an orientation-preserving diffeomorphism  onto
an open neighborhood of $b$ in ${\mathcal{Z}}$ that sends $(b;0)$ to~$b$.
If $V \lra {\mathcal{Z}}$ is a vector bundle, we assume
that such a map $\vt_{\mathcal{Z}}$ comes with a choice of a lift of 
$\vt_{\mathcal{Z}}$ to a bundle identification
$$\vt_{{\mathcal{Z}};V}\co   
\pi_{T{\mathcal{Z}}}^*V\big|_{T{\mathcal{Z}}_{\de}} \lra V,$$
ie an isomorphism of smooth Hermitian vector bundles that
restricts to the identity over~${\mathcal{Z}} \subset T{\mathcal{Z}}$.
If $V$ is given as a direct sum of proper subbundles~$V_i$,
$\vt_{{\mathcal{Z}};V}$ will be assumed to be induced by
the identifications $\vt_{{\mathcal{Z}};V_i}$.
If ${\mathcal{M}}$ is a smooth oriented manifold,
$E \lra {\mathcal{M}}$ is a vector bundle,
and $\ka \in \Ga({\mathcal{M}};E)$ is a section transversal to the zero~set,
we assume that the smooth oriented submanifold 
${\mathcal{Z}} = \ka^{-1}(0)$ comes with a \textit{normal-neighborhood model},
ie a smooth map
$$\vt_{\ka}\co E_{\de}\lra{\mathcal{M}},$$
for some $\de \in C({\mathcal{Z}};\Bbb{R}^+)$, which is an orientation-preserving 
diffeomorphism onto an open neighborhood of ${\mathcal{Z}}$ in ${\mathcal{M}}$ such that
$\vt_{\ka}|{\mathcal{Z}}$ is the identity~map.
If $V \lra {\mathcal{M}}$ is a vector bundle, $\vt_{\ka}$
comes with a choice of a bundle identification
$$\vt_{\ka;V}\co  \pi_E^*V\big|_{E_{\de}}\lra V\big|_{\Im\vt_{\ka}}$$
which respects vector-bundle splittings as above.
Furthermore,
$$\ka\big(\vt_{\ka}(b;w)\big)=\vt_{\ka;E}(b;w)$$
for all $(b;w) \in E_{\de}$.
We will often only imply these identifications 
in equations involving vector-bundle sections.

\begin{dfn}
\label{normal_model_dfn}
{\rm Suppose ${\overM}$ is an ms-manifold as in Definition~\ref{space1_dfn},
${\mathcal{S}}$ is a pseudovariety in ${\overM}$ as in part (2) of
Definition~\ref{space2_dfn}, and
${\mathcal{Z}}$ is a smooth submanifold of ${\mathcal{M}}_k$ 
for some \hbox{$k = 0,\ldots,n$}.
\begin{enumerate}
\item[(1)] A \emph{regularization of ${\mathcal{Z}}$
in ${\overM}$} is a tuple
$\big(U_{\mathcal{Z}},E,\ka,({\mathcal{F}}_{k;j})_{j\in
J({\mathcal{Z}})}\big)$,
where
\begin{enumerate}
\item[(1a)] $U_{\mathcal{Z}}$ is a neighborhood of ${\mathcal{Z}}$
in ${\mathcal{M}}_k$, $E \lra U_{\mathcal{Z}}$ is a vector bundle,
and $\ka \in \Ga(U_{\mathcal{Z}};V)$ is a section transversal to the
zero set such that ${\mathcal{Z}} = \ka^{-1}(0)$;
\item[(1b)] ${\mathcal{F}}_{k;j} \lra {\mathcal{U}}_{\mathcal{Z}}$
is a non-zero subbundle of ${\mathcal{F}}_k|U_{\mathcal{Z}}$ for each
$j \in J({\mathcal{Z}})$
such that 
%\begin{gather*}
${\mathcal{F}}_k|U_{\mathcal{Z}}=   
\bigoplus_{j\in J({\mathcal{Z}})}    {\mathcal{F}}_{k;j}$,
$\phi_k\big( U_k\cap Y({\mathcal{F}}_k|U_{\mathcal{Z}};J({\mathcal{Z}}))\big)
 \subset \partial{\overM}$, and
$\phi_k\big((U_k|U_{\mathcal{Z}}) -
 Y({\mathcal{F}}_k|U_{\mathcal{Z}};J({\mathcal{Z}}))\big)
 \subset {\mathcal{M}}$.
%\end{gather*}
\end{enumerate}
\item[(2)] A \emph{model for ${\mathcal{Z}}$ in ${\mathcal{S}}$} is a tuple 
$\big(U_{\mathcal{Z}},E,\ka,({\mathcal{F}}_{k;j})_{j\in J({\mathcal{Z}})};
{\mathcal{O}}_{\mathcal{Z}},\psi_{\mathcal{Z}}\big)$, where
\begin{enumerate}
\item[(2a)] $(U_{\mathcal{Z}},E,\ka,({\mathcal{F}}_{k;j})_{j\in
  J({\mathcal{Z}})})$ is a regularization of ${\mathcal{Z}}$ in ${\overM}$;
\item[(2b)] ${\mathcal{O}}_{\mathcal{Z}} \lra U_{\mathcal{Z}}$
is a vector bundle of rank $\frac{1}{2}\big(\dim{\mathcal{M}} -
\dim{\mathcal{S}}\big)$;
\item[(2c)] $\psi_{\mathcal{Z}}\co  {\mathcal{F}}_k \lra
\mathcal{O_{\mathcal{Z}}}$ is a bundle map such~that
\begin{enumerate}
\item[(2c-i)] $\psi_{\mathcal{Z}}$ is smooth outside of $Y({\mathcal{F}}_k;J({\mathcal{Z}}))$;
\item[(2c-ii)] if $\ups \in U_k|U_{\mathcal{Z}}$,
$\phi_k(\ups) \in {\mathcal{S}}$ if and only if $\psi_{\mathcal{Z}}(\ups)
= 0$.
\end{enumerate}
\end{enumerate}
\end{enumerate}}
\end{dfn}

Suppose $V \lra {\overS}$ is a vector bundle  
and $(U_{\mathcal{Z}},E,\ka,({\mathcal{F}}_{k;j})_{j\in J({\mathcal{Z}})})$
is a regularization of ${\mathcal{Z}}$ in ${\overM}$.
In such a case, we assume that the tuple 
$$(U_{\mathcal{Z}},E,\ka,({\mathcal{F}}_{k;j})_{j\in J({\mathcal{Z}})})$$
implicitly encodes a Hermitian vector-bundle isomorphism
$$\vt_{U_{\mathcal{Z}};V}\co 
 \pi_{{\mathcal{F}}_k}^*V\big|_{\phi_k^{-1}(U_V)|(U_V\cap U_{\mathcal{Z}})}
\lra V\big|_{U_V\cap\phi_k(U_k|(U_V\cap U_{\mathcal{Z}}))}$$
that covers the map $\phi_k$ and 
restricts to the identity over~$U_V\cap U_{\mathcal{Z}}$.
This isomorphism is to be smooth over the complement of
$Y({\mathcal{F}}_k|U_{\mathcal{Z}};J({\mathcal{Z}}))$.
Along with the map $\vt_{k;V}$, we then obtain an identification
$$\tilde{\vt}_{{\mathcal{Z}};V}\co 
\pi_{E\oplus{\mathcal{F}}_k}^*V\big|_{
(E_{\de} \times {\mathcal{F}}_{k;\de})|(U_V\cap{\mathcal{Z}})}
\lra V\big|_{U_V\cap\phi_k\vt_{\ka;{\mathcal{F}}_k}
 ((E_{\de} \times {\mathcal{F}}_{k;\de})|(U_V\cap{\mathcal{Z}}))},$$
covering the map $\phi_k \circ \vt_{\ka;{\mathcal{F}}_k}$
for $\de \in C({\mathcal{Z}}\cap U_V;\Bbb{R}^+)$ sufficiently small.

\begin{dfn}
\label{negligible_dfn}
{\rm Suppose ${\mathcal{M}}$ is a smooth manifold, 
$F,V \lra {\mathcal{M}}$ are vector bundles,
and $\Om$ is an open subset of~$F$.

(1)\qua A smooth bundle map $\ve\co \Om \lra V$ is 
\emph{$C^0$--negligible} if $\lim_{\ups\lra0}\ve(\ups) = 0$.

(2)\qua A $C^0$--negligible map $\ve\co \Om \lra V$ is 
\emph{$C^1$--negligible} if $\lim_{\ups\lra0}D_{\mathcal{M}}\ve(\ups) = 0$,
where $D_{\mathcal{M}}\ve$ denotes the differentiation of $\ve$
along~${\mathcal{M}}$ with respect to some connections in $F$~and~$V$.

(3)\qua If $\rho$ is a positive or negative monomials map on 
$F = \bigoplus_{i\in I}F_i$ and $\tilde{\rho}$ 
is a monomials map on $F$ with values in a line bundle~$L$,
a smooth bundle map $\ve\co \Om \lra V$ is 
\emph{$(\rho,\tilde{\rho})$--controlled} if there exist 
monomials maps \hbox{$\rho_1,\ldots,\rho_N$} on $F$ with values in~$L$
such that
$$\tilde{\rho}>_{\rho}\rho_j\quad\mbox{for all}\ j = 1,\ldots,N,
\quad\and\ \lim_{\stackrel{\ups\in\Om-Y(F;I)}{\ups\lra0}}
\Big(\sum_{j=1}^N|\rho_j(\ups)|\Big)^{-1}\big|\ve(\ups)\big|<\i.$$}
\end{dfn}

\begin{dfn}
\label{reg_dfn1a}
{\rm Suppose ${\overM}$, ${\mathcal{S}}$, and ${\mathcal{Z}} \subset
{\mathcal{M}}_k$ are as in Definition~\ref{normal_model_dfn},
$V \lra {\overS}$ is a vector bundle, and $s \in \Ga({\mathcal{S}};V)$.

(1)\qua A \emph{semi-regularization of $s$ near ${\mathcal{Z}}$} is a~tuple 
$$\big(U_{\mathcal{Z}},E,\ka,({\mathcal{F}}_{k;j})_{j\in J({\mathcal{Z}})};
{\mathcal{O}}_{\mathcal{Z}}^- \oplus {\mathcal{O}}_{\mathcal{Z}}^+,
\psi_{\mathcal{Z}}^- \oplus \psi_{\mathcal{Z}}^+;
\tilde{\mathcal{F}},\rho,\al_+,\al_V,\nu^*\big),$$
where
\begin{enumerate}
\item[(1a)] $\big(U_{\mathcal{Z}},E,\ka,
  ({\mathcal{F}}_{k;j})_{j\in J({\mathcal{Z}})};
  {\mathcal{O}}_{\mathcal{Z}}^- \oplus {\mathcal{O}}_{\mathcal{Z}}^+,
  \psi_{\mathcal{Z}}^- \oplus \psi_{\mathcal{Z}}^+\big)$
is model for  ${\mathcal{Z}}$ in ${\mathcal{S}}$
such that $\rk E \ge \rk{\mathcal{O}}_{\mathcal{Z}}^-$;

\item[(1b)] $\tilde{\mathcal{F}} =
  \bigoplus_{i\in\tilde{I}({\mathcal{Z}})} \tilde{\mathcal{F}}_i
  \lra U_{\mathcal{Z}}$ is a vector bundle and 
$$\rho\co  {\mathcal{F}}_k|{\mathcal{U}}_{\mathcal{Z}}
  - Y({\mathcal{F}}_k|{\mathcal{U}}_{\mathcal{Z}};J({\mathcal{Z}}))
  \lra \tilde{\mathcal{F}}$$
is a smooth bundle~map;

\item[(1c)] $\nu^* \in \Ga(U_{\mathcal{Z}};{\mathcal{O}}_{\mathcal{Z}}^+
  \oplus V)$ and $\al_+ \oplus \al_V \in 
  \Ga\big(U_{\mathcal{Z}};\hbox{Hom}(\tilde{\mathcal{F}};
                                   {\mathcal{O}}^+_{\mathcal{Z}} \oplus V)\big)$
is an injective linear map such that $\al_+$ is onto along~${\mathcal{Z}}$;
\
\item[(1d)] there exist $\de \in C^{\i}({\mathcal{Z}}\cap U_V;\Bbb{R}^+)$ and
a $C^0$--negligible map
$$\tilde{\ve}_{+,V} \co {\mathcal{F}}_k|{\mathcal{U}}_{\mathcal{Z}} - 
Y({\mathcal{F}}_k|{\mathcal{U}}_{\mathcal{Z}};J({\mathcal{Z}}))\lra\hbox{Hom}(\tilde{\mathcal{F}};
 {\mathcal{O}}_{\mathcal{Z}}^+ \oplus V)$$
such that
$$\big(\psi_+(b;\ups),s(b;\ups)\big)=
\big(\al_+\rho(b;\ups),\al_V\rho(b;\ups)\big)
+\nu^*(b)+\tilde{\ve}_{+,V}(b;\ups)\rho(b;\ups)$$
for all $(b;\ups) \in {\mathcal{F}}_{k;\de}$
such that $\phi_k(b;\ups) \in {\mathcal{S}}$.
\end{enumerate}

(2)\qua A semi-regularization 
$$\big(U_{\mathcal{Z}},E,\ka,({\mathcal{F}}_{k;j})_{j\in J({\mathcal{Z}})};
{\mathcal{O}}_{\mathcal{Z}}^- \oplus {\mathcal{O}}_{\mathcal{Z}}^+,
\psi_{\mathcal{Z}}^- \oplus \psi_{\mathcal{Z}}^+;
\tilde{\mathcal{F}},\rho,\al_+,\al_V,\nu^*\big)$$
is \emph{hollow} if 
$\rk\tilde{\mathcal{F}}_k <  \rk{\mathcal{F}}_k$ and either $\nu^* = 0$ or 
the bundle~map
$$\ti{\mathcal{F}}\lra{\mathcal{O}}^+ \oplus V,\qquad
(b;\ti{\ups})\lra \big\{\al_+ \oplus \al_V\big\}(b;\ti{\ups})+
\nu^*(b),$$
does not vanish over~${\mathcal{Z}}$.

(3)\qua  A semi-regularization 
$$\big(U_{\mathcal{Z}},E,\ka,({\mathcal{F}}_{k;j})_{j\in J({\mathcal{Z}})};
{\mathcal{O}}_{\mathcal{Z}}^- \oplus {\mathcal{O}}_{\mathcal{Z}}^+,
\psi_{\mathcal{Z}}^- \oplus \psi_{\mathcal{Z}}^+;\ti{\mathcal{F}},\rho,\al_+,\al_V,\nu^*\big)$$
is \emph{neutral}  if $\rho$ is a neutral monomials map,
$\al_+|Y(\ti{\mathcal{F}};\{i\})$ is onto over~${\mathcal{Z}}$
for all $i \in \tilde{I}({\mathcal{Z}})$,
and $\nu^* = 0$.
}
\end{dfn}

\begin{dfn}
\label{reg_dfn1b}
{\rm
Suppose ${\overM}$, ${\mathcal{S}}$, ${\mathcal{Z}} \subset {\mathcal{M}}_k$,
$V \lra {\overS}$, and $s \in \Ga({\mathcal{S}};V)$
are as in Definition~\ref{reg_dfn1a}.
A \emph{regularization of $s$ near ${\mathcal{Z}}$} is a~tuple 
\begin{multline*}
\big(U_{\mathcal{Z}},(E_i)_{i\in I({\mathcal{Z}})},
(\ka_i)_{i\in I({\mathcal{Z}})},({\mathcal{F}}_{k;j})_{j\in J({\mathcal{Z}})};\\
{\mathcal{O}}_{\mathcal{Z}}^- \oplus {\mathcal{O}}_{\mathcal{Z}}^+,
\psi_{\mathcal{Z}}^- \oplus \psi_{\mathcal{Z}}^+;
\tilde{\mathcal{F}},\rho,\al_+,\al_V;
(\tilde{E}_i)_{i\in I({\mathcal{Z}})},\tilde{\rho},\al_-\big),
\end{multline*}
where

(1)\qua $E_i \lra U_{\mathcal{Z}}$ is a vector bundle for each
$i \in I({\mathcal{Z}})$;

(2)\qua with $E = \bigoplus_{i\in I({\mathcal{Z}})}E_i$
and $\ka = (\ka_i)_{i\in I({\mathcal{Z}})}$,
$$\big(U_{\mathcal{Z}},E,\ka,({\mathcal{F}}_{k;j})_{j\in J({\mathcal{Z}})};
{\mathcal{O}}_{\mathcal{Z}}^- \oplus {\mathcal{O}}_{\mathcal{Z}}^+,
\psi_{\mathcal{Z}}^- \oplus \psi_{\mathcal{Z}}^+;
\ti{\mathcal{F}},\rho,\al_+,\al_V,0)$$
is a semi-regularization of $s$ near ${\mathcal{Z}}$ such that
$\rho$ is a positive or negative monomials map and
$\al_+|Y(\ti{\mathcal{F}};\{i\})$ is onto over~${\mathcal{Z}}$ for all 
$i \in \ti{I}({\mathcal{Z}})$;

(3)\qua $\ti{\rho} = (\ti{\rho}_i)_{i\in I({\mathcal{Z}})}$
is a degree--${\mathcal{A}}$ monomials map on $E \oplus {\mathcal{F}}_k$
with values in~$\tilde{E} \equiv \bigoplus_{i\in I({\mathcal{Z}})}\tilde{E}_i$, 
for some integer-valued function ${\mathcal{A}}$ on 
$I({\mathcal{Z}}) \times J({\mathcal{Z}})$;

(4)\qua $\al_- \in \Ga\big({\mathcal{Z}};
\hbox{Hom}(\tilde{E},{\mathcal{O}}_{\mathcal{Z}}^-)\big)$
is an isomorphism on every fiber;

(5)\qua there exist $\de \in C({\mathcal{Z}}\cap U_V;\Bbb{R}^+)$
and for each $i \in I({\mathcal{Z}})$
a $C^1$--negligible map and a $(\rho,\tilde{\rho}_i)$--controlled map,
\begin{gather*}
\tilde{\ve}^-_i\co {\mathcal{F}}_k|{\mathcal{U}}_{\mathcal{Z}} 
- Y({\mathcal{F}}_k|{\mathcal{U}}_{\mathcal{Z}};J({\mathcal{Z}})) \lra 
\hbox{Hom}(\tilde{F}_i(\tilde{\rho}),\al_i(\tilde{E}_i))
\qquad\hbox{and}\\
\ve^-_i\co {\mathcal{F}}_k|{\mathcal{U}}_{\mathcal{Z}} 
- Y({\mathcal{F}}_k|{\mathcal{U}}_{\mathcal{Z}};J({\mathcal{Z}})) \lra \al_i(\tilde{E}_i)
\end{gather*}
such that
$$\psi^-_i(b;\ups)=
\al_-\tilde{\rho}_i(b;\ka_i(b),\ups)+
\tilde{\ve}^-_i(b;\ups)\tilde{\rho}_i(\ups)+\ve^-_i(b;\ups)
\quad\mbox{for all}\ (b;\ups) \in {\mathcal{F}}_{k;\de},$$
where $\psi^-_i$ denotes the $i$th
component of~$\psi^-$ under the decomposition
${\mathcal{O}}^-_{\mathcal{Z}}|{\mathcal{Z}} = 
\bigoplus_{i\in I({\mathcal{Z}})}\al_-(\tilde{E}_i)$.}
\end{dfn}

\begin{remark}
{\rm Proposition~\ref{str_prp} ensures that $\psi^-$ admits an expansion as 
in (5) of Definition~\ref{reg_dfn1b}
in most cases one would encounter in counting rational curves
in projective spaces.
However, this expansion is not needed if ${\mathcal{Z}}$ is
$s$--hollow or $s$--neutral; see Definition~\ref{reg_dfn2}
and Proposition~\ref{euler_prp}.}
\end{remark}

If ${\mathcal{F}},E \lra {\mathcal{Z}}$ are vector bundles
and $\al \in \Ga\big({\mathcal{Z}};\hbox{Hom}({\mathcal{F}},E)\big)$,
let
$$\tilde{\al} \in \Ga\big(\Bbb{P}{\mathcal{F}};
\hbox{Hom}(\ga_{\mathcal{F}},\pi_{\Bbb{P}{\mathcal{F}}}^*E)\big)$$
denote the section induced by~$\al$.
Here $\ga_{\mathcal{F}} \lra \Bbb{P}{\mathcal{F}}$ is the tautological line bundle,
while $\pi_{\Bbb{P}{\mathcal{F}}}\co \Bbb{P}{\mathcal{F}} \lra {\mathcal{Z}}$
is the bundle projection map.

If ${\mathcal{F}},E \lra {\mathcal{Z}}$ and $\al$ are as above,
a \textit{closure of $({\mathcal{Z}},\al)$} is a tuple 
$({\overM}',{\mathcal{F}}',E')$, where
${\overM}'$ is an ms-manifold containing ${\mathcal{Z}}$ as
a pseudovariety, and
${\mathcal{F}}'$ and $E'$ are vector bundles 
over ${\overZ}$ that restrict to ${\mathcal{F}}$ and $E$,
respectively, over~${\mathcal{Z}}$.

\begin{dfn}
\label{reg_dfn2a}
{\rm
Let ${\overM}$, ${\mathcal{S}}$, ${\mathcal{Z}}$, $V$, and $s$ 
be as in Definition~\ref{reg_dfn1a}.
A regularization 
\begin{multline*}
\big(U_{\mathcal{Z}},(E_i)_{i\in I({\mathcal{Z}})},
(\ka_i)_{i\in I({\mathcal{Z}})},({\mathcal{F}}_{k;j})_{j\in J({\mathcal{Z}})};\\
{\mathcal{O}}_{\mathcal{Z}}^- \oplus {\mathcal{O}}_{\mathcal{Z}}^+,
\psi_{\mathcal{Z}}^- \oplus \psi_{\mathcal{Z}}^+;
\tilde{\mathcal{F}},\rho,\al_+,\al_V;
(\tilde{E}_i)_{i\in I({\mathcal{Z}})},\tilde{\rho},\al_-\big)
\end{multline*}
of $s$ at ${\mathcal{Z}}$ is \emph{closable} if
$({\mathcal{Z}},\al_V)$ admits a closure
$({\overM}',\tilde{\mathcal{F}}',{\mathcal{O}}')$
and $\Bbb{P}\tilde{\mathcal{F}}'$ admits an ms-orbifold  structure
such~that

(1)\qua ${\overM}' \subset {\overM}$,
$\pi_{\Bbb{P}\tilde{\mathcal{F}}'}\co \Bbb{P}\tilde{\mathcal{F}}' \lra {\overM}$
is an ms-map, and
$\tilde{\al}_+^{-1}(0)$ is a pseudovariety in~$\Bbb{P}\tilde{\mathcal{F}}'$;

(2)\qua $\partial\ov{\tilde{\al}_+^{-1}(0)}$ is a union of subsets ${\mathcal{Z}}_i$
such that $\tilde{\al}_V$ admits a semi-regularization at
each~${\mathcal{Z}}_i$.}
\end{dfn}

\begin{dfn}
\label{reg_dfn2}
{\rm Let ${\overM}$, ${\mathcal{S}}$, $V$, and $s$ 
be as in Definition~\ref{reg_dfn1a}.

(1)\qua ${\mathcal{Z}} \subset {\mathcal{M}}_k$ is \emph{$s$--hollow} (neutral, regular) 
if $s$ admits a hollow semi-reg\-u\-lar\-i\-za\-tion
(neutral semi-regularization, closable regularization) near~${\mathcal{Z}}$.

(2)\qua A section $s$ is \under{re}g\under{ular} if $s$ 
is transversal to the zero set and
there exists a finite collection $\{{\mathcal{Z}}_i\}_{i\in I_s}$
of smooth disjoint manifolds in ${\overM}$ such that
\begin{enumerate}
\item[(2a)] $\partial{\overS}\subset\bigsqcup_{i\in I_s} {\mathcal{Z}}_i$ and 
${\overZ}_i - {\mathcal{Z}}_i\subset
\bigcup_{\dim{\mathcal{Z}}_j<\dim{\mathcal{Z}}_i} {\mathcal{Z}}_j$
for all $i \in I_s$;
\item[(2b)] ${\mathcal{Z}}_i$ is either $s$--hollow, $s$--neutral, or 
$s$--regular for every $i \in I_s$.
\end{enumerate}}
\end{dfn}

Suppose $s \in \Ga({\mathcal{S}};V)$ is a regular section as in (2)
of Definition~\ref{reg_dfn2}.
Let $\{{\mathcal{Z}}_i\}_{i\in I_s^*}$  be the subcollection of
$s$--regular subsets.
To each $i \in I_s^*$, we associate the~tuple
$$\vr_i=({\overM}_i,{\mathcal{S}}_i,\ga_i,{\mathcal{O}}_i,\al_i;\deg\rho_i),$$
where, with notation as in Definitions~\ref{reg_dfn1b}
and~\ref{reg_dfn2a},
\begin{gather*}
{\overM}_i = \Bbb{P}\tilde{\mathcal{F}}', \quad
{\mathcal{S}}_i = \tilde{\al}_+^{-1}(0) \subset{\overM}_i, \quad
\ga_i = \ga_{\tilde{\mathcal{F}}'} \lra {\overS}_i, \quad
{\mathcal{O}}_i = \pi_{\tilde{\mathcal{F}}'}^*V\lra{\overS}_i,\\
\al_i = \tilde{\al}_{V} \in
\Ga\big({\mathcal{S}}_i,\hbox{Hom}(\ga_i,{\mathcal{O}}_i)\big),
\end{gather*}
are the objects corresponding to~${\mathcal{Z}}_i$.
We will write $\deg\vr_i$ for~$\deg\rho_i$.

\subsection{Contributions to the Euler class}
\label{top_comp_subs}

In this subsection, we describe 
the topological part of the computational method of this paper:

\stepcounter{thm}

\setcounter{alprp}0
\begin{alprp}
\label{zeros_prp}
Suppose ${\mathcal{S}}$ is an $m$--pseudovariety in an ms-orbifold ${\overM}$
and $E,{\mathcal{O}} \lra {\overS}$ are vector bundles such that
$\rk E = 1$ and $\rk{\mathcal{O}} = m + 1$.
If $\al \in \Ga\big({\mathcal{S}};\hbox{Hom}(E,{\mathcal{O}})\big)$ 
is a regular section, for a dense open subset 
$\Ga_{\al}({\overS};{\mathcal{O}})$ of sections
$\bar{\nu} \in \Ga({\overS};{\mathcal{O}})$,
the affine map 
$$\psi_{\al,\bar{\nu}} \equiv \al + \bar{\nu}\co  E|{\mathcal{S}}\lra{\mathcal{O}},
\qquad \psi_{\al,\bar{\nu}}(b;\ups)=\al_b(\ups)+\bar{\nu}_b,$$
is transversal to the zero set.
The set $\psi_{\al,\bar{\nu}}^{-1}(0)$ is finite, and
its signed cardinality $^{\pm} |\psi_{\al,\bar{\nu}}^{-1}(0)|$ 
is dependent only on~$\al$ and is given~by
$$N(\al)\equiv ^{\pm} |\psi_{\al,\bar{\nu}}^{-1}(0)|
=\big|\al^{\perp-1}(0)\big|,$$
where $\al^{\perp} \in 
\Ga\big({\mathcal{S}};\hbox{Hom}(E,{\mathcal{O}}/\Bbb{C})\big)$
is the composition of $\al$ with projection map onto
the quotient of ${\mathcal{O}}$ by a generic trivial line subbundle.
\end{alprp}

\begin{alprp}
\label{euler_prp}
Suppose ${\mathcal{S}}$ is an $m$--pseudovariety in ms-orbifold ${\overM}$,
$V \lra {\overS}$ is vector bundle of rank~$m$, 
such that $e(V)$ is the restriction of a cohomology class on~${\overM}$,
and $s \in \Ga(S;V)$ is a regular section.
Then, $s^{-1}(0)$ is a finite set, and
$$^{\pm} \big|s^{-1}(0)\big|
=\big\lan e(V),{\overS}\big\ran-
\sum_{i\in I^*_s}\deg\vr_i\cdot N(\vr_i) \equiv
\big\lan e(V),{\overS}\big\ran-{\mathcal{C}}_{\partial{\overS}}(s),$$
where $I^*_s$ is a complete collection of 
effective regularizations of $s$ on ${\overS}$,
as in the last paragraph of Subsection~\ref{top_dfn_subs}
and $N(\vr_i) \equiv N(\al_i)$.
\end{alprp}

\begin{remarks}
{\rm (1)\qua Together the two propositions give a reductive procedure
for counting the number of zeros of a section 
over the main stratum of a pseudovariety,
provided that the section is ``reasonably nice.''
Indeed, one application of both propositions reduces 
the rank of the target bundle by~one.
In the holomorphic category, every section is in fact ``reasonably nice.''
By Proposition~\ref{str_prp}, many sections 
of interest to us also have the needed properties.

(2)\qua In Proposition~\ref{zeros_prp}, 
$E$ can be a vector bundle of arbitrary rank,
provided the rank of ${\mathcal{O}}$ is adjusted appropriately and 
the section $\tilde{\al}$ is regular.
In such a case, $N(\al) = N(\tilde{\al})$.
In fact, one can obtain such a reduction even if
the original map~$\al$ is a polynomial;
see Subsection~3.3 in~\cite{Z1} for details.
In addition, it is not necessary to assume that $\al$ does not vanish over~${\mathcal{S}}$.
However, in practical applications, the boundary of ${\overM}$
can be enlarged to absorb the zero set of~$\al$.

(3)\qua If $E,E' \lra {\overM}$ are vector bundles,
$\rho \in \Ga\big({\mathcal{S}};\hbox{Hom}(E,E')\big)$ is an isomorphism
on every fiber,
and \hbox{$\al \in \Ga\big({\mathcal{S}};\hbox{Hom}(E,{\mathcal{O}})\big)$},
then $N(\al) = N(\al \circ \rho^{-1})$, provided
both numbers are defined.
Note that the isomorphism $\rho$ is assumed to be defined
over~${\mathcal{S}}$, and not over~${\overS}$.
We will call the replacement of $\al$ by $\al \circ \rho^{-1}$
a \textit{rescaling of the linear~map}.
A good choice of the isomorphism~$\rho$ may greatly simplify
the computation of the number $N(\al)$
via Propositions~\ref{zeros_prp} and~\ref{euler_prp}.
In actual applications, our isomorphisms $\rho$ will be such that
$N(\al)$ is defined if and only if $N(\al \circ \rho^{-1})$
is defined.

(4)\qua If $\al_+$ and  $\al_V$ are as in the last paragraph
of Subsection~\ref{top_dfn_subs},
$$N(\tilde{\al}_V)=N(\al_+ \oplus \al_V).$$
In particular, if ${\mathcal{Z}}_i$ is a finite set and
thus $\al_+ \oplus \al_V$ is an isomorphism over every point 
of ${\mathcal{Z}}_i$, then \hbox{$N(\tilde{\al}_V) = ^{\pm}
|{\mathcal{Z}}_i|$}.}
\end{remarks}

For computational purposes, it is useful to observe that 
if $E \lra {\overS}$ is a vector bundle of rank~$n$ such that 
$c(E)$ is the restriction of an element of $H^*({\overM})$,
\begin{equation}\label{zeros_main_e}\begin{split}
&\la_E^n+ \sum_{k=1}^{k=n} c_k(E)\la_E^{n-k} = 0 \in  
H^{2n}(\Bbb{P}E;\Bbb{Q})
\qquad\hbox{and}\\
&\quad \big\lan \mu\la_E^{n-1},\Bbb{P}E|{\overS}\big\ran = 
\big\lan\mu,{\overS}\big\ran
\quad\mbox{for all}\ \mu \in H^{2m}({\overM};\Bbb{Q});
\end{split}\end{equation}
the same formula in a more standard setting can be found in
\cite[Section~20]{BT}, for example.
In the situation of Proposition~\ref{zeros_prp},
but with $E$ of an arbitrary rank,
Propositions~\ref{zeros_prp} and~\ref{euler_prp}
and equation~\eqref{zeros_main_e} give
\begin{equation}\label{zeros_e}
N(\al) = \big\lan c({\mathcal{O}})c(E)^{-1},{\overS}\big\ran
-{\mathcal{C}}_{\Bbb{P}E|\partial{\overS}}
(\tilde{\al}^{\perp}).
\end{equation}
The last term above is zero if $\al$ extends to a section
of $\hbox{Hom}(E,{\mathcal{O}})$ over ${\overS}$ that has full rank
over every point of~${\overS}$.
In the computational sections of this papers,
we view formulas~\eqref{zeros_main_e} and~\eqref{zeros_e} 
as parts of Propositions~\ref{zeros_prp} and~\ref{euler_prp}.

In the rest of this subsection, we prove 
Propositions~\ref{zeros_prp} and~\ref{euler_prp}.
Before proceeding, we first comment on the topology on 
$\Ga({\overS};{\mathcal{O}})$ to which the first proposition makes an
implicit reference.
There are many topologies in which the statement of the proposition is valid.
One of them is defined via convergence of sequences 
on compact subsets in the $C^0$--norm on $U_{\mathcal{O}}$ and 
the $C^2$--norm on compact subset of~$U_{\mathcal{O}}\cap{\mathcal{M}}_k$;
see \cite[Subsection~3.2]{Z1} for more details.

The proof of Proposition~\ref{zeros_prp} is essentially the same as
the proof of \cite[Lemma~3.14]{Z1}.
The finiteness claim is proved as follows.
Suppose $(b_r,v_r) \in E|{\mathcal{S}}$ is a sequence such that
$\psi_{\al,\bar{\nu}}(b_r,v_r) = 0$ and $\{b_r\}$ converges to
some $b^* \in \partial{\overS}$.
Let 
$$\big(U_{\mathcal{Z}},E',\ka,({\mathcal{F}}_{k;j})_{j\in J({\mathcal{Z}})};
{\mathcal{O}}_{\mathcal{Z}}^- \oplus {\mathcal{O}}_{\mathcal{Z}}^+,
\psi_{\mathcal{Z}}^- \oplus \psi_{\mathcal{Z}}^+;
\tilde{\mathcal{F}},\rho,\al_+,\al_{E^*\otimes{\mathcal{O}}},\nu^*\big),$$
be a semi-regularization of $\al$ at a submanifold ${\mathcal{Z}}$ of
${\mathcal{M}}_k$ containing~$b^*$, 
as provided by (2) of Definition~\ref{reg_dfn2}.
By replacing ${\mathcal{Z}}$ with ${\mathcal{Z}}\cap U_E\cap U_{\mathcal{O}}$, if necessary, 
it can be assumed that \hbox{${\mathcal{Z}}_i \subset U_E\cap U_{\mathcal{O}}$}.
As in the proof of Lemma~3.12 in~\cite{Z1},
from the sequence $\{(b_r,v_r)\}$ we can obtain a zero of the~map
$$\psi_{\al,\bar{\nu};{\mathcal{Z}}}\co 
E \otimes \tilde{\mathcal{F}}\lra E \otimes {\mathcal{O}}^+\oplus{\mathcal{O}}, 
\psi_{\al,\bar{\nu};{\mathcal{Z}}}(b;\ups)=
\big(\al_+(\ups),\al_{E^*\otimes{\mathcal{O}}}(\ups)+\bar{\nu}\big),$$
if $\nu^* = 0$.
This is a bundle map over~${\mathcal{Z}}$.
By (1c) of Definition~\ref{reg_dfn1a}, the first-component map is surjective.
Thus, for a dense open subset of elements of $\Ga({\overS};{\mathcal{O}})$,
the map $\psi_{\al,\bar{\nu};{\mathcal{Z}}}$ is transversal to the zero set.
Assumptions on the dimension and the ranks imply that 
$\psi_{\al,\bar{\nu};{\mathcal{Z}}}(b;\ups)$ does not vanish 
if $\bar{\nu}$ lies in this open dense subset;
see (2) of Definition~\ref{normal_model_dfn} and
(1a)--(1c) of Definition~\ref{reg_dfn1a}.
On the other hand, if $\bar{\nu}^* \neq 0$, we can obtain a zero of the~map
\begin{gather*}
\psi_{\al,\bar{\nu};{\mathcal{Z}}}\co 
E \otimes \tilde{\mathcal{F}} \oplus E\lra 
E \otimes {\mathcal{O}}^+\oplus{\mathcal{O}}, \\
\psi_{\al,\bar{\nu};{\mathcal{Z}}}(b;w,\ups)=
\big(\al_+(\ups),\al_{E^*\otimes{\mathcal{O}}}(\ups)+\bar{\nu}\big)
+\nu^*w.
\end{gather*}
In this case, $\rk\tilde{\mathcal{F}} < \rk{\mathcal{F}}_k$ and thus
the map $\psi_{\al,\bar{\nu};{\mathcal{Z}}}$ again has no zeros
if $\bar{\nu}$ is generic.
We conclude that $\psi_{\al,\bar{\nu}}^{-1}(0)$ is a finite set.
The independence of 
$^{\pm} \big|\psi_{\al,\bar{\nu}}^{-1}(0)\big|$ of 
the choice of~$\bar{\nu}$ is shown by constructing 
a cobordism between $\psi_{\al,\bar{\nu}_1}^{-1}(0)$ 
and~$\psi_{\al,\bar{\nu}_2}^{-1}(0)$;
see part (5) of the proof of
Lemma~\ref{euler_lmm4} below for a similar argument.
The final claim of Proposition~\ref{zeros_prp} is 
nearly immediate from the definition of~$\psi_{\al,\bar{\nu}}$;
see Subsection~3.3 in~\cite{Z1}.
The trivial subbundle mentioned in the statement of the proposition
is simply $\Bbb{C}\bar{\nu}$, if $\bar{\nu}\in\Ga({\overS};{\mathcal{O}})$
is generic and thus does not vanish.

We next prove Proposition~\ref{euler_prp}.
The first step is to construct a section 
\hbox{$\tilde{s} \in \Ga({\overS};V)$}
such that $\tilde{s} = s$ outside of a small neighborhood
of $\partial{\overS}$ that contains no zeros of~$s$.
This is achieved by cutting~$s$ off near~$\partial{\overS}$,
so that the new section extends over $\partial{\overS}$ by~zero.
This procedure changes the estimate of (1d) of 
Definition~\ref{reg_dfn1a} in a well-controlled manner.
We then add a small perturbation $t\nu$ to $\tilde{s}$ such that
$\tilde{s} + t\nu$ has transverse zeros on ${\mathcal{S}}$ and
no zeros on~$\partial{\overS}$.
The total number of zeros of this section is then the euler class of $V$ over~${\overS}$,
$\lr{e(V),{\overS}}$.
On the other hand, for each element of $s^{-1}(0)$ there will be
a nearby zero of the perturbed section.
All the remaining zeros will lie near~$\partial{\overS}$.
The final step is to show that all such zeros lie near 
the $s$--regular subsets of $\partial{\overS}$, 
for a good choice of~$\nu$,
and can be expressed in terms of the zeros of affine~maps.

Let $\{{\mathcal{Z}}_i\}_{i\in I_s}$ be a collection of smooth manifolds
in ${\overM}$
as in Definition~\ref{reg_dfn2} with corresponding semi-regularizations
$$\big(U_i,E_i,\ka_i,({\mathcal{F}}_{k_i;j})_{j\in J({\mathcal{Z}}_i)};
{\mathcal{O}}_i^- \oplus {\mathcal{O}}_i^+,
\psi_i^- \oplus \psi_i^+;
\tilde{\mathcal{F}}_i,\rho_i,\al_{i;+},\al_{i;V},\nu_i^*\big).$$
By replacing ${\mathcal{Z}}_i$ with ${\mathcal{Z}}_i\cap U_V$,
if necessary, it can be assumed that \hbox{${\mathcal{Z}}_i\subset U_V$}.
Let $\tilde{\de}_i \in C({\mathcal{Z}}_i;\Bbb{R}^+)$ be such that
$\vt_{\ka_i}$ is defined on $E_{\tilde{\de}_i}$,
$$\vt_{\ka_i;{\mathcal{F}}_{k_i}}\big(
E_{\tilde{\de}_i} \times_{{\mathcal{Z}}_i} {\mathcal{F}}_{{k_i};\tilde{\de}_i}\big)
\subset U_{k_i},\quad\hbox{and}\quad
\phi_{k_i}\vt_{\ka_i;{\mathcal{F}}_{k_i}}
\big(E_{\tilde{\de}_i} \times_{{\mathcal{Z}}_i} 
{\mathcal{F}}_{{k_i};\tilde{\de}_i}\big)
\subset U_V.$$
If $\de\in C({\mathcal{Z}}_i;\Bbb{R}^+)$ is less than~$\tilde{\de}_i$
and ${\mathcal{K}}$ is a subset of~${\mathcal{Z}}_i$, let
\begin{gather*}
{\mathcal{U}}_i(\de) =
   E_{i;\de} \times _{{\mathcal{Z}}_i}{\mathcal{F}}_{k_i;\de}
	 - Y(E_i \oplus {\mathcal{F}}_{k_i};J({\mathcal{Z}}_i)), \qquad\hbox{and}\\
W_i(\de;{\mathcal{K}}) = \phi_{k_i}\big(\vt_{\ka_i;{\mathcal{F}}_{k_i}}
(E_{i;\de} \times _{\mathcal{K}}{\mathcal{F}}_{k_i;\de})\big)
\subset{\overM}.
\end{gather*}
We denote $W_i(\de;{\mathcal{Z}}_i)$ by $W_i(\de)$.
Since all the zeros of $s$ are transverse, the next lemma
implies that $s^{-1}(0)$ is a finite~set.

\begin{lmm}
\label{euler_lmm1}
There exists $\tilde{\de}_i^* \in C({\mathcal{Z}}_i;\Bbb{R}^+)$ such
 that $W_i(\tilde{\de}_i^*)\cap s^{-1}(0) = \eset$.
\end{lmm}

\begin{proof}
Let $\tilde{\al}_i = \al_{i;+} \oplus \al_{i;V}$.
By (1c), (2), and (3) of Definition~\ref{reg_dfn1a} and 
(2) of Definition~\ref{reg_dfn1b},  there exists 
$\ve_{\tilde{\al}_i} \in C({\mathcal{Z}}_i;\Bbb{R}^+)$ such~that
\begin{equation}\label{euler_lmm1e1}
\big|\tilde{\al}_i(b;\tilde{\ups})+\nu_i^*(b)\big|
\ge \ve_{\tilde{\al}_i}(b)\big|\tilde{\ups}\big|
\qquad\mbox{for all}\  (b;\tilde{\ups}) \in \tilde{\mathcal{F}}_i.
\end{equation}
If $\tilde{\ve}_{i;+,V}$ is a $C^0$--negligible map as in (1d) of 
Definition~\ref{reg_dfn1a} corresponding to~${\mathcal{Z}}_i$, let
\begin{equation*}\begin{split}
\bar{\ve}_{i;+}(b,r) = \sup\big\{
\big|\tilde{\ve}_{i;+,V}(b;w,\ups)\big| :
(b;w,\ups) \in E_{\tilde{\de}_i} \times_{{\mathcal{Z}}_i}
 {\mathcal{F}}_{{k_i};\tilde{\de}_i}; \qquad\qquad\qquad&\\
|w|,|\ups| < r,~(b;\ups) \not\in Y({\mathcal{F}}_{k_i},J({\mathcal{Z}}_i))\big\}.&
\end{split}\end{equation*}
Then, $\bar{\ve}_{i;+}$ is continuous and
$\lim_{r\lra0}\bar{\ve}_{i;+}(b,r) = 0$.
Suppose 
$$(b;w,\ups) \in {\mathcal{U}}_i(\de),
\quad \phi_{k_i}\vt_{\ka_i;{\mathcal{F}}_{k_i}}(b;w,\ups) \in {\mathcal{S}},
\quad\hbox{and}\quad
s\big(\phi_{k_i}\vt_{\ka_i;{\mathcal{F}}_{k_i}}(b;w,\ups)\big)=0,$$
where $\de \le \tilde{\de}_i$.
Then, by Definition~\ref{reg_dfn1a},
\begin{gather}\label{euler_lmm1e2}
\big|\tilde{\al}_i\rho_i(b;\ups)+\nu_i^*(b)\big|\le
\bar{\ve}_{i;+}\big(b,\de(b)\big)\big|\rho_i(b;\ups)\big|.
\end{gather}
By~\eqref{euler_lmm1e1} and~\eqref{euler_lmm1e2},
if $\bar{\ve}_{i;+}\big(b,\de(b)\big) < 
\ve_{\tilde{\al}_i}(b)$ for all $b \in {\mathcal{Z}}_i$, then
$W_i(\de)\cap s^{-1}(0)$ is empty.
\end{proof}

Choose $\nu \in \Ga({\overS};V)$ with the following properties:

(A1)\qua for all $i \in I_s$, $\bar{\nu}_i \equiv \nu|_{{\mathcal{Z}}_i}$
has no zeros, and the~map
$$\psi_{i;\bar{\nu}_i}\co \tilde{\mathcal{F}}_i\lra{\mathcal{O}}_i^+ \oplus V,
\qquad
\psi_{i;\bar{\nu}_i}(b;\ups) = \big(\al_{i;+}(\ups),
\bar{\nu}_i(b) + \al_{i;V}(\ups)\big),$$
is transversal to the zero set.
Furthermore, if ${\mathcal{Z}}_i$ is $s$--hollow, 
the sets
$\Bbb{R}\nu^*_i$
and
$\Bbb{R}^*\cdot\Im\psi_{i;\bar{\nu}_i}(b;\ups)$ are disjoint.

(A2)\qua If ${\mathcal{Z}}_i$ is $s$--regular or $s$--neutral, 
$\psi_{i;\bar{\nu}_i}|Y(\tilde{\mathcal{F}}_i;j)$ does not vanish
for all $j \in \tilde{I}({\mathcal{Z}}_i)$;

(A3)\qua If ${\mathcal{Z}}_i$ is $s$--regular, 
$\psi_{i;\bar{\nu}_i}^{-1}(0)$ is a finite set.

Note that (A1) and (A2) are just transversality assumptions,
due to (1c) and (3) of Definition~\ref{reg_dfn1a}
and (2) of Definition~\ref{reg_dfn1b}.
Conditions~(A3) holds if $\bar{\nu}_i$ is
the restriction of a section $\bar{\nu}' \in \Ga({\overZ}_i;V)$
such that 
$$\pi_{\tilde{\mathcal{F}}_i}^*\bar{\nu}' \in 
\Ga_{\tilde{\al}_{i;V}}\big(\ov{\tilde{\al}_{i;+}^{-1}(0)};
\pi_{\tilde{\mathcal{F}}_i}^*V\big).$$

\begin{lmm}
\label{euler_lmm2}
If ${\mathcal{Z}}_i$ is $s$--hollow, there exists $\de_i^* \in C({\mathcal{Z}}_i;\Bbb{R}^+)$
with the following property.
If $W$ is an open neighborhood of $\partial{\overS}$
in ${\overM}$, there exists $\ep > 0$ such that
for all $\eta \in C({\mathcal{S}};\Bbb{R})$, $t \in \Bbb{R}^+$, and
$\nu' \in \Ga({\mathcal{S}};V)$, satisfying 
$\nu'|_W = \nu|_W$ and $\|\nu - \nu'\|_{C^2({\mathcal{S}}-W)} < \ep$,
$$W_i(\de_i^*)\cap\big\{\eta s + t\nu'\big\}^{-1}(0)=\eset.$$
\end{lmm}

\begin{proof}
By assumption (A1), the map
$$\tilde{\mathcal{F}}_i\oplus\Bbb{R}\lra{\mathcal{O}}_i^+ \oplus V,
\qquad
(b;\ups)\lra\big(\al_{i;+}(b;\ups),\bar{\nu}_i(b) + \al_{i;V}(\ups)\big)
+\tau \nu^*_i(b),$$
has no zeros over ${\mathcal{Z}}_i$. 
Thus, there exists 
$\de_{\tilde{\al}_i,\bar{\nu}_i} \in C({\mathcal{Z}}_i;\Bbb{R}^+)$
such~that 
\begin{equation}\label{euler_lmm2e1}
\big|(0,\bar{\nu}_i(b))+\tilde{\al}_i(b;\tilde{\ups})+\tau \nu^*_i(b)\big|\ge
\de_{\tilde{\al}_i,\bar{\nu}_i}(b)\big(1+ 
\big|\tilde{\al}_i(b;\tilde{\ups})\big|\big)  
\quad\mbox{for}\  (b;\tilde{\ups}) \in \tilde{\mathcal{F}}_i.
\end{equation}
By continuity of $\nu$, there exists 
$\ve_{i,\nu} \in C^0({\mathcal{Z}}_i \times \Bbb{R};\Bbb{R})$
such that 
$$\lim_{r\lra0}\ve_{i,\nu}(b,r) = 0\quad\mbox{and}\quad
\big|\nu\big(\phi_{k_i}\vt_{\ka_i;{\mathcal{F}}_{k_i}}(b;w,\ups)\big) - 
\bar{\nu}_i(b)\big| \le  \ve_{i,\nu}\big(b,\de(b)\big)$$
for all $(b;w,\ups) \in {\mathcal{U}}_i(\de)$ such that
$\phi_{k_i}\vt_{\ka_i;{\mathcal{F}}_{k_i}}(b;w,\ups) \in {\mathcal{S}}$.
If $\nu' \in \Ga({\mathcal{S}};V)$ is such that 
$\nu'|_W = \nu|_W$ and $\|\nu - \nu'\|_{C^2({\mathcal{S}}-W)} < \ep$, then
$$\big|\big\{\nu - \nu'\}
\big(\phi_{k_i}\vt_{\ka_;{\mathcal{F}}_{k_i}}(b;w,\ups)\big)\big|
\le C_{i,W}(b)\de(b)\ep$$
for all $(b;w,\ups) \in {\mathcal{U}}_i(\de)$ such that
$\phi_{k_i}\vt_{\ka_i;{\mathcal{F}}_{k_i}}(b;w,\ups) \in {\mathcal{S}}$,
where $C_{i,W} \in C({\mathcal{Z}}_i;\Bbb{R})$ depends only on~$W$.
Thus, by the above and (1d) of Definition~\ref{reg_dfn1a},
\begin{multline}
\label{euler_lmm2e3}
\big|\big(0,t\bar{\nu}_i(b)\big)+
\tilde{\al}_i\big(\eta(\phi_{k_i}\vt_{\ka_i;{\mathcal{F}}_{k_i}}(b;w,\ups))
\rho_i(b;\ups)\big)
+\eta(\phi_{k_i}\vt_{\ka_i;{\mathcal{F}}_{k_i}}(b;w,\ups))\nu_i^*(b)\big| \\
\le \ve_{\tilde{\al}_i}(b)^{-1}\bar{\ve}_+(b,\de(b))
\big|\tilde{\al}_i\big(\eta(\phi_{k_i}\vt_{\ka_i;{\mathcal{F}}_{k_i}}(b;w,\ups))
                            \rho_i(b;\ups)\big)\big|\\
+t\big(\ve_{i,\nu}(b,\de(b)) + C_{i,W}(b)\de(b)\ep\big)
\end{multline}
for all $(b;w,\ups) \in {\mathcal{U}}_i(\de)$ such that
$\phi_{k_i}\vt_{\ka_i;{\mathcal{F}}_{k_i}}(b;w,\ups) \in {\mathcal{S}}$
and
$$\big\{\eta s +
t\nu'\big\}\phi_{k_i}\vt_{\ka_i;{\mathcal{F}}_{k_i}}(b;w,\ups) = 0,$$
where $\ve_{\tilde{\al}_i}$ is as in the proof of Lemma~\ref{euler_lmm1}.
By \eqref{euler_lmm2e1} and~\eqref{euler_lmm2e3},
$$W_i(\de)\cap\big\{\eta s + t\nu'\big\}^{-1}(0)=\eset$$
if
$\ve_{\tilde{\al}_i}(b)^{-1}\bar{\ve}_+(b,\de(b)),
\ve_{i,\nu}(b,\de(b)),C_{i,W}(b)\de(b)\ep<
\frac{1}{4}\de_{\tilde{\al}_i,\bar{\nu}_i}(b)$
for all $b \in {\mathcal{Z}}_i$.
\end{proof}

\begin{lmm}
\label{euler_lmm3}
If ${\mathcal{Z}}_i$ is $s$--neutral, there exists $\de_i^* \in C({\mathcal{Z}}_i;\Bbb{R}^+)$
with the following property.
If $W$ is an open neighborhood of $\partial{\overS}$
in ${\overM}$, there exists $\ep > 0$ such that
for all $\eta \in C({\mathcal{S}};\Bbb{R})$, $t \in \Bbb{R}^+$, and
$\nu' \in \Ga({\mathcal{S}};V)$, satisfying 
$\nu'|_W = \nu|_W$ and $\|\nu - \nu'\|_{C^2({\mathcal{S}}-W)} < \ep$,
$$W_i(\de_i^*)\cap\big\{\eta s + t\nu'\big\}^{-1}(0)=\eset.$$
\end{lmm}

\begin{proof}
By assumptions (A1) and (A2),
$\psi_{i,\bar{\nu}_i}^{-1}(0)$ is a discrete subset of 
$\tilde{\mathcal{F}}_i - Y(\tilde{\mathcal{F}}_i;\tilde{I}({\mathcal{Z}}_i))$
and contains at most one point of each fiber.
For each $\tilde{\ups} \in \psi_{i,\bar{\nu}_i}^{-1}(0)$,
let ${\mathcal{K}}_{\tilde{\ups}}$ be a neighborhood of 
$\tilde{\ups}$ in 
\hbox{$\tilde{\mathcal{F}}_i - Y(\tilde{\mathcal{F}}_i;\tilde{I}({\mathcal{Z}}_i))$}
such that the closure of ${\mathcal{K}}_{\tilde{\ups}}$ in
$\tilde{\mathcal{F}}_i - Y(\tilde{\mathcal{F}}_i;\tilde{I}({\mathcal{Z}}_i))$
is compact.
For simplicity, it can be assumed that all 
the sets ${\mathcal{K}}_{\tilde{\ups}}$ are disjoint.
By Proposition~\ref{monomials_prp1}, there exists 
$\de_{\rho_i,\bar{\nu}_i} \in C^{\i}({\mathcal{S}};\Bbb{R}^+)$ such~that
$${\mathcal{F}}_{k_i;\de_{\rho_i,\bar{\nu}_i}}\cap\rho_i^{-1}
\Bigl( \bigcup_{\tilde{\ups}\in\psi_{i,\bar{\nu}_i}^{-1}(0)}
   \Bbb{R}^+\cdot{\mathcal{K}}_{\tilde{\ups}}\Bigr)=\eset.$$
Then, there exists $\de_{\bar{\nu}_i} \in C({\mathcal{Z}}_i;\Bbb{R}^+)$
such~that for all $\tau \in \Bbb{R}$
\begin{equation}\label{euler_lmm3e1}
\big|(0,\bar{\nu}_i(b)) + 
\tilde{\al}_i\big(\tau\rho_i(b;\ups)\big)\big|\ge
\de_{\bar{\nu}_i}(b)\big(1 + 
\big|\tilde{\al}_i\big(\tau\rho_i(b;\ups)\big)\big|\big)
\end{equation} 
for all $(b;\ups) \in {\mathcal{F}}_{k_i;\de_{\rho_i,\bar{\nu}_i}}  - 
               Y(\tilde{\mathcal{F}}_i;\tilde{I}({\mathcal{Z}}_i))$.
On the other hand, by the same argument as in the proof of 
Lemma~\ref{euler_lmm2},
\begin{multline}\label{euler_lmm3e3}
\big|\big(0,t\bar{\nu}_i(b)\big)+
\tilde{\al}_i\big(\eta(b;w,\ups) \rho_i(b;\ups)\big)\big| \\
 \le \ve_{\tilde{\al}_i}(b)^{-1}\bar{\ve}_+(b,\de(b))
\big|\tilde{\al}_i\big(\eta(b;w,\ups)\rho_i(b;\ups)\big)\big| \\
+t\big(\ve_{i,\nu}(b,\de(b)) + C_{i,W}(b)\de(b)\ep\big)
\end{multline}
for all $(b;w,\ups) \in {\mathcal{U}}_i(\de)$ such that
$\phi_{k_i}\vt_{\ka_i;{\mathcal{F}}_{k_i}}(b;w,\ups) \in {\mathcal{S}}$
and
$$\big\{\eta s + t\nu'\big\}\phi_{k_i}(b;w,\ups) = 0,$$
if $\nu'|_W = \nu|_W$ and $\|\nu - \nu'\|_{C^2({\mathcal{S}}-W)} < \ep$.
Thus, 
$W_i(\de)\cap\big\{\eta s + t\nu'\big\}^{-1}(0) = \eset$ if
$\ve_{\tilde{\al}_i}(b)^{-1}\bar{\ve}_+(b,\de(b)),
\ve_{i,\nu}(b,\de(b)),C_{i,W}(b)\de(b)\ep<
\frac{1}{4}\de_{\bar{\nu}_i}(b)$
for all $b \in {\mathcal{Z}}_i$ and $\de\le\de_{\rho_i,\bar{\nu}_i}$.
\end{proof}

\begin{lmm}
\label{euler_lmm4}
If ${\mathcal{Z}}_i$ is $s$--regular,
there exist $\de_i^* \in C({\mathcal{Z}}_i;\Bbb{R}^+)$,
a compact subset ${\mathcal{K}}_{\tilde{\al}_i;\bar{\nu}_i}$ of
${\mathcal{Z}}_i$, and $\de_i,\ep_i \in \Bbb{R}^+$ with the following
property.  If $W$ is an open neighborhood of $\partial{\overS}$ 
in ${\overM}$, there exists $\ep_W \in \Bbb{R}^+$ such that
for every $\eta \in C({\mathcal{S}};\Bbb{R})$, satisfying
$$\eta\big(\phi_{k_i}\vt_{\ka_i;{\mathcal{F}}_{k_i}}(w,\ups)\big) =
\be_{\de_i}(|\ups|)$$
for all $(w,\ups) \in E_{i;\de_i^*} 
\times_{{\mathcal{K}}_{\tilde{\al}_i;\bar{\nu}_i}}{\mathcal{F}}_{k_i;\de_i^*}$
such that $\phi_{k_i}\vt_{\ka_i;{\mathcal{F}}_{k_i}}(w,\ups) \in
{\mathcal{S}}$,
$t \in (0,\ep_i)$, and $\nu' \in \Ga({\mathcal{S}};V)$, satisfying
$\nu'|_W = \nu|_W$, $\|\nu - \nu'\|_{C^2({\mathcal{S}}-W)} < \ep_W$,
and $\eta s + t\nu'$ is transversal to the zero set on~${\mathcal{S}}$,
then
$$^{\pm}\big|W_i(\de_i^*)\cap\big\{\eta s + t\nu'\big\}^{-1}(0)\big|
=\deg\rho_i\cdot ^{\pm}\big|\psi_{i,\bar{\nu}}^{-1}(0)\big|$$
and
$$W_i(\de_i^*)\cap\big\{\eta s + t\nu'\big\}^{-1}(0)
 \subset  W_i(\de_i^*;K_{\tilde{\al}_i;\bar{\nu}_i}).$$
Furthermore, if $W^*$ is a neighborhood of ${\mathcal{Z}}_i$
in ${\overM}$, $\de_i^*$ can be chosen so that $W_i(\de_i) \subset W^*$.
\end{lmm}

\begin{proof}
(1)\qua Let ${\mathcal{K}}$ be an open neighborhood
of the finite subset $\psi_{i,\bar{\nu}_i}^{-1}(0)$ of 
$\tilde{\mathcal{F}}_i - Y(\tilde{\mathcal{F}}_i;\tilde{I}({\mathcal{Z}}_i))$
such that the closure of ${\mathcal{K}}$ in 
$\tilde{\mathcal{F}}_i - Y(\tilde{\mathcal{F}}_i;\tilde{I}({\mathcal{Z}}_i))$ is compact.
Let $\tilde{\mathcal{K}}_{\tilde{\al}_i,\bar{\nu}_i} = \rho_i^{-1}
                                           (\Bbb{R}\cdot{\overK})$.
Note that by (1) of Proposition~\ref{monomials_prp2}, 
$\tilde{\mathcal{K}}_{\tilde{\al}_i,\bar{\nu}_i}$ is a closed subset
of ${\mathcal{F}}_{k_i} - {\mathcal{Z}}_i$.
We take
$${\mathcal{K}}_{\tilde{\al}_i,\bar{\nu}_i}=
\pi_{{\mathcal{F}}_{k_i}}(\tilde{\mathcal{K}}_{\tilde{\al}_i,\bar{\nu}_i})
=\pi_{\tilde{\mathcal{F}}_i}({\overK}).$$
This is a compact subset of ${\mathcal{Z}}_i$.

(2)\qua As in the proofs of Lemmas~\ref{euler_lmm2} and~\ref{euler_lmm3},
there exists $\de_{\bar{\nu}_i} \in C({\mathcal{Z}}_i;\Bbb{R}^+)$
such~that 
\begin{equation}\label{euler_lmm4e1}
\big|(0,\bar{\nu}_i(b))+
\tilde{\al}_i(b;\tilde{\ups})\big|>
\de_{\bar{\nu}_i}(b)\big(1 + 
\big|\tilde{\al}_i(b;\tilde{\ups})\big|\big)
\qquad\mbox{for all}\  
(b;\tilde{\ups}) \in \tilde{\mathcal{F}}_i - \Bbb{R}\cdot{\mathcal{K}}.
\end{equation} 
On the other hand, as before, by the estimate (1d) 
of Definition~\ref{reg_dfn1a}, 
\begin{multline}\label{euler_lmm4e2}
\big|\big(0,t\bar{\nu}_i(b)\big)+
\tilde{\al}_i\big(\eta(b;w,\ups)\rho_i(b;\ups)\big)\big|\\
\qquad\qquad\le \tilde{\ve}_+(b,\de_i^*(b))
\big|\tilde{\al}_i\big(\eta(b;w,\ups)\rho_i(b;\ups)\big)\big|
+t\big(\ve_{i,\nu}(b,\de_i^*(b)) + \tilde{C}_{i,W}(b)\ep_W\big)
\end{multline}
for all $(b;w,\ups) \in {\mathcal{U}}_i(\de_i^*)$ such that
$\phi_{k_i}\vt_{\ka_i;{\mathcal{F}}_{k_i}}(b;w,\ups) \in {\mathcal{S}}$
and
$$\big\{\eta s + t\nu'\big\}(b;w,\ups) = 0,$$
where $\tilde{\ve}_+(b,\de_i^*(b)) =
         \ve_{\tilde{\al}_i}(b)^{-1}\bar{\ve}_+(b,\de_i^*(b))$
and $\tilde{C}_{i,W}(b) = C_{i,W}(b)\de_i^*(b)$.
By \eqref{euler_lmm4e1} and \eqref{euler_lmm4e2},
\begin{multline}\label{euler_lmm4e3}
W_i(\de_i^*)\cap\big\{\eta s + t\nu'\big\}^{-1}(0)\\
\subset \phi_{k_i}\vt_{\ka_i;{\mathcal{F}}_{k_i}}\big(
\big\{(w,\ups) \in E_i \oplus {\mathcal{F}}_{k_i} : 
\ups \in \tilde{\mathcal{K}}_{\tilde{\al}_i,\bar{\nu}_i},
\ |w|,|\ups| < \de_i^*(b)\big\}\big) 
\end{multline}
if $\tilde{\ve}_+(b,\de_i^*(b)),\ve_{i,\nu}(b,\de_i^*(b)),
\tilde{C}_{i,W}(b)\ep_W <\frac{1}{4}\de_{\bar{\nu}_i}(b)$
for all $b \in {\mathcal{Z}}_i.$

(3)\qua Let $\tilde{\rho}_i \equiv
          (\tilde{\rho}_{i,i'})_{i'\in I({\mathcal{Z}}_i)}$ and 
$\al_{i;-}$
be the monomials map and the vector-bundle isomorphism
as in (3) and (4) of Definition~\ref{reg_dfn1b}, respectively,
 corresponding to ${\mathcal{Z}}_i$.
Similarly, for each $i' \in I({\mathcal{Z}}_i)$, 
let $\tilde{\ve}_{i,i'}^-$ and $\ve_{i,i'}^-$ 
be a $C^1$--negligible map and a $(\rho_i,\tilde{\rho}_{i,i'})$--controlled
map as in (5) of Definition~\ref{reg_dfn1b}
corresponding to~${\mathcal{Z}}_i$.
Since $\al_{i;-}$ is an isomorphism on every fiber
and the image of $\tilde{\ve}_{i,i'}^-$
lies in the subbundle $\al_{i;-}(E_{i,i'})$ of ${\mathcal{O}}^-$,
we can define a $C^1$--negligible map  
$$\tilde{\ve}_{i,i'}\co 
{\mathcal{F}}_{k_i}|{\mathcal{U}}_i 
- Y({\mathcal{F}}_{k_i}|{\mathcal{U}}_i;J({\mathcal{Z}}_i)) \lra E_i$$
by $\tilde{\ve}_{i,i'}^-(w,\ups)=
\al_{i;-}\big(\tilde{\ve}_{i,i'}(w,\ups)\otimes
\tilde{\rho}_{i,i'}(\ups)\big)$.
By the Contraction Principle, there exists $\de_- \in \Bbb{R}^+$
such that for all
$$(w',\ups) \in E_{i;\de_-} 
 \times_{{\mathcal{K}}_{\tilde{\al}_i,\bar{\nu}_i}} \mathcal{FT}_{k_i,\de_-}$$
the~equation
$$w+\sum_{i'\in I({\mathcal{Z}}_i)}\tilde{\ve}_{i,i'}(w,\ups)=w',
\qquad w \in E_{i;2\de_-},$$
has a unique solution.
Furthermore, this solution satisfies $|w| \le 2|w'|$;
see the proof of \cite[Lemma~3.18]{Z1}, for example.
Let $\ep_- \in \Bbb{R}^+$ be such that
\begin{equation}\label{euler_lmm4e5}
\big|\al_{i;-}(\tilde{w})\big|\ge\ep_-|\tilde{w}|
\qquad\mbox{for all}\  
\tilde{w} \in \tilde{E}|{\mathcal{K}}_{\tilde{\al}_i,\bar{\nu}_i}.
\end{equation}
By (3) of Definition~\ref{negligible_dfn} and 
                             (3) of Proposition~\ref{monomials_prp2},
there exists $\de_+ \in \Bbb{R}^+$ such~that
$$\big|\ve_{i,i'}^-(w,\ups)\big|\le
\frac{\ep_-\de_-}{4|I({\mathcal{Z}}_i)|}\big|\tilde{\rho}_{i,i'}(\ups)\big|$$
for all $i' \in I({\mathcal{Z}}_i)$ and
$\ups \in \tilde{\mathcal{K}}_{\tilde{\al}_i,\bar{\nu}_i}$
such that $|\ups| \le \de_+$.
It can be assumed that $\de_+ \le \de_-$.
Then by (5) of Definition~\ref{reg_dfn1b},
\begin{equation}\label{euler_lmm4e9}
\big|\al_{i;-}\big((w + \tilde{\ve}_{i,i'}(w,\ups))\otimes
\tilde{\rho}_{i,i'}(\ups)\big)\big|
<\frac{\ep_-\de_-}{2|I({\mathcal{Z}}_i)|}\big|\tilde{\rho}_{i,i'}(\ups)\big|
\end{equation}
for all  $(w,\ups) \in E_{i;\de_i^*}
\times_{{\mathcal{Z}}_i}{\mathcal{F}}_{k_i;\de_i^*}$ where
$\ups \in \tilde{\mathcal{K}}_{\tilde{\al}_i,\bar{\nu}_i},~
|\ups| \le \de_+,~~
\phi_{k_i}\vt_{\ka_i;{\mathcal{F}}_{k_i}}(w,\ups) \in {\mathcal{S}}$.

(4)\qua Let $\de_i = \frac{1}{4}\de_+^2$.
If $\de_+ < \inf_{{\mathcal{K}}_{\tilde{\al}_i,\bar{\nu}_i}}\de_i^*$,
by (2) of Proposition~\ref{monomials_prp2}, there exists 
$\ep_i \in \Bbb{R}^+$ such that for all $t \in (0,\ep_i)$,
$${\mathcal{F}}_{k_i;\de_i^*}\cap
\big\{t^{-1}\eta\rho_i\big\}^{-1}({\mathcal{K}})
\subset {\mathcal{F}}_{k_i,\de_+},$$
\hbox{and}
$t^{-1}\eta\rho_i\co {\mathcal{F}}_{k_i;\de_i^*}\cap
\big\{t^{-1}\eta\rho_i\big\}^{-1}({\mathcal{K}})\lra{\mathcal{K}}$
is a covering projection of oriented degree $\deg\rho_i$.
Then, by~\eqref{euler_lmm4e3}--\eqref{euler_lmm4e9}
and the assumption on~$\de_-$,
\begin{multline*}
W_i(\de_i^*)\cap\big\{\eta s+t\nu'\big\}^{-1}(0)\\
\qquad\qquad\subset
\phi_{k_i}\vt_{\ka_i;{\mathcal{F}}_{k_i}}\big(
\big\{(b;w,\ups) \in E_i \oplus {\mathcal{F}}_{k_i} : 
b \in {\mathcal{K}}_{\tilde{\al}_i,\bar{\nu}_i},~
|w|,|\ups| \le\frac{1}{2}\de_i^*(b)\big\}\big)
\end{multline*}
for all $t \in (0,\ep_i)$,
if $\de_{\pm}\le\frac{1}{2}\de_i^*(b)$
for all $b \in {\mathcal{K}}_{\tilde{\al}_i,\bar{\nu}_i}$ and
$$\tilde{\ve}_+(b,\de_i^*(b)),\ve_{i,\nu}(b,\de_i^*(b)),
\tilde{C}_{i,W}(b)\ep_W
<\frac{1}{4}\de_{\bar{\nu}_i}(b)$$
for all $b \in {\mathcal{Z}}_i$.

(5)\qua By Definition~\ref{reg_dfn1a}, the set
$W_i(\de_i^*)\cap\big\{\eta s + t\nu'\big\}^{-1}(0)$
consists of the solutions of the system
\begin{equation*}
\left\{
\begin{array}{r}
\psi_i^-(w,\ups) = 0\in{\mathcal{O}}^-\\
\psi_i^+(w,\ups) = 0\in{\mathcal{O}}^+\\
\big\{\eta s + t\nu'\big\}(w,\ups)
 = 0\in V~~{}
\end{array}\right. \quad
(w,\ups) \in {\mathcal{U}}_i(\de_i^*).
\end{equation*} 
By our assumptions, the zeros of this system of equations are transverse.
By (1)--(4) and Definitions~\ref{reg_dfn1a} and~\ref{reg_dfn1b},
these zeros are the same as the zeros of the~system
\begin{equation}\label{euler_lmm4e25b}
\left\{   
\begin{array}{r}
\al_{i;-}\tilde{\rho}_i(w,\ups)+
\sum_{i'\in I({\mathcal{Z}}_i)}\big(
\tilde{\ve}^-_{i,i'}(w,\ups)\tilde{\rho}_{i,i'}(\ups)
+\ve^-_{i,i'}(w,\ups)\big) = 0\in{\mathcal{O}}^-\\
\be_{\de_i}(|\ups|)
\big(\al_{i;+} + \tilde{\ve}_{i;+}(w,\ups)\big)\rho_i(\ups)
 = 0\in{\mathcal{O}}^+\\
\be_{\de_i}(|\ups|)\big(\al_{i;V} + 
\tilde{\ve}_{i;V}(w,\ups)\big)\rho_i(\ups) + 
t\nu'(w,\ups) = 0\in V~~{}\\
\end{array}\right. 
\end{equation} 
where $\tilde{\ve}_{i;+}$ and $\tilde{\ve}_{i;V}$
denote the ${\mathcal{O}}^+$ and $V$--components of~$\tilde{\ve}_{i;+,V}$
and $(w,\ups) \in {\mathcal{U}}_i(\de_i^*)$ as before.
We will show that the zeros of~\eqref{euler_lmm4e25b}
are cobordant to the zeros of the~system
\begin{equation}\label{euler_lmm4e25c}
\left\{
\begin{array}{r}
\al_{i;-}\tilde{\rho}_i(b;w,\ups) = 0\in{\mathcal{O}}^-\\
\al_{i;V}\big(\be_{\de_i}(|\ups|)\rho_i(b;\ups)\big) = 0\in{\mathcal{O}}^+\\
\al_{i;+}\big(\be_{\de_i}(|\ups|)\rho_i(b;\ups)\big)+t\bar{\nu}_i(b)
 = 0\in V~~{}
\end{array}\right. \quad
(b;w,\ups) \in {\mathcal{U}}_i(\de_i^*).
\end{equation} 
We construct a cobordism between the two zero sets as follows.
If $h_-$, $h_+$, and $h_V$ are bundle maps 
from ${\mathcal{X}} \equiv [0,1] \times {\mathcal{U}}_i(\de_i^*)$
to ${\mathcal{O}}^-$, ${\mathcal{O}}^+$, and $V$, respectively, 
and $\under{h} = (h_-,h_+,h_V)$, let ${\mathcal{X}}(\under{h})$
be the set ot tuples $(\tau,(b;w,\ups)) \in {\mathcal{X}}$ such that
$$\begin{cases}
\be_{\de_i}(|\ups|)
\big(\al_{i;+} + \tau\tilde{\ve}_{i;+}(b;w,\ups)\big)\rho_i(b;\ups)
+h_+(\tau,b;w,\ups)=0,\\
\begin{aligned}
\be_{\de_i}(|\ups|)
\big(\al_{i;V} + \tau\tilde{\ve}_{i;V}(b;w,\ups)\big)\rho_i(b;\ups)
+t\big(\tau\nu'(b;w,\ups) + (1 - \tau)\bar{\nu}_i(b)\big) \qquad&\\
+h_V(\tau,b;w,\ups)=0, &
\end{aligned}\\
\begin{aligned}
\al_{i;-}\tilde{\rho}_i\big(b;w,\ups)
+\tau   \sum_{i'\in I({\mathcal{Z}}_i)}   \big(
\tilde{\ve}^-_{i,i'}(b;w,\ups)\tilde{\rho}_{i,i'}(b;\ups) 
+\ve^-_{i,i'}(b;w,\ups)\big) \qquad&\\
+h_-(\tau,b;w,\ups)=0.&
\end{aligned}\end{cases}$$
Since the zeros of \eqref{euler_lmm4e25b} and~\eqref{euler_lmm4e25c}
are transverse,
for a generic choice of $\under{h}$ with the boundary condition
$\under{h}_{\tau = 0,1} \equiv 0$,
${\mathcal{X}}(\under{h})$ is a smooth oriented manifold such that
\hbox{$\partial{\mathcal{X}}(\under{h}) = {\mathcal{X}}_1(\under{h}) -
 {\mathcal{X}}_0(\under{h})$}.
By the same argument as in (1)--(4),
\begin{multline*}
{\mathcal{X}}(\under{h})\subset [0,1] 
~ \times\big\{(b;w,\ups) \in E_i \oplus {\mathcal{F}}_{k_i} : 
(b;\ups) \in \tilde{\mathcal{K}}_{\tilde{\al}_i,\bar{\nu}_i},\\
\de_i \le |\ups| \le \frac{1}{2}\de_i^*(b),
~|w| < \frac{1}{2}\de_i^*(b)\big\},
\end{multline*}
if $\big|h_V(\tau,b;w,\ups)\big|, 
\big|h_+(\tau,b;w,\ups)\big|\le \frac{1}{3}t\de_{\bar{\nu}_i}(b)$
for all $b \in {\mathcal{Z}}_i$ and
$\big|h_-(\tau,b;w,\ups)\big|\le
\frac{\ep_-\de_-}{4|I({\mathcal{Z}}_i)|}\big|\tilde{\rho}_{i,i'}(\ups)\big|$
for all $\ups \in \tilde{\mathcal{K}}_{\tilde{\al}_i,\bar{\nu}_i}$
such that $\de_i \le |\ups| \le \de_i^*(b),~i' \in I({\mathcal{Z}}_i)$.
The lower bound on $|\ups|$ above follows from the fact that
$\be_{\de_i}(|\ups|)$ is zero if $|\ups| < \de_i$.
Thus, ${\mathcal{X}}(\under{h})$ is a compact space if $\under{h}$ is sufficiently
small.
We conclude that 
\begin{equation}\label{euler_lmm4e27}
^{\pm} \big|{\mathcal{X}}_1(\under{h})\big|
=^{\pm} \big|{\mathcal{X}}_0(\under{h})\big|
=\deg\rho_i\cdot ^{\pm}\big|\psi_{i,\bar{\nu}}^{-1}(0)\big|.
\end{equation}
The second equality is immediate from (A1), (A2), (A3),
and (2) of Proposition~\ref{monomials_prp2}.
This concludes the proof of the main claim of the lemma.
The other claim is clear.
\end{proof}

Proposition~\ref{euler_prp} is essentially proved;
it only remains to construct a cutoff function $\eta$
that has good properties.
Let $W_1 = \bigsqcup_{i\in I_s}W_i(\de_i^*/2)$.
This is an open neighborhood of $\partial{\overS}$ in~${\overM}$.
By Lemmas~\ref{euler_lmm1}--\ref{euler_lmm4}, it can be assumed
that $s^{-1}(0)\cap\bar{W}_1 = \eset$ and that
$$\ov{\phi_{k_i}\big(\vt_{\ka_i;{\mathcal{F}}_{k_i}}
(E_{i;\de_i^*} 
  \times _{{\mathcal{K}}_{\tilde{\al}_i,\bar{\nu}_i}}{\mathcal{F}}_{k_i;\de_i^*})\big)}
\cap
\ov{\phi_{k_j}\big(\vt_{\ka_j;{\mathcal{F}}_{k_j}}
(E_{j;\de_j^*} 
  \times _{{\mathcal{K}}_{\tilde{\al}_j,\bar{\nu}_j}}{\mathcal{F}}_{k_j;\de_j^*})\big)}
=\eset$$
for all $i,j \in I_s^*$ such that $i \neq j$.
Let $W_0$ be an open neighborhood of $\partial{\overS}$ in~${\overM}$
such that $\bar{W}_0 \subset W_1$
and let  $\eta'\co {\mathcal{M}} \lra [0,1]$ be a smooth
function such that $\eta'|W_0 = 0$ and $\eta'|{\mathcal{M}} - W_1 = 1$.
For each $i \in I_s^*$, let 
${\mathcal{K}}_{\tilde{\al}_i,\bar{\nu}_i}'$ be a compact subset of~${\mathcal{Z}}_i$
such~that
${\mathcal{K}}_{\tilde{\al}_i,\bar{\nu}_i}
\subset\hbox{Int}_{{\mathcal{Z}}_i}{\mathcal{K}}_{\tilde{\al}_i,\bar{\nu}_i}'$ and
$$\ov{\phi_{k_i}\big(\vt_{\ka_i;{\mathcal{F}}_{k_i}}
(E_{i;\de_i^*} 
  \times _{{\mathcal{K}}_{\tilde{\al}_i,\bar{\nu}_i}'}{\mathcal{F}}_{k_i;\de_i^*})\big)}
\cap
\ov{\phi_{k_j}\big(\vt_{\ka_j;{\mathcal{F}}_{k_j}}
(E_{j;\de_j^*} 
  \times _{{\mathcal{K}}_{\tilde{\al}_j,\bar{\nu}_j}'}{\mathcal{F}}_{k_j;\de_j^*})\big)}
=\eset$$
for all $i,j \in I_s^*$ such that $i \neq j$.
Choose a smooth function 
$\eta_i'\co E_i \oplus {\mathcal{F}}_{k_i} \lra [0,1]$
such that  
$$\eta_i'\big|E_{i;\de_i^*/2} 
  \times _{{\mathcal{K}}_{\tilde{\al}_i,\bar{\nu}_i}}{\mathcal{F}}_{k_i;\de_i^*/2}=1
\quad\hbox{and}\quad
\eta_i'\big|
\big(E_i \oplus {\mathcal{F}}_{k_i}-
E_{i;\de_i^*} 
  \times _{{\mathcal{K}}_{\tilde{\al}_i,\bar{\nu}_i}'}{\mathcal{F}}_{k_i;\de_i^*}
                                        \big)=0.$$
If $\de_i \in \Bbb{R}^+$ is as in Lemma~\ref{euler_lmm4},
we define $\eta_i\co {\mathcal{S}} \lra [0,1]$ by
$$\eta_i(b^*)=\eta_i'(b;w,\ups)\be_{\de_i}(|\ups|)
+\big(1 - \eta_i'(b;w,\ups)\big)\eta'(b^*)$$
if $b^*=\phi_{k_i}\vt_{\ka_i;{\mathcal{F}}_{k_i}}(b;w,\ups)$
for some $(b;w,\ups) \in E_{i;\de_i^*} 
\times_{{\mathcal{K}}_{\tilde{\al}_i,\bar{\nu}_i}'}{\mathcal{F}}_{k_i;\de_i^*}$,
and
$$\eta_i(b^*)=\eta'(b^*)$$
otherwise.
This function is smooth on ${\mathcal{S}}$, since it is the restriction
of a smooth function on~${\mathcal{M}}$.
Let
$$\eta=\eta'+\sum_{i\in I_s^*}\eta_i.$$
Since $\eta$ vanishes on a neighborhood $W$ of $\partial{\overS}$
in ${\overM}$, the section $\tilde{s} \equiv \eta s$
of $V$ over ${\mathcal{S}}$ extends by zero to a continuous section
over~${\overS}$.
Since $\nu$ does not vanish on ${\overS}$,
we can assume that  $\nu$ does not vanish on~$\bar{W}$ as~well.
Furthermore, $\eta$ satisfies the requirements of Lemma~\ref{euler_lmm4}
with $\nu$ replaced by $(1 - \eta')\nu$ and 
$\de_i^*$ replaced \hbox{by $\de_i^*/2$}.
Since $\eta$ does not vanish on the complement of~$W_1$,
it follows that $\eta s + t(1 - \eta')\nu$ is transverse to the zero set
on $({\mathcal{S}} - W_1)\cup\ov{W}$.
Thus, for all $t,\ep \in \Bbb{R}^+$, there exists 
$\nu'\in\Ga({\overS};V)$ such that
$$\nu'\big|_{{\mathcal{S}} - W_1} = (1 - \eta')\nu\big|_{{\mathcal{S}} - W_1},  
\quad \nu'\big|_W = (1 - \eta')\nu\big|_W,  \quad
\|\nu' - (1 - \eta')\nu\|_{C^2({\mathcal{S}}-W)}<\ep,$$
and $\eta s + t\nu'$ is transversal to the zero~set on~${\mathcal{S}}$.
Since $\eta s + t\nu'$ does not vanish on~$\partial{\overS}$
and is a positive multiple of $s$ outside of~$W_1$, 
\begin{equation}\label{euler_prp_e1}\begin{split}
\blr{e(V),{\overS}}
&=\, ^{\pm} \big|\big\{\eta s + t\nu'\big\}^{-1}(0)\big\}\big| \\
&=\, ^{\pm} \big|s^{-1}(0)\big| + 
\sum_{i\in I_s}
{^{\pm} \big|\big\{\eta s + t\nu'\big\}^{-1}(0)\cap W_i(\de_i^*/2)\big|}.
\end{split}\end{equation}
On the other hand, by Lemmas~\ref{euler_lmm2}--\ref{euler_lmm4},
\begin{equation}\label{euler_prp_e2}
\sum_{i\in I_s}
{^{\pm} \big|\big\{\eta s + t\nu'\big\}^{-1}(0)\cap W_i(\de_i^*/2)\big|}
=\sum_{i\in I_s^*}
\deg\rho_i\cdot ^{\pm}\big|\psi_{i,\bar{\nu}}^{-1}(0)\big|,
\end{equation}
provided $t$ and $\ep$ are sufficiently small.
Proposition~\ref{euler_prp} follows immediately  from
equations~\eqref{euler_prp_e1} and~\eqref{euler_prp_e2},
assumptions (A1) and (A3), and 
the second remark after the statement of the proposition.

\section{Spaces of stable maps}
\label{stable_maps_sec}

\subsection{Notation}
\label{notation_subs}

In this subsection, we describe our notation for spaces of tuples of stable rational 
maps and for important vector bundles over~them.
The zeros of certain sections of these bundles can be identified 
with rational curves with prescribed singularities.
In many cases, using the topological method of Section~\ref{topology_sec}
and the analytic estimates of Subsection~\ref{analysis_subs},
one can express the number of such zeros in terms of
intersection numbers of certain tautological cohomology classes,
that are also defined in this subsection.
The notation described below is a generalization on that of
\cite[Subsection~1.3]{Z1} and
\cite[Section~2]{Z2}. Thus, we omit some details.

\begin{dfn}
\label{comb_dfn1}
{\rm A finite partially ordered set $I$ is a \emph{linearly ordered set}
if for all \hbox{$i_1,i_2,h \in  I$} such that $i_1,i_2 < h$, 
either $i_1 \le  i_2$ \hbox{or $i_2 \le  i_1$.}}
\end{dfn}

The term \textit{linearly ordered set} is sometimes used with
a different meaning in combinatorics, as the referee pointed out.
We continue to use this term with the meaning of Definition~\ref{comb_dfn1}
to be consistent with earlier papers.

If $I$ is a linearly ordered set, let $\hat{I}$ be
the subset of the non-minimal elements of~$I$.
For every $h \in \hat{I}$,  denote by $\io_h \in I$
the largest element of $I$ which is smaller than~$h$, ie
\begin{equation}\label{notation_e1}
\io_h=\max\big\{i \in I:i < h\big\}.
\end{equation}

\begin{dfn}
\label{comb_dfn2}
{\rm A linearly ordered set $I$ is \emph{graded} by $I = I^-\sqcup
I^+$ if $I - \hat{I} \subset I^-$, 
$\io_h \in I^{\mp}$ for all \hbox{$h \in I^{\pm}\cap\hat{I}$}, and
for every $i \in I^-$ there exists $h \in I^+$ such that
$\io_h = i$.}
\end{dfn}

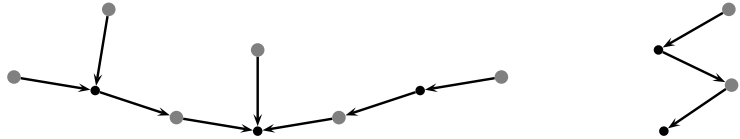
\begin{figure}[ht!]\small
\begin{pspicture}(-1.2,-1.5)(10,1)
\psset{unit=.36cm}
\pscircle[linecolor=black,fillstyle=solid,fillcolor=black](10,-3){.175}
\pscircle[linecolor=gray,fillstyle=solid,fillcolor=gray](7,-2.5){.25}
\psline[linewidth=.09]{->}(7.25,-2.543)(9.825,-2.97)
\pscircle[linecolor=black,fillstyle=solid,fillcolor=black](4,-1.5){.175}
\psline[linewidth=.09]{->}(4.17,-1.56)(6.76,-2.416)
\pscircle[linecolor=gray,fillstyle=solid,fillcolor=gray](1,-1){.25}
\psline[linewidth=.09]{->}(1.25,-1.043)(3.825,-1.47)
\pscircle[linecolor=gray,fillstyle=solid,fillcolor=gray](4.5,1.5){.25}
\psline[linewidth=.09]{->}(4.455,1.25)(4.03,-1.33)
\pscircle[linecolor=gray,fillstyle=solid,fillcolor=gray](13,-2.5){.25}
\psline[linewidth=.09]{->}(12.75,-2.543)(10.175,-2.97)
\pscircle[linecolor=black,fillstyle=solid,fillcolor=black](16,-1.5){.175}
\psline[linewidth=.09]{->}(15.83,-1.56)(13.24,-2.416)
\pscircle[linecolor=gray,fillstyle=solid,fillcolor=gray](19,-1){.25}
\psline[linewidth=.09]{->}(18.75,-1.043)(16.175,-1.47)
\pscircle[linecolor=gray,fillstyle=solid,fillcolor=gray](10,0){.25}
\psline[linewidth=.09]{->}(10,-.252)(10,-2.822)
\pscircle[linecolor=black,fillstyle=solid,fillcolor=black](25,-3){.175}
\pscircle[linecolor=gray,fillstyle=solid,fillcolor=gray](27.5,-1.3){.25}
\psline[linewidth=.09]{->}(27.29,-1.442)(25.147,-2.9)
\pscircle[linecolor=black,fillstyle=solid,fillcolor=black](24.8,0){.175}
\psline[linewidth=.09]{->}(24.96,-.079)(27.27,-1.19)
\pscircle[linecolor=gray,fillstyle=solid,fillcolor=gray](27.4,1.5){.25}
\psline[linewidth=.09]{->}(27.18,1.37)(24.955,.09)
\end{pspicture}
\caption{A graded linearly ordered set}
\label{glos_fig}
\end{figure}

A graded linearly ordered set can be represented by an oriented graph.
In Figure~\ref{glos_fig}, 
the small black and large gray dots denote the elements of $I^-$
and~$I^+$, respectively.
The arrows specify the partial ordering of the linearly ordered set~$I$.
We use graded linearly ordered sets
to encode the structure of an $I^+$--tuple of stable maps.
The set $I^-$ will describe the nodes of the domain of the map.
For example, the domain of an $I^+$--tuple of stable maps
with the structure depicted by Figure~\ref{glos_fig} 
will have two connected components.
One of these components will consist of six irreducible components
and contain a double point and two triple points.
We give more details below.

If $I = I^-\sqcup I^+$ is a graded linearly ordered set, 
$I^-$ and $I^+$ are linearly ordered sets.
We denote by $\hat{I}^-$ and $\hat{I}^+$  
the subsets of the non-minimal elements of $I^-$ and~$I^+$, respectively.
If $h \in \hat{I}^{\pm}$, we define $\io_h^{\pm}$
as in~\eqref{notation_e1}, but with $I$ replaced by~$\hat{I}^{\pm}$.
If $h_1,h_2 \in I$, let
\begin{gather*}
[h_1,h_2]=\big\{i \in I^+ :h_1 \le i \le h_2\big\},\qquad
[h_1,h_2)=\big\{i \in I^+ :h_1 \le i < h_2\big\},\\
(h_1,h_2]=\big\{i \in I^+ :h_1 < i \le h_2\big\},\qquad
(h_1,h_2)=\big\{i \in I^+ :h_1 < i < h_2\big\}.\\
\end{gather*}
If $I = I^-\sqcup I^+$ is a graded linearly ordered set and 
$i^*$ is an element of~$I^+$, we define a new graded linearly ordered set
$$I(i^*)=I^-(i^*)\sqcup I^+(i^*)$$
as follows.
We take $I^{\pm}(i^*) = I^{\pm}\sqcup\{i^*_{\pm}\}$,
where $i^*_+$ and $i^*_-$ are new elements.
We define a partial ordering $\succ$ on the set $I^{\pm}(i^*)$ by
\begin{gather*}
h \succ i\hbox{~~if~~}h,i \in I\hbox{~and~} h > i; \qquad
i\succ i_{\pm}^*\hbox{~~if~~}i \in I\hbox{~and~}i > i^*;\\
i_{\pm}^*\succ i \hbox{~~if~~}i \in I\hbox{~and~}i^* \ge i;\quad
i_+^*\succ i_-^*.
\end{gather*}
It is easy to see that $I(i^*)$ is indeed a graded linearly ordered set.

We denote the south pole of the 2-sphere $S^2 \subset \Bbb{R}^3$
by $\i$ and identify $S^2 - \{\i\}$ with $\Bbb{C}$ via
the standard stereographic projection~$q_N$ mapping the origin in $\Bbb{C}$
to the north pole of~$S^2$.
If $M$ is a finite set, 
a \textit{$\P$--valued bubble map with $M$--marked points} is a tuple
$b=\big(M,I;x,(j,y),u\big)$,
where $I = I^-\sqcup I^+$ is a graded linearly ordered set, and
\begin{gather*}
x\co \hat{I}^- \lra S^2 - \{\i\}, \quad
j\co M \lra I^+, \quad
y\co M \lra S^2 - \{\i\},\\
\hbox{and}\qquad  u\co I^+ \lra C^{\i}(S^2;\P)
\end{gather*}
are maps such that 
$$\big(\io_{h_1},x_{h_1}\big) \neq \big(\io_{h_2},x_{h_2}),\quad
\big(j_{l_1},y_{l_1}\big) \neq \big(j_{l_2},y_{l_2}\big),\quad
\big(\io_h,x_h\big) \neq \big(j_l,y_l\big)$$
for all $h,h_1,h_2 \in \hat{I}^-$ and $l,l_1,l_2 \in M$;
$$u_{h_1}(\i) = u_{h_2}(\i)$$
for all $h_1,h_2 \in I^+$
such that $\io_{h_1} = \io_{h_2}$; and
$$u_h(\i) = u_{\io_h^+}(x_{\io_h})$$
for all $h \in \hat{I}^+$.
We associate such a tuple with Riemann surface
$$\Si_b=\Big(\bigsqcup_{h\in I^+}\Si_{b,h}\Big)\Big/ \sim,$$
where
$\Si_{b,h}=\{h\} \times S^2$, and
$(h_1,\i) \sim (h_2,\i)$ for all $h_1,h_2 \in I^+$ such that
$\io_{h_1} = \io_{h_2}$, and
$(h,\i) \sim \big(\io_h^+,x_{\io_h}\big)$ for all $h \in \hat{I}^+$,
with marked points $(j_l,y_l) \in \Si_{b,j_l}$,
and continuous map $u_b\co \Si_b \lra \P$,
given by $u_b|\Si_{b,h} = u_h$ for \hbox{all $h \in I^+$}.
We require that $\Si_{b,h}$ contain at least two singular and/or
marked points of~$\Si_b$ other than~$(h,\i)$
if $u_{h*}[S^2] = 0 \in  H_2(\P;\Bbb{Z})$.
In addition, we implicitly consider each point $(h,\i)$
to be a special marked point.
Figure~\ref{glos_fig} is basically the dual graph of~$\Si_b$.
The black dots simply specify which of the special marked points 
are identified and thus are mapped to the same point in~$\P$.
If 
$$b_1=\big(M,I_1;x^{(1)},(j^{(1)},y^{(1)}),u^{(1)}\big) 
\quad\hbox{and}\quad
b_2=\big(M,I_2;x^{(2)},(j^{(2)},y^{(2)}),u^{(2)}\big)$$
are two bubble maps and $\tilde{I}$ is a subset of $I_1^+$ and of~$I_2^+$, 
we say $b_1$ and $b_2$ are {\it$\tilde{I}$--equivalent} if there
exists a homeomorphism $\phi\co \Si_{b_1} \lra \Si_{b_2}$
such that $\phi|_{\Si_{b_1,i}}$ is holomorphic for all $i \in I^+$,
\hbox{$\phi(\Si_{b_1,i}) \subset \Si_{b_2,i}$} for all $i \in \tilde{I}$,
$u_{b_1} = u_{b_2} \circ\phi$,
$$\phi(j^{(1)}_l,y^{(1)}_l)=(j^{(2)}_l,y^{(2)}_l) \qquad\mbox{for all}\ l \in M,$$
and for every $i \in I_1^+$ there exists $i' \in I_2^+$ 
such that $\phi(i,\i) = (i',\i)$.

The general structure of bubble maps is described
by tuples ${\mathcal{T}} = (M,I;j,\under{d})$,
with \hbox{$d_i \in \Bbb{Z}$} 
specifying the degree of the map $u_b$ on~$\Si_{b,i}$.
The above equivalence relation on the set of bubble maps
induces an equivalence relation on the set of bubble types.
If $\tilde{I} \subset I^+$, we denote by 
${\mathcal{A}}({\mathcal{T}}|\tilde{I})$ the group of $\tilde{I}$--automorphisms
of~${\mathcal{T}}$; let ${\mathcal{A}}({\mathcal{T}}) = {\mathcal{A}}({\mathcal{T}}|\eset)$.
For each $i \in I^+$, let
$$H_i{\mathcal{T}}=\{h \in \hat{I}^- :\io_h = i\}
\qquad\hbox{and}\qquad
M_i{\mathcal{T}}=\{l \in  M :j_l = i\}.$$
If $i \in I^-$, we put
$$H_i{\mathcal{T}}=\big\{h \in I^+  :\io_h = i\big\}\cup
\begin{cases}
\{i\},&\hbox{if}~i \in \hat{I}^- ;\\
\eset,&\hbox{if}~i \not\in \hat{I}^- .
\end{cases}$$
Let ${\mathcal{H}}_{\mathcal{T}}$ denote the space of all holomorphic bubble maps
with structure~${\mathcal{T}}$.
This is a smooth complex manifold; see \cite[Chapter~3]{MS}, for example.

We denote by ${\mathcal{U}}_{{\mathcal{T}}|\tilde{I}}$ the set of
$\tilde{I}$--equivalence classes of bubble maps in~${\mathcal{H}}_{\mathcal{T}}$.
Then there exists a smooth submanifold 
${\mathcal{B}}_{\mathcal{T}}$ of ${\mathcal{H}}_{\mathcal{T}}$
such that ${\mathcal{U}}_{{\mathcal{T}}|\tilde{I}}$ is the quotient of 
${\mathcal{B}}_{\mathcal{T}}$ by a natural action of the~group
$$G_{{\mathcal{T}}|\tilde{I}}\equiv
{\mathcal{A}}({\mathcal{T}}|\tilde{I})\ltimes G_{\mathcal{T}},
\qquad\hbox{where}\qquad
G_{\mathcal{T}}=\big(S^1\big)^{I^+}.$$
For any $i \in \tilde{I}$, denote by ${\mathcal{U}}_{{\mathcal{T}}|\tilde{I}}^{(i)}$
the quotient of  ${\mathcal{B}}_{\mathcal{T}}$ by the group
$$G_{{\mathcal{T}}|\tilde{I}}^{(i)}\equiv
{\mathcal{A}}({\mathcal{T}}|\tilde{I})\ltimes G_{\mathcal{T}}^{(i)},
\qquad\hbox{where}\qquad
G_{\mathcal{T}}^{(i)}=\big(S^1\big)^{I^+-\{i\}}.$$
Then, ${\mathcal{U}}_{{\mathcal{T}}|\tilde{I}}$ is the quotient of 
${\mathcal{U}}_{{\mathcal{T}}|\tilde{I}}^{(i)}$
by the residual $S^1$--action.
If $i \in I^+$ is fixed by every element of 
the group ${\mathcal{A}}({\mathcal{T}}|\tilde{I})$, corresponding to the first quotient
we obtain a line orbi-bundle 
$L_i{\mathcal{T}}  \lra {\mathcal{U}}_{{\mathcal{T}}|\tilde{I}}$.
In general, the direct sum of the line bundles $L_i{\mathcal{T}}$
taken over all elements of the orbit ${\mathcal{A}}({\mathcal{T}}|\tilde{I}) \cdot i$
is well-defined.
If $h \in \hat{I}^+$, let 
\hbox{${\mathcal{F}}_h{\mathcal{T}} = L_{\io_h^+}^*{\mathcal{T}} \otimes L_h{\mathcal{T}}$}.

The Gromov-convergence topology on the space of all holomorphic maps
induces a partial ordering on the set of bubble types
and their equivalence classes such that the~spaces
$${\overU}_{\tilde{\mathcal{T}}}^{(i)}\equiv
{\overU}_{\tilde{\mathcal{T}}|\tilde{I}^+}^{(i)}=
\bigcup_{{\mathcal{T}}\le\tilde{\mathcal{T}}}{\mathcal{U}}_{{\mathcal{T}}|\tilde{I}^+}^{(i)}
\qquad\hbox{and}\qquad
{\overU}_{\tilde{\mathcal{T}}}\equiv
{\overU}_{\tilde{\mathcal{T}}|\tilde{I}^+}=
\bigcup_{{\mathcal{T}}\le\tilde{\mathcal{T}}}{\mathcal{U}}_{{\mathcal{T}}|\tilde{I}^+}$$
are compact and Hausdorff.
Here $\tilde{\mathcal{T}} \equiv (M,\tilde{I};\tilde{j},\tilde{\under{d}})$,
and the unions are disjoint if taken over \hbox{$\tilde{I}^+$--equivalence}
classes of bubble~types.
If ${\mathcal{T}} \le \tilde{\mathcal{T}}$, let
$${\overU}_{{\mathcal{T}}|\tilde{\mathcal{T}}}^{(i)}
={\overU}_{{\mathcal{T}}|\tilde{I}^+}^{(i)}
\qquad\hbox{and}\qquad
{\overU}_{{\mathcal{T}}|\tilde{\mathcal{T}}}
={\overU}_{{\mathcal{T}}|\tilde{I}^+}.$$
The residual $S^1$--action on ${\mathcal{U}}_{\tilde{\mathcal{T}}}^{(i)}$ extends
to an action on ${\overU}_{\tilde{\mathcal{T}}}^{(i)}$,
and thus the line orbi-bundle 
$L_i{\tilde{\mathcal{T}}} \lra {\mathcal{U}}_{\tilde{\mathcal{T}}}$
extends over~${\overU}_{\mathcal{T}}$ as the line bundle~$L_i{\mathcal{T}}$.

For each $l \in M$, $h \in I^+$,  and $i \in I^-$, 
we define evaluation maps 
\hbox{${\mathcal{H}}_{\mathcal{T}} \lra \P$ by}
$$\ev_l\big((M,I;x,(j,y),u)\big)=u_{j_l}(y_l),\qquad
\ev_h\big((M,I;x,(j,y),u)\big)=u_h(\i),$$
and $\ev_i = \ev_h$ if $h \in I^+$ and $\io_h = i$.
These maps descend to all the quotients defined above and induce continuous
maps on~${\overU}_{\mathcal{T}}$.
If $\tilde{M} \subset M\sqcup I^-$ and
$\mu = \mu_{\tilde{M}}$ is an $\tilde{M}$--tuple 
of submanifolds of~$\P$, let
$${\mathcal{U}}_{\mathcal{T}}(\mu)=
\big\{b \in {\mathcal{U}}_{\mathcal{T}} :\ev_l(b) \in \mu_l~\mbox{for
all}\  l \in \tilde{M}\big\},$$
and define the spaces ${\overU}_{\mathcal{T}}(\mu)$ and
${\mathcal{U}}_{{\mathcal{T}}|\tilde{\mathcal{T}}}(\mu)$ analogously.

If ${\mathcal{T}} = (M,I;j,\under{d})$ is a bubble type and $k \in I^+$, 
we define the bubble type
${\mathcal{T}}_k \equiv (M_k,I_k;j^{(k)},\under{d}^{(k)})$ by
\begin{gather*}
M_k=M_k{\mathcal{T}}\sqcup H_k{\mathcal{T}};\quad
I_k = \{k,\io_k\}\subset I;\quad
j^{(k)}_l = k~~~\mbox{for all}\  l \in M_k;\quad
d^{(k)}_k = d_k. 
\end{gather*}
Let ${\mathcal{U}}_{{\mathcal{T}},{\mathcal{T}}} = \prod_{k\in I^+}
{\mathcal{U}}_{{\mathcal{T}}_k}$ and 
${\overU}_{{\mathcal{T}},{\mathcal{T}}} = \prod_{k\in I^+}
{\overU}_{{\mathcal{T}}_k}$.
Note that the spaces ${\mathcal{U}}_{\mathcal{T}}$ and ${\overU}_{\mathcal{T}}$ 
are contained in ${\mathcal{U}}_{{\mathcal{T}},{\mathcal{T}}}$
and ${\overU}_{{\mathcal{T}},{\mathcal{T}}}$, respectively.

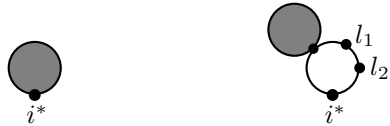
\begin{figure}[ht!]\small
\begin{pspicture}(0.0,-1.5)(10,.5)
\psset{unit=.36cm}
\pscircle[fillstyle=solid,fillcolor=gray](12,-2.5){1}
\pscircle*(12,-3.5){.22}\rput(12.1,-4.2){$i^*$}
\pscircle(23,-2.5){1}
\pscircle*(23,-3.5){.22}\rput(23.1,-4.2){$i^*$}
\pscircle*(24,-2.5){.2}\rput(24.75,-2.5){$l_2$}
\pscircle*(23.5,-1.63){.2}\rput(24.2,-1.3){$l_1$}
\pscircle[fillstyle=solid,fillcolor=gray](21.59,-1.09){1}
\pscircle*(22.29,-1.79){.19}
\end{pspicture}
\caption{The domains of elements of ${\mathcal{U}}_{\mathcal{T}}$
and ${\mathcal{U}}_{{\mathcal{T}}(M_0)}$}
\label{psi_fig}
\end{figure}

Suppose ${\mathcal{T}} = (M,I;j,\under{d})$ is a bubble type, 
$i^*$ is an element of $I^+$ such that $d_{i^*} \neq 0$ and
$M_0$ is nonempty subset of~$M_{i^*}{\mathcal{T}}$.
We define bubble type 
$${\mathcal{T}}(M_0) \equiv (M,I(i^*);j',\under{d}')$$  
by
$$j'_l=\begin{cases}
i^*,&\hbox{if~}l \in  M_0;\\
i^*_+,&\hbox{if~}l \in  M_{i^*}{\mathcal{T}} - M_0;\\
j_l,&\hbox{otherwise};
\end{cases}\qquad
d'_i=\begin{cases}
0,&\hbox{if~}i = i^*;\\
d_{i^*},&\hbox{if~}i = i^*_+;\\
d_i,&\hbox{otherwise}.
\end{cases}$$
If $l \in M$, we will write ${\mathcal{T}}(l)$ for ${\mathcal{T}}(\{l\})$.
In Figure~\ref{psi_fig}, we show the domain of an element
of the space ${\mathcal{U}}_{\mathcal{T}}$, where $I = \{i^*\}$
is a single-element set, and 
the domain of an element of the space ${\mathcal{U}}_{{\mathcal{T}}(M_0)}$, 
where $M_0 = \{l_1,l_2\}$  is a two-element set.

In Figure~\ref{psi_fig}, as well as in later figures,
we denote each component of the domain by a disk and shade
the component(s) on which the map into~$\P$ is nonconstant.
We indicate marked points on the ghost components, 
ie the components on which the map is constant,
by putting small dots on the boundary of the corresponding disk.
The point labeled by~$i^*$, ie the same way as the component,
is the special marked point~$(i^*,\i)$.
Lemma~\ref{str_lmm} and Proposition~\ref{str_prp},
as well as the decomposition~\eqref{cart_split},
show that it is crucial to clearly distinguish between 
ghost and non-ghost components.

If $[N] \equiv \{1,\ldots,N\}$ is a subset of $M$ such that the set
$M - [N]$ contains no positive integers, we~put
\begin{equation}\label{psi_class1}
c_1({\mathcal{L}}_{i^*}^*{\mathcal{T}})\equiv c_1(L_{i^*}^*{\mathcal{T}})-
 \sum_{\eset\neq M_0\subset M_{i^*}{\mathcal{T}}\cap[N]}         
PD_{{\overU}_{{\mathcal{T}},{\mathcal{T}}}}
\big[{\overU}_{{\mathcal{T}}(M_0),{\mathcal{T}}(M_0)}\big]
\in H^2\big({\overU}_{{\mathcal{T}},{\mathcal{T}}}\big).
\end{equation} 
By Proposition~\ref{str_prp}, ${\overU}_{{\mathcal{T}},{\mathcal{T}}}$ 
is an ms-orbifold, while
${\overU}_{{\mathcal{T}}(M_0),{\mathcal{T}}(M_0)}$ is an ms-suborbifold 
of~${\overU}_{{\mathcal{T}},{\mathcal{T}}}$.
Thus, the cohomology class on the left-hand side of~\eqref{psi_class1} 
is well-defined.
We illustrate definition~\eqref{psi_class1} in Figure~\ref{psi_fig2}
in the case $I^+ = \{i^*\}$ is a single-element set.
In this figure, as well in the future ones,
we denote spaces of tuples of stable maps by
drawing a picture of the domain of a typical element of such a~space.
On the other hand, let 
$$\pi_i\co  {\overU}_{{\mathcal{T}},{\mathcal{T}}}\lra
{\overU}_{{\mathcal{T}}_i'},
\quad\hbox{where}\quad
{\mathcal{T}}_i'=\big(M_i - [N],I_i;j^{(i)},\under{d}^{(i)}\big),$$
be the composition of the projection onto the $i$th factor
with the appropriate forgetful map.
By \cite[Lemma~2.2.2]{P},
$$c_1({\mathcal{L}}_{i^*}^*{\mathcal{T}})=\pi_{i^*}^*\psi_{i^*},$$
where $\psi_{i^*}$ is the first chern class of 
the universal cotangent line bundle at the marked point $(i^*,\i)$
over the moduli space
$$\ov{\frak M}_{0,(M_i - [N])\sqcup\{i^*\}}(d_{i^*},\P)$$
of stable rational degree--$d_{i^*}$ maps into $\P$
with marked points labeled by the set $(M_i - [N])\sqcup{(i^*,\i)}$.
In particular, $c_1({\mathcal{L}}_{i^*}^*{\mathcal{T}})$ is 
the first chern class of a line orbi-bundle 
over~${\overU}_{{\mathcal{T}},{\mathcal{T}}}$.
Whenever the bubble type ${\mathcal{T}}$ is clear from context,
we will write $c_1(L_i^*)$ and $c_1({\mathcal{L}}_i^*)$ for
$c_1(L_i^*{\mathcal{T}})$ and $c_1({\mathcal{L}}_i^*{\mathcal{T}})$,
respectively.
If $\tilde{M}$ is a subset of $M\sqcup I^-$ that contains~$[N]$
and $\mu$ is an $\tilde{M}$--tuple of constraints in $\P$
such that $\mu_{l_1} \cap \mu_{l_2} = \eset$
for all distinct elements $l_1,l_2$ of~${N}$,
\begin{equation}\label{psi_class2}\begin{split}
\big[{\overU}_{{\mathcal{T}}(l)}(\mu)\big] 
                        \cap c_1({\mathcal{L}}_{i^*}^*{\mathcal{T}})
&= \big[{\overU}_{{\mathcal{T}}(l)}(\mu)\big]
                        \cap c_1(L_{i^*_+}^*{\mathcal{T}}(l)) \\
&= \big[{\overU}_{{\mathcal{T}}(l)}(\mu)\big]
                    \cap c_1({\mathcal{L}}_{i^*_+}^*{\mathcal{T}}(l))
\quad\mbox{for all}\  l \in M_{i^*}{\mathcal{T}}.
\end{split}\end{equation} 
Note that by Lemma~\ref{str_lmm},
${\mathcal{U}}_{{\mathcal{T}}(l)}(\mu)$ is a pseudocycle in
${\overU}_{{\mathcal{T}},{\mathcal{T}}}$
and thus induces a homology class.
The first equality in~\eqref{psi_class2} can be deduced from
\cite[Subsection~3.2]{P}.

\begin{figure}[ht!]\small
\begin{pspicture}(-0.3,-1.5)(10,.5)
\psset{unit=.36cm}
\rput(4,-2.5){$c_1({\mathcal{L}}_{i^*}^*)\cap~$}
\pscircle[fillstyle=solid,fillcolor=gray](7,-2.5){1}
\pscircle*(7,-3.5){.22}\rput(7.1,-4.2){$i^*$}
\rput(10,-2.5){$=$}
\rput(14,-2.5){$c_1(L_{i^*}^*)\cap~$}
\pscircle[fillstyle=solid,fillcolor=gray](17,-2.5){1}
\pscircle*(17,-3.5){.22}\rput(17.1,-4.2){$i^*$}
\rput(20,-2.5){$-$}
\rput(23,-3){$~~\sum\limits_{\eset\neq M_0\subset[N]}$}
\pscircle(27.5,-2.5){1}
\pscircle*(27.5,-3.5){.22}\rput(27.6,-4.2){$i^*$}
\pscircle[fillstyle=solid,fillcolor=gray](26.09,-1.09){1}
\pscircle*(26.79,-1.79){.19}
\pscircle*(28.5,-2.5){.2}\pscircle*(28,-1.63){.2}
\psline[linewidth=.05]{<-}(28.3,-1.63)(30,-1.63)
\psline[linewidth=.05]{<-}(28.8,-2.47)(30,-1.7)
\rput(30.7,-1.68){$M_0$}
\end{pspicture}
\caption{An example of Definition~\eqref{psi_class1}}
\label{psi_fig2}
\end{figure}

We are now ready to formally explain the notation involved 
in the statement of Theorem~\ref{cusps_thm}.
Let $n$, $d$, $N$ be positive integers and let $\mu$ be an $N$--tuple
of constraints in~$\P$.
If $k \ge 1$, denote by ${\overV}_k(\mu)$ the quotient of 
the disjoint union of the spaces ${\overU}_{\mathcal{T}}(\mu)$
taken over all bubble types
\hbox{${\mathcal{T}} = ([N],I_k;j,\under{d})$} 
such that
$$I_k^+=\{\tilde{1},\ldots,\tilde{k}\}, \qquad \sum d_i = d,$$ 
and $I_k^- = \{\hat{0}\}$ is a one-element set,
by the natural action of the symmetric group~$S_k$.
We define the spaces ${\mathcal{V}}_k(\mu)$ similarly.
Denote by 
$$\eta_{\hat{0},l},\tilde{\eta}_{\hat{0},l}
\in H^{2l}\big({\overV}_k(\mu)\big)$$ 
the cohomology classes such that 
$\pi^*\eta_{\hat{0},l}$  and $\pi^*\tilde{\eta}_{\hat{0},l}$
are the sum of all degree--$l$ monomials in 
$$\big\{c_1({\mathcal{L}}_{\tilde{1}}^*),\ldots,
c_1({\mathcal{L}}_{\tilde{k}}^*)\big\}
\qquad\hbox{and}\qquad
\big\{c_1(L_{\tilde{1}}^*),\ldots,
c_1(L_{\tilde{k}}^*)\big\},$$
respectively, where 
\hbox{$\pi\co \bigcup_{\mathcal{T}}{\overU}_{\mathcal{T}} \lra {\overV}_k(\mu)$}
is the quotient projection map.
For example, if $k = 2$,
$$\pi^*\eta_{\hat{0},3}=
c_1^3({\mathcal{L}}_{\tilde{1}}^*)+
c_1^3({\mathcal{L}}_{\tilde{2}}^*)+
c_1^2({\mathcal{L}}_{\tilde{1}}^*)c_1({\mathcal{L}}_{\tilde{2}}^*)+
c_1({\mathcal{L}}_{\tilde{1}}^*)c_1^2({\mathcal{L}}_{\tilde{2}}^*).$$
Let $a_{\hat{0}} = \ev_{\hat{0}}^*c_1(\ga_{\P}^*) \in 
H^2\big({\overV}_k(\mu)\big)$, where 
$\ga_{\P} \lra \P$ is the tautological line bundle.

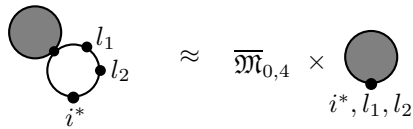
\begin{figure}[ht!]\small
\begin{pspicture}(-0.3,-1.5)(10,.5)
\psset{unit=.36cm}
\pscircle(12,-2.5){1}
\pscircle*(12,-3.5){.22}\rput(12.1,-4.2){$i^*$}
\pscircle*(13,-2.5){.2}\rput(13.75,-2.5){$l_2$}
\pscircle*(12.5,-1.63){.2}\rput(13.2,-1.3){$l_1$}
\pscircle[fillstyle=solid,fillcolor=gray](10.59,-1.09){1}
\pscircle*(11.29,-1.79){.19}
\rput(16.3,-2){$\approx$}
\rput(19,-2.3){$\ov{\frak M}_{0,4}$}
\rput(21,-2.2){$\times$}
\pscircle[fillstyle=solid,fillcolor=gray](23,-2){1}
\pscircle*(23,-3){.22}\rput(23,-3.8){$i^*,l_1,l_2$}
\end{pspicture}
\caption{An example of the decomposition~\eqref{cart_split}}
\label{decomp_fig}
\end{figure}

Suppose $\tilde{\mathcal{T}} = (M,\tilde{I};\tilde{j},\tilde{\under{d}})$ and
${\mathcal{T}} = (M,I;j,\under{d})$ are bubble types,
such that ${\mathcal{T}} < \tilde{\mathcal{T}}$, 
$\tilde{M}$ is a subset of $M\sqcup\tilde{I}^-$,
and $\mu$ is an $\tilde{M}$--tuple of constraints in~$\P$.
Let
$$I_0=\{i \in I^+  :d_i = 0\}.$$
Suppose $I_0 \subset I^+ - \hat{I}^+$ and
for every $i \in I_0$ there exists $h \in I^+$ such that $i < h$.
We can then construct a decomposition of the spaces 
${\mathcal{U}}_{{\mathcal{T}}|\tilde{\mathcal{T}}}(\mu)$
and ${\overU}_{{\mathcal{T}}|\tilde{\mathcal{T}}}(\mu)$
which is useful in computations as follows.
Let ${\overT} = (\bar{M},\bar{I};\bar{j},\bar{d})$ 
be the bubble type given~by
$$\bar{M}=M-\bigcup_{i\in I_0} M_i{\mathcal{T}}
\quad\mbox{and}\quad
\bar{I}=\Big(I\sqcup\bigcup_{i\in I_0}M_i{\mathcal{T}}\Big)\big/ \sim,$$
where $\io_i\sim h$ if $i \in I_0$ and
$h \in \big(\{i\}\cup H_i{\mathcal{T}}\big)\sqcup M_i{\mathcal{T}}$.
The set $\bar{I}$ is a graded linearly ordered set,
with its structure induced from~$I$.
Let $\bar{d}_i = d_i$ and $\bar{j}_l = j_l$
whenever \hbox{$i \in \bar{I}^+ \subset I^+$} and 
$l \in \bar{M} \subset M$.
Let $\tilde{M}'$ be the image of $\tilde{M}$ under the quotient projection map
\hbox{$I\sqcup M \lra \bar{I}\sqcup\bar{M}$}.
We identify the $\tilde{M}$--tuple $\mu$ of constraints 
with an  $\tilde{M}'$--tuple $\bar{\mu}$ 
of constraints in~$\P$~by
$$\bar{\mu}_{\bar{l}}=\bigcap_{[l]=\bar{l}}\mu_l.$$
Since every degree-zero holomorphic map is constant, we obtain
\begin{equation}\label{cart_split}\begin{split}
{\mathcal{U}}_{{\mathcal{T}}|\tilde{\mathcal{T}}}(\mu)
\approx&\Big( 
\prod_{i\in I_0}{\frak M}_{(i\cup H_i{\mathcal{T}})\sqcup M_i{\mathcal{T}}}
\times{\mathcal{U}}_{{\overT}}(\mu)\Big)
\Big/{\mathcal{A}}({\mathcal{T}}|\tilde{I}^+)\\
\subset&
\Big(\prod_{i\in I_0}\ov{\frak M}_{(i\cup H_i{\mathcal{T}})\sqcup M_i{\mathcal{T}}}
\times{\overU}_{{\overT}}(\mu)\Big)
\Big/{\mathcal{A}}({\mathcal{T}}|\tilde{I}^+). 
\end{split}\end{equation}
Here $\ov{\frak M}_{(i\cup H_i{\mathcal{T}})\sqcup M_i{\mathcal{T}}}$ denotes
the Deligne--Mumford moduli space of rational curves with marked
points labeled by the set $(i\cup H_i{\mathcal{T}})\sqcup M_i{\mathcal{T}}$,
and also ${\frak M}_{(i\cup H_i{\mathcal{T}})\sqcup M_i{\mathcal{T}}}$
denotes the main stratum 
of~$\ov{\frak M}_{(i\cup H_i{\mathcal{T}})\sqcup M_i{\mathcal{T}}}$.
If $i \in I_0 \subset \tilde{I}^+$, then by definition
the line bundle $L_i\tilde{\mathcal{T}}$ restricts to
the universal tangent line bundle at the marked point~$i$
over~$\ov{\frak M}_{(i\cup H_i{\mathcal{T}})\sqcup M_i{\mathcal{T}}}$.
We will denote this bundle by~$\ga_{{\mathcal{T}};i}$.
If $\tilde{d}_i \neq 0$ for all $i \in \tilde{I}^+$
and $M_0 \subset [N] \cap M_i{\mathcal{T}}$ for some
$i \in \tilde{I}^+$,
we will write $\tilde{\mathcal{T}}/M_0$ for the bubble type ${\overT}$
corresponding to ${\mathcal{T}} = \tilde{\mathcal{T}}(M_0)$
under the construction of this paragraph.
The decomposition~\eqref{cart_split} for the bubble 
${\mathcal{T}}(M_0)$ of Figure~\ref{psi_fig} is illustrated 
in~Figure~\ref{decomp_fig}.

\subsection{Structural descriptions}
\label{analysis_subs}

In this subsection, we define certain bundle sections over the spaces
${\mathcal{U}}_{{\mathcal{T}}|\tilde{\mathcal{T}}}$.
These sections are central to this paper,
as the zeros of these and closely related sections
count rational curves with pre-specified singularities.
We state a basic transversality lemma 
that implies that these sections are well-behaved 
over~${\mathcal{U}}_{{\mathcal{T}}|\tilde{\mathcal{T}}}$ in most cases.
The general structure of the 
spaces~${\overU}_{\tilde{\mathcal{T}}}(\mu)$ and
the behavior of the sections near the boundary strata 
are described by Proposition~\ref{str_prp}.

Let $q_S$ denote the standard stereographic projection $\Bbb{C} \lra S^2$
mapping the origin in $\Bbb{C}$ to the south pole of~$S^2$.
Suppose 
$$b = \big(M,I;x,(j,y),u\big) \in{\mathcal{B}}_{\mathcal{T}}$$ 
and $m \in \Bbb{Z}^+$.
If $i \in I^+$, let
$${\mathcal{D}}_{{\mathcal{T}},i}^{(m)}b=
\frac{1}{m!}
\frac{D^{m-1}}{ds^{m-1}}\frac{d}{ds}(u_i\circ q_S)\Big|_{(s,t)=0},$$
where $(s,t)$ are the real and imaginary coordinates on~$\Bbb{C}$
and $\frac{D^{m-1}}{ds^{m-1}}$ denotes the $(m - 1)$st covariant
derivative with respect to the Levi-Civita connection of
some metric $g_{b,i}$~on~$\P$. 
If $h \in \hat{I}^-$ and $l \in M$, we similarly define
$${\mathcal{D}}_{{\mathcal{T}},h}^{(m)}b=
\frac{1}{m!}\frac{D^{m-1}}{ds^{m-1}}
\frac{d}{ds}u_{\io_h}(x_h + s\big)\Big|_{s=0}
\quad\!\mbox{and}\quad\!
{\mathcal{D}}_{{\mathcal{T}},l}^{(m)}b=
\frac{1}{m!}\frac{D^{m-1}}{ds^{m-1}}
\frac{d}{ds}u_{j_l}\big(y_l + s\big)\Big|_{s=0}$$
Here we take covariant derivatives with respect to 
some metrics $g_{b,h}$ and $g_{b,l}$ on~$\P$, respectively.
If $\tilde{I}$ is a subset of~$I^+$, 
each metric $g_{b,i}$ is Kahler near~$\ev_i(b)$,
and the family $\{g_{b,i}\}$ is invariant under the action
of the group of $G_{{\mathcal{T}}|\tilde{I}}$,  
${\mathcal{D}}_{{\mathcal{T}},i}^{(m)}$ induces a section of 
$\hbox{Hom}(L_i^{\otimes m}{\mathcal{T}},\ev_i^*T\P)$ over 
${\mathcal{U}}_{{\mathcal{T}}|\tilde{I}}$  given~by
$${\mathcal{D}}_{{\mathcal{T}},i}^{(m)}[b,c_i]
=c_i{\mathcal{D}}_{{\mathcal{T}},i}^{(m)}b,
\quad\hbox{if}~~b \in {\mathcal{B}}_{\mathcal{T}},~c_i \in \Bbb{C}.$$ 
Under analogous circumstances,
${\mathcal{D}}_{{\mathcal{T}},h}^{(m)}$ and ${\mathcal{D}}_{{\mathcal{T}},l}^{(m)}$
induce sections of 
$$\hbox{Hom}(L_{\io_h}^{*\otimes m}{\mathcal{T}},\ev_h^*T\P)
\qquad\hbox{and}\qquad 
\hbox{Hom}(L_{j_l}^{*\otimes m}{\mathcal{T}},\ev_l^*T\P),$$
respectively.
In a certain sense, the choice of the metrics does not matter,
since ${\mathcal{D}}_{{\mathcal{T}},i}^{(m)}b$,
${\mathcal{D}}_{{\mathcal{T}},h}^{(m)}b$, and
${\mathcal{D}}_{{\mathcal{T}},l}^{(m)}b$ are well-defined 
modulo the image of the lower-order derivatives, and 
only these quotients have a geometric meaning.
However, the method of \cite[Subsection~2.5]{Z1}
for proving the explicit estimates of Proposition~\ref{str_prp}
makes use of special properties of the metric near the point
where the derivatives are taken.
Thus, we put
$$g_{b,i} = g_{\P,\ev_i(b)},\qquad
g_{b,h} = g_{\P,\ev_h(b)},\quad\hbox{and}\quad
g_{b,l} = g_{\P,\ev_l(b)},$$
where $g_{\P,\cdot}$ is as in Lemma~\ref{flat_metrics},
which is exactly \cite[Lemma~2.1]{Z1}.

\begin{lmm}
\label{flat_metrics}
There exist $r_{\P} > 0$ and a smooth family of Kahler metrics
$$\{g_{\P,q} : q \in \P\}$$ 
on~$\P$ with the following property.
If $B_q(q',r) \subset \P$ denotes the $g_{\P,q}$--geodesic ball about~$q'$,
the triple $(B_q(q,r_{\P}),J,g_{\P,q})$ is isomorphic
to a ball in $\Bbb{C}^n$ for \hbox{all $q \in \P$}.
\end{lmm}

Suppose $\tilde{\mathcal{T}} = (M,\tilde{I};\tilde{j},\under{\tilde{d}})$
and ${\mathcal{T}} \equiv (M,I;j,\under{d})$ are 
bubble types such that ${\mathcal{T}} < \tilde{\mathcal{T}}$.
For each $k \in \tilde{I}^+$, let 
$\tilde{\mathcal{T}}_k \equiv (M_k,\tilde{I}_k;\tilde{j}^{(k)},
                                 \under{\tilde{d}}^{(k)})$ be
the one-component bubble type defined as in Subsection~\ref{notation_subs}
and let
$${\mathcal{T}}_k(\tilde{\mathcal{T}}) \equiv (M_k,I_k;j^{(k)},\under{d}^{(k)})
<\tilde{\mathcal{T}}_k$$
be the bubble type given~by
\begin{gather*}
I_k=\{\io_k\}\cup\big\{h \in I :h \ge k;~
h \not\ge h'~\mbox{for all}\  h' \in \tilde{I}~\mbox{such that}~h' > k\big\},\\
j^{(k)}_l = j_l~~\mbox{for all}\  l \in M_k; \qquad
d^{(k)}_h = d_h~~\mbox{for all}\  h \in I_k^+.
\end{gather*}
Let
${\mathcal{U}}_{\tilde{\mathcal{T}},{\mathcal{T}}} = 
\prod_{k\in\tilde{I}^+}{\mathcal{U}}_{{\mathcal{T}}_k(\tilde{\mathcal{T}})|\tilde{\mathcal{T}}_k}
\subset {\overU}_{\tilde{\mathcal{T}},\tilde{\mathcal{T}}}$.
If $\tilde{M}$ is a subset of $\tilde{I}^-\cup M$, we define an evaluation~map
\begin{gather*}
\ev_{\tilde{\mathcal{T}},\tilde{M}}\co 
{\overU}_{\tilde{\mathcal{T}},\tilde{\mathcal{T}}} \equiv 
 \prod_{k\in\tilde{I}^+}  {\overU}_{{\mathcal{T}}_k}\lra
{\mathcal{X}}_{\tilde{\mathcal{T}}}(\tilde{M}) \equiv 
(\P)^{\hat{\tilde{I}}^-\sqcup\tilde{I}^+\sqcup(\tilde{M}\cap M)}\\
\hbox{by}\qquad
\pi_l\big(\ev_{\tilde{\mathcal{T}},\tilde{M}}(b)\big)=
\begin{cases}
\ev_l(b_k),&\hbox{if}~l \in M_k;\\
\ev_l(b_l),&\hbox{if}~l \in \tilde{I}^+;\\
\ev_l(b_{\io_l}),&\hbox{if}~l \in \hat{\tilde{I}}^-.
\end{cases}
\end{gather*}
If $\mu$ is an $\tilde{M}$--tuple of submanifolds in $\P$, let
\begin{multline*}
\De_{\tilde{\mathcal{T}}}(\mu)
=\prod_{i\in\tilde{I}^-}  
\biggl\{(x_h)_{h\in H_i\tilde{\mathcal{T}}} \in (\P)^{H_i\tilde{\mathcal{T}}}  :
\begin{array}{ll}
x_h = x_l &\mbox{for all}\ h,l \in H_i\tilde{\mathcal{T}},\\
x_h \in \mu_i & \mbox{if}\ i \in \tilde{M}\end{array}\biggr\} \\
\times  \prod_{l\in(\tilde{M}\cap M)} \mu_l
\subset {\mathcal{X}}_{\tilde{\mathcal{T}}}(\tilde{M}).
\end{multline*}
We denote by
${\mathcal{N}}\De_{\tilde{\mathcal{T}}}(\mu) \subset  
T{\mathcal{X}}_{\tilde{\mathcal{T}}}(\tilde{M})\big|_{\De_{\tilde{\mathcal{T}}}(\mu)}$
the normal bundle of 
$\De_{\tilde{\mathcal{T}}}(\mu)$ in ${\mathcal{X}}_{\tilde{\mathcal{T}}}(\tilde{M})$
as well as an extension of this normal bundle to a neighborhood
of $\De_{\tilde{\mathcal{T}}}(\mu)$ in~${\mathcal{X}}_{\tilde{\mathcal{T}}}(\tilde{M})$.
By definition, 
$${\mathcal{U}}_{\ti{\mathcal{T}}}(\mu)=
\ev_{\ti{\mathcal{T}},\ti{M}}^{-1}\big(\De_{\ti{\mathcal{T}}}(\mu)\big)
\cap {\mathcal{U}}_{\ti{\mathcal{T}},\ti{\mathcal{T}}}
\quad\hbox{and}\quad
{\mathcal{U}}_{{\mathcal{T}}|\ti{\mathcal{T}}}(\mu) = 
\ev_{\ti{\mathcal{T}},\ti{M}}^{-1}\big(\De_{\ti{\mathcal{T}}}(\mu)\big)
\cap {\mathcal{U}}_{\ti{\mathcal{T}},{\mathcal{T}}}.$$
The following lemma and now-standard arguments, such as in
\cite[Chapter~3]{MS},
imply that ${\mathcal{U}}_{\tilde{\mathcal{T}},{\mathcal{T}}}(\mu)$
is a smooth orbifold,
if $\mu$ is a tuple of constraints in general position:

\begin{lmm}
\label{str_lmm}
Suppose $u\co S^2 \lra \P$ is a holomorphic map of degree~$d$,
$$z_1,\ldots,z_k \in S^2, \quad 
v_i \in T_{z_i}S^2 - \{0\} ~\hbox{for}~ 
i = 1,\ldots,k, \quad\hbox{and}\quad m_1,\ldots,m_k \in \Bbb{Z}^+.$$
If $d + 1 \ge \sum_{i}m_i$, the map
$$\phi\co \big\{\xi \in \Ga(S^2;u^*T\P) :\bar{\partial}\xi = 0\big\}\lra
\bigoplus_{i=1}^{i=k}\bigoplus_{j=1}^{j=m_i}T_{u(z_i)}\P,\quad
\phi_{i,j}(\xi) = D_{v_i}^{j-1}\xi\big|_{z_i},$$
is surjective.
\end{lmm}

In the statement of this lemma, $D_{v_i}^j\xi\big|_{z_i}$ denotes 
the $j$th covariant derivative of $\xi$ along $u$ in the direction of~$v$.
The meaning of the lemma is the same 
no matter what connection is used near each point.
The proof is a very slight generalization of that of 
\cite[Corollary~6.3]{Z1}.
Lemma~\ref{str_lmm} implies that if $d$ is sufficiently high, 
a certain section over the space of degree--$d$ holomorphic maps
is transverse to the zero~set.
In many actual computations, low, but positive-degree, cases will not appear
for simple geometric reasons.
For example, the space of degree-one maps whose image is
a cuspidal curve is empty.
Thus, if $k = 1$ and $m_1 = 2$, the relevant implication of 
the lemma is valid as long \hbox{as $d > 0$}.

Let
$${\mathcal{F}}_{\tilde{\mathcal{T}}}{\mathcal{T}}=
\bigoplus_{h\in I^+-\tilde{I}}    {\mathcal{F}}_h{\mathcal{T}}
=\bigoplus_{k\in\tilde{I}^+}{\mathcal{F}}\{{\mathcal{T}}_k(\tilde{\mathcal{T}})\}
\lra{\mathcal{U}}_{\tilde{\mathcal{T}},{\mathcal{T}}}.$$
If $\tilde{\mathcal{T}}$ is a bubble type such that $\hat{\tilde{I}}^- = \eset$,
we will write $\mathcal{FT}$ for~${\mathcal{F}}_{\tilde{\mathcal{T}}}{\mathcal{T}}$.
By (2) of Proposition~\ref{str_prp},
${\mathcal{F}}_{\tilde{\mathcal{T}}}{\mathcal{T}}$ is
the ``normal bundle'' of ${\mathcal{U}}_{\tilde{\mathcal{T}},{\mathcal{T}}}$
in ${\overU}_{\tilde{\mathcal{T}},\tilde{\mathcal{T}}}$.
Part (3) of Proposition~\ref{str_prp} describes the behavior
of various evaluation maps and bundle sections over
${\mathcal{U}}_{\tilde{\mathcal{T}},\tilde{\mathcal{T}}}$ 
near the stratum ${\mathcal{U}}_{\tilde{\mathcal{T}},{\mathcal{T}}}$
of the boundary of ${\overU}_{\tilde{\mathcal{T}},\tilde{\mathcal{T}}}$.
However, before we can state the relevant expansions,
we need to introduce more notation.

If $k \in \tilde{I}^+$, $h_1,h_2 \in I_k^+$, and $l \in M_k$, let
\begin{gather*}
i_{\mathcal{T}}(h_1,h_2) = \max
\big\{i \in I_k^+ : i \le h_1,~i \le h_2\big\},\qquad
i_{\mathcal{T}}(l,h_2) = i_{\mathcal{T}}(j_l,h_2);\\
\chi_{{\mathcal{T}},h_1}(h_2)=
\begin{cases}
0,&\hbox{if~} d_i = 0~\mbox{for all}\  
i \in [i_{\mathcal{T}}(h_1,h_2),h_1]\cup[i_{\mathcal{T}}(h_1,h_2),h_2];\\
1,&
\begin{aligned}
\hbox{if~} &d_i = 0~\mbox{for all}\  
i \in [i_{\mathcal{T}}(h_1,h_2),h_1]\cup[i_{\mathcal{T}}(h_1,h_2),h_2] - \{h_2\},\\
&\hbox{but~} d_{h_2} \neq 0;
\end{aligned}\\
2,&\hbox{otherwise}.
\end{cases}
\end{gather*}
Put $\chi_{h_1}({\mathcal{T}}) = \big\{h_2 \in I_k^+ :
\chi_{{\mathcal{T}},h_1}(h_2) = 1\big\}$ and
$\chi_l({\mathcal{T}}) = \big\{h_2 \in \chi_{j_l}({\mathcal{T}}) :
h_2 \not\le j_l\big\}$.
If 
$$\ups = \big[(M,I;x,(j,y),u),(v_h)_{h\in I^+-\tilde{I}}\big]
\in {\mathcal{F}}_{\tilde{\mathcal{T}}}{\mathcal{T}}=
\bigoplus_{k\in I^+-\tilde{I}}   {\mathcal{F}}_k{\mathcal{T}}$$
and $m,m' \in \Bbb{Z}$, let
\begin{align*}
x_{h_1;h_2}(\ups)&=
\sum_{h\in(i_{\mathcal{T}}(h_1,h_2),h_2]}      x_{\io_h}
  \prod_{i\in(i_{\mathcal{T}}(h_1,h_2),h)}      v_i~\in~ 
L_{i_{\mathcal{T}}(h_1,h_2)}^*{\mathcal{T}};\\
y_{h_1;l}(\ups)&=
\sum_{h\in(i_{\mathcal{T}}(h_1,l),j_l]}      x_{\io_h}
  \prod_{i\in(i_{\mathcal{T}}(h_1,l),h)}      v_i
+y_l  \prod_{i\in(i_{\mathcal{T}}(h_1,l),j_l]}      v_i
\in L_{i_{\mathcal{T}}(h_1,l)}^*{\mathcal{T}}.
\end{align*}
\begin{gather}\label{str_subs_e1}
\rho_{{\mathcal{T}},h_1;h_2}^{(m;m')}(\ups)=
\Big(  \bigotimes_{i\in(k,i_{\mathcal{T}}(h_1,h_2)]}
          v_i\Big)^{\otimes (-m)}  \otimes
\Big(\bigotimes_{i\in(i_{\mathcal{T}}(h_1,h_2),h_2]}     v_i\Big)^{\otimes m'}
\in \tilde{\mathcal{F}}_{h_1;h_2}^{(m;m')}{\mathcal{T}},\\
\hbox{where}\qquad
\tilde{\mathcal{F}}_{h_1;h_2}^{(m;m')}{\mathcal{T}}=
\Big(\bigotimes_{i\in(k,i_{\mathcal{T}}(h_1,h_2)]}        
{\mathcal{F}}_i{\mathcal{T}}
\Big)^{\otimes(-m)} \otimes
\Big(  \bigotimes_{i\in(i_{\mathcal{T}}(h_1,h_2),h_2]}        
{\mathcal{F}}_i{\mathcal{T}} \Big)^{\otimes m'}  \notag.
\end{gather}
The map $\rho_{{\mathcal{T}},h_1;h_2}^{(m;m')}$ is defined if 
$h_1 = k$ or $\ups \not\in Y({\mathcal{F}}_{\tilde{\mathcal{T}}}{\mathcal{T}};
I^+ - \tilde{I})$.
If $i \in I_k^- - \tilde{I}$, we define the map 
$\rho_{{\mathcal{T}},i;h_2}^{(m;m')}$
by \eqref{str_subs_e1}, but with
$i_{\mathcal{T}}(h_1,h_2)$ replaced \hbox{by $\io_i \in I_k^+$}.
Furthermore, we define the map $\rho_{{\mathcal{T}},i;i}^{(m;m')}$
by replacing $h_2$ in \eqref{str_subs_e1}
with the unique element $h(i) \in I_k^+$
such that $\io_{h(i)} = i$.
We will write 
$x_{l;h}(\ups)$, $\rho_{{\mathcal{T}},l;h}$, and ${\tilde F}_{l;h}{\mathcal{T}}$ for 
$x_{j_l;h}(\ups)$, $\rho_{{\mathcal{T}},j_l;h}$, and ${\tilde F}_{j_l;h}{\mathcal{T}}$, 
respectively, whenever $l \in M_k$.
Finally, if $d_i \neq 0$ for some $i \in [k,j_l]$, 
let $\si_{\tilde{\mathcal{T}},{\mathcal{T}}}(l) = 1$; 
otherwise, let $\si_{\tilde{\mathcal{T}},{\mathcal{T}}}(l) = 0$. 
In the former case, we put
\begin{align*}
j_l^*({\mathcal{T}})&=
\begin{cases}
j_l \in M_k,&\hbox{if}~d_{j_l} \neq 0;\\
\min\big\{i \in \hat{I}_k^- : i \le j_l,~
d_h = 0~\mbox{for all}\  h \in (i;j_l]\big\},&\hbox{if}~d_{j_l} = 0;
\end{cases}\\
y_{l;{\mathcal{T}}}(\ups)&=\begin{cases}
0,&\hbox{if}~d_{j_l} \neq 0;\\
y_{h;l}(\ups),&\hbox{if}~d_{j_l} = 0,\\
\end{cases}
\end{align*}
where $h \in I_k^+$ is given by $\io_h = j_l^*({\mathcal{T}})$
if $j_l^*({\mathcal{T}})\in \hat{I}_k^-$.

\begin{prp}
\label{str_prp}
Suppose $\tilde{\mathcal{T}} = (M,\tilde{I};\tilde{j},\under{\tilde{d}})$
and ${\mathcal{T}} = (M,I;j,\under{d})$
are bubble types such that \hbox{${\mathcal{T}} < \tilde{\mathcal{T}}$}.

{\rm(1)}\qua The spaces
${\mathcal{U}}_{\tilde{\mathcal{T}},\tilde{\mathcal{T}}}$ and
${\mathcal{U}}_{\tilde{\mathcal{T}},{\mathcal{T}}}$
are smooth orbifolds, while ${\overU}_{\tilde{\mathcal{T}},\tilde{\mathcal{T}}}$
is an ms-orbifold.

{\rm(2)}\qua There exist 
$\de \in C({\mathcal{U}}_{\tilde{\mathcal{T}},{\mathcal{T}}};\Bbb{R}^+)$
and a map 
$\phi_{\tilde{\mathcal{T}},{\mathcal{T}}}\co 
{\mathcal{F}}_{\tilde{\mathcal{T}}}{\mathcal{T}}_{\de} \lra 
{\overU}_{\tilde{\mathcal{T}},\tilde{\mathcal{T}}}$ such that
$\phi_{\tilde{\mathcal{T}},{\mathcal{T}}}$ is a homeomorphism onto an open neighborhood 
$W_{\tilde{\mathcal{T}},{\mathcal{T}}}$ of 
${\mathcal{U}}_{\tilde{\mathcal{T}},{\mathcal{T}}}$ in 
${\overU}_{\tilde{\mathcal{T}},\tilde{\mathcal{T}}}$,
\begin{gather*}
\phi_{\tilde{\mathcal{T}},{\mathcal{T}}}
\big({\mathcal{F}}_{\tilde{\mathcal{T}}}{\mathcal{T}}_{\de}\cap 
Y({\mathcal{F}}_{\tilde{\mathcal{T}}}{\mathcal{T}};I^+ - \tilde{I})\big)
\subset \partial{\overU}_{\tilde{\mathcal{T}},\tilde{\mathcal{T}}},\\
\hbox{and}\qquad
\phi_{\tilde{\mathcal{T}},{\mathcal{T}}}\co 
{\mathcal{F}}_{\tilde{\mathcal{T}}}{\mathcal{T}}_{\de} -  
Y({\mathcal{F}}_{\tilde{\mathcal{T}}}{\mathcal{T}};I^+ - \tilde{I})\lra
{\mathcal{U}}_{\tilde{\mathcal{T}},\tilde{\mathcal{T}}}\cap W_{\mathcal{T}}
\end{gather*}
is an orientation-preserving diffeomorphism.

{\rm(3)}\qua Furthermore, there exist normal-neighborhood models
and collections of trivializations such that
the following identities are satisfied
by all elements $k \in \tilde{I}^+$ and
$(b;\ups) \in {\mathcal{F}}_{\tilde{\mathcal{T}}}{\mathcal{T}}_{\de} - 
Y({\mathcal{F}}_{\tilde{\mathcal{T}}}{\mathcal{T}};I^+ - \tilde{I})$:

{\rm(3a)}\qua
if $\tilde{M}$ is a subset of $\tilde{I}^-\sqcup M$,
$$\ev_{\tilde{\mathcal{T}},\tilde{M}}
\big(\phi_{\tilde{\mathcal{T}},{\mathcal{T}}}(\ups)\big)=
\ev_{\tilde{\mathcal{T}},\tilde{M}}(b)+
\ve_{\tilde{\mathcal{T}},\tilde{M}}(\ups),$$
where 
$\ve_{\tilde{\mathcal{T}},\tilde{M}}\co 
{\mathcal{F}}_{\tilde{\mathcal{T}}}{\mathcal{T}}_{\de} - 
Y({\mathcal{F}}_{\tilde{\mathcal{T}}}{\mathcal{T}};I^+ - \tilde{I}) \lra 
\ev_{\tilde{\mathcal{T}},\tilde{M}}^*T{\mathcal{X}}_{\tilde{\mathcal{T}}}(\tilde{M})$
is a $C^1$--negligible map;

%\item[{\rm(3b)}]
{\rm(3b)}\qua
if $k \in \tilde{I}^+$ and $m \in \Bbb{Z}^+$,
\begin{equation*}\begin{split}
&{\mathcal{D}}_{\tilde{\mathcal{T}},k}^{(m)}\phi_{\tilde{\mathcal{T}},{\mathcal{T}}}(\ups)\\
&\qquad\quad= \sum_{m'=1}^m\binom{m - 1}{m' - 1} 
\sum_{h\in\chi_k({\mathcal{T}})}  x_{k;h}(\ups)^{m-m'}
\Big\{{\mathcal{D}}_{{\mathcal{T}},h}^{(m')} + 
\ve_{{\mathcal{T}},h}^{(m')}(\ups)\Big\}
\rho_{{\mathcal{T}},k;h}^{(m;m')}(\ups),
\end{split}\end{equation*}
where each map 
$$\ve_{{\mathcal{T}},h}^{(m')}\co 
{\mathcal{F}}_{\tilde{\mathcal{T}}}{\mathcal{T}}_{\de} - 
Y({\mathcal{F}}_{\tilde{\mathcal{T}}}{\mathcal{T}};I^+ - \tilde{I}) \lra
\hbox{Hom}(L_h^{\otimes m'}{\mathcal{T}},\ev_k^*T\P)$$
is $C^1$--negligible;

%\item[{\rm(3c)}]
{\rm(3c)}\qua
if $l \in M$ and $m \in \Bbb{Z}^+$,
\begin{equation*}\begin{split}
&{\mathcal{D}}_{\tilde{\mathcal{T}},l}^{(m)}\phi_{\tilde{\mathcal{T}},{\mathcal{T}}}(\ups)\\
&\qquad=
\si_{\tilde{\mathcal{T}},\cal T}(l)\sum_{m'=m+1}^{\i}   \binom{m'}{m}
y_{l;{\mathcal{T}}}(\ups)^{m'-m}
\Big\{{\mathcal{D}}_{{\mathcal{T}},j_l^*({\mathcal{T}})}^{(m')} + 
\ve_{{\mathcal{T}},j_l^*({\mathcal{T}})}^{(m')}(\ups)\Big\}
\rho_{{\mathcal{T}},j_l^*({\mathcal{T}});j_l}^{(m;m')}(\ups)\\
&\qquad\quad+ (-1)^m\sum_{m'=1}^{\i}\Bigg\{\binom{m + m' - 1}{m}\\
&\qquad\qquad\qquad
\sum_{h\in\chi_l({\mathcal{T}})}     
\big(y_{h;l}(\ups) - x_{l;h}(\ups)\big)^{-(m+m')}
\Big\{{\mathcal{D}}_{{\mathcal{T}},h}^{(m')}+\ve_{{\mathcal{T}},h}^{(m')}(\ups)\Big\}
\rho_{{\mathcal{T}},l;h}^{(m;m')}(\ups)\Bigg\},
\end{split}\end{equation*}
where $\ve_{{\mathcal{T}},j_l^*({\mathcal{T}})}^{(m')}(\ups) = 0$ if $m' \neq m$
and
$\ve_{{\mathcal{T}},j_l^*({\mathcal{T}})}^{(m)}$ is a $C^1$--negligible map on 
$${\mathcal{F}}_{\tilde{\mathcal{T}}}{\mathcal{T}}_{\de}-
Y({\mathcal{F}}_{\tilde{\mathcal{T}}}{\mathcal{T}};I^+ - \tilde{I}).$$
%\end{enumerate}
\end{prp}

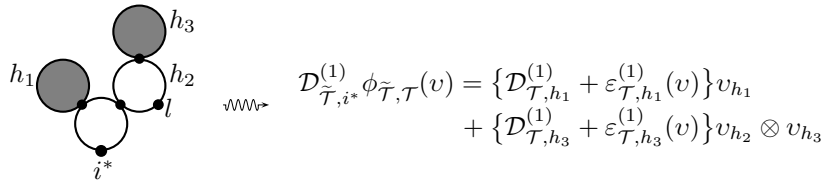
\begin{figure}[ht!]\small
\begin{pspicture}(-1.8,-1.7)(10,1.5)
\psset{unit=.36cm}
\pscircle(1,-2.5){1}
\pscircle*(1,-3.5){.22}\rput(1.1,-4.2){$i^*$}
\pscircle[fillstyle=solid,fillcolor=gray](-.41,-1.09){1}
\pscircle*(.29,-1.79){.19}\rput(-1.9,-.7){$h_1$}
\pscircle(2.41,-1.09){1}\pscircle*(1.71,-1.79){.19}\rput(4,-.7){$h_2$}
\pscircle[fillstyle=solid,fillcolor=gray](2.41,.91){1}
\pscircle*(2.41,-.09){.19}\rput(4,1.3){$h_3$}
\pscircle*(3.11,-1.79){.2}\rput(3.5,-1.9){$l$}
\pscoil[linewidth=.03,coilarmA=0.1,coilarmB=0.4,coilaspect=0,coilheight=.7,
coilwidth=.4]{->}(5.5,-1.8)(7.2,-1.8)
\rput(18,-1.8){\sm{\begin{tabular}{l}
${\mathcal{D}}_{\tilde{\mathcal{T}},i^*}^{(1)}
\phi_{\tilde{\mathcal{T}},{\mathcal{T}}}(\ups)=
\big\{{\mathcal{D}}_{{\mathcal{T}},h_1}^{(1)} + \ve_{{\mathcal{T}},h_1}^{(1)}(\ups)
\big\}\ups_{h_1}$\\
${}\qquad\qquad\qquad
+\big\{{\mathcal{D}}_{{\mathcal{T}},h_3}^{(1)} + \ve_{{\mathcal{T}},h_3}^{(1)}(\ups)
\big\}\ups_{h_2} \otimes \ups_{h_3}$\end{tabular}}}
\end{pspicture}
\caption{An example of the expansion (3b) of Proposition~\ref{str_prp}}
\label{exp_fig}
\end{figure}

The expansions (3b) and (3c) above look quite complicated.
However, it is clear from the construction that they involve
monomials maps between vector bundles.
Figure~\ref{exp_fig} illustrates the expansion~(3b) in a case 
when $\tilde{I}^+ = \{i^*\}$ is a single-element set and $m = 1$.
Note that, while the stratum ${\mathcal{U}}_{{\mathcal{T}}|\tilde{\mathcal{T}}}$
of Figure~\ref{exp_fig} has codimension three 
in~${\overU}_{\tilde{\mathcal{T}}}$,
the section ${\mathcal{D}}_{\tilde{\mathcal{T}},i^*}^{(1)}$
depends only on two parameters of the normal bundle,
$\ups_{h_1}$ \hbox{and $\ups_{h_2} \otimes \ups_{h_3}$},
at least up to negligible~terms.
Such bubble types~${\mathcal{T}}$ will always be hollow in the sense of 
Definition~\ref{reg_dfn1a} and will not effect our~computations.

One of the crucial points about the expansions (3b) and (3c) of Proposition~\ref{str_prp} 
is that the terms that appear between the curly brackets depend on~$h$, 
$j_l^*({\mathcal{T}})$, and $m'$, and not~$k$, $l$, or~$m$.
Thus, by subtracting expansions of lower-derivatives 
with appropriate coefficients from  expansions of higher-derivatives, 
we can get rather good estimates on the latter along the subspaces
of the main stratum of a moduli space on which the former vanish;
see the proofs of Lemmas~\ref{n3p3_contr_lmm1b} and \ref{n3p3_contr_lmm1c},
for example.

Statements (1) and (2) of Proposition~\ref{str_prp} are basically 
special cases of Theorem~2.8 in~\cite{Z1}.
The map $\phi_{\ti{\mathcal{T}},{\mathcal{T}}}$ of Proposition~\ref{str_prp} 
is the product of the gluing maps, as constructed in Subsection~3.6 of~\cite{Z2},  
corresponding to each pair of bubble types ${\mathcal{T}}_k(\ti{\mathcal{T}}) < \ti{\mathcal{T}}_k$.
Claims (3a) and (3b) are proved in Subsection~4.1 of~\cite{Z2}
and in Subsection~2.5 of~\cite{Z1}, though only a slightly weaker version of 
the $m = 1$ case of (3b) is stated as part of Theorem~2.8 in~\cite{Z1}.
The proof of (3c) uses essentially the same trick as the proof of~(3b) in~\cite{Z1}.
The difference is that we make $(j_l,y_l)$, instead of~$(k,\i)$, a node and then use 
the explicit nature of the gluing map~$\phi_{\ti{\mathcal{T}}_k,{\mathcal{T}}_k(\ti{\mathcal{T}})}$
to do integration by parts as before.
The appropriate normal-neighborhood models and collections of trivializations referred 
to in (3) are described as follows.
In (3a), we use the product of exponential maps in every component taken with respect to 
the family of metrics of Lemma~\ref{flat_metrics}.
In (3b) and (3c), all the relevant bundles have the form 
$\hbox{Hom}(L_{\tilde{\mathcal{T}},i}^{\otimes\pm m},\ev_i^*T\P)$.
We use parallel transport in the metric $g_{b,i}$ to identify
$\phi_{\ti{\mathcal{T}},{\mathcal{T}}}^*\ev_i^*T\P$ with 
$\pi_{{\mathcal{F}}_{\ti{\mathcal{T}}}{\mathcal{T}}}^*\ev_i^*T\P$.
On the other hand, the map~$\phi_{\ti{\mathcal{T}}_k,{\mathcal{T}}_k(\ti{\mathcal{T}})}$, 
constructed in Subsection~3.6 of~\cite{Z2}, is descendant from an $S^1$--equivariant map 
$\ti{\phi}_{\ti{\mathcal{T}}_k,{\mathcal{T}}_k(\tilde{\mathcal{T}})}$
from a bundle over ${\mathcal{B}}_{{\mathcal{T}}_k(\ti{\mathcal{T}})}$ 
to~${\overU}_{\ti{\mathcal{T}}_k}^{(k)}$.
This map $\ti{\phi}_{\ti{\mathcal{T}}_k,{\mathcal{T}}_k(\ti{\mathcal{T}})}$ 
induces an identification of the line bundles 
$\phi_{\ti{\mathcal{T}},{\mathcal{T}}}^*L_k\ti{\mathcal{T}}$
and~$\pi_{{\mathcal{F}}_{\ti{\mathcal{T}}}{\mathcal{T}}}^*L_k{\mathcal{T}}$ 
over~${\mathcal{F}}_{\ti{\mathcal{T}}}{\mathcal{T}}_{\de}$.

\begin{remarks}
{\rm
(1)\qua By the construction of the map~$\phi_{\ti{\mathcal{T}}_k,{\mathcal{T}}_k(\ti{\mathcal{T}})}$ 
in \cite[Subsection\,3.6]{Z2}, if $\ups$ is an element of 
${\mathcal{F}}_{\ti{\mathcal{T}}}{\mathcal{T}}_{\de} - Y({\mathcal{F}}_{\ti{\mathcal{T}}}{\mathcal{T}};I^+ - \ti{I})$,
$$\phi_{\ti{\mathcal{T}},{\mathcal{T}}}(\ups)=
\big(M,\ti{I};x(\ups),(\ti{j},y(\ups)),u_{\ups}\big),$$
where $y_l(\ups)=y_{k;l}(\ups)$ for all $l \in M_k$.
We use this fact in Section~\ref{n3p3_sec}.

(2)\qua Even in rather simple cases, the bundle 
${\mathcal{F}}_{\ti{\mathcal{T}}}{\mathcal{T}} \lra {\overU}_{\ti{\mathcal{T}},{\mathcal{T}}}$
is not the normal bundle of ${\overU}_{\tilde{\mathcal{T}},{\mathcal{T}}}$
in~${\overU}_{\ti{\mathcal{T}},\ti{\mathcal{T}}}$,
as can be seen from Subsection~3.2 in~\cite{P}.
Statement~(2) of Proposition~\ref{str_prp} only implies that
the restrictions of ${\mathcal{F}}_{\tilde{\mathcal{T}}}{\mathcal{T}}$
of the normal bundle of ${\overU}_{\tilde{\mathcal{T}},{\mathcal{T}}}$
in~${\overU}_{\tilde{\mathcal{T}},\tilde{\mathcal{T}}}$ to
${\mathcal{U}}_{\tilde{\mathcal{T}},\tilde{\mathcal{T}}}$ are isomorphic.}
\end{remarks}

\section{Example 1: Rational triple-pointed curves in $\PPP$}
\label{n3p3_sec}

\subsection{Summary}
\label{n3p3_sum_subs}

In this section, we illustrate our computational method by proving Theorem~\ref{n3p3_thm}.
We first describe the set ${\mathcal{V}}_1^{(2)}(\mu)$ that appears
in the statement of the theorem.
If $N = p + q$, we view $\mu$ as an $[N]$--tuple of
constraints in $\PPP$, where 
$$[N]=\big\{1,\ldots,N\big\}.$$
Let \hbox{$I_1 = I_1^- \sqcup I_1^+$} be a two-element graded linearly ordered set. 
We denote the unique elements of $I_1^-$ and $I_1^+$
by~$\hat{0}$ and~$\tilde{1}$, respectively. 
Put
$${\mathcal{T}}_0=\big(M_0,I_1,;j^{(0)},d^{(0)}\big),\quad
{\mathcal{T}}_1=\big(M_1,I_1;j^{(1)},d^{(0)}\big),\quad
{\mathcal{T}}_2=\big(M_2,I_1;j^{(2)},d^{(0)}\big),$$
where
$$M_0=[N],\  M_1=M_0 \sqcup \{\hat{1}\},
\  M_2=M_1 \sqcup \{\hat{2}\}; \ 
j^{(k)}_l = \tilde{1}~~\mbox{for all}\  l \in M_k;\ 
d^{(0)}_{\tilde{1}} = d.$$
The tuples ${\mathcal{T}}_0$, ${\mathcal{T}}_1$, and ${\mathcal{T}}_2$
are bubble types, and we set
\begin{equation*}\begin{split}
{\mathcal{V}}_1^{(1)}(\mu) &=\big\{b \in {\mathcal{U}}_{{\mathcal{T}}_1}(\mu) :
\ev_{\tilde{1}}(b) = \ev_{\hat{1}}(b)\big\},\\
{\mathcal{U}}_{{\mathcal{T}}_2}^{(1)}(\mu) &= \big\{b \in {\mathcal{U}}_{{\mathcal{T}}_2}(\mu) :
\ev_{\tilde{1}}(b) = \ev_{\hat{1}}(b)\big\},\\
{\mathcal{V}}_1^{(2)}(\mu) &=\big\{b \in {\mathcal{U}}_{{\mathcal{T}}_2}^{(1)}(\mu) :
\ev_{\tilde{1}}(b) = \ev_{\hat{2}}(b)\big\}.
\end{split}\end{equation*}
The cardinality of the last set is clearly six times
the number of rational curves that have a triple point and 
pass through the tuple $\mu$ of points and lines in~$\PPP$.

Let ${\overV}_1^{(1)}(\mu)$ and ${\overU}_{{\mathcal{T}}_2}^{(1)}(\mu)$
denote the closures of the space ${\mathcal{V}}_1^{(1)}(\mu)$ 
in ${\overU}_{{\mathcal{T}}_1}(\mu)$ and
of the space ${\mathcal{U}}_{{\mathcal{T}}_2}^{(1)}(\mu)$ in
${\overU}_{{\mathcal{T}}_2}(\mu)$, respectively.
In the next subsection, we describe the boundary of the set
${\overU}_{{\mathcal{T}}_2}^{(1)}(\mu)$ and conclude that 
${\mathcal{U}}_{{\mathcal{T}}_2}^{(1)}(\mu)$ is a 3-pseudovariety 
in~${\overU}_{{\mathcal{T}}_2}$.
Thus, the map
\begin{equation}\label{n3p3_sume1}
\ev_{\tilde{1}} \times \ev_{\hat{2}}\co 
{\mathcal{U}}_{{\mathcal{T}}_2}^{(1)}(\mu)\lra\PPP \times \PPP,\quad
\big\{\ev_{\tilde{1}} \times \ev_{\hat{2}}\big\}(b)=
\big(\ev_{\tilde{1}}(b),\ev_{\hat{2}}(b)\big),
\end{equation}
is a 6-pseudocycle in $\PPP \times \PPP$ in 
the sense of \cite[Chapter~7]{MS} and \cite[Section~1]{RT}, 
ie the map~\eqref{n3p3_sume1} defines an element
in $H_6(\PPP \times \PPP;\Bbb{Z})$.
In particular, there is a well-defined homology-intersection number
\begin{equation*}\begin{split}
\LlRr{ {\mathcal{V}}_1^{(2)}(\mu)}
&\equiv \LlRr{
\{\ev_{\tilde{1}} \times \ev_{\hat{2}}\}^{-1}
(\De_{\PPP\times\PPP}),{\overU}_{{\mathcal{T}}_2}^{(1)}(\mu)}\\
&=\sum_{r+s=3}
\LlRr{\{\ev_{\tilde{1}} \times \ev_{\hat{2}}\}^{-1}
(H^r \times H^s),{\overU}_{{\mathcal{T}}_2}^{(1)}(\mu)},
\end{split}\end{equation*}
where $\De_{\PPP\times\PPP}$ and $H^r$ 
denote the diagonal in $\PPP \times \PPP$ and 
a linear subspace of complex dimension~$r$ in~$\PPP$, respectively.
However, this number is \textit{not} $|{\mathcal{V}}_1^{(2)}(\mu)|$ in general.
Since the map $\ev_{\tilde{1}} \times \ev_{\hat{2}}$
is transversal to a generic submanifold $H^r \times H^s$
of $\PPP \times \PPP$, by definition,
\begin{equation}\label{n3p3_sume7}
\LlRr{ {\mathcal{V}}_1^{(2)}(\mu)}=
\big|{\mathcal{V}}_1^{(1)}(\mu + H^0)\big|
+\big\lan a_{\hat{0}},{\overV}_1^{(1)}(\mu + H^1)\big\ran
+d\big\lan a_{\hat{0}}^2,{\overV}_1^{(1)}(\mu)\big\ran,
\end{equation}
where the spaces ${\mathcal{V}}_1^{(1)}(\mu + H^0)$ 
and ${\overV}_1^{(1)}(\mu + H^1)$
are defined as above, but with $p$ replaced by $p + 1$
in the first case and with $q$ replaced by $q + 1$
in the second~case.

The number $\llrr{{\mathcal{V}}_1^{(2)}(\mu)}$ can also be obtained
by perturbing the map \hbox{$\ev_{\tilde{1}} \times \ev_{\hat{2}}$}.
If 
$$\th\co{\overU}_{{\mathcal{T}}_2}^{(1)}(\mu)\lra \PPP \times \PPP$$
is a small perturbation of $\ev_{\tilde{1}} \times \ev_{\hat{2}}$
such that 
$$\th^{-1}(\De_{\PPP\times\PPP})\cap\partial{\overU}_{{\mathcal{T}}_2}^{(1)}(\mu)
=\eset$$ 
and
$\th|{\mathcal{U}}_{{\mathcal{T}}_2}^{(1)}(\mu)$ is smooth and transversal to
$\De_{\PPP\times\PPP}$, then
$$\llrr{ {\mathcal{V}}_1^{(2)}(\mu)} =^{\pm} |\th^{-1}(\De_{\PPP\times\PPP})|.$$
Since ${\mathcal{U}}_{{\mathcal{T}}_2}^{(1)}(\mu) = \eset$ if $d = 1$, 
by Lemma~\ref{str_lmm} the map $\ev_{\tilde{1}} \times \ev_{\hat{2}}$
is transversal to $\De_{\PPP\times\PPP}$ on~${\mathcal{U}}_{{\mathcal{T}}_2}^{(1)}(\mu)$.
Thus, $\th$ can be chosen so that
\hbox{$\th = \ev_{\tilde{1}} \times \ev_{\hat{2}}$} outside of
a very small neighborhood~$W$ 
of~$\partial{\overU}_{{\mathcal{T}}_2}^{(1)}(\mu)$.
Then,
\begin{equation}\label{n3p3_sume9}\begin{split}
\LlRr{{\mathcal{V}}_1^{(2)}(\mu)}
&=\, ^{\pm} \big|\th^{-1}(\De_{\PPP\times\PPP}\big|\\
&=\big|{\mathcal{V}}_1^{(2)}(\mu)\big|+
\, ^{\pm} \big|\{\ev_{\tilde{1}} \times \ev_{\hat{2}}\}^{-1}
(\De_{\PPP\times\PPP})\cap W\big|\\
&\equiv \big|{\mathcal{V}}_1^{(2)}(\mu)\big|+
{\mathcal{C}}_{\partial{\overU}_{{\mathcal{T}}_2}^{(1)}(\mu)}
\big(\ev_{\tilde{1}} \times \ev_{\hat{2}};\De_{\PPP\times\PPP}\big).
\end{split}\end{equation}
The last term above is the contribution to 
$\llrr{{\mathcal{V}}_1^{(2)}(\mu)}$ from the boundary of 
${\overU}_{{\mathcal{T}}_2}^{(1)}(\mu)$;
it is the analogue of the term
${\mathcal{C}}_{\partial{\overS}}(s)$ in Proposition~\ref{euler_prp}.
If $\ev_{\tilde{1}} \times \ev_{\hat{2}}$ maps a stratum ${\mathcal{Z}}$
of $\partial{\overU}_{{\mathcal{T}}_2}^{(1)}(\mu)$ into~$\De_{\PPP\times\PPP}$,
near ${\mathcal{Z}}$ the map $\ev_{\tilde{1}} \times \ev_{\hat{2}}$
can be modeled on a section of a bundle over~${\mathcal{Z}}$.
In Subsection~\ref{n3p3_corr_subs}, we use
the topological approach of Section~\ref{topology_sec}
to compute this contribution.
Theorem~\ref{n3p3_thm} follows from equations~\eqref{n3p3_sume7} 
and~\eqref{n3p3_sume9} and Corollary~\ref{n3p3_contr_crl}.

Before concluding this subsection, we formally define the space
${\overV}_2^{(1,1)}(\mu)$.
We do not need this space in this section, but it is used in the next section and 
it is natural to describe its structure along 
with the structure of the space~${\overV}_1^{(1)}(\mu)$.
Let $I_2 = I_2^-\sqcup I_2^+$ be the graded linearly ordered set such that
$I_2^- = \{\hat{0}\}$ is a one-element set and
$I_2^+ = \{\tilde{1},\tilde{2}\}$ is a two-element set.
If ${\mathcal{T}} = (M_2,I_2;j,\under{d})$ is a bubble type such that
$j_{\hat{1}} = \tilde{1}$ and  $j_{\hat{2}} = \tilde{2}$, put
$${\mathcal{U}}_{\mathcal{T}}^{(1)}= 
\big\{b \in {\mathcal{U}}_{\mathcal{T}} : 
\ev_{\hat{1}}(b) = \ev_{\hat{2}}(b)\big\}.$$
Let ${\overU}_{\mathcal{T}}^{(1)}$ be the closure of 
${\mathcal{U}}_{\mathcal{T}}^{(1)}$ in the space 
${\overU}_{\mathcal{T}}$ or equivalently 
in~${\overU}_{{\mathcal{T}},{\mathcal{T}}}$. 
We define ${\mathcal{V}}_2^{(1,1)}(\mu)$ and ${\overV}_2^{(1,1)}(\mu)$ 
to be the disjoint unions of the spaces 
${\mathcal{U}}_{\mathcal{T}}^{(1)}(\mu)$  and ${\overU}_{\mathcal{T}}^{(1)}(\mu)$,
respectively, taken over all bubble types ${\mathcal{T}}$ as above
such that $d_{\tilde{1}},d_{\tilde{2}} > 0$
and $d_{\tilde{1}} + d_{\tilde{2}} = d$.

\subsection[On the structure of UT2,VT2, and similar spaces]{On the structure of ${\overU}_{{\mathcal{T}}_2}^{(1)}(\mu)$,
${\overV}_1^{(1)}(\mu)$, and similar spaces}
\label{n3p3_str_subs}

In this subsection, we describe the closure
${\overU}_{{\mathcal{T}}_2}^{(1)}(\mu)$
of the space ${\mathcal{U}}_{{\mathcal{T}}_2}^{(1)}(\mu)$ 
in ${\overU}_{{\mathcal{T}}_2}(\mu)$, or equivalently, 
in~${\overU}_{{\mathcal{T}}_2}$.
The tuple~$\mu$ can be arbitrary, and, in fact, 
$\PPP$ can be replaced by any other projective space.
Lemmas~\ref{n3p3_str_lmm1} and~\ref{n3p3_str_lmm2} imply that
${\mathcal{U}}_{{\mathcal{T}}_2}^{(1)}(\mu)$ is a pseudovariety in 
${\overU}_{{\mathcal{T}}_2}$
if $\mu$ is as in the previous subsection
and is a pseudocycle in~general.
The two lemmas in particular describe the kinds of curves
that can appear in the limit of rational one-component nodal curves,
reproducing a  known result in algebraic geometry, but
in a fairly direct way.
More importantly, we obtain a description
of what happens in the limit on the finer level of stable~maps.
The analytic expansion (3b) of Proposition~\ref{str_prp}
plays a crucial role in the proof of the second lemma.
We conclude this subsection with Lemma~\ref{n3p3_str_lmm3},
which describes the structure of the space~${\overV}_2^{(1,1)}(\mu)$.

If ${\mathcal{T}} = (M_2,I;j,\under{d})$ is a bubble type such that
${\mathcal{T}} < {\mathcal{T}}_2$ and $i,l \in I^+ \cup M_2$, let
$$\chi_{\mathcal{T}}(i,l)=
\max\big(\chi_{{\mathcal{T}},i}(l),\chi_{{\mathcal{T}},l}(i)\big);$$
see Subsection~\ref{notation_subs}.
Note that by continuity of the map $\ev_{\tilde{1}} \times \ev_{\hat{1}}$,
$${\overU}_{{\mathcal{T}}_2}^{(1)}(\mu) \subset
\{\ev_{\tilde{1}} \times \ev_{\hat{1}}\}^{-1}(\De_{\PPP\times\PPP}).$$

\begin{figure}[ht!]\small
\begin{pspicture}(-1.1,-1.7)(10,.5)
\psset{unit=.4cm}
% 1st diagram
\pscircle[fillstyle=solid,fillcolor=lightgray](1,-2.5){1}
\cnode*(1,-3.5){0}{A1}\pscircle*(1,-3.5){.22}\rput(1.1,-4.4){$\tilde{1}$}
\pscircle[fillstyle=solid,fillcolor=gray](-.41,-1.09){1}
\pscircle*(.29,-1.79){.19}
\cnode*(.29,-.39){0}{B1}\pscircle*(.29,-.39){.2}\rput(.75,.05){$\hat{1}$}
\ncarc[nodesep=.35,arcangleA=-105,arcangleB=-50,ncurv=2]{<->}{A1}{B1}
% 2nd diagram
\pscircle(7,-2.5){1}
\cnode*(7,-3.5){0}{A2}\pscircle*(7,-3.5){.22}\rput(7.1,-4.4){$\tilde{1}$}
\pscircle[fillstyle=solid,fillcolor=gray](5.59,-1.09){1}
\pscircle[fillstyle=solid,fillcolor=gray](8.41,-1.09){1}
\pscircle*(6.29,-1.79){.19}\pscircle*(7.71,-1.79){.19}
\cnode*(9.11,-.39){0}{B2}\pscircle*(9.11,-.39){.2}\rput(9.65,.05){$\hat{1}$}
\ncarc[nodesep=.35,arcangleA=-60,arcangleB=-90,ncurv=1.2]{<->}{A2}{B2}
% 3rd diagram
\pscircle(17,-2.5){1}
\pscircle*(17,-3.5){.22}\rput(17.1,-4.4){$\tilde{1}$}
\pscircle[fillstyle=solid,fillcolor=gray](15.59,-1.09){1}
\cnode*(16.29,-1.79){0}{B3}\pscircle*(16.29,-1.79){.19}
\pscircle*(17.71,-1.79){.2}\rput(18.15,-1.7){$l$}
\pnode(19.5,-2.5){A3}\rput(19.5,-2.5){\footnotesize cusp}
\ncarc[nodesep=.3,arcangleA=90,arcangleB=30,ncurv=1]{->}{A3}{B3}
% 4th diagram
\pscircle(24,-2.5){1}
\pscircle*(24,-3.5){.22}\rput(24.1,-4.4){$\tilde{1}$}
\pscircle[fillstyle=solid,fillcolor=gray](22.59,-1.09){1}
\pscircle[fillstyle=solid,fillcolor=gray](25.41,-1.09){1}
\pscircle*(23.29,-3.21){.2}\rput(22.8,-3.1){$\hat{1}$}
\cnode*(23.29,-1.79){0}{B4a}\cnode*(24.71,-1.79){0}{B4b}
\pscircle*(23.29,-1.79){.19}\pscircle*(24.71,-1.79){.19}
\pnode(28,-2.5){A4}\rput(28,-2.5){\footnotesize tacnode}
\ncarc[nodesep=.35,arcangleA=90,arcangleB=50,ncurv=.8]{->}{A4}{B4a}
\ncarc[nodesep=.35,arcangleA=90,arcangleB=30,ncurv=1]{->}{A4}{B4b}
\end{pspicture}
\caption{Some boundary strata of ${\overV}_1^{(1)}(\mu)$}
\label{n3p3_str_fig1}
\end{figure}
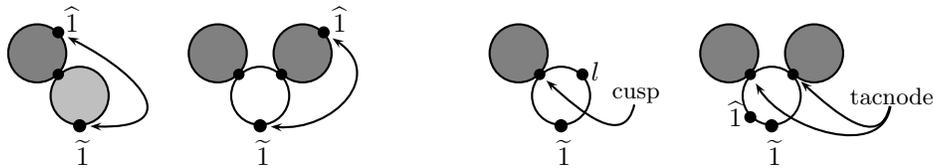

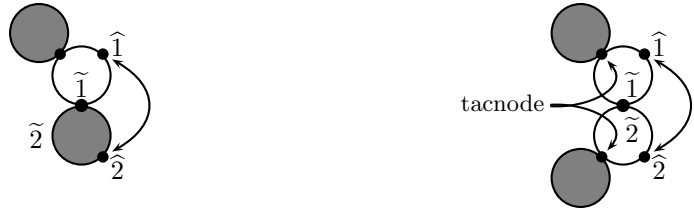
\begin{figure}[ht!]\small
\begin{pspicture}(.2,-3)(10,.5)
\psset{unit=.4cm}
\pscircle(8,-2.5){1}\rput(8,-2.8){$\tilde{1}$}
\pscircle[fillstyle=solid,fillcolor=gray](6.59,-1.09){1}
\pscircle*(7.29,-1.79){.19}
\cnode*(8.71,-1.79){0}{A1}\pscircle*(8.71,-1.79){.19}\rput(9.2,-1.45){$\hat{1}$}
                \pscircle[fillstyle=solid,fillcolor=gray](8,-4.5){1}
\rput(6.5,-4.5){$\tilde{2}$}
\pscircle*(8,-3.5){.22}
\cnode*(8.71,-5.21){0}{B1}\pscircle*(8.71,-5.21){.19}\rput(9.2,-5.55){$\hat{2}$}
\ncarc[nodesep=.35,arcangleA=60,arcangleB=60,ncurv=1.2]{<->}{A1}{B1}
\pscircle(26,-2.5){1}
\pscircle[fillstyle=solid,fillcolor=gray](24.59,-1.09){1}
\cnode*(25.29,-1.79){0}{B3a}\pscircle*(25.29,-1.79){.19}
\cnode*(26.71,-1.79){0}{A2}\pscircle*(26.71,-1.79){.19}\rput(27.2,-1.45){$\hat{1}$}
\pscircle(26,-4.5){1}
\pscircle[fillstyle=solid,fillcolor=gray](24.59,-5.91){1}
\cnode*(25.29,-5.21){0}{B3b}\pscircle*(25.29,-5.21){.19}
\pscircle*(26,-3.5){.22}\rput(26.3,-2.8){$\tilde{1}$}\rput(26.3,-4.3){$\tilde{2}$}
\cnode*(26.71,-5.21){0}{B2}\pscircle*(26.71,-5.21){.19}\rput(27.2,-5.55){$\hat{2}$}
\ncarc[nodesep=.35,arcangleA=60,arcangleB=60,ncurv=1.2]{<->}{A2}{B2}
\rput(22,-3.4){\footnotesize tacnode}
\cnode*(23.3,-3.5){0}{A3}
\ncarc[nodesep=.28,arcangleA=-45,arcangleB=-90,ncurv=1]{->}{A3}{B3a}
\ncarc[nodesep=.28,arcangleA=45,arcangleB=90,ncurv=1]{->}{A3}{B3b}
\end{pspicture}
\caption{Some boundary strata of ${\overV}_2^{(1,1)}(\mu)$}
\label{n3p3_str_fig2}
\end{figure}

Figures~\ref{n3p3_str_fig1} and \ref{n3p3_str_fig2} summarize
the three lemmas below.
All other boundary strata are either empty or
will be hollow with respect to all sections that we encounter.
The map may be constant or not on the disk shaded light~gray.
The lines connecting two marked points indicate that the map
has the same value at the two~points.

\begin{lmm}
\label{n3p3_str_lmm1}
If ${\mathcal{T}} = (M_2,I;j,\under{d})$ is a bubble type such that
${\mathcal{T}} < {\mathcal{T}}_2$ and  $\chi_{\mathcal{T}}(\tilde{1},\hat{1}) > 0$,
the~map 
$$\ev_{\tilde{1}} \times \ev_{\hat{1}}\co 
{\mathcal{U}}_{{\mathcal{T}}|{\mathcal{T}}_2}(\mu) \lra \PPP \times \PPP$$
is transversal to the diagonal~$\De_{\PPP\times\PPP}$.
Thus, 
$${\mathcal{U}}_{{\mathcal{T}}|{\mathcal{T}}_2}^{(1)}(\mu)\equiv
\{\ev_{\tilde{1}} \times \ev_{\hat{1}}\}^{-1}(\De_{\PPP\times\PPP})
\cap{\mathcal{U}}_{{\mathcal{T}}_2,{\mathcal{T}}}(\mu)$$ 
is a smooth submanifold of ${\mathcal{U}}_{{\mathcal{T}}|{\mathcal{T}}_2}(\mu)$
of dimension less than the dimension of ${\mathcal{U}}_{{\mathcal{T}}_2}^{(1)}(\mu)$
with normal bundle isomorphic to~$\ev_{\tilde{1}}^*T\PPP$.
\end{lmm}

\begin{proof}
The first statement is immediate from Lemma~\ref{str_lmm}.
The second claim follows from the first, since the dimension of
${\mathcal{U}}_{{\mathcal{T}}|{\mathcal{T}}_2}(\mu)$ is less than the
dimension of~${\mathcal{U}}_{{\mathcal{T}}_2}(\mu)$.
\end{proof}

\begin{lmm}
\label{n3p3_str_lmm2}
If ${\mathcal{T}}$ is as in Lemma~\ref{n3p3_str_lmm1}, but
$\chi_{\mathcal{T}}(\tilde{1},\hat{1}) = 0$,
for every $h \in \hat{I}^+$ and $k \in \Bbb{Z}$
there exists a $C^1$--negligible map
$$\ve_{{\mathcal{T}},\hat{1};h}^{(k)}\co 
\mathcal{FT}_{\de} - Y(\mathcal{FT};\hat{I}^+)\lra
\hbox{Hom}(L_h^{\otimes k}{\mathcal{T}},\ev_{\tilde{1}}^*T\PPP),$$
where $\de$ is as in Proposition~\ref{str_prp},
such that with notation as in Proposition~\ref{str_prp}
and with appropriate identifications,
\begin{equation*}\begin{split}
&\big\{\ev_{\tilde{1}} \times \ev_{\hat{1}}\big\}
\big(\phi_{{\mathcal{T}}_2,{\mathcal{T}}}(\ups)\big) \\
&\qquad\qquad\quad
=\sum_{k=1}^{\i} \sum_{{}~h\in\chi_{\tilde{1}}({\mathcal{T}})}  
    \big(y_{h;\hat{1}}(\ups) - x_{\hat{1};h}(\ups)\big)^{-k}
\Big\{{\mathcal{D}}_{{\mathcal{T}},h}^{(k)} +
 \ve_{{\mathcal{T}},\hat{1};h}^{(k)}(\ups)\Big\}
\rho_{{\mathcal{T}},\hat{1};h}^{(0;k)}(\ups)
\end{split}\end{equation*}
for all $\ups \in \mathcal{FT}_{\de} - Y(\mathcal{FT};\hat{I}^+)$.
Thus, ${\overU}_{{\mathcal{T}}_2}^{(1)}(\mu)\cap{\mathcal{U}}_{{\mathcal{T}}_2,{\mathcal{T}}}
\subset{\mathcal{S}}_{{\mathcal{T}}|{\mathcal{T}}_2}(\mu)$, where
\begin{multline*}
{\mathcal{S}}_{{\mathcal{T}}|{\mathcal{T}}_2}(\mu)
= \biggl\{b \in {\mathcal{U}}_{{\mathcal{T}}|{\mathcal{T}}_2}(\mu) :
\sum_{h\in\chi_{\tilde{1}}({\mathcal{T}})}      
{\mathcal{D}}_{{\mathcal{T}},h}\ups_h=0
\begin{array}{l}
\mbox{for some}~ \ups_h \in L_h{\mathcal{T}}\\
\mbox{with}~\big(\ups_h\big)_{h\in\chi_{\tilde{1}}({\mathcal{T}})}
 \neq 0\end{array}\biggr\}.
\end{multline*}
In particular, ${\overU}_{{\mathcal{T}}_2}^{(1)}(\mu)
\cap{\mathcal{U}}_{{\mathcal{T}}_2,{\mathcal{T}}}$
is contained in a finite union of
smooth submanifolds of ${\mathcal{U}}_{{\mathcal{T}}_2,{\mathcal{T}}}$ 
of dimension less than the dimension of~${\mathcal{U}}_{{\mathcal{T}}_2}^{(1)}(\mu)$.
\end{lmm}

\begin{proof}
In this case, we choose a specific identification
of small neighborhoods of $\De_{\PPP\times\PPP}$ in
$\pi_{\tilde{1}}^*T\PPP$ and in $\PPP \times \PPP$:
$$\big((x,x),(0,w)\big)\lra \big(x,\exp_{x,x}w),$$
where $\exp_{x,\cdot}$ denotes the exponential map with respect
to the metric~$g_{\PPP,x}$; see Lemma~\ref{flat_metrics}.
Since $\chi_{\mathcal{T}}(\tilde{1},\hat{1}) = 0$,
$\ev_{\tilde{1}}(b) = \ev_{\hat{1}}(b)$ for all
$b \in {\mathcal{U}}_{{\mathcal{T}}_2,{\mathcal{T}}}$, and thus
$g_{b,\tilde{1}} = g_{b,\hat{1}}$ for all
$b \in {\mathcal{U}}_{{\mathcal{T}}_2,{\mathcal{T}}}$.
The above expression for $\big\{\ev_{\tilde{1}} \times \ev_{\hat{1}}\big\}
\big(\phi_{{\mathcal{T}}_2,{\mathcal{T}}}(\ups)\big)$ is then simply the ``difference''
between the values of $\ev_{\tilde{1}}$ and $\ev_{\hat{1}}$
at~$\phi_{{\mathcal{T}}_2,{\mathcal{T}}}(\ups)$, which is computable 
from (3b) of Proposition~\ref{str_prp}:
\begin{alignat}{1}
&\big\{\ev_{\tilde{1}} \times \ev_{\hat{1}}\big\}
\big(\phi_{{\mathcal{T}}_2,{\mathcal{T}}}(\ups)\big)
=\sum_{m=1}^{\i}y_{\tilde{1};\hat{1}}(\ups)^{-m}
{\mathcal{D}}_{{\mathcal{T}}_2,\tilde{1}}^{(m)}\phi_{\tilde{\mathcal{T}},{\mathcal{T}}}(\ups) \notag\\
&~=\sum_{k=1}^{\i}\!\sum_{{}~h\in\chi_{\tilde{1}}({\mathcal{T}})}   
\!\!\!\bigg(\sum_{m=k}^{\i}
\!\binom{m{-}1}{k{-}1}
  y_{\tilde{1};\hat{1}}(\ups)^{-m}x_{\tilde{1};h}(\ups)^{m-k}\bigg)
\!\Big\{{\mathcal{D}}_{{\mathcal{T}},h}^{(k)} + \ve_{{\mathcal{T}},\hat{1};h}^{(k)}
(\ups)\Big\}\rho_{{\mathcal{T}},\tilde{1};h}^{(0;k)}(\ups) \notag\\
\label{n3p3_str_lmm2e3}
&~=\sum_{k=1}^{\i} \sum_{{}~h\in\chi_{\tilde{1}}({\mathcal{T}})}  
    \big(y_{\tilde{1};\hat{1}}(\ups) - x_{\tilde{1};h}(\ups)\big)^{-k}
\Big\{{\mathcal{D}}_{{\mathcal{T}},h}^{(k)} +
 \ve_{{\mathcal{T}},\hat{1};h}^{(k)}(\ups)\Big\}
\rho_{{\mathcal{T}},\tilde{1};h}^{(0;k)}(\ups).
\end{alignat}
Note that $\rho_{{\mathcal{T}},\tilde{1};h}^{(m;k)} = 
\rho_{{\mathcal{T}},\tilde{1};h}^{(0;k)}$
for all $m \in \Bbb{Z}$.
The last expression in~\eqref{n3p3_str_lmm2e3} is the same 
as the right-hand side of the expansion in the statement of the lemma;
see Subsection~\ref{notation_subs}.
This sum is absolutely convergent for all $\de$ sufficiently small,
since there exists
$C \in C({\mathcal{U}}_{{\mathcal{T}}|{\mathcal{T}}_2};\Bbb{R}^+)$
such~that
\begin{gather*}
\big|y_{h;\hat{1}}(b;\ups) - x_{\hat{1};h}(b;\ups)\big|^{-1}\le C(b)
\qquad\hbox{and}\\
\big|\rho_{{\mathcal{T}},\hat{1};h}^{(0;k)}(\ups)\big|\le|\ups|
\qquad\mbox{for all}\  h \in \chi_{\tilde{1}}({\mathcal{T}}),~~
(b;\ups) \in \mathcal{FT}.
\end{gather*}
The first inequality is immediate from the definitions of
$y_{h;\hat{1}}$ and $x_{\hat{1};h}$, while the second follows
from the assumption $\chi_{\mathcal{T}}(\tilde{1},\hat{1}) = 0$.
The above expansion of 
$\big\{\ev_{\ti{1}} \times \ev_{\hat{1}}\big\} \circ 
\phi_{{\mathcal{T}}_2,{\mathcal{T}}}$ immediately implies that
$${\overU}_{{\mathcal{T}}_2}^{(1)}(\mu)\cap{\mathcal{U}}_{{\mathcal{T}}_2,{\mathcal{T}}}
\subset {\mathcal{S}}_{{\mathcal{T}}|{\mathcal{T}}_2}(\mu).$$
In fact, the opposite inclusion also holds, as can be seen from Lemma~\ref{str_lmm}
and the Contraction Principle.
The remaining claim of the lemma is obtained by simple dimension-counting
from Lemma~\ref{str_lmm}.
\end{proof}

\begin{lmm}
\label{n3p3_str_lmm3}
Suppose $\tilde{\mathcal{T}} = (M_2,I_2;\tilde{j},\tilde{d})$ and
${\mathcal{T}} = (M_2,I;j,d)$ are bubble types such that
$\tilde{j}_{\hat{1}} = \tilde{1}$, 
$\tilde{j}_{\hat{2}} = \tilde{2}$, and
${\mathcal{T}} < \tilde{\mathcal{T}}$.

{\rm(1)}\qua If $d_h > 0$ for some $h \in I^+$ such that
$h \le j_{\hat{1}}$ or $h \le j_{\hat{2}}$,
the~map
$$\ev_{\tilde{1}} \times \ev_{\hat{1}}\co 
{\mathcal{U}}_{{\mathcal{T}}|\tilde{\mathcal{T}}}(\mu) \lra \PPP \times \PPP$$
is transversal to the diagonal~$\De_{\PPP\times\PPP}$.
Thus, 
$${\mathcal{U}}_{{\mathcal{T}}|\tilde{\mathcal{T}}}^{(1)}(\mu)\equiv
\{\ev_{\tilde{1}} \times \ev_{\hat{1}}\}^{-1}(\De_{\PPP\times\PPP})
\cap{\mathcal{U}}_{{\mathcal{T}}|\tilde{\mathcal{T}}}(\mu)$$ 
is a smooth submanifold of ${\mathcal{U}}_{{\mathcal{T}}|\tilde{\mathcal{T}}}(\mu)$
with normal bundle isomorphic to~$\ev_{\tilde{1}}^*T\PPP$.

{\rm(2)}\qua If $d_h > 0$ for all $h \in I^+$ such that
$h \le j_{\hat{1}}$ or $h \le j_{\hat{2}}$,
\begin{multline*}
\big\{\ev_{\tilde{1}} \times \ev_{\hat{1}}\big\}
\big(\phi_{\tilde{\mathcal{T}},{\mathcal{T}}}(\ups)\big)\\
=\sum_{i=1,2}  (-1)^i\sum_{k=1}^{\i} 
\sum_{{}~h\in\chi_{\tilde{i}}({\mathcal{T}})}  
    \big(y_{h;\hat{i}}(\ups) - x_{\hat{i};h}(\ups)\big)^{-k}
\Big\{{\mathcal{D}}_{{\mathcal{T}},h}^{(k)} +
 \ve_{{\mathcal{T}},\hat{i};h}^{(k)}(\ups)\Big\}
\rho_{{\mathcal{T}},\hat{i};h}^{(0;k)}(\ups)
\end{multline*}
for all $\ups \in \mathcal{FT}_{\de} - Y(\mathcal{FT};\hat{I}^+)$.
Thus, ${\overU}_{\tilde{\mathcal{T}}}^{(1)}(\mu)
\cap{\mathcal{U}}_{\tilde{\mathcal{T}},{\mathcal{T}}}
\subset{\mathcal{S}}_{{\mathcal{T}}|\tilde{\mathcal{T}}}(\mu)$, where
\begin{multline*}
{\mathcal{S}}_{{\mathcal{T}}|\tilde{\mathcal{T}}}(\mu)
= \biggl\{b \in {\mathcal{U}}_{{\mathcal{T}}|\tilde{\mathcal{T}}}(\mu) :
\sum_{h\in\chi_{\tilde{1}}({\mathcal{T}})
\cup\chi_{\tilde{2}}({\mathcal{T}})}           
\!\!\!\!\!\!{\mathcal{D}}_{{\mathcal{T}},h}\ups_h=0\\
\mbox{for some}\ \ups_h \in L_h{\mathcal{T}}
\ \mbox{such that}\ \big(\ups_h\big)_{h\in\chi_{\tilde{1}}\cup
\chi_{\tilde{2}}({\mathcal{T}})} \neq 0\biggr\}.
\end{multline*}
In either case, 
${\overU}_{\tilde{\mathcal{T}}}^{(1)}(\mu)\cap
{\mathcal{U}}_{\tilde{\mathcal{T}},{\mathcal{T}}}$
is contained in a finite union of 
smooth submanifolds of ${\mathcal{U}}_{\tilde{\mathcal{T}},{\mathcal{T}}}$
of dimension less than the dimension 
of~${\mathcal{U}}_{\tilde{\mathcal{T}}}^{(1)}(\mu)$.
\end{lmm}

\begin{proof}
The proof is essentially the same as that of Lemmas~\ref{n3p3_str_lmm1}
and~\ref{n3p3_str_lmm2}.
The only change is that in the second case 
we first obtain expansions for $\ev_{\hat{1}} - \ev_{\tilde{1}}$ and
$\ev_{\hat{2}} - \ev_{\tilde{2}}$
and then take their difference.
\end{proof}

\subsection[Behavior of the map ev1xev2 near dT2]{Behavior of the map
$\mathrm{ev}_{\tilde{1}}\times\mathrm{ev}_{\hat{2}}$ near
$\partial\bar{{\mathcal{U}}}_{{{\mathcal{T}}}_2}^{(1)}(\mu)$}
%$\ev_{\tilde{1}}\times\ev_{\hat{2}}$ near}
% $\partial\bar{{\mathcal{U}}}_{{{\mathcal{T}}}_2}^{(1)}(\mu)$}
\label{n3p3_exp_subs}

In this subsection, we use Lemma~\ref{str_lmm} and 
Proposition~\ref{str_prp} to describe the behavior of 
$\ev_{\tilde{1}} \times \ev_{\hat{2}}$
near the boundary of the space~${\overU}_{{\mathcal{T}}_2}^{(1)}(\mu)$.
We assume that $\mu$ is a tuple of points and lines in $\PPP$
as in Subsection~\ref{n3p3_sum_subs}.

\begin{lmm}
\label{n3p3_exp_lmm1}
If ${\mathcal{T}} = (M_2,I;j,\under{d})$ is a bubble type such that
${\mathcal{T}} < {\mathcal{T}}_2$,
$\chi_{\mathcal{T}}(\tilde{1},\hat{2}) > 0$, and
$\chi_{\mathcal{T}}(\hat{1},\hat{2}) > 0$,
$$\big\{\ev_{\tilde{1}} \times \ev_{\hat{2}}\big\}^{-1}
(\De_{\PPP\times\PPP})\cap
\big({\overU}_{{\mathcal{T}}_2}^{(1)}(\mu)
\cap{\mathcal{U}}_{{\mathcal{T}}_2,{\mathcal{T}}}\big)
=\eset.$$
\end{lmm}

\proof
(1)\qua If $\chi_{\mathcal{T}}(\tilde{1},\hat{1}) > 0$, by 
Lemma~\ref{n3p3_str_lmm1},
$${\overU}_{{\mathcal{T}}_2}^{(1)}(\mu)\cap{\mathcal{U}}_{{\mathcal{T}}_2,{\mathcal{T}}}
\subset{\mathcal{U}}_{{\mathcal{T}}|{\mathcal{T}}_2}^{(1)}(\mu).$$
Since every degree-one map into $\PPP$ is injective, the~map
$$\ev_{\tilde{1}} \times \ev_{\hat{2}}\co 
{\mathcal{U}}_{{\mathcal{T}}|{\mathcal{T}}_2}^{(1)}(\mu)\lra
\PPP \times \PPP$$
is transversal to $\De_{\PPP\times\PPP}$ by Lemma~\ref{str_lmm}.
Since the complex dimension of ${\mathcal{U}}_{{\mathcal{T}}|{\mathcal{T}}_2}^{(1)}(\mu)$
is less than three by Lemma~\ref{n3p3_str_lmm1}, it follows~that
$$\big\{\ev_{\tilde{1}} \times \ev_{\hat{2}}\big\}^{-1}
(\De_{\PPP\times\PPP})\cap
\big({\overU}_{{\mathcal{T}}_2}^{(1)}(\mu)
\cap{\mathcal{U}}_{{\mathcal{T}}_2,{\mathcal{T}}}\big)=\eset.$$

(2)\qua If $\chi_{\mathcal{T}}(\tilde{1},\hat{1}) = 0$, by
Lemma~\ref{n3p3_str_lmm2},
$${\overU}_{{\mathcal{T}}_2}^{(1)}(\mu)\cap{\mathcal{U}}_{{\mathcal{T}}|{\mathcal{T}}_2}
\subset{\mathcal{S}}_{{\mathcal{T}}_2,{\mathcal{T}}}(\mu).$$
Since every degree-one map into $\PPP$ is an immersion, the~map
$$\ev_{\tilde{1}} \times \ev_{\hat{2}}\co 
{\mathcal{S}}_{{\mathcal{T}}|{\mathcal{T}}_2}(\mu)\lra
\PPP \times \PPP$$
is transversal to $\De_{\PPP\times\PPP}$ by Lemma~\ref{str_lmm}.
Since the complex dimension of ${\mathcal{S}}_{{\mathcal{T}}|{\mathcal{T}}_2}(\mu)$
is less than three by Lemma~\ref{n3p3_str_lmm2}, it follows~that
$$\big\{\ev_{\tilde{1}} \times \ev_{\hat{2}}\big\}^{-1}
(\De_{\PPP\times\PPP})\cap 
\big({\overU}_{{\mathcal{T}}_2}^{(1)}(\mu)\cap{\mathcal{U}}_{{\mathcal{T}}_2,{\mathcal{T}}}\big)
=\eset.\eqno{\qed}$$

\begin{lmm}
\label{n3p3_exp_lmm2}
If ${\mathcal{T}} \equiv (M_2,I;j,\under{d}) < {\mathcal{T}}_2$ 
is a bubble type such that
$\chi_{\mathcal{T}}(\tilde{1},\hat{2}) = 0$, 
$\chi_{\mathcal{T}}(\hat{1},\hat{2}) > 0$, and
\hbox{${\mathcal{U}}_{{\mathcal{T}}|{\mathcal{T}}_2}^{(1)}(\mu) \neq \eset$}, then
$${\overU}_{{\mathcal{T}}_2}^{(1)}(\mu)\cap{\mathcal{U}}_{{\mathcal{T}}_2,{\mathcal{T}}}
 \subset \big\{\ev_{\tilde{1}} \times \ev_{\hat{2}}\big\}^{-1}
(\De_{\PPP\times\PPP}),\quad
\chi_{\mathcal{T}}(\tilde{1},\hat{1}) > 0, \quad
|\chi_{\tilde{1}}({\mathcal{T}})| \in \{1,2\}.$$
Furthermore, there exist a rank-$|\chi_{\tilde{1}}({\mathcal{T}})|$
vector bundle 
$\tilde{\mathcal{F}}{\mathcal{T}}  \lra {\mathcal{U}}_{{\mathcal{T}}_2,{\mathcal{T}}}$,
a monomials map \hbox{$\rho\co \mathcal{FT}  \lra \tilde{\mathcal{F}}{\mathcal{T}}$},
a section $\al \in \Ga\big({\mathcal{U}}_{{\mathcal{T}}_2,{\mathcal{T}}};
\hbox{Hom}(\tilde{\mathcal{F}}{\mathcal{T}};\ev_{\tilde{1}}^*T\PPP)\big)$,
and a $C^0$--negligible map 
$$\ve\co \mathcal{FT} - Y(\mathcal{FT};\hat{I}^+)\lra
\hbox{Hom}(\tilde{\mathcal{F}}{\mathcal{T}},\ev_{\tilde{1}}^*T\PPP)$$
such~that
$$\big\{\ev_{\tilde{1}} \times \ev_{\hat{2}}\big\}
\big(\phi_{{\mathcal{T}}_2,{\mathcal{T}}}(\ups)\big)
=\{\al + \ve(\ups)\big\}\rho(\ups)
\qquad\mbox{for all}\ \ups \in \mathcal{FT}_{\de} - Y(\mathcal{FT};\hat{I}^+).$$
The vector-bundle map~$\al$ is injective over
${\mathcal{U}}_{{\mathcal{T}}|{\mathcal{T}}_2}^{(1)}(\mu)$.
Finally, if $\hat{I}^+ = \chi_{\tilde{1}}({\mathcal{T}})$,
$\rho$ is the identity map,~and
$$\al(\ups)=\sum_{h\in\chi_{\tilde{1}}({\mathcal{T}})}    
\big(y_{\hat{2}} - x_h\big)^{-1}
 \otimes{\mathcal{D}}_{{\mathcal{T}},h}^{(1)}\ups_h
\qquad\mbox{for all}\ \ups = (\ups_h)_{h\in\chi_{\tilde{1}}({\mathcal{T}})}\in\mathcal{FT}.$$
\end{lmm}

\begin{proof}
The first two claims of the first sentence are clear;
the third follows by dimension-counting from Lemma~\ref{str_lmm}.
On the other hand, 
equation~\eqref{n3p3_str_lmm2e3} with $\hat{1}$ replaced by~$\hat{2}$
gives
\begin{equation*}\begin{split}
& \big\{\ev_{\tilde{1}} \times \ev_{\hat{2}}\big\}
\big(\phi_{{\mathcal{T}}_2,{\mathcal{T}}}(\ups)\big) \\
&\qquad\quad=\sum_{{}~h\in\chi_{\tilde{1}}({\mathcal{T}})}    
\big(y_{h;\hat{2}}(\ups) - x_{\hat{2};h}(\ups)\big)^{-1}
\Big\{{\mathcal{D}}_{{\mathcal{T}},h}^{(1)} +
 \tilde{\ve}_{{\mathcal{T}},\hat{2};h}(\ups)\Big\}
\rho_{{\mathcal{T}},\hat{2};h}^{(0;1)}(\ups).
\end{split}\end{equation*}
Thus, we define the monomials map 
$\rho = (\rho_h)_{h\in\chi_{\tilde{1}}({\mathcal{T}})}$
and the linear map~$\al$~by:
$$\rho_h(\ups)=  
\prod_{i\in(i_{\mathcal{T}}(\hat{2},h),h]}       \ups_i,
\qquad
\al\big((\tilde{\ups}_h)_{h\in\chi_{\tilde{1}}({\mathcal{T}})}\big)
=  \sum_{h\in\chi_{\tilde{1}}({\mathcal{T}})}    
\big(y_{h;\hat{2}} - x_{\hat{2};h}\big)^{-1}
 \otimes{\mathcal{D}}_{{\mathcal{T}},h}^{(1)}\tilde{\ups}_h,$$
where
\begin{gather*}
x_{\hat{2};h} = x_{h'}\quad\hbox{if~~}
h' \in (i_{\mathcal{T}}(\hat{2},h);h]\hbox{~~and~~}
\io_{h'}^+=i_{\mathcal{T}}(\hat{2},h);\\
y_{h;\hat{2}}=
\begin{cases}
y_{\hat{2}},&\hbox{if}~j_{\hat{2}} = i_{\mathcal{T}}(\hat{2},h);\\
x_{h'},&\hbox{if}~h' \in (i_{\mathcal{T}}(\hat{2},h);j_{\hat{2}}]
\hbox{~and~}\io_{h'}^+=i_{\mathcal{T}}(\hat{2},h).
\end{cases}
\end{gather*}
We write $b = (M,I;x,(j,y),u)$ as before;
then $x_h$ and $y_l$ are sections of a bundle 
over~${\mathcal{U}}_{{\mathcal{T}}_2,|{\mathcal{T}}}$.
By Lemma~\ref{str_lmm}, the map $\al$ is injective 
over~${\mathcal{U}}_{{\mathcal{T}}|{\mathcal{T}}_2}^{(1)}(\mu)$.
\end{proof}

\begin{lmm}
\label{n3p3_exp_lmm3}
If ${\mathcal{T}} = (M_2,I;j,\under{d})$ is a bubble type such that
${\mathcal{T}} < {\mathcal{T}}_2$,
$$\chi_{\mathcal{T}}(\tilde{1},\hat{1}) = \chi_{\mathcal{T}}(\tilde{1},\hat{2})
= \chi_{\mathcal{T}}(\hat{1},\hat{2}) = 0,$$
and ${\mathcal{S}}_{{\mathcal{T}}|{\mathcal{T}}_2}(\mu) \neq \eset$,
$|\chi_{\tilde{1}}({\mathcal{T}})| \in \{1,2\}$.
Furthermore, the following properties~hold.

{\rm(1)}\qua There exist a rank-$|\chi_{\tilde{1}}({\mathcal{T}})|$
vector bundle 
$\tilde{\mathcal{F}}{\mathcal{T}} \lra {\mathcal{U}}_{{\mathcal{T}}|{\mathcal{T}}_2}$,
section 
$$\al \in \Ga\big({\mathcal{U}}_{{\mathcal{T}}|{\mathcal{T}}_2};
\hbox{Hom}(\tilde{\mathcal{F}}{\mathcal{T}} ,\ev_{\tilde{1}}^*T\PPP)\big),$$
and a monomials map and $C^0$--negligible map
$$\rho,\ve\co \mathcal{FT} - Y(\mathcal{FT};\hat{I}^+)\lra\tilde{\mathcal{F}}{\mathcal{T}}, \,
\hbox{Hom}(\tilde{\mathcal{F}}{\mathcal{T}},\ev_{\tilde{1}}^*T\PPP)$$
such~that 
$$\big\{\ev_{\tilde{1}} \times \ev_{\hat{2}}\big\}
\big(\phi_{{\mathcal{T}}_2,{\mathcal{T}}}(\ups)\big)
=\{\al + \ve(\ups)\}\rho(\ups)$$
for all $\ups \in \mathcal{FT}_{\de}$ such that
$\phi_{{\mathcal{T}}_2,\cal T}(\ups) \in
{\mathcal{U}}_{{\mathcal{T}}_2}^{(1)}$.

{\rm(2)}\qua If $\chi_{\tilde{1}}({\mathcal{T}}) = \{h\}$ is a
single-element set, $\al$ is injective over
${\mathcal{S}}_{{\mathcal{T}}|{\mathcal{T}}_2}(\mu)$.
If, in addition,  $\hat{I}^+ = \chi_{\tilde{1}}({\mathcal{T}})$, then
$$\rho(\ups_h)=\ups \otimes \ups$$
for all $\ups \in \mathcal{FT}$ and
$$\al(\tilde{\ups}_h) =
\big(y_{\hat{2}} - x_h\big)^{-2} \otimes 
\big(y_{\hat{1}} - x_h\big)^{-1} \otimes 
\big(y_{\hat{2}} - y_{\hat{1}}\big) \otimes 
{\mathcal{D}}_{{\mathcal{T}},h}^{(2)}\tilde{\ups}_h$$
for all $\tilde{\ups} \in \tilde{\mathcal{F}}{\mathcal{T}} \equiv 
\mathcal{FT}^{\otimes2}$.

{\rm(3)}\qua If $|\chi_{\tilde{1}}({\mathcal{T}})| > 1$,
there exist a line bundle ${\mathcal{L}} \lra {\mathcal{U}}_{{\mathcal{T}}|{\mathcal{T}}_2}$,
a section
$$\al_+ \in \Ga\big({\mathcal{U}}_{{\mathcal{T}}|{\mathcal{T}}_2};
\hbox{Hom}(\tilde{\mathcal{F}}{\mathcal{T}},
{\mathcal{L}}^* \otimes \ev_{\tilde{1}}^*T\PPP)\big),$$
and  a monomials map and a $C^0$--negligible map
$$\rho_+,\ve_+\co \mathcal{FT} - Y(\mathcal{FT};\hat{I}^+) \lra 
{\mathcal{L}},\hbox{Hom}(\tilde{\mathcal{F}}{\mathcal{T}},
{\mathcal{L}}^* \otimes \ev_{\tilde{1}}^*T\PPP)$$
such~that 
$\al_+ \oplus \al$ is injective over ${\mathcal{S}}_{{\mathcal{T}}|{\mathcal{T}}_2}(\mu)$,
$$\al_+|Y(\tilde{\mathcal{F}}{\mathcal{T}},\{h_1\})  \qquad\hbox{and}\qquad
\al_+|Y(\tilde{\mathcal{F}}{\mathcal{T}},\{h_2\})$$ 
are onto~$\hbox{Im}(\al_+)$
over ${\mathcal{S}}_{{\mathcal{T}}|{\mathcal{T}}_2}(\mu)$, and
$$\big\{\ev_{\tilde{1}} \times \ev_{\hat{1}}\big\}
(\phi_{{\mathcal{T}}_2,{\mathcal{T}}}(\ups))
=\rho_+(\ups)\otimes\big\{\al_+ + \ve_+(\ups)\big\}\rho(\ups)$$
for all $\ups \in \mathcal{FT}_{\de} - Y(\mathcal{FT};\hat{I}^+)$.
If, in addition, $\hat{I}^+ = \chi_{\tilde{1}}({\mathcal{T}}) = \{h_1,h_2\}$,
$\rho$ is the identity map,~and
\begin{gather*}
\al_+(\ups)=\big(y_{\hat{1}} - x_{h_1}\big)^{-1}\otimes
{\mathcal{D}}_{{\mathcal{T}},h_1}^{(1)}\ups_{h_1}+
\big(y_{\hat{1}} - x_{h_2}\big)^{-1}\otimes
{\mathcal{D}}_{{\mathcal{T}},h_2}^{(1)}\ups_{h_2};\\
\begin{split}
\al(\ups)=\big(y_{\hat{1}} - x_{h_2}\big)^{-1} \otimes 
\big(y_{\hat{2}} - x_{h_1}\big)^{-1} \otimes 
\big(y_{\hat{2}} - x_{h_2}\big)^{-1}
&\otimes  \big(y_{\hat{2}} - y_{\hat{2}}\big)\\
&\otimes  \big(x_{h_1} - x_{h_2}\big) \otimes 
{\mathcal{D}}_{{\mathcal{T}},h_2}^{(1)}\ups_{h_2}
\end{split} \end{gather*}
for all $\ups = (\ups_{h_1},\ups_{h_2}) \in \mathcal{FT}_{\de}$.
\end{lmm}

\begin{proof}
(1)\qua The first statement of this lemma follows 
from Lemma~\ref{str_lmm} by dimension-counting.
If $\chi_{\tilde{1}}({\mathcal{T}}) = \{h\}$ is a single-element set,
the remaining claims are obtained by subtracting the expansion
of $\{\ev_{\tilde{1}} \times \ev_{\hat{1}}\}
 \circ \phi_{{\mathcal{T}}_2,{\mathcal{T}}}$
given in Lemma~\ref{n3p3_str_lmm2} times
$$\big(y_{\tilde{1};\hat{1}}(\ups) - x_{\tilde{1};h}(\ups)\big)
\big(y_{\tilde{1};\hat{2}}(\ups) - x_{\tilde{1};h}(\ups)\big)^{-1}$$
from the corresponding expression for 
$\{\ev_{\tilde{1}} \times \ev_{\hat{2}}\} \circ 
\phi_{{\mathcal{T}}_2,{\mathcal{T}}}$.

(2)\qua If $|\chi_{\tilde{1}}({\mathcal{T}})| > 1$,
$h \in \chi_{\tilde{1}}({\mathcal{T}})$, and $l = \hat{1},\hat{2}$, we put
\begin{gather*}
i_{\mathcal{T}}^*(h)=\max\big\{i_{\mathcal{T}}(h,\hat{1}),i_{\mathcal{T}}(h,\hat{2})\big\},\\
i_{\mathcal{T}}^+(l)=
\max\big\{i_{\mathcal{T}}(l,h) :h \in \chi_{\tilde{1}}({\mathcal{T}})\big\},\quad
i_{\mathcal{T}}^-(l)=
\min\big\{i_{\mathcal{T}}(l,h) :h \in \chi_{\tilde{1}}({\mathcal{T}})\big\}.
\end{gather*}
If $h_1 \in \chi_{\tilde{1}}({\mathcal{T}})$ is such that
either $i_{\mathcal{T}}(h_1,\hat{1}) = i_{\mathcal{T}}^*(h_1)$ or
$i_{\mathcal{T}}(h_1,\hat{2}) = i_{\mathcal{T}}^-(\hat{2})$,
we subtract the expansion
of $\{\ev_{\tilde{1}} \times \ev_{\hat{1}}\} \circ 
\phi_{{\mathcal{T}}_2,{\mathcal{T}}}$
given in Lemma~\ref{n3p3_str_lmm2} times
$$\big(y_{\tilde{1};\hat{1}}(\ups) - x_{\tilde{1};h_1}(\ups)\big)
\big(y_{\tilde{1};\hat{2}}(\ups) - x_{\tilde{1};h_1}(\ups)\big)^{-1}$$
from the corresponding expression for 
$\{\ev_{\tilde{1}} \times \ev_{\hat{2}}\} \circ \phi_{{\mathcal{T}}_2,{\mathcal{T}}}$
and then take the leading term.
\end{proof}

\subsection[Computation of the number C(evxev;D)]{Computation of the number 
${\mathcal{C}}_{\partial\overline{{\mathcal{U}}}_{{\mathcal{T}}_2}^{(1)}(\mu)}
  \big(\mathrm{ev}_{\tilde{1}}\times\mathrm{ev}_{\hat{2}};
  \De_{\PPP\times\PPP}\big)$}
%${\mathcal{C}}_{{\partial{\overU}}_{{\mathcal{T}}_2}^{(1)}(\mu)}$
%\big(\ev_{\tilde{1}} \times \ev_{\hat{2}};\De_{\PPP\times\PPP}\big)$}
\label{n3p3_corr_subs}

We are now ready to compute the term 
${\mathcal{C}}_{\partial{\overU}_{{\mathcal{T}}_2}^{(1)}(\mu)}
\big(\ev_{\tilde{1}} \times \ev_{\hat{2}};\De_{\PPP\times\PPP}\big)$
appearing in equation~\eqref{n3p3_sume9}.
We perturb the map $\ev_{\tilde{1}} \times \ev_{\hat{2}}$
to a new continuous map~$\th$ on ${\overU}_{{\mathcal{T}}_2}$
such that the image of 
$\partial{\overU}_{{\mathcal{T}}_2}^{(1)}(\mu)$
under~$\th$ is disjoint from~$\De_{\PPP\times\PPP}$
and $\th|{\mathcal{U}}_{{\mathcal{T}}_2}^{(1)}(\mu)$ is smooth and transversal
to~$\De_{\PPP\times\PPP}$.
In order to achieve these requirements,
it is sufficient to perturb $\ev_{\tilde{1}} \times \ev_{\hat{2}}$
very slightly on a small neighborhood~$W$ of 
$$\partial{\overU}_{{\mathcal{T}}_2}^{(1)}(\mu)\cap
\{\ev_{\tilde{1}} \times \ev_{\hat{2}}\}^{-1}(\De_{\PPP\times\PPP}).$$
Then, 
$${\mathcal{C}}_{\partial{\overU}_{{\mathcal{T}}_2}^{(1)}(\mu)}
\big(\ev_{\tilde{1}} \times \ev_{\hat{2}};\De_{\PPP\times\PPP}\big)
=\, ^{\pm} \big|\th^{-1}(\De_{\PPP\times\PPP})\cap W\big|.$$
Along the set $\partial{\overU}_{{\mathcal{T}}_2}^{(1)}(\mu)\cap
\{\ev_{\tilde{1}} \times \ev_{\hat{2}}\}^{-1}(\De_{\PPP\times\PPP})$,
the maps $\ev_{\tilde{1}} \times \ev_{\hat{2}}$ and $\th$
can be viewed as sections of the bundle~$\ev_{\tilde{1}}^*T\PPP$.
Thus, we can apply the terminology and the computational method of
Section~\ref{topology_sec} to determine the number of zeros
of a small perturbation of $\ev_{\tilde{1}} \times \ev_{\hat{2}}$
near
$$\partial{\overU}_{{\mathcal{T}}_2}^{(1)}(\mu)\cap
\{\ev_{\tilde{1}} \times \ev_{\hat{2}}\}^{-1}(\De_{\PPP\times\PPP}).$$
Of course, we cannot ``cut off'' the map near the entire 
boundary of~${\overU}_{{\mathcal{T}}_2}^{(1)}(\mu)$,
as was done for vector-bundle sections in Subsection~\ref{top_comp_subs}.
However, the entire approach of Subsection~\ref{top_comp_subs}
goes through, since 
$\{\ev_{\tilde{1}} \times \ev_{\hat{2}}\}^{-1}(\De_{\PPP\times\PPP})$
is well-defined on all of the space 
${\overU}_{{\mathcal{T}}_2}^{(1)}(\mu)$,
and not just on~${\mathcal{U}}_{{\mathcal{T}}_2}^{(1)}(\mu)$.

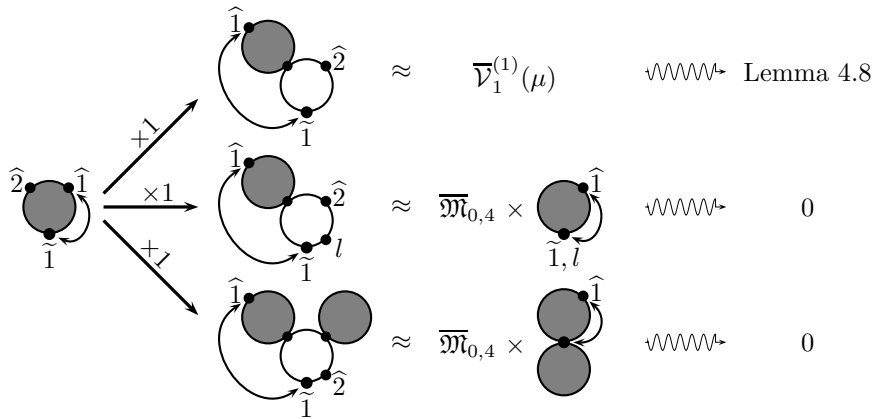
\begin{figure}[ht!]\small
\begin{pspicture}(-.7,-4)(10,2.7)
\psset{unit=.36cm}
\pscircle[fillstyle=solid,fillcolor=gray](1,-2){1}
\pscircle*(1,-3){.22}\rput(1,-3.8){$\tilde{1}$}
\pscircle*(.29,-1.29){.2}\pscircle*(1.71,-1.29){.2}
\rput(2.2,-1){$\hat{1}$}\rput(-.2,-1){$\hat{2}$}
\pnode(1,-3){A1}\pnode(1.71,-1.29){B1}
\ncarc[nodesep=.35,arcangleA=-90,arcangleB=-75,ncurv=1.2]{<->}{A1}{B1}
%2nd column starts here
\psline[linewidth=.12]{->}(3,-1.5)(6.5,2)
\rput{45}(4.5,.7){$\times1$}
\psline[linewidth=.12]{->}(3,-2)(6.5,-2)
\rput(5,-1.5){$\times1$}
\psline[linewidth=.12]{->}(3,-2.5)(6.5,-6)
\rput{-45}(4.9,-3.7){$\times1$}
%3rd column starts here
\pscircle(10.5,2.5){1}
\pscircle[fillstyle=solid,fillcolor=gray](9.08,3.92){1}
\pscircle*(10.5,1.5){.22}\rput(10.5,.7){$\tilde{1}$}
\pscircle*(9.79,3.21){.19}\pscircle*(11.21,3.21){.2}
\rput(11.7,3.5){$\hat{2}$}
\pscircle*(8.37,4.63){.2}\rput(7.88,4.92){$\hat{1}$}
\pnode(10.5,1.5){A2}\pnode(8.37,4.63){B2}
\ncarc[nodesep=.35,arcangleA=90,arcangleB=90,ncurv=1.2]{<->}{A2}{B2}
%2nd row starts here
\pscircle(10.5,-2.5){1}
\pscircle[fillstyle=solid,fillcolor=gray](9.08,-1.08){1}
\pscircle*(10.5,-3.5){.22}\rput(10.5,-4.3){$\tilde{1}$}
\pscircle*(9.79,-1.79){.19}
\pscircle*(11.21,-1.79){.2}\rput(11.7,-1.5){$\hat{2}$}
\pscircle*(8.37,-.37){.2}\rput(7.88,-.08){$\hat{1}$}
\pscircle*(11.21,-3.21){.2}\rput(11.7,-3.5){$l$}
\pnode(10.5,-3.5){A3}\pnode(8.37,-.37){B3}
\ncarc[nodesep=.35,arcangleA=90,arcangleB=90,ncurv=1.2]{<->}{A3}{B3}
%3rd row starts here
\pscircle(10.5,-7.5){1}
\pscircle[fillstyle=solid,fillcolor=gray](9.08,-6.08){1}
\pscircle*(10.5,-8.5){.22}\rput(10.5,-9.3){$\tilde{1}$}
\pscircle*(9.79,-6.79){.19}
\pscircle[fillstyle=solid,fillcolor=gray](11.92,-6.08){1}
\pscircle*(11.21,-6.79){.19}
\pscircle*(8.37,-5.37){.2}\rput(7.88,-5.08){$\hat{1}$}
\pscircle*(11.21,-8.21){.2}\rput(11.7,-8.5){$\hat{2}$}
\pnode(10.5,-8.5){A4}\pnode(8.37,-5.37){B4}
\ncarc[nodesep=.35,arcangleA=90,arcangleB=90,ncurv=1.2]{<->}{A4}{B4}
\rput(14,3){$\approx$}\rput(14,-2){$\approx$}\rput(14,-7){$\approx$}
\rput(18.2,2.9){${\overV}_1^{(1)}(\mu)$}
\rput(17,-2.1){$\ov{\frak M}_{0,4}~\times$}
\pscircle[fillstyle=solid,fillcolor=gray](20,-2){1}
\pscircle*(20,-3){.22}\rput(20,-3.8){$\tilde{1},l$}
\pscircle*(20.71,-1.29){.2}\rput(21.2,-1){$\hat{1}$}
\pnode(20,-3){A6}\pnode(20.71,-1.29){B6}
\ncarc[nodesep=.3,arcangleA=-90,arcangleB=-70,ncurv=1.2]{<->}{A6}{B6}
\rput(17,-7.1){$\ov{\frak M}_{0,4}~\times$}
\pscircle[fillstyle=solid,fillcolor=gray](20,-6){1}
\pscircle[fillstyle=solid,fillcolor=gray](20,-8){1}
\pscircle*(20,-7){.22}
\pscircle*(20.71,-5.29){.2}\rput(21.2,-5){$\hat{1}$}
\pnode(20.05,-7){A7}\pnode(20.71,-5.29){B7}
\ncarc[nodesep=.3,arcangleA=-70,arcangleB=-70,ncurv=1.2]{<->}{A7}{B7}
\pscoil[linewidth=.03,coilarmA=0.1,coilarmB=0.4,coilaspect=0,coilheight=.72,
coilwidth=.6]{->}(23,3)(26,3)
\pscoil[linewidth=.03,coilarmA=0.1,coilarmB=0.4,coilaspect=0,coilheight=.72,
coilwidth=.6]{->}(23,-2)(26,-2)
\pscoil[linewidth=.03,coilarmA=0.1,coilarmB=0.4,coilaspect=0,coilheight=.72,
coilwidth=.6]{->}(23,-7)(26,-7)
\rput(29,3){\small Lemma~\ref{n3p3_contr_lmm1b}}\rput(29,-2){$0$}\rput(29,-7){$0$}
\end{pspicture}
\caption{An outline of the proof of Lemma~\ref{n3p3_contr_lmm1}}
\label{n3p3_contr_fig1}
\end{figure}

\begin{figure}[ht!]\small
\begin{pspicture}(-1.1,-4)(10,2.7)
\psset{unit=.36cm}
\pscircle[fillstyle=solid,fillcolor=gray](1,-2){1}
\pscircle*(1,-3){.22}\rput(1,-3.8){$\tilde{1}$}
\pscircle*(.29,-1.29){.2}\pscircle*(1.71,-1.29){.2}
\rput(2.2,-1){$\hat{1}$}\rput(-.2,-1){$\hat{2}$}
\pnode(1,-3){A1}\pnode(1.71,-1.29){B1}
\ncarc[nodesep=.35,arcangleA=-90,arcangleB=-75,ncurv=1.2]{<->}{A1}{B1}
% 2nd column
\psline[linewidth=.12]{->}(3,-1.5)(6.5,2)
\rput{45}(4.5,.7){$\times2$}
\psline[linewidth=.12]{->}(3,-2)(6.5,-2)
\rput(5,-1.5){$\times2$}
\psline[linewidth=.12]{->}(3,-2.5)(6.5,-6)
\rput{-45}(4.9,-3.7){$\times1$}
% 3rd column
\pscircle(10,2.5){1}
\pscircle[fillstyle=solid,fillcolor=gray](8.58,3.92){1}
\pscircle*(10,1.5){.22}\rput(10,.7){$\tilde{1}$}
\pscircle*(9.29,3.21){.19}\pscircle*(10.71,3.21){.2}
\rput(11.2,3.5){$\hat{2}$}
\pscircle*(10.71,1.79){.2}\rput(11.2,1.5){$\hat{1}$}
%2nd row starts here
\pscircle(10,-2.5){1}
\pscircle[fillstyle=solid,fillcolor=gray](8.58,-1.08){1}
\pscircle*(10,-3.5){.22}\rput(10,-4.3){$\tilde{1}$}
\pscircle*(9.29,-1.79){.19}
\pscircle*(10.71,-1.79){.2}\rput(11.2,-1.5){$\hat{2}$}
\pscircle*(10.71,-3.21){.2}\rput(11.2,-3.5){$\hat{1}$}
\pscircle*(9.29,-3.21){.2}\rput(8.8,-3.5){$l$}
%3rd row starts here
\pscircle(10,-7.5){1}
\pscircle[fillstyle=solid,fillcolor=gray](8.58,-6.08){1}
\pscircle*(10,-8.5){.22}\rput(10,-9.3){$\tilde{1}$}
\pscircle*(9.29,-6.79){.19}
\pscircle[fillstyle=solid,fillcolor=gray](11.42,-6.08){1}
\pscircle*(10.71,-6.79){.19}
\pscircle*(10.71,-8.21){.2}\rput(11.2,-8.7){$\hat{2}$}
\pscircle*(9.29,-8.21){.2}\rput(8.8,-8.5){$\hat{1}$}
%4th column starts now
\rput(14,3){$\approx$}\rput(14,-2){$\approx$}\rput(14,-7){$\approx$}
\rput(14,1.5){\footnotesize cusp}\rput(14,-3.5){\footnotesize cusp}
\rput(14,-8.5){\footnotesize tacnode}
%1st row starts here
\rput(17.5,2.9){$\ov{\frak M}_{0,4}~\times$}
\pscircle[fillstyle=solid,fillcolor=gray](20.5,3){1}
\pscircle*(20.5,2){.22}\rput(20.5,1.2){$\tilde{1}$}
\pnode(13.3,1.5){A2a}\pnode(14.7,1.6){A2b}
\pnode(9.29,3.21){B2a}\pnode(20.5,2){B2b}
\ncarc[nodesep=.3,arcangleA=10,arcangleB=10,ncurv=.8]{->}{A2a}{B2a}
\ncarc[nodesep=.3,arcangleA=-10,arcangleB=-20,ncurv=.7]{->}{A2b}{B2b}
%2nd row starts here
\rput(17.5,-2.1){$\ov{\frak M}_{0,5}~\times$}
\pscircle[fillstyle=solid,fillcolor=gray](20.5,-2){1}
\pscircle*(20.5,-3){.22}\rput(20.6,-3.8){$\tilde{1},l$}
\pnode(13.3,-3.5){A3a}\pnode(14.7,-3.4){A3b}
\pnode(9.29,-1.79){B3a}\pnode(20.3,-3){B3b}
\ncarc[nodesep=.3,arcangleA=10,arcangleB=10,ncurv=.8]{->}{A3a}{B3a}
\ncarc[nodesep=.3,arcangleA=-10,arcangleB=-20,ncurv=.6]{->}{A3b}{B3b}
%3rd row starts here
\rput(17.5,-7.1){$\ov{\frak M}_{0,5}\times$}
\pscircle[fillstyle=solid,fillcolor=gray](20.5,-6){1}
\pscircle[fillstyle=solid,fillcolor=gray](20.5,-8){1}
\pscircle*(20.5,-7){.22}
\pnode(12.3,-8.5){A4a1}\pnode(10,-7.5){A4a2}\pnode(15.3,-8.6){A4b}
\pnode(9.44,-6.94){B4a1}\pnode(10.56,-6.94){B4a2}\pnode(20.4,-7){B4b}
\ncarc[nodesep=0,arcangleA=-10,arcangleB=0,ncurv=.1]{-}{A4a1}{A4a2}
\ncline[nodesep=0]{->}{A4a2}{B4a1}\ncline[nodesep=0]{->}{A4a2}{B4a2}
\ncarc[nodesep=.3,arcangleA=-35,arcangleB=15,ncurv=.2]{->}{A4b}{B4b}
%3rd row ends here
\pscoil[linewidth=.03,coilarmA=0.1,coilarmB=0.4,coilaspect=0,coilheight=.72,
coilwidth=.6]{->}(23,3)(26,3)
\pscoil[linewidth=.03,coilarmA=0.1,coilarmB=0.4,coilaspect=0,coilheight=.72,
coilwidth=.6]{->}(23,-2)(26,-2)
\pscoil[linewidth=.03,coilarmA=0.1,coilarmB=0.4,coilaspect=0,coilheight=.72,
coilwidth=.6]{->}(23,-7)(26,-7)
\rput(27,3){$0$}\rput(27,-2){$0$}\rput(27,-7){$0$}
\end{pspicture}
\caption{An outline of the proof of Lemma~\ref{n3p3_contr_lmm3}}
\label{n3p3_contr_fig2}
\end{figure}
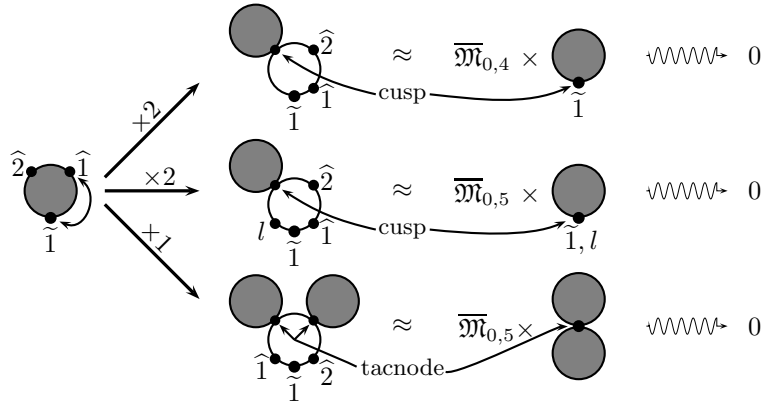

We prove Corollary~\ref{n3p3_contr_crl},
which expresses the boundary contribution
$${\mathcal{C}}_{\partial{\overU}_{{\mathcal{T}}_2}^{(1)}(\mu)}
\big(\ev_{\tilde{1}} \times \ev_{\hat{2}};\De_{\PPP\times\PPP}\big)$$
in terms of level~1 numbers, by computing contributions
from the individual strata~${\mathcal{U}}_{{\mathcal{T}}_2,{\mathcal{T}}}$.
We split the computation into four cases, depending
on whether $\chi_{\mathcal{T}}(\tilde{1},\hat{2})$ and 
$\chi_{\mathcal{T}}(\hat{1},\hat{2})$ are zero or~not.
By Lemma~\ref{n3p3_exp_lmm1}, if
$\chi_{\mathcal{T}}(\tilde{1},\hat{2}) \neq 0$ and
$\chi_{\mathcal{T}}(\hat{1},\hat{2}) \neq 0$,
the space ${\mathcal{U}}_{{\mathcal{T}}_2,{\mathcal{T}}}$ makes
no contribution.
The case $\chi_{\mathcal{T}}(\tilde{1},\hat{2}) = 0$ and
$\chi_{\mathcal{T}}(\hat{1},\hat{2}) \neq 0$ is handled in
Lemma~\ref{n3p3_contr_lmm1}.
Figure~\ref{n3p3_contr_fig1} shows the three possibilities 
for non-hollow spaces~${\mathcal{U}}_{{\mathcal{T}}_2,{\mathcal{T}}}$.
In all three cases, we express the contribution from the stratum
in terms of the number $N(\al_1)$
of zeros of an affine map between vector bundles.
However, in two of the cases, this number is zero,
basically for dimensional reasons;
the remaining number is computed in Lemma~\ref{n3p3_contr_lmm1b}.
The  case $\chi_{\mathcal{T}}(\tilde{1},\hat{2}) \neq 0$ and
$\chi_{\mathcal{T}}(\hat{1},\hat{2}) = 0$ is symmetric
to the one just considered and no separate computation is~needed.
The remaining case is dealt with in Lemma~\ref{n3p3_contr_lmm3}.
Figure~\ref{n3p3_contr_fig2} shows the three possibilities 
for non-hollow spaces~${\mathcal{U}}_{{\mathcal{T}}_2,{\mathcal{T}}}$,
but in all three cases the  corresponding number $N(\al_1)$
is zero for dimensional reasons.
In both figures, the numbers above the arrows show the multiplicity
with which the corresponding number~$N(\al_1)$
enters into 
${\mathcal{C}}_{\partial{\overU}_{{\mathcal{T}}_2}^{(1)}(\mu)}
\big(\ev_{\tilde{1}} \times \ev_{\hat{2}};\De_{\PPP\times\PPP}\big)$.

Before proceeding with the actual proofs, we formally define more spaces
of tuples of stable genus-zero maps that appear in the statements
of Theorems~\ref{n3p3_thm} and~\ref{n3tac_thm}
and describe curves pictured in Figure~\ref{images_fig}.
First, let
$${\mathcal{S}}_1(\mu)= 
\big\{b \in {\mathcal{U}}_{{\mathcal{T}}_0}(\mu) :
{\mathcal{D}}_{{\mathcal{T}}_0,\tilde{1}}^{(1)}b = 0\big\}$$
and let ${\overS}_1(\mu)$ be the closure of 
${\mathcal{S}}_1(\mu)$ in ${\overU}_{{\mathcal{T}}_0}$.
If ${\mathcal{T}} = (M_0,I_2;j,\under{d})$, let
\begin{equation*}\begin{split}
{\mathcal{S}}_{\mathcal{T}}(\mu)
&= \big\{b \in {\mathcal{U}}_{\mathcal{T}}(\mu) :\\
&\qquad\quad {\mathcal{D}}_{{\mathcal{T}},\tilde{1}}^{(1)}\ups_{\tilde{1}} + 
{\mathcal{D}}_{{\mathcal{T}},\tilde{2}}^{(1)}\ups_{\tilde{2}} = 0
~\hbox{for some}~
(b;\ups_{\tilde{1}},\ups_{\tilde{2}}) \in 
L_{\tilde{1}}{\mathcal{T}} \oplus L_{\tilde{2}}{\mathcal{T}} - {\mathcal{U}}_{\mathcal{T}}\big\}.
\end{split}\end{equation*}
We denote by ${\mathcal{S}}_2(\mu)$ the quotient of 
the disjoint union of the spaces ${\mathcal{S}}_{\mathcal{T}}(\mu)$,
taken over all bubble types ${\mathcal{T}}$ as above such that  
$d_{\tilde{1}},d_{\tilde{2}} > 0$ and 
$d_{\tilde{1}}+d_{\tilde{2}} = d$,
by the natural action of the symmetric group~$S_2$.
Finally, if ${\mathcal{T}} = (M_1,I_2;j,\under{d})$ is a bubble type
such that $j_{\hat{1}} = \tilde{1}$, let
$${\mathcal{U}}_{\mathcal{T}}^{(1)}(\mu)=
\big\{b \in {\mathcal{U}}_{\mathcal{T}}(\mu) :
\ev_{\ti{1}}(b) = \ev_{\hat{1}}(b)\big\}.$$
We denote by ${\mathcal{V}}_2^{(1)}(\mu)$
the disjoint union of the spaces  ${\mathcal{U}}_{\mathcal{T}}^{(1)}(\mu)$,
taken over all bubble types ${\mathcal{T}}$ as above such that  
$d_{\ti{1}},d_{\ti{2}} > 0$ and 
$d_{\ti{1}}+d_{\ti{2}} = d$.

\begin{lmm}
\label{n3p3_contr_lmm1}
Suppose ${\mathcal{T}} = (M_2,I;j,\under{d})$ is a bubble type such that
${\mathcal{T}} < {\mathcal{T}}_2$, $\chi_{\mathcal{T}}(\tilde{1},\hat{2}) = 0$, and
\hbox{$\chi_{\mathcal{T}}(\hat{1},\hat{2}) > 0$}.

{\rm(1)}\qua If $|\hat{I}^+| > |\chi_{\tilde{1}}({\mathcal{T}})|$,
${\mathcal{U}}_{\mathcal{T}}^{(1)}(\mu)$
is $(\ev_{\tilde{1}} \times \ev_{\hat{2}},\De_{\PPP\times\PPP})$--hollow,
and~thus 
$${\mathcal{C}}_{{\mathcal{U}}_{{\mathcal{T}}_2,{\mathcal{T}}}}
\big(\ev_{\tilde{1}} \times \ev_{\hat{2}};\De_{\PPP\times\PPP}\big)
 = 0.$$
{\rm(2)}\qua If $|\hat{I}^+| = |\chi_{\tilde{1}}({\mathcal{T}})| = 1$,
$|M_{\tilde{1}}{\mathcal{T}}| \in \{1,2\}$.
If $|M_{\tilde{1}}{\mathcal{T}}| = 2$, 
$${\mathcal{C}}_{{\mathcal{U}}_{{\mathcal{T}}_2,{\mathcal{T}}}}
\big(\ev_{\tilde{1}} \times \ev_{\hat{2}};\De_{\PPP\times\PPP}\big)=0.$$
If $|M_{\tilde{1}}{\mathcal{T}}| = 1$, 
\begin{equation*}\begin{split}
{\mathcal{C}}_{{\mathcal{U}}_{{\mathcal{T}}_2,{\mathcal{T}}}}
\big(\ev_{\tilde{1}} \times \ev_{\hat{2}};\De_{\PPP\times\PPP}\big)
=\big\lan 6a_{\hat{0}}^2 + 4a_{\hat{0}}c_1({\mathcal{L}}_{\tilde{1}}^*)
 + c_1^2({\mathcal{L}}_{\tilde{1}}^*),{\overV}_1^{(1)}(\mu)\big\ran
+\big|{\mathcal{S}}_2(\mu)\big| \qquad&\\
-\big|{\mathcal{V}}_2^{(1)}(\mu)\big|
-\big\lan 8a_{\hat{0}}^2 + 4c_1({\mathcal{L}}_{\tilde{1}}^*),
{\overS}_1(\mu)\big\ran.&
\end{split}\end{equation*}
{\rm(3)}\qua If $|\hat{I}^+| = |\chi_{\tilde{1}}({\mathcal{T}})| = 2$,
${\mathcal{C}}_{{\mathcal{U}}_{{\mathcal{T}}_2,{\mathcal{T}}}}
\big(\ev_{\tilde{1}} \times \ev_{\hat{2}};\De_{\PPP\times\PPP}\big)=0$.
\end{lmm}

\begin{proof}
(1)\qua  By Lemma~\ref{n3p3_str_lmm1},
${\overU}_{{\mathcal{T}}_2}^{(1)}(\mu)\cap{\mathcal{U}}_{{\mathcal{T}}_2,{\mathcal{T}}}
 \subset {\mathcal{U}}_{{\mathcal{T}}|{\mathcal{T}}_2}^{(1)}(\mu)$.
With appropriate identifications, 
${\mathcal{U}}_{{\mathcal{T}}|{\mathcal{T}}_2}^{(1)}(\mu)$
is the zero set of the section
$\ev_{{\mathcal{T}}_2,M_0} \oplus (\ev_{\hat{1}} - \ev_{\tilde{1}})$
of the bundle 
$$\ev_{{\mathcal{T}}_2,M_0}^*{\mathcal{N}}\De_{{\mathcal{T}}_2}(\mu)\oplus
\ev_{\tilde{1}}^*T\PPP$$
over an open neighborhood of 
${\mathcal{U}}_{{\mathcal{T}}|{\mathcal{T}}_2}^{(1)}(\mu)$
in~${\mathcal{U}}_{{\mathcal{T}}_2,{\mathcal{T}}}$.
By Lemma~\ref{str_lmm}, this section is transversal to the zero set,
since the constraints~$\mu$ are assumed to be in general position.
By Proposition~\ref{str_prp}, there exists 
a $C^1$--negligible map 
$$\ve_-\co \mathcal{FT}_{\de} - Y(\mathcal{FT};\hat{I}^+)
\lra\ev_{{\mathcal{T}}_2,M_0}^*{\mathcal{N}}\De_{{\mathcal{T}}_2}(\mu) \oplus 
\ev_{\tilde{1}}^*T\PPP$$ 
such~that
$$\big\{\ev_{{\mathcal{T}}_2,M_0} \times 
\ev_{\tilde{1}} \times \ev_{\hat{1}}\big\}
\big(\phi_{{\mathcal{T}}_2,{\mathcal{T}}}(b;\ups)\big)=
\big\{\ev_{{\mathcal{T}}_2,M_0} \times 
\ev_{\tilde{1}} \times \ev_{\hat{1}}\big\}(b)+\ve_-(b;\ups)$$
for all $(b;\ups) \in \mathcal{FT}_{\de} - Y(\mathcal{FT};\hat{I}^+)$.
On the other hand, by Lemma~\ref{n3p3_exp_lmm2}, 
$$\big\{\ev_{\tilde{1}} \times \ev_{\hat{2}}\big\}
\big(\phi_{{\mathcal{T}}_2,{\mathcal{T}}}(\ups)\big)
=\big\{\al + \ve(\ups)\big\}\rho(\ups)\in\ev_{\tilde{1}}^*T\PPP$$
for all $\ups \in \mathcal{FT}_{\de^*} - Y(\mathcal{FT};\hat{I}^+)$,
where $\al$ is a linear map, which is injective over 
${\mathcal{U}}_{{\mathcal{T}}|{\mathcal{T}}_2}^{(1)}(\mu)$,
and its domain is a vector bundle of rank~$|\chi_{\tilde{1}}({\mathcal{T}})|$.
Thus, if $|\hat{I}^+| > |\chi_{\tilde{1}}({\mathcal{T}})|$, then
${\mathcal{U}}_{{\mathcal{T}}|{\mathcal{T}}_2}^{(1)}(\mu)$ is
$(\ev_{\tilde{1}} \times \ev_{\hat{2}},\De_{\PPP\times\PPP})$--hollow, 
and
$${\mathcal{C}}_{{\mathcal{U}}_{{\mathcal{T}}_2,{\mathcal{T}}}}
\big(\ev_{\tilde{1}} \times \ev_{\hat{2}};\De_{\PPP\times\PPP}\big)
=0$$
by Proposition~\ref{euler_prp}, or Lemma~\ref{euler_lmm2},
and Lemma~\ref{n3p3_str_lmm1}.

(2)\qua On the other hand, if $|\hat{I}^+| =
|\chi_{\tilde{1}}({\mathcal{T}})|$,
by the above and Lemma~\ref{n3p3_exp_lmm2}, 
${\mathcal{U}}_{{\mathcal{T}}|{\mathcal{T}}_2}^{(1)}(\mu)$ 
is $(\ev_{\tilde{1}} \times \ev_{\hat{2}},\De_{\PPP\times\PPP})$--regular,
and by Proposition~\ref{euler_prp} and rescaling of the linear~map,
\begin{gather*}
{\mathcal{C}}_{{\mathcal{U}}_{{\mathcal{T}}_2,{\mathcal{T}}}}
\big(\ev_{\tilde{1}} \times \ev_{\hat{2}};\De_{\PPP\times\PPP}\big)
=N(\al),\qquad\hbox{where}\\
\al \in \Ga\big({\mathcal{U}}_{{\mathcal{T}}|{\mathcal{T}}_2}^{(1)}(\mu);
\hbox{Hom}(\mathcal{FT},\ev_{\tilde{1}}^*T\PPP)\big),\quad
\al(\ups)=\sum_{h\in\chi_{\tilde{1}}({\mathcal{T}})}  
\big(y_{\hat{2}} - x_h\big)^{-1} \otimes
{\mathcal{D}}_{{\mathcal{T}},h}^{(1)}\ups_h.
\end{gather*}
provided $\al$ is a regular section, as is implied by what follows.
Since the map
$$\mathcal{FT}\equiv\bigoplus_{h\in\chi_{\tilde{1}}({\mathcal{T}})}  
L_{\tilde{1}}^*{\mathcal{T}} \otimes L_h{\mathcal{T}}
\lra
{\mathcal{F}}=\bigoplus_{h\in\chi_{\tilde{1}}({\mathcal{T}})}   L_h{\mathcal{T}},\qquad
\ups_h\lra \big(y_{\hat{2}} - x_h\big)^{-1} \otimes\ups_h,$$
is simply a rescaling of factors over~${\mathcal{U}}_{{\mathcal{T}}_2}(\mu)$,
\begin{gather*}
{\mathcal{C}}_{{\mathcal{U}}_{{\mathcal{T}}|{\mathcal{T}}_2}}
\big(\ev_{\tilde{1}} \times \ev_{\hat{2}};\De_{\PPP\times\PPP}\big)
=N(\al'),\qquad\hbox{where}\\
\al' \in \Ga\big({\overU}_{{\mathcal{T}}|{\mathcal{T}}_2}^{(1)}(\mu);
\hbox{Hom}({\mathcal{F}},\ev_{\tilde{1}}^*T\PPP)\big),\quad
\al'(\ups)=\sum_{h\in\chi_{\tilde{1}}({\mathcal{T}})}   
{\mathcal{D}}_{{\mathcal{T}},h}^{(1)}\ups_h.
\end{gather*}
Note that with respect to the decomposition~\eqref{cart_split},
shown in Figure~\ref{n3p3_contr_fig1},
the linear map $\al'$ comes entirely from the second factor.
Thus, if the first factor is positive-dimensional,
\hbox{$N(\al') = 0$, ie}
$${\mathcal{C}}_{{\mathcal{U}}_{{\mathcal{T}}_2,{\mathcal{T}}}}
\big(\ev_{\tilde{1}} \times \ev_{\hat{2}};\De_{\PPP\times\PPP}\big)=0$$
unless $|\chi_{\tilde{1}}({\mathcal{T}})| = 1$ and
$|M_{\tilde{1}}{\mathcal{T}}| = 1$.
If $|\chi_{\tilde{1}}({\mathcal{T}})| = 1$ and
$|M_{\tilde{1}}{\mathcal{T}}| = 1$, we conclude that
\begin{gather*}
{\mathcal{C}}_{{\mathcal{U}}_{{\mathcal{T}}_2,{\mathcal{T}}}}
\big(\ev_{\tilde{1}} \times \ev_{\hat{2}};\De_{\PPP\times\PPP}\big)
=N(\al_1),\qquad\hbox{where}\\
\al_1 \in \Ga\big({\overV}_1^{(1)}(\mu);
\hbox{Hom}(L_{\tilde{1}}{\mathcal{T}}_1,\ev_{\hat{0}}^*T\PPP)\big),\quad
\al_1(\ups)={\mathcal{D}}_{{\mathcal{T}}_1,\tilde{1}}^{(1)}\ups.
\end{gather*}
The number $N(\al_1)$ is computed below.
\end{proof}

\begin{lmm}
\label{n3p3_contr_lmm1b}
If $\al_1 \in \Ga\big({\overV}_1^{(1)}(\mu);
\hbox{Hom}(L_{\tilde{1}},\ev_{\hat{0}}^*T\PPP)\big)$ is given by
$\al_1 = {\mathcal{D}}_{{\mathcal{T}}_1,\tilde{1}}^{(1)}$,
\begin{equation*}\begin{split}
N(\al_1) =
\blr{6a_{\hat{0}}^2 + 4a_{\hat{0}}c_1({\mathcal{L}}_{\tilde{1}}^*)
 + c_1^2({\mathcal{L}}_{\tilde{1}}^*),{\overV}_1^{(1)}(\mu)}
+\big|{\mathcal{S}}_2(\mu)\big|-\big|{\mathcal{V}}_2^{(1)}(\mu)\big|\qquad\quad&\\
-\blr{8a_{\hat{0}}^2 + 4c_1({\mathcal{L}}_{\tilde{1}}^*),
{\overS}_1(\mu)}.&
\end{split}\end{equation*}
\end{lmm}

\begin{proof}
(1)\qua Since $\al_1$ does not vanish on ${\mathcal{V}}_1^{(1)}(\mu)$ 
by Lemma~\ref{str_lmm}, by Propositions~\ref{zeros_prp} and~\ref{euler_prp},
\begin{equation}\label{n3p3_contr_lmm1b_e1}
N(\al_1)=\big\lan 6a_{\hat{0}}^2 + 4a_{\hat{0}}c_1(L_{\tilde{1}}^*)
+c_1^2(L_{\tilde{1}}^*),{\overV}_1^{(1)}(\mu)\big\ran
-{\mathcal{C}}_{\partial{\overV}_1^{(1)}(\mu)}(\al_1^{\perp}),
\end{equation}
where $\al_1^{\perp}$ denotes the composition of $\al_1$
with the projection $\pi_{\bar{\nu}_1}^{\perp}$
onto the quotient ${\mathcal{O}}_1$ of $\ev_{\hat{0}}^*T\PPP$
by a generic trivial line subbundle~$\Bbb{C}\bar{\nu}_1$.
Figure~\ref{n3p3_contr_fig1b} shows the five types of boundary
strata that are not $\al_1^{\perp}$--hollow.
Contributions from the first two are computed in (2) below,
from the following two in~(3), and from the last one in (4)~below.

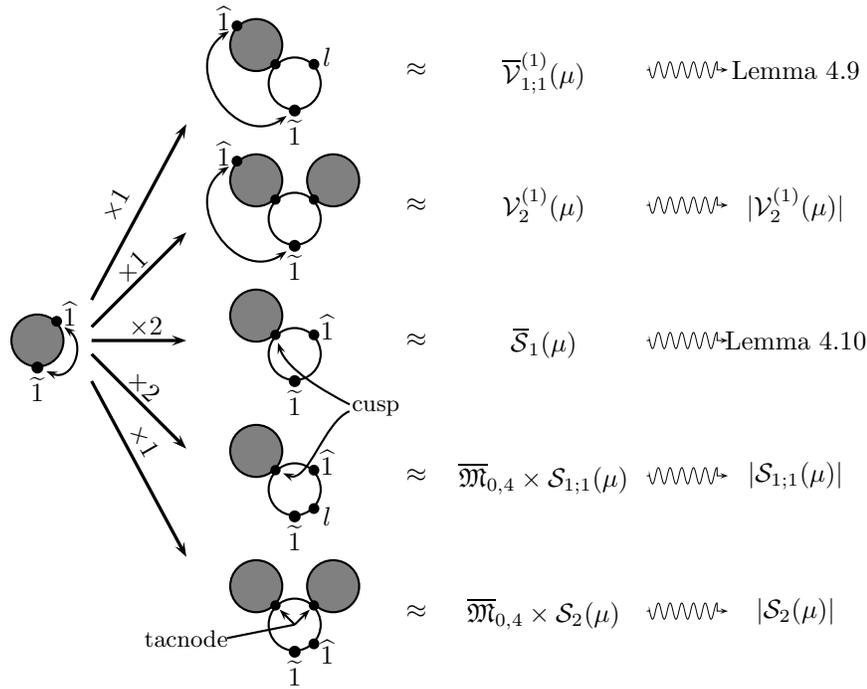
\begin{figure}[ht!]\small
\begin{pspicture}(-.5,-5.7)(10,4)
\psset{unit=.36cm}
\pscircle[fillstyle=solid,fillcolor=gray](1,-2){1}
\pscircle*(1,-3){.22}\rput(1,-3.8){$\tilde{1}$}
\pscircle*(1.71,-1.29){.2}\rput(2.2,-1){$\hat{1}$}
\pnode(1,-3){A0}\pnode(1.71,-1.29){B0}
\ncarc[nodesep=.35,arcangleA=-90,arcangleB=-75,ncurv=1.2]{<->}{A0}{B0}
% 2nd column starts here
\psline[linewidth=.12]{->}(3,-.5)(6.5,6)
\rput{60}(3.9,3){$\times1$}
\psline[linewidth=.12]{->}(3,-1.5)(6.5,2)
\rput{45}(4.5,.7){$\times1$}
\psline[linewidth=.12]{->}(3,-2)(6.5,-2)
\rput(5,-1.5){$\times2$}
\psline[linewidth=.12]{->}(3,-2.5)(6.5,-6)
\rput{-45}(4.9,-3.7){$\times2$}
\psline[linewidth=.12]{->}(3,-3.5)(6.5,-10)
\rput{-60}(4.9,-5.7){$\times1$}
% 3rd column starts here
\pscircle(10.5,7.5){1}
\pscircle[fillstyle=solid,fillcolor=gray](9.08,8.92){1}
\pscircle*(10.5,6.5){.22}\rput(10.5,5.7){$\tilde{1}$}
\pscircle*(9.79,8.21){.19}
\pscircle*(11.21,8.21){.2}\rput(11.7,8.5){$l$}
\pscircle*(8.37,9.63){.2}\rput(7.88,9.92){$\hat{1}$}
\pnode(10.5,6.5){A1}\pnode(8.37,9.63){B1}
\ncarc[nodesep=.35,arcangleA=90,arcangleB=90,ncurv=1.2]{<->}{A1}{B1}
% 2nd row starts here
\pscircle(10.5,2.5){1}
\pscircle[fillstyle=solid,fillcolor=gray](9.08,3.92){1}
\pscircle*(10.5,1.5){.22}\rput(10.5,.7){$\tilde{1}$}
\pscircle*(9.79,3.21){.19}
\pscircle[fillstyle=solid,fillcolor=gray](11.92,3.92){1}
\pscircle*(11.21,3.21){.19}
\pscircle*(8.37,4.63){.2}\rput(7.88,4.92){$\hat{1}$}
\pnode(10.5,1.5){A2}\pnode(8.37,4.63){B2}
\ncarc[nodesep=.35,arcangleA=90,arcangleB=90,ncurv=1.2]{<->}{A2}{B2}
% 3rd row starts here
\pscircle(10.5,-2.5){1}
\pscircle[fillstyle=solid,fillcolor=gray](9.08,-1.08){1}
\pscircle*(10.5,-3.5){.22}\rput(10.5,-4.3){$\tilde{1}$}
\pscircle*(9.79,-1.79){.19}\pscircle*(11.21,-1.79){.2}
\rput(11.7,-1.5){$\hat{1}$}
% 4th row starts here
\pscircle(10.5,-7.5){1}
\pscircle[fillstyle=solid,fillcolor=gray](9.08,-6.08){1}
\pscircle*(10.5,-8.5){.22}\rput(10.5,-9.3){$\tilde{1}$}
\pscircle*(9.79,-6.79){.19}
\pscircle*(11.21,-6.79){.2}\rput(11.7,-6.5){$\hat{1}$}
\pscircle*(11.21,-8.21){.2}\rput(11.7,-8.5){$l$}
\pnode(12.8,-4.5){C}\rput(13.5,-4.5){\footnotesize cusp}
\pnode(9.79,-1.79){C1}\pnode(9.79,-6.79){C2}
\ncarc[nodesep=.3,arcangleA=15,arcangleB=30,ncurv=.6]{->}{C}{C1}
\ncarc[nodesep=.3,arcangleA=-15,arcangleB=60,ncurv=.6]{->}{C}{C2}
% 5th row starts here
\pscircle(10.5,-12.5){1}
\pscircle[fillstyle=solid,fillcolor=gray](9.08,-11.08){1}
\pscircle*(10.5,-13.5){.22}\rput(10.5,-14.3){$\tilde{1}$}
\pscircle*(9.79,-11.79){.19}
\pscircle[fillstyle=solid,fillcolor=gray](11.92,-11.08){1}
\pscircle*(11.21,-11.79){.19}
\pscircle*(11.21,-13.21){.2}\rput(11.7,-13.5){$\hat{1}$}
\pnode(8,-13){T1}\rput(6.5,-13){\footnotesize tacnode}
\pnode(10.5,-12.5){T2}\pnode(9.94,-11.94){T2a}\pnode(11.06,-11.94){T2b}\
\ncarc[nodesep=0,arcangleA=0,arcangleB=0,ncurv=.5]{-}{T1}{T2}
\ncline[nodesep=0]{->}{T2}{T2a}\ncline[nodesep=0]{->}{T2}{T2b}
% 4th column starts here
\rput(15,8){$\approx$}\rput(15,3){$\approx$}\rput(15,-2){$\approx$}
\rput(15,-7){$\approx$}\rput(15,-12){$\approx$}
% 5th column starts here
\rput(19.7,7.9){${\overV}_{1;1}^{(1)}(\mu)$}
\rput(19.7,3){${\mathcal{V}}_2^{(1)}(\mu)$}
\rput(19.7,-2.1){${\overS}_1(\mu)$}
\rput(19.7,-7.1){$\ov{\frak M}_{0,4}\times{\mathcal{S}}_{1;1}(\mu)$}
\rput(19.7,-12.1){$\ov{\frak M}_{0,4}\times{\mathcal{S}}_2(\mu)$}
%6th column starts here
\pscoil[linewidth=.03,coilarmA=0.1,coilarmB=0.4,coilaspect=0,coilheight=.72,
coilwidth=.6]{->}(23.5,8)(26.5,8)
\pscoil[linewidth=.03,coilarmA=0.1,coilarmB=0.4,coilaspect=0,coilheight=.72,
coilwidth=.6]{->}(23.5,3)(26.5,3)
\pscoil[linewidth=.03,coilarmA=0.1,coilarmB=0.4,coilaspect=0,coilheight=.72,
coilwidth=.6]{->}(23.5,-2)(26.5,-2)
\pscoil[linewidth=.03,coilarmA=0.1,coilarmB=0.4,coilaspect=0,coilheight=.72,
coilwidth=.6]{->}(23.5,-7)(26.5,-7)
\pscoil[linewidth=.03,coilarmA=0.1,coilarmB=0.4,coilaspect=0,coilheight=.72,
coilwidth=.6]{->}(23.5,-12)(26.5,-12)
%7th column starts here
\rput(29,8){\small Lemma~\ref{n3p3_contr_lmm1d}}
\rput(29,3){$|{\mathcal{V}}_2^{(1)}(\mu)|$}
\rput(29,-2){\small Lemma~\ref{n3p3_contr_lmm1c}}
\rput(29,-7){$|{\mathcal{S}}_{1;1}(\mu)|$}
\rput(29,-12){$|{\mathcal{S}}_2(\mu)|$}
\end{pspicture}
\caption{An outline of the proof of Lemma~\ref{n3p3_contr_lmm1b}}
\label{n3p3_contr_fig1b}
\end{figure}

(2)\qua  If
${\mathcal{T}} < {\mathcal{T}}_1$ and $\chi_{\mathcal{T}}(\tilde{1},\hat{1}) > 0$,
${\overV}_1^{(1)}(\mu)\cap{\mathcal{U}}_{{\mathcal{T}}_1,{\mathcal{T}}}
 \subset {\mathcal{U}}_{{\mathcal{T}}|{\mathcal{T}}_1}^{(1)}(\mu)$
by Lemma~\ref{n3p3_str_lmm1}.
By Proposition~\ref{str_prp},
\begin{gather*}
{\mathcal{D}}_{{\mathcal{T}}_1,\tilde{1}}^{(1)}
\big(\phi_{{\mathcal{T}}_1,{\mathcal{T}}}(\ups)\big)
=  \sum_{h\in\chi_{\tilde{1}}({\mathcal{T}})}     
\big({\mathcal{D}}_{{\mathcal{T}},h}^{(1)} + \ve_h(\ups)\big)\rho_h(\ups),
\qquad\hbox{where}\\
\rho_h(\ups)=  \prod_{i\in(\tilde{1},h]}    \ups_i,
\quad\mbox{for all}\ 
~\ups \in \mathcal{FT}_{\de^*} - Y(\mathcal{FT};\hat{I}^+),
\end{gather*}
and for some $C^0$--negligible maps 
$$\ve_h\co \mathcal{FT}_{\de^*} - Y(\mathcal{FT};\hat{I}^+) \lra
\hbox{Hom}(L_h{\mathcal{T}},\ev_{\hat{0}}^*T\PPP).$$
Let $\tilde{\mathcal{F}}_h{\mathcal{T}}$ denote the line bundle determined by~$\rho_h$.
By Lemma~\ref{str_lmm}, the~map
$$\tilde{\mathcal{F}}{\mathcal{T}}\equiv\bigoplus_{h\in\chi_{\tilde{1}}({\mathcal{T}})}
    \tilde{\mathcal{F}}_h{\mathcal{T}}\lra 
L_{\tilde{1}}^*{\mathcal{T}} \otimes \ev_{\hat{0}}^*T\PPP,\qquad
\big(\tilde{\ups}_h\big)_{h\in\chi_{\tilde{1}}({\mathcal{T}})}\lra
\sum_{h\in\chi_{\tilde{1}}({\mathcal{T}})}   
{\mathcal{D}}_{{\mathcal{T}},h}^{(1)}\tilde{\ups}_h,$$
is injective over ${\mathcal{U}}_{{\mathcal{T}}|{\mathcal{T}}_1}^{(1)}(\mu)$.
If $\bar{\nu}_1$ is generic, the same is true of the~map
$$\al_2'\co \tilde{\mathcal{F}}{\mathcal{T}}\lra
L_{\tilde{1}}^*{\mathcal{T}} \otimes {\mathcal{O}}_1,\qquad
\big\{\al_2'(\tilde{\ups})\big\}(\ups)=
\pi_{\bar{\nu}_1}^{\perp}\sum_{h\in\chi_{\tilde{1}}({\mathcal{T}})}   
{\mathcal{D}}_{{\mathcal{T}},h}^{(1)}\tilde{\ups}_h \otimes \ups.$$
By the same argument as in (1) of the proof of Lemma~\ref{n3p3_contr_lmm1},
$${\mathcal{C}}_{{\mathcal{U}}_{{\mathcal{T}}_1,{\mathcal{T}}}}(\al_1^{\perp})=0
\qquad\hbox{unless}\quad \hat{I}^+ = \chi_{\tilde{1}}({\mathcal{T}}).$$
If $\hat{I}^+ = \chi_{\tilde{1}}({\mathcal{T}})$, 
by dimension-counting either
\begin{gather*}
|\hat{I}^+| = |\chi_{\tilde{1}}({\mathcal{T}})|=1
\quad\hbox{and}\quad |M_{\tilde{1}}{\mathcal{T}}|=1 \qquad\hbox{or}\\ 
|\hat{I}^+| = |\chi_{\tilde{1}}({\mathcal{T}})|=2
\quad\hbox{and}\quad |M_{\tilde{1}}{\mathcal{T}}|=0.
\end{gather*}
In the second case, the map $\al_2'$ is an isomorphism on every fiber
of $\mathcal{FT}$ over the finite~set 
${\mathcal{U}}_{{\mathcal{T}}|{\mathcal{T}}_1}^{(1)}(\mu)$.
Thus, by Proposition~\ref{euler_prp},
$${\mathcal{C}}_{{\mathcal{U}}_{{\mathcal{T}}_1,{\mathcal{T}}}}(\al_1^{\perp})
=\big|{\mathcal{U}}_{{\mathcal{T}}|{\mathcal{T}}_1}^{(1)}(\mu)\big|.$$
Note that the sum of the numbers
$|{\mathcal{U}}_{{\mathcal{T}}|{\mathcal{T}}_1}^{(1)}(\mu)|$,
taken over all bubble types ${\mathcal{T}} < {\mathcal{T}}_1$
such that
$$|\hat{I}^+| = |\chi_{\tilde{1}}({\mathcal{T}})|=2 \qquad\hbox{and}\qquad
|M_{\tilde{1}}{\mathcal{T}}| = 0,$$
is~$|{\mathcal{V}}_2^{(1)}(\mu)|$.
On the other hand, if 
$\hat{I}^+ = \chi_{\tilde{1}}({\mathcal{T}}) = \{h\}$ is a single-element set
and $|M_{\tilde{1}}{\mathcal{T}}| = 1$, ie
${\mathcal{T}} = {\mathcal{T}}_1(l)$ for some $l \in [N]$,
by Proposition~\ref{euler_prp} and the decomposition~\eqref{cart_split},
$${\mathcal{C}}_{{\mathcal{U}}_{{\mathcal{T}}|{\mathcal{T}}_1}^{(1)}(\mu)}(\al_1^{\perp})
=N(\al_2),\qquad\hbox{where}\qquad
\al_2\in\Ga\big({\overU}_{{\mathcal{T}}_1/l}^{(1)}(\mu);
\hbox{Hom}(L_h,{\mathcal{O}}_1)\big)$$
is the map induced by $\al_2'$, ie
the composition of ${\mathcal{D}}_{{\mathcal{T}}_1/l,h}^{(1)}$
with the projection $\pi_{\bar{\nu}_1}^{\perp}$
onto the quotient ${\mathcal{O}}_1$ of $\ev_{\hat{0}}^*T\PPP$
by a generic line subbundle~$\Bbb{C}\bar{\nu}_1$.
Thus, by Lemma~\ref{n3p3_contr_lmm1d},
$$\sum_{l\in[N]}
{\mathcal{C}}_{{\mathcal{U}}_{{\mathcal{T}}_1(l)|{\mathcal{T}}_1}}(\al_1^{\perp})
=\big\lan 4a_{\hat{0}} + c_1(L_{\tilde{1}}^*),
{\overV}_{1;1}^{(1)}(\mu)\big\ran
-2\big|{\mathcal{S}}_{1;1}(\mu)\big|.$$
Summing up the above contributions, we find that
\begin{equation}\label{n3p3_contr_lmm1b_e11}
\sum_{\chi_{\mathcal{T}}(\tilde{1},\hat{1})>0}     
{\mathcal{C}}_{{\mathcal{U}}_{{\mathcal{T}}_1,{\mathcal{T}}}}(\al_1^{\perp})
=\big|{\mathcal{V}}_2^{(1)}(\mu)\big|+
\big\lan 4a_{\hat{0}} + c_1(L_{\tilde{1}}^*),
{\overV}_{1;1}^{(1)}(\mu)\big\ran
-2\big|{\mathcal{S}}_{1;1}(\mu)\big|.
\end{equation}
In~\eqref{n3p3_contr_lmm1b_e11}, ${\mathcal{S}}_{1;1}(\mu)$
denotes the disjoint union of the sets 
\begin{gather*}
{\mathcal{S}}_{{\mathcal{T}}_{0;l}}(\mu)=\{b \in {\mathcal{U}}_{{\mathcal{T}}_{0;l}}(\mu) :
{\mathcal{D}}_{{\mathcal{T}}_{0;l},\tilde{1}}^{(1)}b = 0\big\},
\qquad\hbox{where}\\
{\mathcal{T}}_{0;l} = \big(M_0 - \{l\},I_1(l);j,d),~~
I_1(l)=\{\hat{0} = l\}\sqcup I_1^+,
\end{gather*}
taken over all $l \in [N]$.
The space ${\mathcal{V}}_{1;1}^{(1)}(\mu)$ is the disjoint union
of the sets
\begin{gather*}
{\mathcal{U}}_{{\mathcal{T}}_{1;l}}^{(1)}(\mu)=
\{b \in {\mathcal{U}}_{{\mathcal{T}}_{1;l}}(\mu) :
\ev_{\tilde{1}}(b) = \ev_{\hat{1}}(b)\big\},
\qquad\hbox{where}\\
{\mathcal{T}}_{1;l} = \big(M_1 - \{l\},I_1(l);j,d),
\end{gather*}
taken over all $l \in [N]$.
As usual, ${\overV}_{1;1}^{(1)}(\mu)$ denotes 
the closure of ${\mathcal{V}}_{1;1}^{(1)}(\mu)$ inside of 
a union of moduli spaces of stable~maps.
More geometrically, the image of every element of 
${\mathcal{S}}_{1;1}(\mu)$ (of ${\mathcal{V}}_{1;1}^{(1)}(\mu)$)
has a cusp (a node) at one of the constraints $\mu_1,\ldots,\mu_N$.

(3)\qua If ${\mathcal{T}} < {\mathcal{T}}_1$ and $\chi_{\mathcal{T}}(\tilde{1},\hat{1}) = 0$,
${\overV}_1^{(1)}(\mu)\cap{\mathcal{U}}_{{\mathcal{T}}_1,{\mathcal{T}}}
 \subset {\mathcal{S}}_{{\mathcal{T}}|{\mathcal{T}}_1}(\mu)$
by Lemma~\ref{n3p3_str_lmm2}.
On the other hand, by Lemma~\ref{str_lmm},
${\mathcal{S}}_{{\mathcal{T}}|{\mathcal{T}}_1}(\mu) = \eset$
unless $|\chi_{\tilde{1}}({\mathcal{T}})| \in \{1,2\}$.
Suppose $\chi_{\tilde{1}}({\mathcal{T}}) = \{h\}$ is a single-element set.
With appropriate identifications,
${\mathcal{S}}_{{\mathcal{T}}|{\mathcal{T}}_1}(\mu)$
is the zero set of the section
$\ev_{{\mathcal{T}}_1,M_0} \oplus {\mathcal{D}}_{{\mathcal{T}},h}^{(1)}$
of the bundle 
$$\ev_{{\mathcal{T}}_1,M_0}^*{\mathcal{N}}\De_{{\mathcal{T}}_1}(\mu)
\oplus L_h^* \otimes \ev_{\hat{0}}^*T\PPP$$
defined over a neighborhood of ${\mathcal{S}}_{{\mathcal{T}}|{\mathcal{T}}_1}(\mu)$
in~${\mathcal{U}}_{{\mathcal{T}}_1,{\mathcal{T}}}$.
By Lemma~\ref{str_lmm}, this section is transverse to the zero set.
By Proposition~\ref{str_prp} and Lemma~\ref{n3p3_str_lmm2}, 
\begin{gather*}
\ev_{{\mathcal{T}}_1,M_0}\big(\phi_{{\mathcal{T}}_1,{\mathcal{T}}}(b;\ups)\big)
=\ev_{{\mathcal{T}}_1,M_0}(b)+\ve_{-;1}(b;\ups),\\
\{\ev_{\tilde{1}} \times \ev_{\hat{1}}\}
\big(\phi_{{\mathcal{T}}_1,{\mathcal{T}}}(\ups)\big)=
\big(y_{h;\hat{1}} - x_{h;\hat{1}}\big)^{-1}\otimes
\big\{{\mathcal{D}}_{{\mathcal{T}},h}^{(1)} + \ve_{-;2}(\ups)\big\}
\rho_{{\mathcal{T}},\hat{1};h}^{(0;1)}(\ups)
\end{gather*}
for all $(b;\ups) \in \mathcal{FT}_{\de} - Y(\mathcal{FT};\hat{I}^+)$
and some $C^1$--negligible maps
$$\ve_{-;1},\ve_{-;2}\co 
\mathcal{FT}_{\de} - Y(\mathcal{FT};\hat{I}^+)
\lra\ev_{{\mathcal{T}}_1,M_0}^*{\mathcal{N}}\De_{{\mathcal{T}}_1}(\mu),
L_h^* \otimes \ev_{\hat{0}}^*T\PPP.$$
On the other hand, subtracting 
$\big(y_{\tilde{1};\hat{1}}(\ups) - x_{\tilde{1};h}(\ups)\big)$
times the expansion of 
$\{\ev_{\tilde{1}} \times \ev_{\hat{1}}\} \circ 
\phi_{{\mathcal{T}}_1,{\mathcal{T}}}$ in Lemma~\ref{n3p3_str_lmm2} from
the expansion of 
${\mathcal{D}}_{{\mathcal{T}}_1,\tilde{1}}^{(1)} \circ \phi_{{\mathcal{T}}_1,{\mathcal{T}}}$
in (3b) of Proposition~\ref{str_prp}, we obtain
\begin{gather*}
{\mathcal{D}}_{{\mathcal{T}}_1,\tilde{1}}^{(1)}\phi_{{\mathcal{T}}_1,{\mathcal{T}}}(\ups)=
-\big(y_{h;\hat{1}} - x_{\hat{1};h}\big) \otimes 
\big\{{\mathcal{D}}_{{\mathcal{T}},h}^{(2)} + \ve(\ups)\big\}\rho(\ups),\\
\hbox{where}\qquad
\rho(\ups) =  
\prod_{i\in(\chi_{\mathcal{T}}(\hat{1},h),h]}      \ups_i^{\otimes2}
\otimes   
\prod_{i\in(\tilde{1},\chi_{\mathcal{T}}(\hat{1},h)]}      \ups_i,
\end{gather*}
for all $\ups \in \mathcal{FT}_{\de}$ such that
$\phi_{{\mathcal{T}}_1,{\mathcal{T}}}(\ups) \in {\mathcal{U}}_{{\mathcal{T}}_1}^{(1)}(\mu)$
and for some  $C^0$--negligible map
$$\ve\co \mathcal{FT}_{\de} - Y(\mathcal{FT};\hat{I}^+)
\lra L_h^{*\otimes2} \otimes\ev_{\hat{0}}^*T\PPP.$$
By Lemma~\ref{str_lmm}, ${\mathcal{D}}_{{\mathcal{T}},h}^{(2)}$ does not vanish
over ${\mathcal{S}}_{{\mathcal{T}}|{\mathcal{T}}_1}(\mu)$,
and neither does the linear~map
$$\al_2'\co {\mathcal{F}}_h{\mathcal{T}}^{\otimes2} \lra 
L_{\tilde{1}}^*{\mathcal{T}} \otimes {\mathcal{O}}_1,\qquad
\big\{\al_2'(\tilde{\ups})\big\}(\ups)
=\pi_{\bar{\nu}_1}^{\perp}{\mathcal{D}}_{{\mathcal{T}},h}^{(2)}
\tilde{\ups} \otimes \ups,$$
provided $\bar{\nu}_1$ is generic.
Thus, ${\mathcal{S}}_{{\mathcal{T}}|{\mathcal{T}}_1}(\mu)$ is
$\al_1^{\perp}$--hollow unless
$\hat{I}^+ = \chi_{\tilde{1}}({\mathcal{T}})$.
On the other hand, if $\hat{I}^+ = \chi_{\tilde{1}}({\mathcal{T}})$,
by Proposition~\ref{euler_prp}, a rescaling of the linear map, 
and the decomposition~\eqref{cart_split},
\begin{gather*}
{\mathcal{C}}_{{\mathcal{U}}_{{\mathcal{T}}_1,{\mathcal{T}}}(\mu)}(\al_1^{\perp})
=2N(\al_2),\qquad\hbox{where}\\
\al_2 \in \Ga(\ov{\frak M}_{\{\tilde{1},h\}\sqcup M_{\tilde{1}}{\mathcal{T}}}
 \times {\mathcal{S}}_{{\overT}}(\mu);
\hbox{Hom}(\pi_1^*L_{\tilde{1}}^* \otimes \pi_2^*L_h^{\otimes2},
\pi_1^*L_{\hat{1}}^* \otimes {\mathcal{O}}_1)\big),\\
\al_2(\tilde{\ups})=\pi_{\bar{\nu}_1}^{\perp}\circ\big\{
{\mathcal{D}}_{{\overT},h}^{(2)}\tilde{\ups}\big\}.
\end{gather*} 
By dimension-counting, $|M_{\tilde{1}}{\mathcal{T}}| \in \{1,2\}$.
If $|M_{\tilde{1}}{\mathcal{T}}| = 2$, 
ie $M_{\tilde{1}}{\mathcal{T}} = \{\hat{2},l\}$ for some $l \in [N]$,
${\mathcal{S}}_{{\overT}}(\mu)$ is a finite set
and ${\mathcal{D}}_{{\overT},h}^{(2)}$ does not vanish.
Thus, by Propositions~\ref{zeros_prp} and~\ref{euler_prp},
$$\sum_{\chi_{\mathcal{T}}(\tilde{1},\hat{1})=0,|M_{\tilde{1}}{\mathcal{T}}|=2}
{\mathcal{C}}_{{\mathcal{U}}_{{\mathcal{T}}|{\mathcal{T}}_1}}(\al_1^{\perp})
=2\big\lan c_1(L_{\tilde{1}}^*),
\ov{\frak M}_{\{\tilde{1},h\}\sqcup M_{\tilde{1}}{\mathcal{T}}}
\big\ran \big|{\mathcal{S}}_{1;1}(\mu)\big|=2\big|{\mathcal{S}}_{1;1}(\mu)\big|.$$
If $|M_{\tilde{1}}{\mathcal{T}}| = 1$, that is, $M_{\tilde{1}}{\mathcal{T}} = \{\hat{2}\}$,
by Propositions~\ref{zeros_prp} and~\ref{euler_prp} and 
Lemma~\ref{n3p3_contr_lmm1c},
$$N(\al_2)=\big\lan 4a_{\hat{0}} + 2c_1({\mathcal{L}}_{\tilde{1}}^*),
{\overS}_1(\mu)\big\ran-\big|{\mathcal{S}}_2(\mu)\big|.$$
Thus, summing up the above contributions, we obtain
\begin{equation}\label{n3p3_contr_lmm1b_e17}
\sum_{\chi_{\mathcal{T}}(\tilde{1},\hat{1})=0,|\chi_{\tilde{1}}({\mathcal{T}})|=1}
{\mathcal{C}}_{{\mathcal{U}}_{{\mathcal{T}}_1,{\mathcal{T}}}}(\al_1^{\perp})=
\big\lan 8a_{\hat{0}} + 4c_1({\mathcal{L}}_{\tilde{1}}^*),
{\overS}_1(\mu)\big\ran
+2\big|{\mathcal{S}}_{1;1}(\mu)\big|-2\big|{\mathcal{S}}_2(\mu)\big|.
\end{equation}
(4)\qua Finally, if ${\mathcal{T}} < {\mathcal{T}}_1$, 
$\chi_{\mathcal{T}}(\tilde{1},\hat{1}) = 0$,
and $\chi_{\tilde{1}}({\mathcal{T}}) = \{h_1,h_2\}$ is a two-element set,
the section ${\mathcal{D}}_{{\mathcal{T}},h_1}^{(1)}$ does not vanish 
over the set~${\mathcal{S}}_{{\mathcal{T}}|{\mathcal{T}}_1}(\mu)$.
We denote by 
$$\pi_{h_1}\co \ev_{\hat{0}}^*T\PPP\lra\Im{\mathcal{D}}_{{\mathcal{T}},h_1}^{(1)}
\qquad\hbox{and}\qquad
\pi_{h_1}^{\perp}\co \ev_{\hat{0}}^*T\PPP\lra E_{h_1}$$
the orthogonal projections, 
defined over a neighborhood of ${\mathcal{S}}_{{\mathcal{T}}|{\mathcal{T}}_1}(\mu)$
in ${\mathcal{U}}_{{\mathcal{T}}_1,{\mathcal{T}}}$, onto $\Im{\mathcal{D}}_{{\mathcal{T}},h_1}^{(1)}$
and its orthogonal complement $E_1$ in~$\ev_{\hat{0}}^*T\PPP$.
With appropriate identifications, ${\mathcal{S}}_{{\mathcal{T}}|{\mathcal{T}}_1}(\mu)$
is the zero set of the section
$\ev_{{\mathcal{T}}_1,M_0} \oplus 
\pi_{h_1}^{\perp} \circ {\mathcal{D}}_{{\mathcal{T}},h_2}^{(1)}$
of the bundle 
$$\ev_{{\mathcal{T}}_1,M_0}^*{\mathcal{N}}\De_{{\mathcal{T}}_1}(\mu)
\oplus L_{h_2}^* \otimes E_{h_1}$$
defined over a neighborhood of ${\mathcal{S}}_{{\mathcal{T}}|{\mathcal{T}}_1}(\mu)$
in~${\mathcal{U}}_{{\mathcal{T}}_1,{\mathcal{T}}}$.
By Lemma~\ref{str_lmm}, this section is transverse to the zero set.
By Proposition~\ref{str_prp} and Lemma~\ref{n3p3_str_lmm2}, 
\begin{gather*}
\ev_{{\mathcal{T}}_1,M_0}\big(\phi_{{\mathcal{T}}_1,{\mathcal{T}}}(b;\ups)\big)
=\ev_{{\mathcal{T}}_1,M_0}(b)+\ve_{-;1}(b;\ups),\\
\begin{split}
&\pi_{h_1}^{\perp}\{\ev_{\tilde{1}} \times \ev_{\hat{1}}\}
\big(\phi_{{\mathcal{T}}_1,{\mathcal{T}}}(\ups)\big)=
\big(y_{h_2;\hat{1}} - x_{\hat{1};h_2}\big)^{-1}  \otimes \\
&\qquad\qquad\qquad\qquad\qquad\qquad\qquad\qquad
\big\{\pi_{h_1}^{\perp}  \circ {\mathcal{D}}_{{\mathcal{T}},h_2}^{(1)}
 + \ve_{-;2}(\ups)\big\}
\rho_{{\mathcal{T}},\hat{1};h_2}^{(0;1)}(\ups_{h_2}),
\end{split}\\
\begin{split}
&\pi_{h_1}\{\ev_{\tilde{1}} \times \ev_{\hat{1}}\}
\big(\phi_{{\mathcal{T}}_1,{\mathcal{T}}}(\ups)\big)=
\big\{y_{h_1;\hat{1}} - x_{\hat{1};h_1}\big)^{-1}\otimes
\big({\mathcal{D}}_{{\mathcal{T}},h_1}^{(1)}
 + \ve_{+;h_1}(\ups)\big\}\rho_{{\mathcal{T}},\hat{1};h_1}^{(0;1)}(\ups_{h_1})\\
&\quad\qquad\qquad\qquad\qquad
+\big(y_{h_2;\hat{1}} - x_{\hat{1};h_2}\big)^{-1}\otimes
\big\{\pi_{h_1}{\mathcal{D}}_{{\mathcal{T}},h_2}^{(1)}
 + \ve_{+;h_2}(\ups)\big\}\rho_{{\mathcal{T}},\hat{1};h_2}^{(0;1)}(\ups_{h_2})
\end{split}
\end{gather*}
for all $(b;\ups) \in \mathcal{FT}_{\de} - Y(\mathcal{FT};\hat{I}^+)$
and some $C^1$--negligible maps
\begin{equation*}\begin{split}
\ve_{-;1},\ve_{-;2},\ve_{+;h}\co & \mathcal{FT}_{\de} - Y(\mathcal{FT};\hat{I}^+)\\
&\qquad \lra
\ev_{{\mathcal{T}}_1,M_0}^*{\mathcal{N}}\De_{{\mathcal{T}}_1}(\mu),
L_{h_2}^* \otimes E_{h_1},
L_h^* \otimes \Im{\mathcal{D}}_{{\mathcal{T}},h_1}^{(1)}.
\end{split}\end{equation*}
However, subtracting the expansion of 
$\{\ev_{\tilde{1}} \times \ev_{\hat{1}}\} \circ 
\phi_{{\mathcal{T}}_1,{\mathcal{T}}}$ in Lemma~\ref{n3p3_str_lmm2} 
multiplied by  \hbox{$\big(y_{\tilde{1};\hat{1}}(\ups)
 - x_{\tilde{1};h_1}(\ups)\big)$} 
from the expansion of 
${\mathcal{D}}_{{\mathcal{T}},\hat{1}}^{(1)} \circ \phi_{{\mathcal{T}}_1,{\mathcal{T}}}$
in (3b) of Proposition~\ref{str_prp}, we find~that
$${\mathcal{D}}_{{\mathcal{T}}_1,\tilde{1}}^{(1)}\phi_{{\mathcal{T}}_1,{\mathcal{T}}}(\ups)
=\big\{\al_2' + \ve(\ups)\big\}\rho(\ups)$$
for all $\ups \in \mathcal{FT}_{\de}$
such that $\phi_{{\mathcal{T}}_1,{\mathcal{T}}}(\ups)
 \in {\mathcal{U}}_{{\mathcal{T}}_1}^{(1)}(\mu)$,
where $\rho$ is a monomials map on $\mathcal{FT}$ 
with values in a rank-two bundle~$\tilde{\mathcal{F}}{\mathcal{T}}$,
$$\al_2'\co \tilde{\mathcal{F}}{\mathcal{T}} \lra
L_{\tilde{1}}^* \otimes \ev_{\hat{0}}^*T\PPP$$
is a linear map, such that $\al_+ \oplus \al$ is injective over
${\mathcal{U}}_{{\mathcal{T}}_1}^{(1)}(\mu)$, and
$$\ve\co  \mathcal{FT} - Y(\mathcal{FT};\hat{I}^+)\lra
\hbox{Hom}(\tilde{\mathcal{F}}{\mathcal{T}},
L_{\tilde{1}}^* \otimes \ev_{\hat{0}}^*T\PPP)$$
is a $C^0$--negligible map.
Explicitly, if $\hat{I} = \chi_{\tilde{1}}({\mathcal{T}})$,
$\rho$ is the identity map, and
$$\al\big(\ups_{h_1},\ups_{h_2}\big)=
\big(y_{\hat{1}} - x_{h_1}\big)^{-1}\otimes
\big(x_{h_2} - x_{h_2}\big)\otimes{\mathcal{D}}_{{\mathcal{T}},h_2}^{(1)}\ups_{h_2}.$$
Thus, if $\bar{\nu}_1$ is generic, 
${\mathcal{S}}_{{\mathcal{T}}|{\mathcal{T}}_1}(\mu)$ is
$\al_1^{\perp}$--hollow unless
$\hat{I}^+ = \chi_{\tilde{1}}({\mathcal{T}})$.
If $\hat{I}^+ = \chi_{\tilde{1}}({\mathcal{T}})$,
by dimension-counting $M_{\tilde{1}}{\mathcal{T}} = \{\hat{2}\}$
if ${\mathcal{S}}_{{\mathcal{T}}|{\mathcal{T}}_1}(\mu) \neq \eset$, and
by Proposition~\ref{euler_prp},
rescaling of the linear map, and the decomposition~\eqref{cart_split},
\begin{gather*}
{\mathcal{C}}_{{\mathcal{S}}_{{\mathcal{T}}|{\mathcal{T}}_1}(\mu)}(\al_1^{\perp})
=N(\al_2),\qquad\hbox{where}\\
\al_2 \in \Ga(\ov{\frak M}_{\{\tilde{1},\hat{1},h_1,h_2\}}
 \times {\mathcal{S}}_{{\overT}}(\mu);
\hbox{Hom}(\pi_1^*L_{\tilde{1}}^* \otimes \pi_2^*L_{h_2},
\pi_1^*L_{\tilde{1}}^* \otimes {\mathcal{O}}_1)\big), \\
\al_2(\tilde{\ups})=\pi_{\bar{\nu}_1}^{\perp}\circ\big\{
{\mathcal{D}}_{{\overT}}^{(1)}\ups\big\}.
\end{gather*} 
The set ${\mathcal{S}}_{{\overT}}(\mu)$ is finite
and ${\mathcal{D}}_{{\overT}}^{(1)}$ does not vanish.
Thus, by Propositions~\ref{zeros_prp} and~\ref{euler_prp},
\begin{equation}\label{n3p3_contr_lmm1b_e21}
\sum_{\chi_{\mathcal{T}}(\tilde{1},\hat{1})=0,|\chi_{\tilde{1}}({\mathcal{T}})|=2}
{\mathcal{C}}_{{\mathcal{U}}_{{\mathcal{T}}|{\mathcal{T}}_1}}(\al_1^{\perp})
=\blr{c_1(L_{\tilde{1}}^*),\ov{\frak M}_{\{\tilde{1},\hat{1},h_1,h_2\}}{\mathcal{T}}}
\big|{\mathcal{S}}_2(\mu)\big|
=\big|{\mathcal{S}}_2(\mu)\big|.
\end{equation}
The claim follows by plugging equations~\eqref{n3p3_contr_lmm1b_e11},
\eqref{n3p3_contr_lmm1b_e17}, and \eqref{n3p3_contr_lmm1b_e21} into
\eqref{n3p3_contr_lmm1b_e1} and using \eqref{psi_class1} and~\eqref{psi_class2}.
\end{proof}

\begin{lmm}
\label{n3p3_contr_lmm1d}
If ${\mathcal{O}}_1 \lra {\overV}_{1;1}^{(1)}(\mu)$ is the quotient of
the bundle $\ev_{\hat{0}}^*T\PPP$ by a generic trivial
line subbundle~$\Bbb{C}\bar{\nu}_1$,
$\pi_{\bar{\nu}_1}^{\perp}\co  \ev_{\hat{0}}^*T\PPP \lra {\mathcal{O}}_1$
is the quotient projection, and 
$\al_2 \in \Ga\big({\overV}_{1;1}^{(1)}(\mu);
\hbox{Hom}(L_{\tilde{1}},{\mathcal{O}}_1)\big)$ is given by
$\al_2 = \pi_{\bar{\nu}_1}^{\perp} \circ 
{\mathcal{D}}_{{\mathcal{T}}_{1;l},\tilde{1}}^{(1)}$
on~${\mathcal{U}}_{{\mathcal{T}}_{1;l}}^{(1)}(\mu)$,
$$N(\al_2)=
\blr{4a_{\hat{0}} + c_1(L_{\tilde{1}}^*),{\overV}_{1;1}^{(1)}(\mu)}
-2\big|{\mathcal{S}}_{1;1}(\mu)\big|.$$
\end{lmm}

\begin{proof}
By Propositions~\ref{zeros_prp} and~\ref{euler_prp},
\begin{equation}\label{n3p3_contr_lmm1d_e1}
N(\al_2)=\big\lan 4a_{\hat{0}} + c_1(L_{\tilde{1}}^*),
{\overV}_{1;1}^{(1)}(\mu)\big\ran
-{\mathcal{C}}_{\partial{\overV}_{1;1}^{(1)}(\mu)}(\al_2^{\perp}),
\end{equation}
where $\al_2^{\perp}$ denotes the composition of $\al_2$
with the projection $\pi_{\bar{\nu}_2}^{\perp}$
onto the quotient ${\mathcal{O}}_2$ of ${\mathcal{O}}_1$
by a generic trivial line subbundle~$\Bbb{C}\bar{\nu}_2$.
Suppose 
$${\mathcal{T}} \equiv (M_1,I;j,\under{d})<{\mathcal{T}}_{1;l}$$
is a bubble type such that ${\mathcal{D}}_{{\mathcal{T}}_{1;l},\tilde{1}}^{(1)}$
vanishes somewhere on 
${\overV}_{1;1}^{(1)}(\mu)\cap{\mathcal{U}}_{{\mathcal{T}}_{1;l},{\mathcal{T}}}$.
Then, by Lemmas~\ref{str_lmm}, \ref{n3p3_str_lmm1}, and~\ref{n3p3_str_lmm2},
$${\mathcal{T}} = {\mathcal{T}}_{1;l}(\hat{1})  \qquad\hbox{and}\qquad 
{\overV}_{1;1}^{(1)}(\mu)\cap{\mathcal{U}}_{{\mathcal{T}}_{1;l},{\mathcal{T}}}
\subset{\mathcal{S}}_{{\mathcal{T}}|{\mathcal{T}}_{1;l}}(\mu).$$
By the same argument as in (3) of the proof of  Lemma~\ref{n3p3_contr_lmm1b},
$${\mathcal{C}}_{{\mathcal{S}}_{{\mathcal{T}}|{\mathcal{T}}_{1;l}}(\mu)}(\al_2^{\perp})
=2\big|{\mathcal{S}}_{{\mathcal{T}}|{\mathcal{T}}_{1;l}}(\mu)\big|.$$
We conclude that
\begin{equation}\label{n3p3_contr_lmm1d_e3}
{\mathcal{C}}_{\partial{\overV}_{1;1}^{(1)}(\mu)}(\al_2^{\perp})
=2\big|{\mathcal{S}}_{1;1}(\mu)\big|.
\end{equation}
The claim follows from \eqref{n3p3_contr_lmm1d_e1}
and~\eqref{n3p3_contr_lmm1d_e3}.
\end{proof}

\begin{lmm}
\label{n3p3_contr_lmm1c}
If ${\mathcal{O}}_1 \lra {\overS}_1(\mu)$ is the quotient of
the bundle $\ev_{\hat{0}}^*T\PPP$ by a generic trivial
line subbundle~$\Bbb{C}\bar{\nu}_1$,
$\pi_{\bar{\nu}_1}^{\perp}\co  \ev_{\hat{0}}^*T\PPP \lra {\mathcal{O}}_1$
is the quotient projection, and 
$\al_2 \in \Ga\big({\overS}_1(\mu);
\hbox{Hom}(L_{\tilde{1}}^{\otimes2},{\mathcal{O}}_1)\big)$ is given by
$\al_2 = \pi_{\bar{\nu}_1}^{\perp} \circ 
{\mathcal{D}}_{{\mathcal{T}}_0,\tilde{1}}^{(2)}$,
$$N(\al_2)=\big\lan 4a_{\hat{0}} + 2c_1({\mathcal{L}}_{\tilde{1}}^*),
{\overS}_1(\mu)\big\ran-\big|{\mathcal{S}}_2(\mu)\big|.$$
\end{lmm}

\begin{proof}
(1)\qua By Propositions~\ref{zeros_prp} and~\ref{euler_prp}, 
\begin{equation}\label{n3p3_contr_lmm1c_e1}
N(\al_2)=\big\lan 4a_{\hat{0}} + 2c_1(L_{\tilde{1}}^*),
{\overS}_1(\mu)\big\ran-
{\mathcal{C}}_{\partial{\overS}_1(\mu)}(\al_2^{\perp}),
\end{equation}
where $\al_2^{\perp}$ denotes the composition of $\al_2$
with the projection $\pi_{\bar{\nu}_2}^{\perp}$
onto the quotient ${\mathcal{O}}_2$ of ${\mathcal{O}}_1$
by a generic trivial line subbundle~$\Bbb{C}\bar{\nu}_2$.
We now use the expansion (3d) of Proposition~\ref{str_prp}
to describe the boundary strata of ${\overS}_1(\mu)$ and 
compute the contribution of each stratum to 
${\mathcal{C}}_{\partial{\overS}_1(\mu)}(\al_2^{\perp})$.
Our description shows that $\partial{\overS}_1(\mu)$ is a finite
set and thus ${\mathcal{S}}_1(\mu)$ is a one-pseudovariety in 
${\overV}_1(\mu)$ and~${\overV}_1$.

(2)\qua If ${\mathcal{T}} = (M_0,I;j,\under{d}) < {\mathcal{T}}_0$,
by Proposition~\ref{str_prp},
$${\mathcal{D}}_{{\mathcal{T}}_0,\tilde{1}}^{(1)}\phi_{{\mathcal{T}}_0,{\mathcal{T}}}(\ups)
=\sum_{h\in\chi_{\tilde{1}}({\mathcal{T}})}  
\big({\mathcal{D}}_{{\mathcal{T}},h}^{(1)} + \ve_{{\mathcal{T}},\tilde{1};h}^{(1)}\big)
\rho_{{\mathcal{T}},\tilde{1}}^{(0;1)}(\ups),
\quad\hbox{where}\quad
\rho_{{\mathcal{T}},\tilde{1}}^{(0;1)}(\ups)=
\prod_{i\in(\tilde{1},h]}   \ups_i,$$
for all $\ups \in \mathcal{FT}$ sufficiently small.
Thus,
${\overS}_1(\mu)\cap{\mathcal{U}}_{{\mathcal{T}}_0,{\mathcal{T}}}$
is contained in the finite set~${\mathcal{S}}_{{\mathcal{T}}|{\mathcal{T}}_0}(\mu)$.
If $d_{\tilde{1}} \neq 0$, section ${\mathcal{D}}_{{\mathcal{T}}_0,\tilde{1}}^{(2)}$
does not vanish on ${\mathcal{S}}_{{\mathcal{T}}|{\mathcal{T}}_0}(\mu)$
and thus ${\mathcal{U}}_{{\mathcal{T}}|{\mathcal{T}}_0}$ does not contribute to
${\mathcal{C}}_{\partial{\overS}_1(\mu)}(\al_2^{\perp})$,
since ${\mathcal{D}}_{{\mathcal{T}}_0,\tilde{1}}^{(1)}$ is defined
everywhere on~${\overS}_1(\mu)$.

(3)\qua If $d_{\tilde{1}} = 0$ and 
${\mathcal{S}}_{{\mathcal{T}}|{\mathcal{T}}_0}(\mu) \neq \eset$,
 $\hat{I}^+ = \chi_{\tilde{1}}({\mathcal{T}})$
and $|\chi_{\tilde{1}}({\mathcal{T}})| \in \{1,2\}$.
Suppose $\chi_{\tilde{1}}({\mathcal{T}}) = \{h\}$ is a single-element set,
ie ${\mathcal{T}} = {\mathcal{T}}_0(l)$ for some $l \in [N]$.
With appropriate identifications, ${\mathcal{S}}_{{\mathcal{T}}|{\mathcal{T}}_0}(\mu)$
is the zero set of the section
$\ev_{{\mathcal{T}}_0,M_0} \oplus {\mathcal{D}}_{{\mathcal{T}},h}^{(1)}$
of the bundle 
$$\ev_{{\mathcal{T}}_0,M_0}^*{\mathcal{N}}\De_{{\mathcal{T}}_0}(\mu)
\oplus L_h^* \otimes \ev_{\hat{0}}^*T\PPP$$
defined over a neighborhood of ${\mathcal{S}}_{{\mathcal{T}}|{\mathcal{T}}_0}(\mu)$
in~${\mathcal{U}}_{{\mathcal{T}}_0,{\mathcal{T}}}$.
By Lemma~\ref{str_lmm}, this section is transverse to the zero set.
By Proposition~\ref{str_prp} and Lemma~\ref{n3p3_str_lmm2}, 
\begin{gather*}
\ev_{{\mathcal{T}}_0,M_0}\big(\phi_{{\mathcal{T}}_0,{\mathcal{T}}}(b;\ups)\big)
=\ev_{{\mathcal{T}}_0,M_0}(b)+\ve_{-;1}(b;\ups),\\
{\mathcal{D}}_{{\mathcal{T}}_0,\tilde{1}}^{(1)}
\big(\phi_{{\mathcal{T}}_0,{\mathcal{T}}}(\ups)\big)=
\big\{{\mathcal{D}}_{{\mathcal{T}},h}^{(1)} + \ve_{-;2}(\ups)\big\}\ups,
\end{gather*}
for all $(b;\ups) \in \mathcal{FT}_{\de} - Y(\mathcal{FT};\hat{I}^+)$
and some $C^1$--negligible maps
$$\ve_{-;1},\ve_{-;2}\co 
\mathcal{FT}_{\de} - Y(\mathcal{FT};\hat{I}^+)\lra
\ev_{{\mathcal{T}}_0,M_0}^*{\mathcal{N}}\De_{{\mathcal{T}}_0}(\mu),
L_h^* \otimes \ev_{\hat{0}}^*T\PPP.$$
On the other hand, subtracting the expansion of 
${\mathcal{D}}_{{\mathcal{T}}_0,\tilde{1}}^{(1)} \circ \phi_{{\mathcal{T}}_0,{\mathcal{T}}}$ 
of (3d) of Proposition~\ref{str_prp} times $x_h$ 
from the expansion of 
${\mathcal{D}}_{{\mathcal{T}}_0,\tilde{1}}^{(2)} \circ \phi_{{\mathcal{T}}_0,{\mathcal{T}}}$ 
of (3d) of Proposition~\ref{str_prp}, we~obtain
$${\mathcal{D}}_{{\mathcal{T}}_0,\tilde{1}}^{(2)}\phi_{{\mathcal{T}}_0,{\mathcal{T}}}(\ups)
=\big\{{\mathcal{D}}_{{\mathcal{T}},h}^{(2)} + 
\ve_{{\mathcal{T}},\tilde{1};h}^{(2)}(\ups)\big\}\ups \otimes \ups
\quad\mbox{for all}\  \ups \in \mathcal{FT}_{\de}~\mbox{with}~
\phi_{{\mathcal{T}}_0,{\mathcal{T}}}(\ups) \in {\mathcal{S}}_1.$$
By Lemma~\ref{str_lmm}, ${\mathcal{D}}_{{\mathcal{T}},h}^{(2)}$ does not vanish
on the finite set ${\mathcal{S}}_{{\mathcal{T}}|{\mathcal{T}}_0}(\mu)$.
Thus, by Proposition~\ref{euler_prp},
\begin{equation}\label{n3p3_contr_lmm1c_e3}
\sum_{|\chi_{\tilde{1}}({\mathcal{T}})|=1}  
{\mathcal{C}}_{{\mathcal{U}}_{{\mathcal{T}}_0,{\mathcal{T}}}(\mu)}(\al_2^{\perp})
=2\big|{\mathcal{S}}_{1;1}(\mu)\big|.
\end{equation}
(4)\qua If $\chi_{\tilde{1}}({\mathcal{T}}) = \{h_1,h_2\}$ is a two-element set,
the section ${\mathcal{D}}_{{\mathcal{T}},h_1}^{(1)}$ does not vanish 
over the set ${\mathcal{S}}_{{\mathcal{T}}|{\mathcal{T}}_0}(\mu)$.
Let $\pi_{h_1}$, $\pi_{h_1}^{\perp}$, and $E_{h_1}$ be
as in (4) of the proof of Lemma~\ref{n3p3_contr_lmm1b}.
Similarly to the previous case, 
${\mathcal{S}}_{{\mathcal{T}}|{\mathcal{T}}_0}(\mu)$ is the zero set of the section
$$\ev_{{\mathcal{T}}_0,M_0} \oplus
\pi_{h_1}^{\perp} \circ {\mathcal{D}}_{{\mathcal{T}},h_2}^{(1)}$$
of the bundle 
$\ev_{{\mathcal{T}}_0,M_0}^*{\mathcal{N}}\De_{{\mathcal{T}}_0}(\mu)
 \oplus L_{h_2}^* \otimes E_{h_1}$
defined over a neighborhood of ${\mathcal{S}}_{{\mathcal{T}}|{\mathcal{T}}_0}(\mu)$
in~${\mathcal{U}}_{{\mathcal{T}}_0,{\mathcal{T}}}$.
By Lemma~\ref{str_lmm}, this section is transverse to the zero set.
By Proposition~\ref{str_prp} and Lemma~\ref{n3p3_str_lmm2}, 
\begin{gather*}
\ev_{{\mathcal{T}}_0,M_0}\big(\phi_{{\mathcal{T}}_0,{\mathcal{T}}}(b;\ups)\big)
=\ev_{{\mathcal{T}}_0,M_0}(b)+\ve_{-;1}(b;\ups),\\
\pi_{h_1}^{\perp}{\mathcal{D}}_{{\mathcal{T}}_0,\tilde{1}}^{(1)}
\big(\phi_{{\mathcal{T}}_0,{\mathcal{T}}}(\ups)\big)=
\big\{\pi_{h_1}^{\perp} \circ {\mathcal{D}}_{{\mathcal{T}},h_2}^{(1)}
 + \ve_{-;2}(\ups)\big\}\ups_{h_2},\\
\pi_{h_1}{\mathcal{D}}_{{\mathcal{T}}_0,\tilde{1}}^{(1)}
\big(\phi_{{\mathcal{T}}_0,{\mathcal{T}}}(\ups)\big)=
\big\{{\mathcal{D}}_{{\mathcal{T}},h_1}^{(1)}
 + \ve_{+;h_1}(\ups)\big\}\ups_{h_1})+
\big\{\pi_{h_1}{\mathcal{D}}_{{\mathcal{T}},h_2}^{(1)}
 + \ve_{+;h_2}(\ups)\big\}\ups_{h_2}
\end{gather*}
for all $(b;\ups) \in \mathcal{FT}_{\de} - Y(\mathcal{FT};\hat{I}^+)$
and some $C^1$--negligible maps $\ve_{-;1},\ve_{-;2},\ve_{+;h}$.
On the other hand, subtracting the expansion of
${\mathcal{D}}_{{\mathcal{T}},\tilde{1}}^{(1)} \circ \phi_{{\mathcal{T}}_0},{\mathcal{T}}$ 
of (3d) of Proposition~\ref{str_prp} times $x_{h_1}$ 
from the expansion of 
${\mathcal{D}}_{{\mathcal{T}}_0,\tilde{1}}^{(2)} \circ \phi_{{\mathcal{T}}_0,{\mathcal{T}}}$ 
of (3d) of Proposition~\ref{str_prp}, we~obtain
\begin{gather*}
{\mathcal{D}}_{{\mathcal{T}}_0,\tilde{1}}^{(2)}\phi_{{\mathcal{T}}_0,{\mathcal{T}}}(\ups)=
(x_{h_2} - x_{h_1}) \otimes 
\big\{{\mathcal{D}}_{{\mathcal{T}},h_2}^{(1)} + \ve(\ups)\big\}\ups_{h_2}\\
\mbox{for all}\  \ups \in \mathcal{FT}_{\de}
~\mbox{such that}~\phi_{{\mathcal{T}}_0,{\mathcal{T}}}(\ups) \in {\mathcal{S}}_1(\mu),
\end{gather*}
where  $\ve$ is a $C^0$--negligible map.
Since ${\mathcal{D}}_{{\mathcal{T}},h}^{(1)}$ does not vanish on
${\mathcal{S}}_{{\mathcal{T}}|{\mathcal{T}}_0}(\mu)$, 
from Proposition~\ref{euler_prp}, we conclude that
\begin{equation}\label{n3p3_contr_lmm1c_e5}
\sum_{|\chi_{\tilde{1}}({\mathcal{T}})|=2}  
{\mathcal{C}}_{{\mathcal{U}}_{{\mathcal{T}}|{\mathcal{T}}_0}(\mu)}(\al_2^{\perp})
=\big|{\mathcal{S}}_2(\mu)\big|.
\end{equation}
The claim follows by plugging equations~\eqref{n3p3_contr_lmm1c_e3} and 
\eqref{n3p3_contr_lmm1c_e5} into \eqref{n3p3_contr_lmm1c_e1} and 
using \eqref{psi_class1} and~\eqref{psi_class2}. 
\end{proof}

\begin{lmm}
\label{n3p3_contr_lmm2}
Suppose ${\mathcal{T}} = (M_2,I;j,\under{d})$ is a bubble type such that
${\mathcal{T}} < {\mathcal{T}}_2$, $\chi_{\mathcal{T}}(\tilde{1},\hat{2}) > 0$, and
\hbox{$\chi_{\mathcal{T}}(\hat{1},\hat{2}) = 0$}.

{\rm(1)}\qua If $|\hat{I}^+| \neq 1$ or $|M_{j_{\hat{1}}}{\mathcal{T}}| \neq 1$,
${\mathcal{C}}_{{\mathcal{U}}_{{\mathcal{T}}_2,{\mathcal{T}}}}
\big(\ev_{\tilde{1}} \times \ev_{\hat{2}};\De_{\PPP\times\PPP}\big)
 = 0$.

{\rm(2)}\qua If $|\hat{I}^+| = 1$ and $|M_{j_{\hat{1}}}{\mathcal{T}}| = 1$,
\begin{equation*}\begin{split}
{\mathcal{C}}_{{\mathcal{U}}_{{\mathcal{T}}_2,{\mathcal{T}}}}
\big(\ev_{\tilde{1}} \times \ev_{\hat{2}};\De_{\PPP\times\PPP}\big)
=\big\lan 6a_{\hat{0}}^2 + 4a_{\hat{0}}c_1({\mathcal{L}}_{\tilde{1}}^*)
 + c_1^2({\mathcal{L}}_{\tilde{1}}^*),{\overV}_1^{(1)}(\mu)\big\ran \qquad\quad&\\
+\big|{\mathcal{S}}_2(\mu)\big|-\big|{\mathcal{V}}_2^{(1)}(\mu)\big|
-\big\lan 8a_{\hat{0}}^2 + 4c_1({\mathcal{L}}_{\tilde{1}}^*),
{\overS}_1(\mu)\big\ran.&
\end{split}\end{equation*}
\end{lmm}

\begin{proof}
This lemma follows from Lemma~\ref{n3p3_contr_lmm1} by symmetry.
\end{proof}

\begin{remark}
{\rm These contributions can be computed directly, 
ie similarly to the proof of Lemma~\ref{n3p3_contr_lmm1}, and in fact 
one finds a somewhat different expression for the contribution in~(2).  
What this means is that we have found a relationship
between certain intersection  numbers:
\begin{equation}\label{rel_e1}
\blr{8a_{\hat{0}}c_1({\mathcal{L}}_{\tilde{1}}^*),{\overV}_1^{(1)}(\mu)}
+2\big|{\mathcal{V}}_2^{(1)}(\mu)\big|
=\blr{4a_{\hat{0}}-\eta_{\hat{0},1},{\overV}_2^{(1,1)}(\mu)}
+\big|{\mathcal{V}}_{2,(0,1)}^{(1;0,1)}(\mu)\big|.
\end{equation}
Here  ${\mathcal{V}}_{2,(0,1)}^{(1;0,1)}(\mu)$ denotes a set of tuples
of stable maps
whose cardinality is six times the number of rational curves that pass through 
the constraints~$\mu$ and 
have the form described by the last picture of Figure~\ref{images_fig}.
Using Lemma~2.2.2 in~\cite{P}, it is possible to restate the relation~\eqref{rel_e1} 
in terms of numbers of rational curves of various~shapes.}
\end{remark}

\begin{lmm}
\label{n3p3_contr_lmm3}
If ${\mathcal{T}} = (M_2,I;j,\under{d})$ is a bubble type such that
${\mathcal{T}} < {\mathcal{T}}_2$ and 
$\chi_{\mathcal{T}}(\tilde{1},\hat{2}) = \chi_{\mathcal{T}}(\hat{1},\hat{2}) = 0$,
$${\mathcal{C}}_{{\mathcal{U}}_{{\mathcal{T}}_2,{\mathcal{T}}}}
\big(\ev_{\tilde{1}} \times \ev_{\hat{2}};\De_{\PPP\times\PPP}\big) = 0.$$
\end{lmm}

\proof
(1)\qua Since $\chi_{\mathcal{T}}(\tilde{1},\hat{1}) = 0$, by 
Lemma~\ref{n3p3_str_lmm2},
$${\overU}_{{\mathcal{T}}_2}^{(1)}(\mu)\cap{\mathcal{U}}_{{\mathcal{T}}_2,{\mathcal{T}}}
\subset {\mathcal{S}}_{{\mathcal{T}}|{\mathcal{T}}_2}(\mu).$$
If ${\mathcal{S}}_{{\mathcal{T}}|{\mathcal{T}}_2}(\mu) \neq \eset$,
\hbox{$|\chi_{\tilde{1}}({\mathcal{T}})| \in \{1,2\}$}.
If $h$ is the unique element of~$\chi_{\tilde{1}}({\mathcal{T}})$,
mixing the argument in (1) of the proof of Lemma~\ref{n3p3_contr_lmm1}
with (3) of the proof of Lemma~\ref{n3p3_contr_lmm1b},
we find that ${\mathcal{S}}_{{\mathcal{T}}|{\mathcal{T}}_2}(\mu)$ is
$(\ev_{\tilde{1}} \times \ev_{\hat{2}},\De_{\PPP\times\PPP})$--hollow
unless $\hat{I}^+ = \chi_{\tilde{1}}({\mathcal{T}})$.
On the other hand, if 
$\hat{I}^+ = \chi_{\tilde{1}}({\mathcal{T}}) = \{h\}$,
by Proposition~\ref{euler_prp}, the decomposition~\eqref{cart_split},
and a rescaling of the linear~map,
\begin{gather*}
{\mathcal{C}}_{{\mathcal{U}}_{{\mathcal{T}}_2,{\mathcal{T}}}}
\big(\ev_{\tilde{1}} \times \ev_{\hat{2}};\De_{\PPP\times\PPP}\big)
=2N(\al_1),\qquad\hbox{where}\\
\al_1 = \pi_2^*
{\mathcal{D}}_{{\overT}}^{(2)}
\in\Ga\big(\ov{\frak M}_{\{\tilde{1},h\}\sqcup M_{\tilde{1}}{\mathcal{T}}} \times 
{\mathcal{S}}_{{\overT}}(\mu);
\pi_2^*\hbox{Hom}(L_{\tilde{1}}^*,\ev_{\hat{0}}^*T\PPP)\big).
\end{gather*}
Since $|M_{\tilde{1}}{\mathcal{T}}| \ge 2$,
the first factor is positive-dimensional, while
the linear map~$\al_1$ comes entirely from the second factor.
Thus,
$${\mathcal{C}}_{{\mathcal{U}}_{{\mathcal{T}}_2,{\mathcal{T}}}}
\big(\ev_{\tilde{1}} \times \ev_{\hat{2}};\De_{\PPP\times\PPP}\big)=0.$$
(2)\qua If $\chi_{\tilde{1}}({\mathcal{T}}) = \{h_1,h_2\}$ is a two-element set,
the section ${\mathcal{D}}_{{\mathcal{T}},h_1}^{(1)}$ does not vanish 
over the set ${\mathcal{S}}_{{\mathcal{T}}|{\mathcal{T}}_2}(\mu)$.
Mixing (1) of the proof of Lemma~\ref{n3p3_contr_lmm1}
with (4) of the proof of Lemma~\ref{n3p3_contr_lmm1b},
we find that ${\mathcal{S}}_{{\mathcal{T}}|{\mathcal{T}}_2}(\mu)$ is
$(\ev_{\tilde{1}} \times \ev_{\hat{2}},\De_{\PPP\times\PPP})$--hollow
unless $\hat{I}^+ = \chi_{\tilde{1}}({\mathcal{T}})$.
On the other hand, if 
$\hat{I}^+ = \chi_{\tilde{1}}({\mathcal{T}}) = \{h_1,h_2\}$,
by Proposition~\ref{euler_prp}, the decomposition~\eqref{cart_split},
and a rescaling of the linear~map,
\begin{gather*}
{\mathcal{C}}_{{\mathcal{U}}_{{\mathcal{T}}_2,{\mathcal{T}}}}
\big(\ev_{\tilde{1}} \times \ev_{\hat{2}};\De_{\PPP\times\PPP}\big)
=N(\al_1),\qquad\hbox{where}\\
\al_1\in\Ga\big(\ov{\frak M}_{\{\tilde{1},h_1,h_2\}\sqcup M_{\tilde{1}}{\mathcal{T}}}
 \times  {\mathcal{S}}_{{\overT}}(\mu);
\pi_2^*\hbox{Hom}(L_{h_1} \oplus L_{h_2},
\Im{\mathcal{D}}_{{\overT}}^{(1)} \oplus 
\ev_{\hat{0}}^*T\PPP)\big),\\
\al_1\big(\ups_{h_1},\ups_{h_2}\big)
 = \big(
{\mathcal{D}}_{{\overT}}^{(1)}\ups_{h_1} + 
{\mathcal{D}}_{{\overT}}^{(1)}\ups_{h_2},
{\mathcal{D}}_{{\overT}}^{(1)}\ups_{h_2}\big). 
\end{gather*}
Since $|M_{\tilde{1}}{\mathcal{T}}| \ge 2$,
the first factor is positive-dimensional, while
the linear map~$\al_1$ comes entirely from the second factor.
Thus,
$${\mathcal{C}}_{{\mathcal{U}}_{{\mathcal{T}}_2,{\mathcal{T}}}}
\big(\ev_{\tilde{1}} \times
\ev_{\hat{2}};\De_{\PPP\times\PPP}\big)=0.\eqno{\qed}$$

\begin{crl}
\label{n3p3_contr_crl}
The contribution from the boundary to the number
$\llrr{{\mathcal{V}}_1^{(2)}\!(\mu)}$ is given~by
\begin{equation*}\begin{split}
{\mathcal{C}}_{\partial{\overU}_{{\mathcal{T}}_2}^{(2)}}
\big(\ev_{\tilde{1}} \times \ev_{\hat{2}};\De_{\PPP\times\PPP}\big)
=\blr{12a_{\hat{0}}^2 + 8a_{\hat{0}}\eta_{\hat{0},1}
 + 2\eta_{\hat{0},1}^2,{\overV}_1^{(1)}(\mu)}
+2\big|{\mathcal{S}}_2(\mu)\big| \qquad& \\
-2\big|{\mathcal{V}}_2^{(1)}(\mu)\big| 
- \blr{16a_{\hat{0}} + 8\eta_{\hat{0},1},{\overS}_1(\mu)}.&
\end{split}\end{equation*}
\end{crl}

\begin{proof}
This corollary follows immediately from
Lemmas~\ref{n3p3_str_lmm1}--\ref{n3p3_exp_lmm1},
\ref{n3p3_contr_lmm1}, \ref{n3p3_contr_lmm2}, and~\ref{n3p3_contr_lmm3}.
\end{proof}

\section{Example 2: Rational tacnodal curves in $\PPP$}
\label{n3tac_sec}

\subsection{Summary}
\label{n3tac_sum_subs}

In this section, we prove Theorem~\ref{n3tac_thm}.
The general approach is the same as in Section~\ref{n3p3_sec}.
If ${\mathcal{T}}_1$ is the bubble type as defined in Section~\ref{n3p3_sec}, let
\begin{gather*}
{\mathcal{S}}_1^{(1)}(\mu)=\big\{(b,[\ups_{\tilde{1}},\ups_{\hat{1}}])
 \in \Bbb{P}(L_{\tilde{1}} \oplus L_{\tilde{1}}^*)
 \to {\overV}_1^{(1)}: b \in {\mathcal{V}}_1^{(1)}(\mu),~
{\mathcal{D}}_{\tilde{1},\hat{1}}(b,[\ups_{\tilde{1}},\ups_{\hat{1}}]) = 0\big\}\\
\hbox{where}\qquad
{\mathcal{D}}_{\tilde{1},\hat{1}} \in \Ga\big(
\Bbb{P}(L_{\tilde{1}} \oplus L_{\tilde{1}}^*)|{\mathcal{V}}_1^{(1)}(\mu);
\ga_{L_{\tilde{1}} \oplus L_{\tilde{1}}^*}^* \otimes 
\ev_{\hat{0}}^*T\PPP\big),\\
{\mathcal{D}}_{\tilde{1},\hat{1}}(\ups_{\tilde{1}},\ups_{\hat{1}})=
{\mathcal{D}}_{{\mathcal{T}}_1,\tilde{1}}^{(1)}\ups_{\tilde{1}}+
{\mathcal{D}}_{{\mathcal{T}}_1,\hat{1}}^{(1)}\ups_{\hat{1}}.
\end{gather*}
The set ${\mathcal{S}}_1^{(1)}(\mu)$ can be identified with the set
of rational one-component tacnodal curves passing through 
the constraints~$\mu$, but with a choice of a branch at each~node.
In particular, the cardinality of the set ${\mathcal{S}}_1^{(1)}(\mu)$
is twice the enumerative number of Theorem~\ref{n3tac_thm}.

Note that section ${\mathcal{D}}_{\tilde{1},\hat{1}}$ does \textit{not} extend continuously 
over all of the boundary of ${\overV}_1^{(1)}(\mu)$.
In fact, this can be seen from (3c) of Proposition~\ref{str_prp}.
Nevertheless, the behavior of this section can be understood everywhere.
By Proposition~\ref{euler_prp} and equation~\eqref{zeros_main_e},
we~have
\begin{equation}\label{n3tac_e1}
\big|{\mathcal{S}}_1^{(1)}(\mu)\big|=
\blr{6a_{\hat{0}}^2 + c_1^2(L_{\tilde{1}}^*),{\overV}_1^{(1)}(\mu)}
-{\mathcal{C}}_{\partial\Bbb{P}(L_{\tilde{1}}\oplus L_{\tilde{1}}^*)}
({\mathcal{D}}_{\tilde{1},\hat{1}}),
\end{equation}
where 
${\mathcal{C}}_{\partial\Bbb{P}(L_{\tilde{1}}\oplus L_{\tilde{1}}^*)}({\mathcal{D}}_{\tilde{1},\hat{1}})$ 
is the ${\mathcal{D}}_{\tilde{1},\hat{1}}$--contribution from the boundary strata of 
$\Bbb{P}(L_{\tilde{1}} \oplus L_{\tilde{1}}^*)$
to the euler class of 
\hbox{$\ga_{L_{\tilde{1}} \oplus L_{\tilde{1}}^*}^* \otimes \ev_{\hat{0}}^*T\PPP$}.
This contribution is computed in the rest of this section.
Theorem~\ref{n3tac_thm} is obtained by plugging the expressions of
Corollaries~\ref{n3tac_contr_crl1} and~\ref{n3tac_contr_crl2} into~\eqref{n3tac_e1} 
and then using identities~\eqref{psi_class1} and~\eqref{psi_class2}.

Before proceeding with our computation of the contributions from various
strata, we observe that the section 
${\mathcal{D}}_{\tilde{1},\hat{1}}$ extends over 
${\overV}_1^{(1)}(\mu)\cap{\mathcal{U}}_{{\mathcal{T}}_1,{\mathcal{T}}}$
if 
$${\mathcal{T}} \equiv (M_1,I;j,\under{d})<{\mathcal{T}}_1$$ 
is a bubble type such that $j_{\hat{1}} = \tilde{1}$, 
as can be seen from Proposition~\ref{str_prp}.
If in addition $d_{\tilde{1}} \neq 0$, by Lemmas~\ref{str_lmm}
and~\ref{n3p3_str_lmm1},
${\mathcal{D}}_{\tilde{1},\hat{1}}$ does not vanish over 
$${\overV}_1^{(1)}(\mu)\cap{\mathcal{U}}_{{\mathcal{T}}_1,{\mathcal{T}}}
\subset {\mathcal{U}}_{{\mathcal{T}}|{\mathcal{T}}_1}^{(1)}(\mu).$$
Thus, such spaces ${\mathcal{U}}_{{\mathcal{T}}_1,{\mathcal{T}}}$
do not contribute to 
${\mathcal{C}}_{\partial\Bbb{P}(L_{\tilde{1}}\oplus L_{\tilde{1}}^*)}
({\mathcal{D}}_{\tilde{1},\hat{1}})$ and will not be considered below.

\subsection[Contributions from the spaces UT1T2 with kT(1,1)=0]{Contributions from the spaces ${\mathcal{U}}_{{\mathcal{T}}_1,{\mathcal{T}}}$
with $\chi_{\mathcal{T}}(\tilde{1},\hat{1}) = 0$}
\label{n3tac_corr_subs1}

In this subsection, we prove Corollary~\ref{n3tac_contr_crl1},
which gives the total contribution to
${\mathcal{C}}_{\partial\Bbb{P}(L_{\tilde{1}}\oplus L_{\tilde{1}}^*)}
({\mathcal{D}}_{\tilde{1},\hat{1}})$
from all the spaces ${\mathcal{U}}_{{\mathcal{T}}_1,{\mathcal{T}}}$,
where 
$${\mathcal{T}} \equiv (M_1,I;j,\under{d}) < {\mathcal{T}}_1$$ 
is a bubble type such that $\chi_{\mathcal{T}}(\tilde{1},\hat{1}) = 0$.
We use Lemma~\ref{n3p3_str_lmm2}, which describes the intersection
${\overV}_1^{(1)}(\mu)\cap{\mathcal{U}}_{{\mathcal{T}}_1,{\mathcal{T}}}$,
along with Proposition~\ref{str_prp}.

\begin{figure}[ht!]\small
\begin{pspicture}(-1.1,-4)(10,2.7)
\psset{unit=.36cm}
\pscircle[fillstyle=solid,fillcolor=gray](1,-2){1}
\pscircle*(1,-3){.22}\rput(1,-3.8){$\tilde{1}$}
\pscircle*(1.71,-1.29){.2}\rput(2.2,-1){$\hat{1}$}
\pnode(1,-3){A0}\pnode(1.71,-1.29){B0}
\ncarc[nodesep=.35,arcangleA=-90,arcangleB=-75,ncurv=1.2]{<->}{A0}{B0}
% 2nd column starts here
\psline[linewidth=.12]{->}(3,-1.5)(6.5,2)
\psline[linewidth=.12]{->}(3,-2)(6.5,-2)
\psline[linewidth=.12]{->}(3,-2.5)(6.5,-6)
% 3rd column starts here
\pscircle(10,2.5){1}
\pscircle[fillstyle=solid,fillcolor=gray](8.58,3.92){1}
\pscircle*(10,1.5){.22}\rput(10,.7){$\tilde{1}$}
\pscircle*(9.29,3.21){.19}\pscircle*(10.71,3.21){.2}
\rput(11.2,3.5){$\hat{1}$}
% 2nd row starts here
\pscircle(10,-2.5){1}
\pscircle[fillstyle=solid,fillcolor=gray](8.58,-1.08){1}
\pscircle*(10,-3.5){.22}\rput(10,-4.3){$\tilde{1}$}
\pscircle*(9.29,-1.79){.19}
\pscircle*(10.71,-1.79){.2}\rput(11.2,-1.5){$\hat{1}$}
\pscircle*(10.71,-3.21){.2}\rput(11.2,-3.5){$l$}
% 3rd row starts here
\pscircle(10,-7.5){1}
\pscircle[fillstyle=solid,fillcolor=gray](8.58,-6.08){1}
\pscircle*(10,-8.5){.22}\rput(10,-9.3){$\tilde{1}$}
\pscircle*(9.29,-6.79){.19}
\pscircle[fillstyle=solid,fillcolor=gray](11.42,-6.08){1}
\pscircle*(10.71,-6.79){.19}
\pscircle*(10.71,-8.21){.2}\rput(11.2,-8.5){$\hat{1}$}
% 4th column starts now
\rput(14,3){$\approx$}\rput(14,-2){$\approx$}\rput(14,-7){$\approx$}
% 5th column starts here
\rput(17.5,3){${\overS}_1(\mu)$}
\rput(18.2,-2.1){$\ov{\frak M}_{0,4}\times{\mathcal{S}}_{1;1}(\mu)$}
%3rd row starts here
\rput(18.2,-7.1){$\ov{\frak M}_{0,4}\times{\mathcal{S}}_2(\mu)$}
% 6th column
\psline{->}(22,3.25)(26.5,4)\psline{->}(22,2.75)(26.5,2)
\psline{->}(22,-1.75)(26.5,-1)\psline{->}(22,-2.25)(26.5,-3)
\psline{->}(22,-6.75)(26.5,-6)\psline{->}(22,-7.25)(26.5,-8)
% 7th column
\rput(27.5,4.25){$\times2$}\rput(27.5,1.75){$\times3$}
\rput(27.5,-.75){$\times2$}\rput(27.5,-3.25){$\times3$}
\rput(27.5,-5.75){$\times1$}\rput(27.5,-8.25){$\times2$}
% and labels now
\pnode(12.3,.5){C}\rput(13,.5){\footnotesize cusp}
\pnode(9.29,3.21){C1}\pnode(9.29,-1.79){C2}
\ncarc[nodesep=.3,arcangleA=15,arcangleB=30,ncurv=.6]{->}{C}{C1}
\ncarc[nodesep=.3,arcangleA=-15,arcangleB=60,ncurv=.6]{->}{C}{C2}
\pnode(7.5,-8){T1}\rput(6,-8){\footnotesize tacnode}
\pnode(10,-7.5){T2}\pnode(9.44,-6.94){T2a}\pnode(10.56,-6.94){T2b}\
\ncarc[nodesep=0,arcangleA=0,arcangleB=0,ncurv=.5]{-}{T1}{T2}
\ncline[nodesep=0]{->}{T2}{T2a}\ncline[nodesep=0]{->}{T2}{T2b}
\end{pspicture}
\caption{An outline of Subsection~\ref{n3tac_corr_subs1}}
\label{n3tac_contr_fig1}
\end{figure}
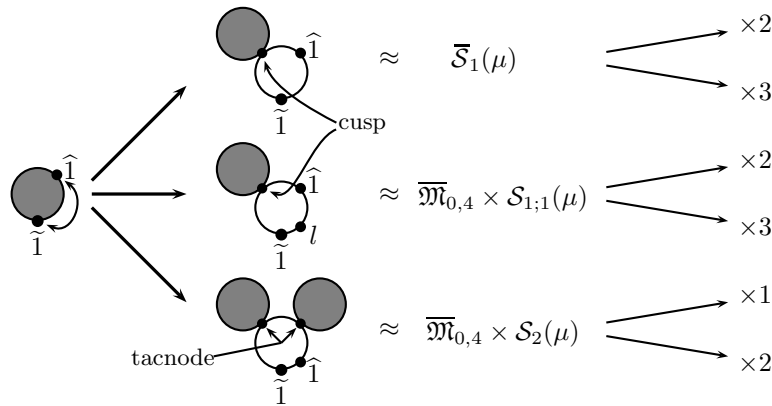

Figure~\ref{n3tac_contr_fig1} shows the three types of boundary strata 
${\overV}_1^{(1)}(\mu)\cap{\mathcal{U}}_{{\mathcal{T}}_1,{\mathcal{T}}}$ such that
$$\Bbb{P}(L_{\tilde{1}} \oplus L_{\tilde{1}}^*)|
{\overV}_1^{(1)}(\mu)\cap{\mathcal{U}}_{{\mathcal{T}}_1,{\mathcal{T}}}$$
is not contained in a finite union of 
${\mathcal{D}}_{\tilde{1},\hat{1}}$--hollow sets.
For such boundary strata, 
$$\Bbb{P}(L_{\tilde{1}} \oplus L_{\tilde{1}}^*)|
{\overV}_1^{(1)}(\mu)\cap{\mathcal{U}}_{{\mathcal{T}}_1,{\mathcal{T}}}$$
is a union of two 
${\mathcal{D}}_{\tilde{1},\hat{1}}$--regular subsets:
a section over the base 
${\overV}_1^{(1)}(\mu)\cap{\mathcal{U}}_{{\mathcal{T}}_1,{\mathcal{T}}}$
and its complement.
Each number in the odd rows of the last column in Figure~\ref{n3tac_contr_fig1}
gives the multiplicity with which the number $N(\al)$ of zeros 
of an affine map over the larger  ${\mathcal{D}}_{\tilde{1},\hat{1}}$--regular
set enters into 
${\mathcal{C}}_{\partial\Bbb{P}(L_{\tilde{1}}\oplus L_{\tilde{1}}^*)}
({\mathcal{D}}_{\tilde{1},\hat{1}})$;
each number in the even rows gives such a multiplicity 
for the smaller set.
Lemma~\ref{n3tac_contr_lmm1} computes the contributions
from the first two types of boundary strata of Figure~\ref{n3tac_contr_fig1};
Lemma~\ref{n3tac_contr_lmm2} deals with remaining~one.

\begin{lmm}
\label{n3tac_contr_lmm1}
Suppose ${\mathcal{T}} = (M_1,I;j,\under{d})$ is a bubble type such that
${\mathcal{T}} < {\mathcal{T}}_1$, $\chi_{\mathcal{T}}(\tilde{1},\hat{1}) = 0$,
and \hbox{$|\chi_{\tilde{1}}({\mathcal{T}})| = 1$}.

{\rm(1)}\qua If $|\hat{I}^+| > |\chi_{\tilde{1}}({\mathcal{T}})|$,
$\Bbb{P}(L_{\tilde{1}} \oplus L_{\tilde{1}}^*)|
{\mathcal{S}}_{{\mathcal{T}}|{\mathcal{T}}_1}(\mu)$
is a finite union of ${\mathcal{D}}_{\tilde{1},\hat{1}}$--hollow subsets 
and~thus 
$${\mathcal{C}}_{\Bbb{P}(L_{\tilde{1}}\oplus L_{\tilde{1}}^*)|
{\mathcal{U}}_{{\mathcal{T}}_1,{\mathcal{T}}}}({\mathcal{D}}_{\tilde{1},\hat{1}})=0.$$
{\rm(2)}\qua If $|\hat{I}^+| = |\chi_{\tilde{1}}({\mathcal{T}})|$
and $M_{\tilde{1}}{\mathcal{T}} = \{\hat{1}\}$,
$${\mathcal{C}}_{\Bbb{P}(L_{\tilde{1}} \oplus L_{\tilde{1}}^*)|
{\mathcal{U}}_{{\mathcal{T}}_1,{\mathcal{T}}}}({\mathcal{D}}_{\tilde{1},\hat{1}})=
\blr{20a_{\hat{0}} + 19c_1({\mathcal{L}}_{\tilde{1}}^*),{\overS}_1(\mu)}
-11\big|{\mathcal{S}}_2(\mu)\big|.$$
{\rm(3)}\qua If $|\hat{I}^+| = |\chi_{\tilde{1}}({\mathcal{T}})|$
and $M_{\tilde{1}}{\mathcal{T}} = \{\hat{1},l\}$ for some $l \in [N]$,
$${\mathcal{C}}_{\Bbb{P}(L_{\tilde{1}}\oplus L_{\tilde{1}}^*)|
{\mathcal{U}}_{{\mathcal{T}}_1,{\mathcal{T}}}}({\mathcal{D}}_{\tilde{1},\hat{1}})=
3\big|{\mathcal{S}}_{{\mathcal{T}}_0/l}(\mu)\big|.$$
\end{lmm}

\begin{proof}
(1)\qua Let $h$ be the unique element of $\chi_{\tilde{1}}({\mathcal{T}})$.
By Lemma~\ref{n3p3_str_lmm2}, 
$${\overV}_1^{(1)}(\mu)\cap{\mathcal{U}}_{{\mathcal{T}}_1,{\mathcal{T}}}
\subset{\mathcal{S}}_{{\mathcal{T}}|{\mathcal{T}}_1}(\mu)\equiv
\big\{b \in {\mathcal{U}}_{{\mathcal{T}}|{\mathcal{T}}_1}(\mu) :
{\mathcal{D}}_{{\mathcal{T}},h}^{(1)}b = 0\big\}.$$
With appropriate identifications, ${\mathcal{S}}_{{\mathcal{T}}|{\mathcal{T}}_1}(\mu)$
is the zero set of the section
$$\ev_{{\mathcal{T}}_1,M_0} \oplus {\mathcal{D}}_{{\mathcal{T}},h}^{(1)}$$
of the bundle 
$$\ev_{{\mathcal{T}}_1,M_0}^*{\mathcal{N}}\De_{{\mathcal{T}}_1}(\mu)
\oplus L_h^* \otimes \ev_{\hat{0}}^*T\PPP$$
defined over a neighborhood of ${\mathcal{S}}_{{\mathcal{T}}|{\mathcal{T}}_1}(\mu)$
in~${\mathcal{U}}_{{\mathcal{T}}_1,{\mathcal{T}}}$.
By Lemma~\ref{str_lmm}, this section is transverse to the zero set.
By Proposition~\ref{str_prp} and Lemma~\ref{n3p3_str_lmm2}, 
\begin{gather*}
\ev_{{\mathcal{T}}_1,M_0}\big(\phi_{{\mathcal{T}}_1,{\mathcal{T}}}(b;\ups)\big)
=\ev_{{\mathcal{T}}_1,M_0}(b)+\ve_{-;1}(b;\ups),\\
\{\ev_{\tilde{1}} \times \ev_{\hat{1}}\}
\big(\phi_{{\mathcal{T}}_1,{\mathcal{T}}}(\ups)\big)=
\big(y_{h;\hat{1}} - x_{\hat{1};h}\big)^{-1}\otimes
\big\{{\mathcal{D}}_{{\mathcal{T}},h}^{(1)} + \ve_{-;2}(\ups)\big\}
\otimes\rho_{{\mathcal{T}},\hat{1};h}^{(0;1)}(\ups)
\end{gather*}
for all $(b;\ups) \in \mathcal{FT}_{\de} - Y(\mathcal{FT};\hat{I}^+)$
and some $C^1$--negligible maps
$$\ve_{-;1},\ve_{-;2}\co 
\mathcal{FT}_{\de} - Y(\mathcal{FT};\hat{I}^+)\lra
\ev_{{\mathcal{T}}_1,M_0}^*{\mathcal{N}}\De_{{\mathcal{T}}_2}(\mu),
L_h^* \otimes \ev_{\hat{0}}^*T\PPP.$$
(2)\qua Subtracting the expansion of 
$\{\ev_{\tilde{1}} \times \ev_{\hat{1}}\} \circ \phi_{{\mathcal{T}}_1,{\mathcal{T}}}$
of Lemma~\ref{n3p3_str_lmm2} multiplied by
$$\big(y_{\tilde{1};\hat{1}}(\ups) - x_{\tilde{1};h}(\ups)\big)
\qquad\hbox{and by}\qquad
-\big(y_{\tilde{1};\hat{1}}(\ups) - x_{\tilde{1};h}(\ups)\big)^{-1}$$
from the expansions of 
${\mathcal{D}}_{{\mathcal{T}}_1,\tilde{1}}^{(1)} \circ \phi_{{\mathcal{T}}_1,{\mathcal{T}}}$
and ${\mathcal{D}}_{{\mathcal{T}}_1,\hat{1}}^{(1)} \circ \phi_{{\mathcal{T}}_1,{\mathcal{T}}}$,
respectively, given by Proposition~\ref{str_prp},
we~obtain
\begin{gather*}
{\mathcal{D}}_{\tilde{1},\hat{1}}\phi_{{\mathcal{T}}_1,{\mathcal{T}}}
([\ups_{\tilde{1}},\ups_{\hat{1}}];\ups)=
\big\{\al + \ve(\ups)\}\rho(\ups)\\
\mbox{for all}\ 
([\ups_{\tilde{1}},\ups_{\hat{1}}];\ups) \in \mathcal{FT}_{\de}
~\mbox{with}~
\phi_{{\mathcal{T}}_1,{\mathcal{T}}}(\ups) \in {\mathcal{V}}_1^{(1)},
\end{gather*}
where $\rho$ is a monomials map on $\mathcal{FT}$
with values in a line bundle $\tilde{\mathcal{F}}{\mathcal{T}}$,
$$\al\co \tilde{\mathcal{F}}{\mathcal{T}} \lra
\ga_{L_{\tilde{1}}\oplus L_{\tilde{1}}^*}^* \otimes \ev_{\hat{0}}^*T\PPP$$
is a linear map, and 
$$\ve\co \mathcal{FT} - Y(\mathcal{FT};\hat{I}^+)\lra
\hbox{Hom}(\tilde{\mathcal{F}}{\mathcal{T}},
\ga_{L_{\tilde{1}}\oplus L_{\tilde{1}}^*}^* \otimes \ev_{\hat{0}}^*T\PPP)$$
is a $C^0$--negligible map.
Explicitly,
\begin{gather*}
\rho(\ups) = \prod_{i\in(i_{\mathcal{T}}(h,\hat{1}),h]}        
\ups_i^{\otimes2}
\otimes \prod_{i\in(\tilde{1},i_{\mathcal{T}}(h,\hat{1})]}       \ups_i^*,\\
\begin{split}
\al\big([\ups_{\tilde{1}},\ups_{\hat{1}}],\tilde{\ups}\big)=
-(y_{h;\hat{1}} - x_{\hat{1};h})^{-3}
&\otimes {\mathcal{D}}_{{\mathcal{T}},h}^{(2)}\tilde{\ups}\\
&\otimes
\begin{cases}
(y_{h;\hat{1}} - x_{\hat{1};h})^2 \otimes \ups_{\tilde{1}}
 + \ups_{\hat{1}},&\hbox{if}~i_{\mathcal{T}}(h,\hat{1}) = \tilde{1};\\
\ups_{\hat{1}},&\hbox{if}~i_{\mathcal{T}}(h,\hat{1}) > \tilde{1}.
\end{cases}\end{split}
\end{gather*}
In particular, $\al$ is an injective linear map outside
of a section ${\mathcal{Z}}_{\mathcal{T}}$ of 
$\Bbb{P}(L_{\tilde{1}} \oplus L_{\tilde{1}}^*)$
over ${\mathcal{S}}_{{\mathcal{T}}|{\mathcal{T}}_1}(\mu)$.
Thus, 
$\Bbb{P}(L_{\tilde{1}} \oplus L_{\tilde{1}}^*)|
{\mathcal{S}}_{{\mathcal{T}}|{\mathcal{T}}_1}(\mu) - {\mathcal{Z}}_{\mathcal{T}}$
is ${\mathcal{D}}_{\tilde{1},\hat{1}}$--hollow
unless $\hat{I}^+ = \chi_{\tilde{1}}({\mathcal{T}})$.
If $\hat{I}^+ = \chi_{\tilde{1}}({\mathcal{T}})$, by Proposition~\ref{euler_prp},
the decomposition~\eqref{cart_split}, and 
a rescaling of the linear~map,
\begin{gather*}
{\mathcal{C}}_{\Bbb{P}(L_{\tilde{1}}\oplus L_{\tilde{1}}^*)|
{\mathcal{S}}_{{\mathcal{T}}|{\mathcal{T}}_1}(\mu)-{\mathcal{Z}}_{\mathcal{T}}}
({\mathcal{D}}_{\tilde{1},\hat{1}})
=2N(\al_1),\qquad\hbox{where}\\
\al_1 \in \Ga\big(\Bbb{P}{\mathcal{F}} \times 
{\mathcal{S}}_{{\overT}}(\mu);
\hbox{Hom}(\ga_{\mathcal{F}}^* \otimes L_h^{\otimes2} , 
\ga_{\mathcal{F}}^* \otimes \ev_{\hat{0}}^*T\PPP)\big),\\
{\mathcal{F}} = L_{\tilde{1}} \oplus L_{\tilde{1}}^* \lra 
\ov{\frak M}_{\{\tilde{1},h\}\sqcup M_{\tilde{1}}{\mathcal{T}}}, \qquad
\al_1 = {\mathcal{D}}_{{\overT},h}^{(2)}.
\end{gather*}
If $M_{\tilde{1}}{\mathcal{T}}=\{\hat{1},l\}$ for some $l \in [N]$,
${\mathcal{D}}_{{\overT},h}^{(2)}$ does not vanish on the finite set 
${\mathcal{S}}_{{\overT}}(\mu)$ and thus
\begin{equation}\label{n3tac_contr_lmm1e5}
{\mathcal{C}}_{\Bbb{P}(L_{\tilde{1}}\oplus L_{\tilde{1}}^*)|
{\mathcal{S}}_{{\mathcal{T}}|{\mathcal{T}}_1}(\mu)-{\mathcal{Z}}_{\mathcal{T}}}
({\mathcal{D}}_{\tilde{1},\hat{1}})=2\big\lan 3\la_{\mathcal{F}}^2 - 
3\la_{\mathcal{F}}^2 + \la_{\mathcal{F}}^2,\Bbb{P}{\mathcal{F}}\big\ran
\big|{\mathcal{S}}_{{\overT}}(\mu)\big|=0.
\end{equation}
If $M_{\tilde{1}}{\mathcal{T}}=\{\hat{1}\}$, 
${\overT} = {\mathcal{T}}_0$, and
by Propositions~\ref{zeros_prp} and~\ref{euler_prp}, and
identity~\eqref{zeros_main_e},
$${\mathcal{C}}_{\Bbb{P}(L_{\tilde{1}}\oplus L_{\tilde{1}}^*)|
{\mathcal{S}}_{{\mathcal{T}}|{\mathcal{T}}_1}(\mu)-{\mathcal{Z}}_{\mathcal{T}}}
({\mathcal{D}}_{\tilde{1},\hat{1}})=
2\Big(\big\lan 4a_{\hat{0}} + 2c_1(L_{\tilde{1}}^*),
{\overS}_1(\mu)\big\ran
-{\mathcal{C}}_{\al_1^{-1}(0)}(\al_1^{\perp})\Big).$$
The zero set of $\al_1$ is precisely 
$\Bbb{P}{\mathcal{F}} \times {\mathcal{D}}_{{\mathcal{T}}_0,\tilde{1}}^{(2)-1}(0)$.
{}From the argument in (3) and~(4) of the proof of Lemma~\ref{n3p3_contr_lmm1c},
we obtain
$${\mathcal{C}}_{\al_1^{-1}(0)}(\al_1^{\perp})=
\big\lan \la_{\mathcal{F}},\Bbb{P}{\mathcal{F}}\big\ran
\big(2|{\mathcal{S}}_{1;1}(\mu)|+|{\mathcal{S}}_2(\mu)|\big)
=2\big|{\mathcal{S}}_{1;1}(\mu)\big|+\big|{\mathcal{S}}_2(\mu)\big|.$$
Putting the last two equations together, we conclude that
if $M_{\tilde{1}}{\mathcal{T}} = \{\hat{1}\}$,
\begin{equation}\label{n3tac_contr_lmm1e7}
{\mathcal{C}}_{\Bbb{P}(L_{\tilde{1}}\oplus L_{\tilde{1}}^*)|
{\mathcal{S}}_{{\mathcal{T}}|{\mathcal{T}}_1}(\mu)-{\mathcal{Z}}_{\mathcal{T}}}
({\mathcal{D}}_{\tilde{1},\hat{1}})=
\big\lan 8a_{\hat{0}} + 4c_1({\mathcal{L}}_{\tilde{1}}^*),
{\overS}_1(\mu)\big\ran-2|{\mathcal{S}}_2(\mu)|.
\end{equation}
(3)\qua In order to compute the contribution from the space~${\mathcal{Z}}_{\mathcal{T}}$,
we keep the two leading terms of the expression 
for~${\mathcal{D}}_{\tilde{1},\hat{1}}$ obtained as in~(2).
We can model a neighborhood of ${\mathcal{Z}}_{\mathcal{T}}$ in
$\Bbb{P}(L_{\tilde{1}} \oplus L_{\tilde{1}}^*)$ by the~map
$$L_{\tilde{1}}^*\otimes L_{\tilde{1}}^*\lra
\Bbb{P}(L_{\tilde{1}} \oplus L_{\tilde{1}}^*),\qquad
\big([\ups_{\tilde{1}},\ups_{\hat{1}}],u\big)\lra 
\big[\ups_{\tilde{1}},\ups_{\hat{1}} + u(\ups_{\tilde{1}})\big].$$
If $i_{\mathcal{T}}(\hat{1},h) = \tilde{1}$, near ${\mathcal{Z}}_{\mathcal{T}}$,
\begin{equation*}\begin{split}
{\mathcal{D}}_{\tilde{1},\hat{1}}\phi_{{\mathcal{T}}_1,{\mathcal{T}}}
([\ups_{\tilde{1}},\ups_{\hat{1}}];u,\ups)=&
-(y_{\tilde{1};\hat{1}} - x_{\tilde{1};h})^{-3} \otimes 
\big({\mathcal{D}}_{{\mathcal{T}},h}^{(2)}+\ve_2(u,\ups)\big)
\rho_{{\mathcal{T}},\tilde{1};h}^{(0;2)}(\ups) \otimes u\\
&-(y_{\tilde{1};\hat{1}} - x_{\tilde{1};h})^{-2}
 \otimes \ups_{\tilde{1}} \otimes 
\big({\mathcal{D}}_{{\mathcal{T}},h}^{(3)}+\ve_3(u,\ups)\big)
\rho_{{\mathcal{T}},\tilde{1};h}^{(0;3)}(\ups).
\end{split}\end{equation*}
Note that by Lemma~\ref{str_lmm}, the images of
${\mathcal{D}}_{{\mathcal{T}},h}^{(2)}$ and ${\mathcal{D}}_{{\mathcal{T}},h}^{(3)}$
are distinct over ${\mathcal{S}}_{{\mathcal{T}}|{\mathcal{T}}_1}(\mu)$.
If $i_{\mathcal{T}}(\hat{1},h) > \tilde{1}$, we similarly find 
${\mathcal{Z}}_{\mathcal{T}}$ is ${\mathcal{D}}_{\tilde{1},\hat{1}}$--hollow unless
$\hat{I}^+ = \chi_{\tilde{1}}({\mathcal{T}})$.
If $\hat{I}^+ = \chi_{\tilde{1}}({\mathcal{T}})$,
by Proposition~\ref{euler_prp},
the decomposition~\eqref{cart_split}, and 
a rescaling of the linear~map,
\begin{gather*}
{\mathcal{C}}_{{\mathcal{Z}}_{\mathcal{T}}}({\mathcal{D}}_{\tilde{1},\hat{1}})
=3N(\al_1),\qquad\hbox{where}\\
\al_1 \in \Ga\big({\overZ}_{\mathcal{T}};
\hbox{Hom}(\ga_{\mathcal{F}}^* \otimes L_h^{\otimes2} \oplus 
\ga_{\mathcal{F}}^* \otimes L_h^{\otimes3};
\ga_{\mathcal{F}}^* \otimes \ev_{\hat{0}}^*T\PPP)\big),\\
{\overZ}_{\mathcal{T}} \subset 
\Bbb{P}{\mathcal{F}} \times 
{\mathcal{S}}_{{\overT}}(\mu),\quad
{\mathcal{F}} = L_{\tilde{1}} \oplus L_{\tilde{1}}^* \lra 
\ov{\frak M}_{\{\tilde{1},h\}\sqcup M_{\tilde{1}}{\mathcal{T}}},\\
\al_1(\ups_2,\ups_3)={\mathcal{D}}_{{\overT},h}^{(2)}\ups_2
 + {\mathcal{D}}_{{\overT},h}^{(3)}\ups_3.
\end{gather*}
If $M_{\tilde{1}}{\mathcal{T}}=\{\hat{1},l\}$ for some $l \in [N]$,
$\al_1$ has full rank over ${\overZ}_{\mathcal{T}}$ and thus
\begin{equation}\label{n3tac_contr_lmm1e9}
\sum_{M_{\tilde{1}}{\mathcal{T}}=\{\hat{1},l\}}     
{\mathcal{C}}_{{\mathcal{Z}}_{\mathcal{T}}}
({\mathcal{D}}_{\tilde{1},\hat{1}})=
3\sum_{l\in[N]}
\big\lan 3\la_{\mathcal{F}} - 2\la_{\mathcal{F}},{\overZ}_{\mathcal{T}}\big\ran
\big|{\mathcal{S}}_{{\overT}}(\mu)\big|=
3\big|{\mathcal{S}}_{1;1}(\mu)\big|.
\end{equation}
If $M_{\tilde{1}}{\mathcal{T}}=\{\hat{1}\}$, 
${\overZ}_{\mathcal{T}} \approx {\overS}_1(\mu)$, and
\begin{equation}\label{n3tac_contr_lmm1e11}
{\mathcal{C}}_{{\mathcal{Z}}_{\mathcal{T}}}
({\mathcal{D}}_{\tilde{1},\hat{1}})=3N(\al_1)=
\big\lan 12a_{\hat{0}} + 15c_1({\mathcal{L}}_{\tilde{1}}^*),
{\overS}_1(\mu)\big\ran -9\big|{\mathcal{S}}_2(\mu)\big|;
\end{equation}
see \cite[Lemma~5.12]{Z1}.
The claim follows from equations
\eqref{n3tac_contr_lmm1e5}--\eqref{n3tac_contr_lmm1e11}.
\end{proof}

\begin{lmm}
\label{n3tac_contr_lmm2}
Suppose ${\mathcal{T}} = (M_1,I;j,\under{d})$ is a bubble type such that
${\mathcal{T}} < {\mathcal{T}}_1$, $\chi_{\mathcal{T}}(\tilde{1},\hat{1}) = 0$,
and \hbox{$|\chi_{\tilde{1}}({\mathcal{T}})| = 2$}.

{\rm(1)}\qua If $|\hat{I}^+| > |\chi_{\tilde{1}}({\mathcal{T}})|$,
$\Bbb{P}(L_{\tilde{1}} \oplus L_{\tilde{1}}^*)|
{\mathcal{S}}_{{\mathcal{T}}|{\mathcal{T}}_1}(\mu)$
is a finite union of ${\mathcal{D}}_{\tilde{1},\hat{1}}$--hollow subsets 
and~thus 
$${\mathcal{C}}_{\Bbb{P}(L_{\tilde{1}}\oplus L_{\tilde{1}}^*)|
{\mathcal{U}}_{{\mathcal{T}}_1,{\mathcal{T}}}}({\mathcal{D}}_{\tilde{1},\hat{1}})
=0.$$
{\rm(2)}\qua If $|\hat{I}^+| = |\chi_{\tilde{1}}({\mathcal{T}})|$,
$M_{\tilde{1}}{\mathcal{T}} = \{\hat{1}\}$ if 
${\mathcal{S}}_{{\mathcal{T}}|{\mathcal{T}}_1}(\mu) \neq \eset$, and
$$\sum_{{\mathcal{T}}<{\mathcal{T}}_1}
{\mathcal{C}}_{\Bbb{P}(L_{\tilde{1}}\oplus L_{\tilde{1}}^*)|
{\mathcal{U}}_{{\mathcal{T}}_1,{\mathcal{T}}}}({\mathcal{D}}_{\tilde{1},\hat{1}})=
2\big|{\mathcal{S}}_2(\mu)\big|,$$
where the sum is taken over all equivalence classes of 
bubble types of the above form.
\end{lmm}

\begin{proof}
The proof is a mixture of the proof of 
Lemma~\ref{n3tac_contr_lmm1}
with (2) of the proof of Lemma~\ref{n3p3_contr_lmm3};
thus, we omit~it.
\end{proof}

\begin{crl}
\label{n3tac_contr_crl1}
The total contribution from the boundary strata 
${\mathcal{U}}_{{\mathcal{T}}_1,{\mathcal{T}}}$ such that 
$\chi_{\mathcal{T}}(\tilde{1},\hat{1}) = 0$ to the number 
${\mathcal{C}}_{\partial\Bbb{P}(L_{\tilde{1}}\oplus L_{\tilde{1}}^*)}
({\mathcal{D}}_{\tilde{1},\hat{1}})$ is given~by
$$\sum_{\chi_{\mathcal{T}}(\tilde{1},\hat{1})=0}     
{\mathcal{C}}_{\Bbb{P}(L_{\tilde{1}}\oplus L_{\tilde{1}}^*)|
{\mathcal{U}}_{{\mathcal{T}}_1,{\mathcal{T}}}}({\mathcal{D}}_{\tilde{1},\hat{1}})=
\big\lan 20a_{\hat{0}} + 19\eta_{\hat{0},1},
{\overS}_1(\mu)\big\ran-9\big|{\mathcal{S}}_2(\mu)\big|
+3\big|{\mathcal{S}}_{1;1}(\mu)\big|.$$
\end{crl}

\begin{proof}
This Corollary follows immediately from 
Lemmas~\ref{n3tac_contr_lmm1} and~\ref{n3tac_contr_lmm2}.
\end{proof}

\subsection[Contributions from the spaces UT|T1 with kT(1,1)>0]{Contributions from the spaces ${\mathcal{U}}_{{\mathcal{T}}|{\mathcal{T}}_1}$
with $\chi_{\mathcal{T}}(\tilde{1},\hat{1}) > 0$}
\label{n3tac_corr_subs2}

In this subsection, we prove Corollary~\ref{n3tac_contr_crl2},
which gives the total contribution to the number
${\mathcal{C}}_{\partial\Bbb{P}(L_{\tilde{1}}\oplus L_{\tilde{1}}^*)}
({\mathcal{D}}_{\tilde{1},\hat{1}})$
from the spaces ${\mathcal{U}}_{{\mathcal{T}}_1,{\mathcal{T}}}$,
where 
$${\mathcal{T}} \equiv (M_1,I;j,\under{d}) < {\mathcal{T}}_1$$ 
is a bubble type such that \hbox{$\chi_{\mathcal{T}}(\tilde{1},\hat{1}) > 0$}.
Note that by the last paragraph of Subsection~\ref{n3tac_sum_subs}
it is sufficient to consider bubble types ${\mathcal{T}}$ such 
\hbox{that $j_{\hat{1}} > \tilde{1}$}.

\begin{figure}[ht!]\small
\begin{pspicture}(-1,-4)(10,2.7)
\psset{unit=.36cm}
\pscircle[fillstyle=solid,fillcolor=gray](1,-2){1}
\pscircle*(1,-3){.22}\rput(1,-3.8){$\tilde{1}$}
\pscircle*(1.71,-1.29){.2}\rput(2.2,-1){$\hat{1}$}
\pnode(1,-3){A0}\pnode(1.71,-1.29){B0}
\ncarc[nodesep=.35,arcangleA=-90,arcangleB=-75,ncurv=1.2]{<->}{A0}{B0}
% 2nd column starts here
\psline[linewidth=.12]{->}(3,-1.5)(6.5,2)
\psline[linewidth=.12]{->}(3,-2)(6.5,-2)
\psline[linewidth=.12]{->}(3,-2.5)(6.5,-6)
% 3rd column starts here
\pscircle(10.5,2.5){1}
\pscircle[fillstyle=solid,fillcolor=gray](9.08,3.92){1}
\pscircle*(10.5,1.5){.22}\rput(10.5,.7){$\tilde{1}$}
\pscircle*(9.79,3.21){.19}
\pscircle*(11.21,3.21){.2}\rput(11.7,3.5){$l$}
\pscircle*(8.37,4.63){.2}\rput(7.88,4.92){$\hat{1}$}
\pnode(10.5,1.5){A1}\pnode(8.37,4.63){B1}
\ncarc[nodesep=.35,arcangleA=90,arcangleB=90,ncurv=1.2]{<->}{A1}{B1}
% 2nd row starts here
\pscircle(10.5,-2.5){1}
\pscircle[fillstyle=solid,fillcolor=gray](9.08,-1.08){1}
\pscircle*(10.5,-3.5){.22}\rput(10.5,-4.3){$\tilde{1}$}
\pscircle*(9.79,-1.79){.19}
\pscircle[fillstyle=solid,fillcolor=gray](11.92,-1.08){1}
\pscircle*(11.21,-1.79){.19}
\pscircle*(8.37,-.37){.2}\rput(7.88,-.08){$\hat{1}$}
\pnode(10.5,-3.5){A2}\pnode(8.37,-.37){B2}
\ncarc[nodesep=.35,arcangleA=90,arcangleB=90,ncurv=1.2]{<->}{A2}{B2}
% 3rd row starts here
\pscircle[fillstyle=solid,fillcolor=gray](10.5,-7.5){1}
\pscircle[fillstyle=solid,fillcolor=gray](9.08,-6.08){1}
\pscircle*(10.5,-8.5){.22}\rput(10.5,-9.3){$\tilde{1}$}
\pscircle*(9.79,-6.79){.19}
\pscircle*(8.37,-5.37){.2}\rput(7.88,-5.08){$\hat{1}$}
\pnode(10.5,-8.5){A3}\pnode(8.37,-5.37){B3}
\ncarc[nodesep=.35,arcangleA=90,arcangleB=90,ncurv=1.2]{<->}{A3}{B3}
% 4th column starts here
\rput(14,3){$\approx$}\rput(14,-2){$\approx$}\rput(14,-7){$\approx$}
% 5th column starts here
\rput(16.7,2.9){${\overV}_{1;1}^{(1)}(\mu)$}
\rput(16.7,-2){${\mathcal{V}}_2^{(1)}(\mu)$}
\rput(16.7,-7.1){${\overV}_2^{(1,1)}(\mu)$}
%6th column starts here
\psline{->}(19,3.25)(24,4)\rput{8}(21.5,4.2){$\times(-1)$}
\psline{->}(19,2.75)(24,2)\rput{-8}(22,2.7){$\times1$}
\psline{->}(19,-1.75)(24,-.5)
\psline{->}(19,-2)(24,-2)\rput(22,-1.6){$\times1$}
\psline{->}(19,-2.25)(24,-3.5)
\psline{->}(19,-6.75)(24,-6)\rput{8}(21.5,-5.8){$\times(-1)$}
\psline{->}(19,-7.25)(24,-8)
%7th column starts here
\rput(27,4.25){\small Lemma~\ref{n3tac_contr_lmm3b}}
\rput(27,2){\small Lemma~\ref{n3tac_contr_lmm3c}}
\rput(27,-.35){${\mathcal{D}}_{\tilde{1},\hat{1}}$--\small{neutral}}
\rput(27,-1.9){$|{\mathcal{V}}_2^{(1)}(\mu)|$}
\rput(27,-3.85){${\mathcal{D}}_{\tilde{1},\hat{1}}$--\small{hollow}}
\rput(27,-5.9){\small Lemma~\ref{n3tac_contr_lmm4b}}
\rput(27,-8.1){${\mathcal{D}}_{\tilde{1},\hat{1}}$--\small{hollow}}
\end{pspicture}
\caption{An outline of Subsection~\ref{n3tac_corr_subs2}}
\label{n3tac_contr_fig2}
\end{figure}
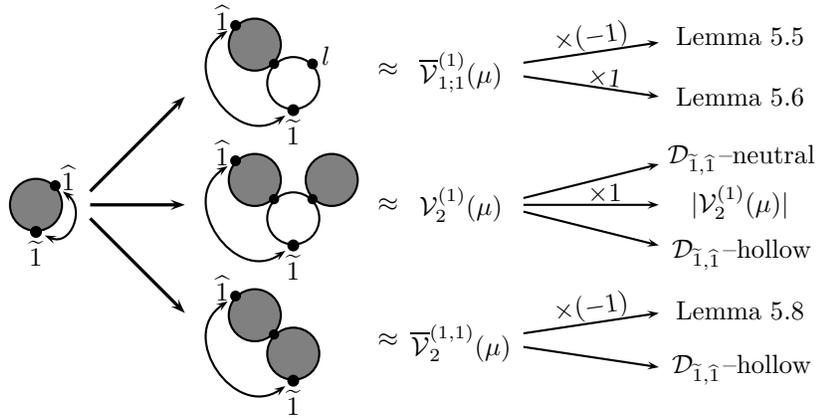

Figure~\ref{n3tac_contr_fig2} shows the three types of boundary strata 
${\overV}_1^{(1)}(\mu)\cap{\mathcal{U}}_{{\mathcal{T}}_1,{\mathcal{T}}}$ such that
$$\Bbb{P}(L_{\tilde{1}} \oplus L_{\tilde{1}}^*)|
{\overV}_1^{(1)}(\mu)\cap{\mathcal{U}}_{{\mathcal{T}}_1,{\mathcal{T}}}$$
is not contained in a finite union of 
${\mathcal{D}}_{\tilde{1},\hat{1}}$--hollow sets.
As in Subsection~\ref{n3tac_corr_subs1}, we have to split each~space
$$\Bbb{P}(L_{\tilde{1}} \oplus L_{\tilde{1}}^*)|
{\overV}_1^{(1)}(\mu)\cap{\mathcal{U}}_{{\mathcal{T}}_1,{\mathcal{T}}}$$
into two or three subspaces, as indicated on the right-hand side
of Figure~\ref{n3tac_contr_fig2}.
Lemma~\ref{n3tac_contr_lmm3} computes the contributions
from the first two types of boundary strata of Figure~\ref{n3tac_contr_fig2};
Lemma~\ref{n3tac_contr_lmm4} deals with remaining~one.

\begin{lmm}
\label{n3tac_contr_lmm3}
Suppose ${\mathcal{T}} = (M_1,I;j,\under{d})$ is a bubble type such that
${\mathcal{T}} < {\mathcal{T}}_1$, $\chi_{\mathcal{T}}(\tilde{1},\hat{1}) > 0$,
and \hbox{$d_{\tilde{1}} = 0$}.
{\rm(1)}\qua If $|\hat{I}^+| > |\chi_{\tilde{1}}({\mathcal{T}})|$,
$\Bbb{P}(L_{\tilde{1}} \oplus L_{\tilde{1}}^*)|
{\mathcal{U}}_{{\mathcal{T}}|{\mathcal{T}}_1}^{(1)}(\mu)$ is a finite union of
${\mathcal{D}}_{\tilde{1},\hat{1}}$--hollow subspaces and thus
$${\mathcal{C}}_{\Bbb{P}(L_{\tilde{1}}\oplus L_{\tilde{1}}^*)|
{\mathcal{U}}_{{\mathcal{T}}_1,{\mathcal{T}}}}({\mathcal{D}}_{\tilde{1},\hat{1}}) = 0.$$
{\rm(2)}\qua The total contribution from the boundary strata 
${\mathcal{U}}_{{\mathcal{T}}_1,{\mathcal{T}}}$
such that $|\hat{I}^+| = |\chi_{\tilde{1}}({\mathcal{T}})| = 1$ 
is given~by
$$\sum_{|\hat{I}^+|=|\chi_{\tilde{1}}({\mathcal{T}})|=1}      
{\mathcal{C}}_{\Bbb{P}(L_{\tilde{1}}\oplus L_{\tilde{1}}^*)|
{\mathcal{U}}_{{\mathcal{T}}_1,{\mathcal{T}}}}({\mathcal{D}}_{\tilde{1},\hat{1}})
=\big\lan c_1(L_{\tilde{1}}^*),{\overV}_{1;1}^{(1)}(\mu)\big\ran
-3\big|{\mathcal{S}}_{1;1}(\mu)\big|.$$
{\rm(3)}\qua The total contribution from the boundary strata 
${\mathcal{U}}_{{\mathcal{T}}_1,{\mathcal{T}}}$
such that $|\hat{I}^+| = |\chi_{\tilde{1}}({\mathcal{T}})| = 2$ 
is given~by
$$\sum_{|\hat{I}^+|=|\chi_{\tilde{1}}({\mathcal{T}})|=2}      
{\mathcal{C}}_{\Bbb{P}(L_{\tilde{1}}\oplus L_{\tilde{1}}^*)|
{\mathcal{U}}_{{\mathcal{T}}_1,{\mathcal{T}}}}({\mathcal{D}}_{\tilde{1},\hat{1}})
=\big|{\mathcal{V}}_2^{(1)}(\mu)\big|.$$
\end{lmm}

\begin{proof}
(1)\qua By Lemma~\ref{n3p3_str_lmm1},
$${\overU}_{{\mathcal{T}}_1}^{(1)}(\mu)\cap{\mathcal{U}}_{{\mathcal{T}}_1,{\mathcal{T}}}
\subset{\mathcal{U}}_{{\mathcal{T}}|{\mathcal{T}}_1}^{(1)}(\mu).$$
With appropriate identifications, 
${\mathcal{U}}_{{\mathcal{T}}|{\mathcal{T}}_1}^{(1)}(\mu)$
is the zero set of the section
$$\ev_{{\mathcal{T}}_1,M_0}\oplus(\ev_{\hat{1}} - \ev_{\tilde{1}})$$
of the bundle 
$$\ev_{{\mathcal{T}}_1,M_0}^*{\mathcal{N}}\De_{{\mathcal{T}}_1}(\mu)\oplus
\ev_{\tilde{1}}^*T\PPP$$ 
over an open neighborhood of 
${\mathcal{U}}_{{\mathcal{T}}|{\mathcal{T}}_1}^{(1)}(\mu)$
in~${\mathcal{U}}_{{\mathcal{T}}_1,{\mathcal{T}}}$.
By Lemma~\ref{str_lmm}, this section is transversal to the zero set.
By Proposition~\ref{str_prp}, there exists 
a $C^1$--negligible map 
$$\ve_-\co \mathcal{FT}_{\de} - Y(\mathcal{FT};\hat{I}^+)\lra
\ev_{{\mathcal{T}}_1,M_0}^*{\mathcal{N}}\De_{{\mathcal{T}}_1}(\mu)
 \oplus \ev_{\hat{0}}^*T\PPP$$ 
such that
$$\big\{\ev_{{\mathcal{T}}_1,M_0} \times 
\ev_{\tilde{1}} \times \ev_{\hat{1}}\big\}
\big(\phi_{{\mathcal{T}}_1,{\mathcal{T}}}(b;\ups)\big)=
\big\{\ev_{{\mathcal{T}}_1,M_0} \times 
\ev_{\tilde{1}} \times \ev_{\hat{1}}\big\}(b)+\ve_-(b;\ups)$$
for all $(b;\ups) \in \mathcal{FT}_{\de} - Y(\mathcal{FT};\hat{I}^+)$.
On the other hand, by Proposition~\ref{str_prp},
$${\mathcal{D}}_{\tilde{1},\hat{1}}\phi_{{\mathcal{T}}_1,{\mathcal{T}}}
([\ups_{\tilde{1}},\ups_{\hat{1}}];\ups)=
\big\{\al + \ve(\ups)\big\}\rho(\ups)
\quad\mbox{for all}\ \ups \in \mathcal{FT}_{\de} - 
Y(\mathcal{FT};\hat{I}^+).$$
In this equation, $\rho$ is the monomials map on $\mathcal{FT}$
defined~by
$$\rho_h(\ups)=
\begin{cases}
\rho_{{\mathcal{T}},\hat{1};\hat{1}}^{(1)}(\ups),&
\hbox{if}~h \in \chi_{j_{\hat{1}}}({\mathcal{T}}) - \chi_{\hat{1}}({\mathcal{T}});\\
\rho_{{\mathcal{T}},\hat{1};h}^{(1;1)}(\ups),&
\hbox{if}~h \in \chi_{\hat{1}}({\mathcal{T}});\\
\rho_{{\mathcal{T}},\tilde{1};h}^{(1;1)}(\ups),&
\hbox{if}~h \in \chi_{\tilde{1}}({\mathcal{T}}),~h \not\le j_{\hat{1}};
\end{cases}$$
with values in the bundle
$\tilde{\mathcal{F}}{\mathcal{T}}=\bigoplus_{h\in\tilde{I}}\tilde{\mathcal{F}}_h{\mathcal{T}}$,
where
\begin{gather*}
\tilde{I}=\chi_{j_{\hat{1}}}({\mathcal{T}})\cup
\big\{h \in \chi_{\tilde{1}}({\mathcal{T}}) :h \not\le j_{\hat{1}}\big\},\\
\tilde{\mathcal{F}}_h{\mathcal{T}}=
\begin{cases}
\tilde{\mathcal{F}}_{{\mathcal{T}},\hat{1};\hat{1}}^{(1)}{\mathcal{T}},&
\hbox{if}~h \in \chi_{j_{\hat{1}}}({\mathcal{T}}) - \chi_{\hat{1}}({\mathcal{T}});\\
\tilde{\mathcal{F}}_{{\mathcal{T}},\hat{1};h}^{(1;1)}{\mathcal{T}},&
\hbox{if}~h \in \chi_{\hat{1}}({\mathcal{T}});\\
\tilde{\mathcal{F}}_{{\mathcal{T}},\tilde{1};h}^{(1;1)}{\mathcal{T}},&
\hbox{if}~h \in \chi_{\tilde{1}}({\mathcal{T}}),~h \not\le j_{\hat{1}}.
\end{cases}\end{gather*}
The linear map 
$\al\co \tilde{\mathcal{F}}{\mathcal{T}}\lra
\ga_{L_{\tilde{1}}\oplus L_{\tilde{1}}^*}^* \otimes \ev_{\hat{0}}^*T\PPP$
is given~by
$$\al\big([\ups_{\tilde{1}},\ups_{\hat{1}}],\tilde{\ups}_h\big)=
\begin{cases}
\ups_{\hat{1}} \otimes 
{\mathcal{D}}_{{\mathcal{T}},j_{\hat{1}}^*({\mathcal{T}})}^{(1)}\tilde{\ups}_h,&
\hbox{if}~h \in \chi_{j_{\hat{1}}}({\mathcal{T}}) - \chi_{\hat{1}}({\mathcal{T}});\\
-\ups_{\hat{1}} \otimes (y_{h;\hat{1}} - x_{\hat{1};h})^{-2}
\otimes{\mathcal{D}}_{{\mathcal{T}},h}^{(1)}\tilde{\ups}_h,&
\hbox{if}~h \in \chi_{\hat{1}}({\mathcal{T}});\\
\ups_{\tilde{1}} \otimes {\mathcal{D}}_{{\mathcal{T}},h}^{(1)}\tilde{\ups}_h,&
\hbox{if}~h \in \chi_{\tilde{1}}({\mathcal{T}}),
~h \not\le j_{\hat{1}}.
\end{cases}$$
In particular, by Lemma~\ref{str_lmm}, $\al$ has full rank over 
${\mathcal{U}}_{{\mathcal{T}}|{\mathcal{T}}_1}^{(1)}(\mu)$ outside of the set
$${\mathcal{Z}}_{\mathcal{T}}\equiv
\begin{cases}
\Bbb{P}L_{\tilde{1}},&\hbox{if}~h \le j_{\hat{1}}~\mbox{for all}\  
h \in \chi_{\tilde{1}}({\mathcal{T}});\\
\Bbb{P}L_{\tilde{1}}\cup\Bbb{P}L_{\tilde{1}}^*, &\hbox{otherwise}.
\end{cases}$$
As usual,
$$\ve\co \mathcal{FT} - Y(\mathcal{FT};\hat{I}^+) \lra
\hbox{Hom}(\tilde{\mathcal{F}}{\mathcal{T}},
\ga_{L_{\tilde{1}}\oplus L_{\tilde{1}}^*}^* \otimes \ev_{\hat{0}}^*T\PPP)$$
is a $C^0$--negligible map.
Thus, 
$$\Bbb{P}(L_{\tilde{1}} \oplus L_{\tilde{1}}^*)|
{\mathcal{U}}_{{\mathcal{T}}|{\mathcal{T}}_1}^{(1)}(\mu) - {\mathcal{Z}}_{\mathcal{T}}$$
is ${\mathcal{D}}_{\tilde{1},\hat{1}}$--hollow unless
$\hat{I}^+ = \tilde{I}$.
Since $\chi_{\mathcal{T}}(\tilde{1},\hat{1}) > 0$ and $d_{\tilde{1}} = 0$, 
if $\hat{I}^+ = \tilde{I}$ and
${\mathcal{U}}_{{\mathcal{T}}|{\mathcal{T}}_1}(\mu) \neq \eset$,
either ${\mathcal{T}} = {\mathcal{T}}_1(l)$ for some $l \in [N]$
or $|\hat{I}^+| = |\chi_{\tilde{1}}({\mathcal{T}})| = 2$
and $M_{\tilde{1}}{\mathcal{T}} = \eset$.
In the second case, the matrix corresponding to the monomials map~$\rho$
is 
$$\left(\begin{array}{cc}-1& 0\\ 0& 1\end{array}\right).$$
Thus, $\rho$ is neutral, and
\begin{equation}\label{n3tac_contr_lmm3e5}
\sum_{l\in[N]}
{\mathcal{C}}_{\Bbb{P}(L_{\tilde{1}}\oplus L_{\tilde{1}}^*)|
{\mathcal{U}}_{{\mathcal{T}}|{\mathcal{T}}_1}^{(1)}(\mu)-{\mathcal{Z}}_{\mathcal{T}}}
({\mathcal{D}}_{\tilde{1},\hat{1}}) =0.
\end{equation}
In the first case, the degree of $\rho$ is~$-1$.
Thus, by Proposition~\ref{euler_prp}, a rescaling of the linear map, and
the decomposition~\eqref{cart_split}, we obtain
\begin{gather*}
\sum_{|\hat{I}^+|=|\chi_{\tilde{1}}({\mathcal{T}})|=1}   
{\mathcal{C}}_{\Bbb{P}(L_{\tilde{1}}\oplus L_{\tilde{1}}^*)|
{\mathcal{U}}_{{\mathcal{T}}_1(l)|{\mathcal{T}}_1}^{(1)}(\mu)-{\mathcal{Z}}_{\mathcal{T}}}
({\mathcal{D}}_{\tilde{1},\hat{1}})=-N(\al_1),\qquad\hbox{where}\\
\al_1 \in \Ga\big(\Bbb{P}^1 \times {\overV}_{1;1}^{(1)}(\mu);
\hbox{Hom}(\ga^* \otimes L_{\tilde{1}}^*,
\ga^* \otimes \ev_{\hat{0}}^*T\PPP)\big),\quad
\al_1\big|{\mathcal{U}}_{{\mathcal{T}}_1/l}^{(1)}(\mu)
 = {\mathcal{D}}_{{\mathcal{T}}_1/l,\hat{1}}^{(1)}.
\end{gather*}
Thus, by Lemma~\ref{n3tac_contr_lmm3b}, 
\begin{equation}\label{n3tac_contr_lmm3e7}\begin{split}
\sum_{l\in[N]}
{\mathcal{C}}_{\Bbb{P}(L_{\tilde{1}}\oplus L_{\tilde{1}}^*)|
{\mathcal{U}}_{{\mathcal{T}}_1(l)|{\mathcal{T}}_1}^{(1)}(\mu)-{\mathcal{Z}}_{\mathcal{T}}}
({\mathcal{D}}_{\tilde{1},\hat{1}})
=&-\big\lan 4a_{\hat{0}}-c_1(L_{\tilde{1}}^*),
{\overV}_{1;1}^{(1)}(\mu)\big\ran\\
&\qquad +2\big|{\mathcal{S}}_{1;1}(\mu)\big|
-\big|{\mathcal{V}}_{2;1}^{(1,1)}(\mu)\big|.
\end{split}\end{equation}
The space  ${\mathcal{V}}_{2;1}^{(1,1)}(\mu)$ is the disjoint union of sets
$${\mathcal{U}}_{\mathcal{T}}^{(1)}(\mu)\equiv
\big\{b \in {\mathcal{U}}_{\mathcal{T}}(\mu) :
\ev_{\hat{1}}(b) = \ev_{\hat{2}}(b)\big\},$$
taken over all bubble types
${\mathcal{T}} = (M_2 - \{l\},I_2(l);j,\under{d})$,
where 
$$I_2(l)=\{\hat{0} = l\}\sqcup I_2^+, \quad 
j_{\hat{1}}=\tilde{1}, \quad j_{\hat{2}}=\tilde{2}, \quad
d_{\tilde{1}},d_{\tilde{2}}>0,  \quad\hbox{and}\quad
d_{\tilde{1}} + d_{\tilde{2}}=d.$$
The image of every element of ${\mathcal{V}}_{2;1}^{(1,1)}(\mu)$
has two components arranged in  a circle, as in the fifth
picture of Figure~\ref{images_fig}, with one of the two nodes
lying on one of the constraints $\mu_1,\ldots,\mu_N$.

(3)\qua We next consider the contribution from the space 
$\Bbb{P}L_{\tilde{1}}|{\mathcal{U}}_{{\mathcal{T}}|{\mathcal{T}}_1}^{(1)}(\mu)$.
We model a neighborhood of  $\Bbb{P}L_{\tilde{1}}$ in 
$\Bbb{P}(L_{\tilde{1}} \oplus L_{\tilde{1}}^*)$ by the~map
$$L_{\tilde{1}}^*\otimes L_{\tilde{1}}^*\lra 
\Bbb{P}(L_{\tilde{1}} \oplus L_{\tilde{1}}^*),\qquad
\big([\ups_{\tilde{1}},\ups_{\hat{1}}],u\big)
\lra \big[\ups_{\tilde{1}},u(\ups_{\tilde{1}})\big].$$
In this case, with notation as in (2) above,
$\tilde{I} = \chi_{j_{\hat{1}}}({\mathcal{T}}) \sqcup
 \chi_{\tilde{1}}({\mathcal{T}})$, 
$$\big(\rho_h(u,\ups),\al(\tilde{\ups}_h)\big)=
\begin{cases} 
\big(u \otimes \rho_{{\mathcal{T}},\hat{1};\hat{1}}^{(1)}(\ups),
{\mathcal{D}}_{{\mathcal{T}},j_{\hat{1}}^*({\mathcal{T}})}^{(1)}\tilde{\ups}_h\big),&
\!\!\!\!\hbox{if}~h \in \chi_{j_{\hat{1}}}({\mathcal{T}}){-}\chi_{\hat{1}}({\mathcal{T}});\\
\begin{aligned}
&\big(u \otimes \rho_{{\mathcal{T}},\hat{1};h}^{(1;1)}(\ups),\\
&\qquad-(y_{h;\hat{1}} - x_{\hat{1};h})^{-2} 
\otimes {\mathcal{D}}_{{\mathcal{T}},h}^{(1)}\tilde{\ups}_h\big),
\end{aligned}&
\!\!\!\!\hbox{if}~h \in \chi_{\hat{1}}({\mathcal{T}});\\
\big(\rho_{{\mathcal{T}},\tilde{1};h}^{(1;1)}(\ups),
{\mathcal{D}}_{{\mathcal{T}},h}^{(1)}\tilde{\ups}_h\big),&
\!\!\!\!\hbox{if}~h \in \chi_{\tilde{1}}({\mathcal{T}}).
\end{cases}$$
As before $\al$ has full rank on 
$\Bbb{P}L_{\tilde{1}}|{\mathcal{U}}_{{\mathcal{T}}|{\mathcal{T}}_1}^{(1)}(\mu)$.
Thus, $\Bbb{P}L_{\tilde{1}}|{\mathcal{U}}_{{\mathcal{T}}|{\mathcal{T}}_1}^{(1)}(\mu)$
is ${\mathcal{D}}_{\tilde{1},\hat{1}}$--hollow unless
either ${\mathcal{T}} = {\mathcal{T}}_1(l)$ for some $l \in [N]$
or $|\hat{I}^+| = |\chi_{\tilde{1}}({\mathcal{T}})| = 2$
and $M_{\tilde{1}}{\mathcal{T}} = \eset$.
In both cases, the degree of $\rho$ is one.
In the second case, $\al$ is an isomorphism on every fiber, and thus
\begin{equation}\label{n3tac_contr_lmm3e9}
\sum_{|\hat{I}^+|=|\chi_{\tilde{1}}({\mathcal{T}})|=2}
{\mathcal{C}}_{\Bbb{P}L_{\tilde{1}}|{\mathcal{U}}_{{\mathcal{T}}|{\mathcal{T}}_1}^{(1)}(\mu)}
({\mathcal{D}}_{\tilde{1},\hat{1}})
=\big|{\mathcal{V}}_2^{(1)}(\mu)\big|.
\end{equation}
In the first case, via the decomposition~\eqref{cart_split}, we obtain
\begin{gather*}
\sum_{l\in[N]}{\mathcal{C}}_{\Bbb{P}L_{\tilde{1}}|
{\mathcal{U}}_{{\mathcal{T}}|{\mathcal{T}}_1(l)}^{(1)}(\mu)}
({\mathcal{D}}_{\tilde{1},\hat{1}})=N(\al_1),\qquad\hbox{where}\\
\al_1 \in \Ga\big({\overV}_{1;1}^{(1)}(\mu);
\hbox{Hom}(L_{\tilde{1}} \oplus L_{\tilde{1}}^*,
\ev_{\hat{0}}^*T\PPP)\big),\\
\big\{\al_1\big|{\mathcal{U}}_{{\mathcal{T}}_{1;l}}^{(1)}(\mu)\big\}
(\ups_{\tilde{1}},\ups_{\hat{1}})
 = {\mathcal{D}}_{{\mathcal{T}}_{1;l},\tilde{1}}^{(1)}\ups_{\tilde{1}}+
{\mathcal{D}}_{{\mathcal{T}}_{1;l},\hat{1}}^{(1)}\ups_{\hat{1}}.
\end{gather*}
Thus, by Lemma~\ref{n3tac_contr_lmm3c},
\begin{equation}\label{n3tac_contr_lmm3e11}
\sum_{l\in[N]}
{\mathcal{C}}_{\Bbb{P}L_{\tilde{1}}|{\mathcal{U}}_{{\mathcal{T}}|{\mathcal{T}}_1(l)}}
({\mathcal{D}}_{\tilde{1},\hat{1}})=
\big\lan 4a_{\hat{0}},{\overV}_{1;1}^{(1)}(\mu)\big\ran
-5\big|{\mathcal{S}}_{1;1}(\mu)\big|+
\big|{\mathcal{V}}_{2;1}^{(1,1)}(\mu)\big|.
\end{equation}
(4)\qua Finally, it is easy to see that the set
$\Bbb{P}L_{\tilde{1}}^*|{\mathcal{U}}_{{\mathcal{T}}|{\mathcal{T}}_1}^{(1)}$
is ${\mathcal{D}}_{\tilde{1},\hat{1}}$--hollow.
Indeed, in this case, the target bundle $\tilde{\mathcal{F}}{\mathcal{T}}$
has the same rank as the target in case~(2),
but the domain of $\rho$ is 
$L_{\tilde{1}} \otimes L_{\tilde{1}} \oplus \mathcal{FT}$,
instead of~$\mathcal{FT}$.
Thus, the claim follows from 
equations~\eqref{n3tac_contr_lmm3e5}--\eqref{n3tac_contr_lmm3e11}.
\end{proof}

\begin{lmm}
\label{n3tac_contr_lmm3b}
If $\al_1 \in \Ga\big(\Bbb{P}^1 \times {\overV}_{1;1}^{(1)}(\mu);
\hbox{Hom}(\ga^* \otimes L_{\tilde{1}}^*,
\ga^* \otimes \ev_{\hat{0}}^*T\PPP)\big)$ is given by
$\al_1\big|{\mathcal{U}}_{{\mathcal{T}}_{1;l}}^{(1)}(\mu)
 = {\mathcal{D}}_{{\mathcal{T}}_{1;l},\hat{1}}^{(1)}$,
$$N(\al_1)= \big\lan 4a_{\hat{0}} - c_1(L_{\tilde{1}}^*),
{\overV}_{1;1}^{(1)}(\mu)\big\ran
-2\big|{\mathcal{S}}_{1;1}(\mu)\big|
+\big|{\mathcal{V}}_{2;1}^{(1)}(\mu)\big|.$$
\end{lmm}

\begin{proof}
(1)\qua By Propositions~\ref{zeros_prp} and~\ref{euler_prp},
\begin{equation}\label{n3tac_contr_lmm3b_e1}
N(\al_1)=\big\lan 4a_{\hat{0}} - c_1(L_{\tilde{1}}^*),
{\overV}_{1;1}^{(1)}(\mu)\big\ran-
{\mathcal{C}}_{\Bbb{P}^1\times\partial{\overV}_{1;1}^{(1)}(\mu)}(\al_1^{\perp}),
\end{equation}
where $\al_1^{\perp}$ denotes the composition of $\al_1$ with the projection map 
onto the quotient ${\mathcal{O}}_1$ of \hbox{$\ga^* \otimes \ev_{\hat{0}}^*T\PPP$}
by generic trivial line subbundle~$\Bbb{C}\bar{\nu}_1$.
Suppose ${\mathcal{T}} < {\mathcal{T}}_{1;l}$ is a bubble type such that
$${\overV}_{1;1}^{(1)}(\mu)\cap
{\mathcal{U}}_{{\mathcal{T}}_{1;l},{\mathcal{T}}} \neq \eset.$$
(2)\qua If $\chi_{\mathcal{T}}(\tilde{1},\hat{1}) > 0$, by
Lemma~\ref{n3p3_str_lmm1},
$${\overV}_{1;1}^{(1)}(\mu)\cap
{\mathcal{U}}_{{\mathcal{T}}_{1;l},{\mathcal{T}}}\subset
{\mathcal{U}}_{{\mathcal{T}}|{\mathcal{T}}_{1;l}}^{(1)}(\mu),$$ 
and thus $\hat{I}^+ = \{h\}$ is a single-element set.
Furthermore, $d_{\tilde{1}} \neq 0$, since our constraints~$\mu$ are disjoint.
Thus, if $j_{\hat{1}} = \tilde{1}$, $\al_1$ extends over
$\Bbb{P}^1 \times {\mathcal{U}}_{{\mathcal{T}}|{\mathcal{T}}_{1;l}}^{(1)}(\mu)$,
and this extension does not vanish by Lemma~\ref{str_lmm}.
It follows that $\Bbb{P}^1 \times {\mathcal{U}}_{{\mathcal{T}}_{1;l},{\mathcal{T}}}$
does not contribute to  ${\mathcal{C}}_{\Bbb{P}^1
               \times\partial{\overV}_{1;1}^{(1)}(\mu)}(\al_1^{\perp})$.
If $j_{\hat{1}} = h$ and $d_h = 0$, by Proposition~\ref{str_prp}, 
$\al_1$ again has a nonvanishing extension over
$\Bbb{P}^1 \times {\mathcal{U}}_{{\mathcal{T}}|{\mathcal{T}}_{1;l}}^{(1)}(\mu)$.
Thus, we only need to consider the case 
$j_{\hat{1}} = h$ and $d_{\tilde{1}},d_h > 0$.
By Proposition~\ref{str_prp},
$$\al_1\big(\phi_{{\mathcal{T}}_{1;l},{\mathcal{T}}}(\ups)\big)=
\{{\mathcal{D}}_{{\mathcal{T}},\hat{1}}^{(1)} + \ve(\ups)\big\}\ups^*
\quad\mbox{for all}\  \ups \in \mathcal{FT}_{\de} - Y(\mathcal{FT};\{h\}).$$
Since ${\mathcal{D}}_{{\mathcal{T}},\hat{1}}^{(1)}$ does not vanish 
on ${\mathcal{U}}_{{\mathcal{T}}|{\mathcal{T}}_{1;l}}^{(1)}(\mu)$,
from Proposition~\ref{euler_prp}, we conclude that
\begin{gather*}
\sum_{\chi_{\mathcal{T}}(\tilde{1},\hat{1})>0}    
{\mathcal{C}}_{\Bbb{P}^1\times{\mathcal{U}}_{{\mathcal{T}}_{1;l},{\mathcal{T}}}}(\al_1^{\perp})
=-N(\al_2)\cdot\big|{\mathcal{V}}_{2;1}^{(1,1)}(\mu)\big|,
\qquad\hbox{where}\\
\al_2\in\Ga\big(\Bbb{P}^1;\hbox{Hom}(\Bbb{C},\Bbb{C}^3/\ga)\big)
\end{gather*}
is a nonvanishing section.
Thus, by Proposition~\ref{zeros_prp},
\begin{equation}\label{n3tac_contr_lmm3b_e3}
\sum_{\chi_{\mathcal{T}}(\tilde{1},\hat{1})>0}    
{\mathcal{C}}_{\Bbb{P}^1\times{\mathcal{U}}_{{\mathcal{T}}_{1;l},{\mathcal{T}}}}(\al_1^{\perp})
=-\big|{\mathcal{V}}_{2;1}^{(1,1)}(\mu)\big|.
\end{equation}
(3)\qua If $\chi_{\mathcal{T}}(\tilde{1},\hat{1}) = 0$, by
Lemma~\ref{n3p3_str_lmm2},
$${\overV}_{1;1}^{(1)}(\mu)\cap {\mathcal{U}}_{{\mathcal{T}}_{1;l},{\mathcal{T}}} 
\subset {\mathcal{S}}_{{\mathcal{T}}|{\mathcal{T}}_{1;l}}(\mu).$$
Thus, $\hat{I}^+ = \{h\}$ is again a single-element set.
Adding the expansion of $\ev_{\tilde{1}} \times \ev_{\hat{1}}$
of Lemma~\ref{n3p3_str_lmm2} times $(y_{\hat{1}} - x_h)^{-1}$
to the expansion of ${\mathcal{D}}_{{\mathcal{T}}_{1;l},\hat{1}}^{(1)}$
of~(3c) of Proposition~\ref{str_prp}, we obtain
\begin{gather*}
\al_1\big(\phi_{{\mathcal{T}}_{1;l},{\mathcal{T}}}(\ups)\big)=
-(y_{\hat{1}} - x_h)^{-2}\otimes
\big\{{\mathcal{D}}_{{\mathcal{T}},h}^{(2)} + \ve(\ups)\big\}\ups \otimes \ups\\
\mbox{for all}\ \ups \in \mathcal{FT}_{\de}
~\mbox{such that}~
\phi_{{\mathcal{T}}_{1;l},{\mathcal{T}}}(\ups) \in {\mathcal{V}}_{1;1}^{(1)}.
\end{gather*}
Thus, as in the second half of (2) of the proof of
Lemma~\ref{n3p3_contr_lmm1b}, we can conclude that
\begin{gather*}
\sum_{\chi_{\mathcal{T}}(\tilde{1},\hat{1})=0}    
{\mathcal{C}}_{\Bbb{P}^1\times{\mathcal{U}}_{{\mathcal{T}}(l),{\mathcal{T}}}}(\al_1^{\perp})
=2N(\al_2)\cdot\big|{\mathcal{S}}_{1;1}(\mu)\big|,
\qquad\hbox{where}\\
\al_2 \in \Ga\big(\Bbb{P}^1;\hbox{Hom}(\Bbb{C},\Bbb{C}^3/\ga)\big)
\end{gather*}
is a nonvanishing section.
Thus, by Proposition~\ref{zeros_prp},
\begin{equation}\label{n3tac_contr_lmm3b_e5}
\sum_{\chi_{\mathcal{T}}(\tilde{1},\hat{1})=0}    
{\mathcal{C}}_{\Bbb{P}^1\times{\mathcal{U}}_{{\mathcal{T}}_{1;l},{\mathcal{T}}}}(\al_1^{\perp})
=2\big|{\mathcal{S}}_{1;1}(\mu)\big|.
\end{equation}
The claim follows from 
equations~\eqref{n3tac_contr_lmm3b_e1}--\eqref{n3tac_contr_lmm3b_e5}.
\end{proof}

\begin{lmm}
\label{n3tac_contr_lmm3c}
If $\al_1 \in \Ga\big({\overV}_{1;1}^{(1)}(\mu);
\hbox{Hom}(L_{\tilde{1}} \oplus L_{\tilde{1}}^*,
\ev_{\hat{0}}^*T\PPP)\big)$ is given by
$$\al_1|_{{\mathcal{U}}_{{\mathcal{T}}_{1;l}}^{(1)}(\mu)}
(\ups_{\tilde{1}},\ups_{\hat{1}})
={\mathcal{D}}_{{\mathcal{T}}_{1;l},\tilde{1}}^{(1)}\ups_{\tilde{1}}+
{\mathcal{D}}_{{\mathcal{T}}_{1;l},\hat{1}}^{(1)}\ups_{\hat{1}}$$
for all $l \in [N]$,
$$N(\al_1)= \blr{4a_{\hat{0}},{\overV}_{1;1}^{(1)}(\mu)}
-5\big|{\mathcal{S}}_{1;1}(\mu)\big| +\big|{\mathcal{V}}_{2;1}^{(1,1)}(\mu)\big|.$$
\end{lmm}

\begin{proof}
(1)\qua Let ${\mathcal{F}} = L_{\tilde{1}} \oplus L_{\tilde{1}}^*$.
By Propositions~\ref{zeros_prp} and~\ref{euler_prp},
\begin{equation}\label{n3tac_contr_lmm3c_e1}
N(\al_1)=\blr{4a_{\hat{0}},{\overV}_{1;1}^{(1)}(\mu)}
- {\mathcal{C}}_{\Bbb{P}{\mathcal{F}}|\partial{\overV}_{1;1}^{(1)}(\mu)}(\al_1^{\perp}),
\end{equation}
where $\tilde{\al}_1^{\perp}$ denotes the composition of  the section 
$$\tilde{\al}_1\in\Ga\big(\Bbb{P}{\mathcal{F}};
\hbox{Hom}(\ga_{\mathcal{F}},\pi_{\Bbb{P}{\mathcal{F}}}^*\ev_{\hat{0}}^*T\PPP)\big),$$ 
induced by~$\al$, with the projection map onto the quotient ${\mathcal{O}}_1$ 
of $\pi_{\Bbb{P}{\mathcal{F}}}^*\ev_{\hat{0}}^*T\PPP$
by generic trivial line subbundle~$\Bbb{C}\bar{\nu}_1$.
Suppose ${\mathcal{T}} < {\mathcal{T}}_{1;l}$ is a bubble type such~that
$${\overV}_{1;1}^{(1)}(\mu)\cap
{\mathcal{U}}_{{\mathcal{T}}_{1;l},{\mathcal{T}}} \neq \eset.$$
(2)\qua  If $\chi_{\mathcal{T}}(\tilde{1},\hat{1}) > 0$, then
$${\overV}_{1;1}^{(1)}(\mu)\cap {\mathcal{U}}_{{\mathcal{T}}_{1;l},{\mathcal{T}}}
\subset {\mathcal{U}}_{{\mathcal{T}}|{\mathcal{T}}_{1;l}}^{(1)}(\mu),$$ 
$\hat{I}^+ = \{h\}$ is a single-element set, and 
$d_{\tilde{1}} \neq 0$.
As before, we only need to consider the case
$j_{\hat{1}} = h$ and $d_{\tilde{1}},d_h > 0$.
By Proposition~\ref{str_prp},
$$\tilde{\al}_1\big(\phi_{{\mathcal{T}}_{1;l},{\mathcal{T}}}
([\ups_{\tilde{1}},\ups_{\hat{1}}],\ups)\big)
=\big\{\ups_{\hat{1}} \otimes {\mathcal{D}}_{{\mathcal{T}},\hat{1}}^{(1)} +
 \ve(\ups)\big\}\ups^*
\quad\mbox{for all}\ \ups \in \mathcal{FT}_{\de} - Y(\mathcal{FT};\{h\}),$$
where $\ve$ is a $C^0$--negligible map.
Since the linear map 
$\ups_{\hat{1}} \otimes {\mathcal{D}}_{{\mathcal{T}},\hat{1}}^{(1)}$
does not vanish outside of the set 
${\mathcal{Z}}_{\mathcal{T}} = \Bbb{P}L_{\tilde{1}}$,
by Proposition~\ref{str_prp} and a rescaling of the linear~map, 
\begin{gather*}
\sum_{\chi_{\mathcal{T}}(\tilde{1},\hat{1})>0}    
{\mathcal{C}}_{\Bbb{P}{\mathcal{F}}|{\mathcal{U}}_{{\mathcal{T}}|{\mathcal{T}}_{1;l}}^{(1)}(\mu)
-{\mathcal{Z}}_{\mathcal{T}}}(\al_1^{\perp})
=-N(\al_2)\cdot\big|{\mathcal{V}}_{2;1}^{(1,1)}(\mu)\big|,
\qquad\hbox{where}\\
\al_2 \in \Ga\big(\Bbb{P}^1;
\hbox{Hom}(\ga^*,\ga^* \otimes \Bbb{C}^2)\big)
\end{gather*}
is a nonvanishing section.
Thus, by Proposition~\ref{zeros_prp},
$$\sum_{\chi_{\mathcal{T}}(\tilde{1},\hat{1})>0}    
{\mathcal{C}}_{\Bbb{P}{\mathcal{F}}|{\mathcal{U}}_{{\mathcal{T}}|{\mathcal{T}}_{1;l}}^{(1)}(\mu)
-{\mathcal{Z}}_{\mathcal{T}}}(\tilde{\al}_1^{\perp})
=-\big|{\mathcal{V}}_{2;1}^{(1,1)}(\mu)\big|.$$
On the other hand, we can model a neighborhood of ${\mathcal{Z}}_{\mathcal{T}}$
in $\Bbb{P}{\mathcal{F}}$ by the map
$$L_{\tilde{1}}^*\otimes L_{\tilde{1}}^*\lra\Bbb{P}{\mathcal{F}}, \qquad 
\big([\ups_{\tilde{1}},\ups_{\hat{1}}],u\big)
\lra\big[\ups_{\tilde{1}},u(\ups_{\tilde{1}})\big].$$
Since by Proposition~\ref{str_prp}
\begin{gather*}
\tilde{\al}_1\big(\phi_{{\mathcal{T}}_{1;l},{\mathcal{T}}}
([\ups_{\tilde{1}},\ups_{\hat{1}}],u,\ups)\big)
=\big\{{\mathcal{D}}_{{\mathcal{T}},\hat{1}}^{(1)} +
 \ve_1(\ups)\big\}u \otimes \ups^*
+{\mathcal{D}}_{{\mathcal{T}},\tilde{1}} + \ve_1(\ups)\\
\mbox{for all}\ \ups \in \mathcal{FT}_{\de} - Y(\mathcal{FT};\{h\}),
\end{gather*}
it follows that ${\mathcal{Z}}_{\mathcal{T}}|
{\mathcal{U}}_{{\mathcal{T}}|{\mathcal{T}}_{1;l}}^{(1)}(\mu)$ is
$\tilde{\al}_1^{\perp}$--hollow.
Thus,
\begin{equation}\label{n3tac_contr_lmm3c_e3}
\sum_{\chi_{\mathcal{T}}(\tilde{1},\hat{1})>0}    
{\mathcal{C}}_{\Bbb{P}{\mathcal{F}}|{\mathcal{U}}_{{\mathcal{T}}_{1;l},{\mathcal{T}}}}
(\tilde{\al}_1^{\perp})=
-\big|{\mathcal{V}}_{2;1}^{(1,1)}(\mu)\big|.
\end{equation}
(3)\qua  If $\chi_{\mathcal{T}}(\tilde{1},\hat{1}) = 0$, 
$${\overV}_{1;1}^{(1)}(\mu)\cap {\mathcal{U}}_{{\mathcal{T}}_{1;l},{\mathcal{T}}}
\subset {\mathcal{S}}_{{\mathcal{T}}|{\mathcal{T}}_{1;l}}(\mu)$$
and $\hat{I}^+ = \{h\}$ is a single-element set.
Subtracting the expansion of $\ev_{\tilde{1}} \times \ev_{\hat{1}}$
of Lemma~\ref{n3p3_str_lmm2} times 
$$(y_{\hat{1}} - x_h)\qquad\hbox{and times}\qquad
(y_{\hat{1}} - x_h)^{-1}$$
from the expansions of ${\mathcal{D}}_{{\mathcal{T}}_{1;l},\tilde{1}}^{(1)}$
and  ${\mathcal{D}}_{{\mathcal{T}}_{1;l},\hat{1}}^{(1)}$
in Proposition~\ref{str_prp}, we obtain
\begin{equation*}\begin{split}
\tilde{\al}_1\big(\phi_{{\mathcal{T}}_{1;l},{\mathcal{T}}}
([\ups_{\tilde{1}},\ups_{\hat{1}}],\ups)\big)
=-(y_{\hat{1}} - x_h)^{-1} \otimes  
\big(\ups_{\tilde{1}} \otimes 
{\mathcal{D}}_{{\mathcal{T}},h}^{(2)} + \ve_{\tilde{1}}(\ups)\big)\ups^{\otimes2}
\qquad\qquad&\\
- (y_{\hat{1}} - x_h)^{-3} \otimes \big(\ups_{\hat{1}} \otimes 
{\mathcal{D}}_{{\mathcal{T}},h}^{(2)} + \ve_{\hat{1}}(\ups)\big)
\ups^{\otimes2}&
\end{split}\end{equation*}
for all $\ups \in \mathcal{FT}_{\de}$
such that $\phi_{{\mathcal{T}}_{1;l},{\mathcal{T}}}(\ups) \in 
{\mathcal{V}}_{1;1}^{(1)}$.
Let 
$${\mathcal{Z}}_{\mathcal{T}}=\big\{
[\ups_{\tilde{1}},\ups_{\hat{1}}] \in \Bbb{P}{\mathcal{F}} :
\ups_{\tilde{1}} + (y_{\hat{1}} - x_h)^{-2}\ups_{\hat{0}} = 0
 \in L_{\tilde{1}}\big\}.$$
Similarly to the argument in (2) of the proof of Lemma~\ref{n3tac_contr_lmm3},
from the above we can conclude that
\begin{gather*}
\sum_{\chi_{\mathcal{T}}(\tilde{1},\hat{1})=0}    
{\mathcal{C}}_{\Bbb{P}{\mathcal{F}}|{\mathcal{S}}_{{\mathcal{T}}|{\mathcal{T}}_{1;l}}(\mu)
-{\mathcal{Z}}_{\mathcal{T}}}(\tilde{\al}_1^{\perp})=
2N(\al_2)\cdot\big|{\mathcal{S}}_{1;1}(\mu)\big|,
\quad\hbox{where}\\
\al_2\in\Ga\big(\Bbb{P}^1;
\hbox{Hom}(\ga^*,\ga^* \otimes \Bbb{C}^2)\big)
\end{gather*}
is a nonvanishing section.
Thus, by Lemma~\ref{zeros_prp},
$$\sum_{\chi_{\mathcal{T}}(\tilde{1},\hat{1})=0}    
{\mathcal{C}}_{\Bbb{P}{\mathcal{F}}|{\mathcal{S}}_{{\mathcal{T}}|{\mathcal{T}}_{1;l}}(\mu)
-{\mathcal{Z}}_{\mathcal{T}}}(\tilde{\al}_1^{\perp})=
2\big|{\mathcal{S}}_{1;1}(\mu)\big|.$$
On the other hand, the same argument as in (3)
of the proof of Lemma~\ref{n3tac_contr_lmm3} shows that
$${\mathcal{C}}_{{\mathcal{Z}}_{\mathcal{T}}|{\mathcal{S}}_{{\mathcal{T}}|{\mathcal{T}}_{1;l}}(\mu)}
(\tilde{\al}_1)=3\big|{\mathcal{Z}}_{\mathcal{T}}|.$$
Thus, we conclude that
\begin{equation}\label{n3tac_contr_lmm3c_e5}
\sum_{\chi_{\mathcal{T}}(\tilde{1},\hat{1})=0}    
{\mathcal{C}}_{\Bbb{P}{\mathcal{F}}|{\mathcal{U}}_{{\mathcal{T}}_{1;l},{\mathcal{T}}}}
(\tilde{\al}_1^{\perp})=5\big|{\mathcal{S}}_{1;1}(\mu)\big|.
\end{equation}
The claim follows from
equations~\eqref{n3tac_contr_lmm3c_e1}--\eqref{n3tac_contr_lmm3c_e5}.
\end{proof}

\begin{lmm}
\label{n3tac_contr_lmm4}
Suppose ${\mathcal{T}} = (M_1,I;j,\under{d})$ is a bubble type such that
${\mathcal{T}} < {\mathcal{T}}_1$, $\chi_{\mathcal{T}}(\tilde{1},\hat{1}) > 0$,
and \hbox{$d_{\tilde{1}} > 0$}.

{\rm(1)}\qua If $j_{\hat{1}} = \tilde{1}$, $|\hat{I}^+| \neq 1$, 
or $d_{j_{\hat{1}}} = 0$,
$\Bbb{P}(L_{\tilde{1}} \oplus L_{\tilde{1}}^*)|
{\mathcal{U}}_{{\mathcal{T}}|{\mathcal{T}}_1}^{(0)}$ is a finite union of
${\mathcal{D}}_{\tilde{1},\hat{1}}$--hollow subspaces and thus
$${\mathcal{C}}_{\Bbb{P}(L_{\tilde{1}}\oplus L_{\tilde{1}}^*)|
{\mathcal{U}}_{{\mathcal{T}}_1,{\mathcal{T}}}}({\mathcal{D}}_{\tilde{1},\hat{1}}) = 0.$$
{\rm(2)}\qua The total contribution from the boundary spaces 
${\mathcal{U}}_{{\mathcal{T}}_1,{\mathcal{T}}}$
such that $j_{\hat{1}} > \tilde{1}$, $|\hat{I}^+| = 1$, 
and $d_{j_{\hat{1}}} > 0$,
is given~by
$$\sum_{d_{\tilde{1}}>0}
{\mathcal{C}}_{\Bbb{P}(L_{\tilde{1}}\oplus L_{\tilde{1}}^*)|
{\mathcal{U}}_{{\mathcal{T}}_1,{\mathcal{T}}}}({\mathcal{D}}_{\tilde{1},\hat{1}})
=-\big\lan 4a_{\hat{0}} + c_1({\mathcal{L}}_{\tilde{1}}^*),
{\overV}_2^{(1,1)}(\mu)\big\ran
+\big|{\mathcal{V}}_2^{(1)}(\mu)\big|
+2\big|{\mathcal{S}}_2(\mu)\big|.$$
\end{lmm}

\begin{proof}
(1)\qua We proceed as in (1) of the proof of Lemma~\ref{n3tac_contr_lmm3}.
In particular, we have
$${\mathcal{D}}_{\tilde{1},\hat{1}}\phi_{{\mathcal{T}}_1,{\mathcal{T}}}
([\ups_{\tilde{1}},\ups_{\hat{1}}];\ups)=
\big\{\al + \ve(\ups)\big\}\rho(\ups)
\quad\mbox{for all}\ \ups \in \mathcal{FT}_{\de} - 
Y(\mathcal{FT};\hat{I}^+).$$
In this equation, $\rho$ is the monomials map on $\mathcal{FT}$
defined~by
$$\rho_h(\ups)=
\begin{cases}
\rho_{{\mathcal{T}},\hat{1};\hat{1}}^{(1)}(\ups),&
\hbox{if}~h \in \chi_{j_{\hat{1}}}({\mathcal{T}}) - \chi_{\hat{1}}({\mathcal{T}});\\
\rho_{{\mathcal{T}},\hat{1};h}^{(1;1)}(\ups),&
\hbox{if}~h \in \chi_{\hat{1}}({\mathcal{T}})
\end{cases}$$
with values in the bundle
$$\tilde{\mathcal{F}}{\mathcal{T}}=
\bigoplus_{h\in\chi_{j_{\hat{1}}}({\mathcal{T}})}\tilde{\mathcal{F}}_h{\mathcal{T}},
\quad\hbox{where}\quad
\tilde{\mathcal{F}}_h{\mathcal{T}}=
\begin{cases}
\tilde{\mathcal{F}}_{{\mathcal{T}},\hat{1};\hat{1}}^{(1)}{\mathcal{T}},&
\hbox{if}~h \in \chi_{j_{\hat{1}}}({\mathcal{T}}) - \chi_{\hat{1}}({\mathcal{T}});\\
\tilde{\mathcal{F}}_{{\mathcal{T}},\hat{1};h}^{(1;1)}{\mathcal{T}},&
\hbox{if}~h \in \chi_{\hat{1}}({\mathcal{T}}).
\end{cases}$$
The linear map 
$\al\co \tilde{\mathcal{F}}{\mathcal{T}}\lra
\ga_{L_{\tilde{1}}\oplus L_{\tilde{1}}^*}^* \otimes \ev_{\hat{0}}^*T\PPP$
is given~by
$$\al\big([\ups_{\tilde{1}},\ups_{\hat{1}}],\tilde{\ups}_h\big)=
\ups_{\hat{1}}\otimes
\begin{cases}
{\mathcal{D}}_{{\mathcal{T}},j_{\hat{1}}^*({\mathcal{T}})}^{(1)}\tilde{\ups}_h,&
\hbox{if}~h \in \chi_{j_{\hat{1}}}({\mathcal{T}}) - \chi_{\hat{1}}({\mathcal{T}});\\
-(y_{h;\hat{1}} - x_{\hat{1};h})^{-2}
\otimes{\mathcal{D}}_{{\mathcal{T}},h}^{(1)}\tilde{\ups}_h,&
\hbox{if}~h \in \chi_{\hat{1}}({\mathcal{T}}).
\end{cases}$$
In particular, by Lemma~\ref{str_lmm}, $\al$ has full rank over 
${\mathcal{U}}_{{\mathcal{T}}|{\mathcal{T}}_1}^{(1)}(\mu)$ outside of the set
${\mathcal{Z}}_{\mathcal{T}} \equiv \Bbb{P}L_{\tilde{1}}$.
Thus, 
$$\Bbb{P}(L_{\tilde{1}} \oplus L_{\tilde{1}}^*)|
{\mathcal{U}}_{{\mathcal{T}}|{\mathcal{T}}_1}^{(1)}(\mu) - {\mathcal{Z}}_{\mathcal{T}}$$
is ${\mathcal{D}}_{\tilde{1},\hat{1}}$--hollow unless
$\hat{I}^+ = \{h\}$ is a one-element set, $j_{\hat{1}} = h$,
and \hbox{$d_h \neq 0$}.
If $\hat{I}^+ = \{h\}$, $j_{\hat{1}} = h$, and $d_h \neq 0$,
the degree of the map $\rho$ is~$-1$, and 
from Proposition~\ref{euler_prp} and a rescaling of the linear map,
\begin{gather*}
{\mathcal{C}}_{\Bbb{P}(L_{\tilde{1}}\oplus L_{\tilde{1}}^*)|
{\mathcal{U}}_{{\mathcal{T}}|{\mathcal{T}}_1}^{(1)}(\mu)-{\mathcal{Z}}_{\mathcal{T}}}
({\mathcal{D}}_{\tilde{1},\hat{1}})=-N(\al_1'),\qquad\hbox{where}\\
\al_1'\in\Ga
\big(\Bbb{P}{\mathcal{F}};\hbox{Hom}(\ga_{\mathcal{F}}^* \otimes L_h^*,
\ga_{\mathcal{F}}^* \otimes \ev_{\hat{0}}^*T\PPP)\big),\\
{\mathcal{F}} = L_{\tilde{1}} \oplus L_{\tilde{1}}^*
 \lra {\mathcal{U}}_{{\mathcal{T}}|{\mathcal{T}}_1}^{(1)}(\mu), \quad
\al_1'={\mathcal{D}}_{{\mathcal{T}},\hat{1}}^{(1)}.
\end{gather*}
Thus, using Lemma~\ref{n3tac_contr_lmm4b},
another rescaling of the linear map, and an obvious symmetry, 
we obtain
\begin{equation}\label{n3tac_contr_lmm4e3}\begin{split}
\sum_{\mathcal{T}}{\mathcal{C}}_{\Bbb{P}(L_{\tilde{1}}\oplus L_{\tilde{1}}^*)|
{\mathcal{U}}_{{\mathcal{T}}_1|{\mathcal{T}}_1}^{(1)}(\mu)-{\mathcal{Z}}_{\mathcal{T}}}
({\mathcal{D}}_{\tilde{1},\hat{1}})
=&-\big\lan 4a_{\hat{0}}+c_1({\mathcal{L}}_{\tilde{1}}^*),
{\overV}_2^{(1,1)}(\mu)\big\ran\\
&\qquad\qquad +\big|{\mathcal{V}}_2^{(1)}(\mu)\big| +2\big|{\mathcal{S}}_2(\mu)\big|.
\end{split}\end{equation}
Finally, it is easy to see that the set
${\mathcal{Z}}_{\mathcal{T}}|{\mathcal{U}}_{{\mathcal{T}}|{\mathcal{T}}_1}^{(1)}(\mu)$
is ${\mathcal{D}}_{\tilde{1},\hat{1}}$--hollow for all bubble types~${\mathcal{T}}$.
\end{proof}

\begin{lmm}
\label{n3tac_contr_lmm4b}
Suppose $\tilde{\mathcal{T}} = (M_2,I_2;\tilde{j},\tilde{d})$ is
a bubble type such that 
$$\tilde{j}_{\hat{1}}=\tilde{1}, \quad j_{\hat{2}}=\tilde{2}, \quad 
\tilde{d}_{\tilde{1}},\tilde{d}_{\tilde{2}}>0, \quad\hbox{and}\quad
\tilde{d}_{\tilde{1}} + \tilde{d}_{\tilde{2}}=d,$$
and ${\mathcal{F}} \lra {\overU}_{\tilde{\mathcal{T}},\tilde{\mathcal{T}}}$ 
is a rank-two vector bundle such that 
$c_1({\mathcal{F}})|{\overU}_{\tilde{\mathcal{T}}}^{(1)}(\mu) = 0$.
If 
$$\al_1 \in \Ga\big(\Bbb{P}{\mathcal{F}}|
{\mathcal{U}}_{\tilde{\mathcal{T}}}^{(1)}(\mu);
\hbox{Hom}(\ga^* \otimes L_{\tilde{1}};
\ga^* \otimes \ev_{\hat{0}}^*T\PPP)\big)$$
is given by  $\al_1(\ups_{\tilde{1}}) = 
{\mathcal{D}}_{\tilde{\mathcal{T}},\tilde{1}}^{(1)}\ups_{\tilde{1}}$,
$$N(\al_1)=\blr{4a_{\hat{0}}+c_1({\mathcal{L}}_{\tilde{1}}^*),
{\overU}_{\tilde{\mathcal{T}}}^{(1)}(\mu)}
-\big|{\mathcal{U}}_{\tilde{\mathcal{T}}/\hat{1}}^{(1)}(\mu)\big|
-\big|{\mathcal{U}}_{\tilde{\mathcal{T}}/\{\hat{1},\hat{2}\}}(\mu)\big|.$$
In other words, the sum of the numbers $N(\al_1)$ taken over all
bubble types $\tilde{\mathcal{T}}$ over the above form is given~by
$$\sum_{\tilde{\mathcal{T}}}N(\al_1)=
\blr{4a_{\hat{0}}+c_1({\mathcal{L}}_{\tilde{1}}^*),{\overV}_2^{(1,1)}(\mu)}
-\big|{\mathcal{V}}_2^{(1)}(\mu)\big|-2\big|{\mathcal{S}}_2(\mu)\big|.$$
\end{lmm}

\begin{proof}
By Proposition~\ref{zeros_prp} and the assumption
$c_1({\mathcal{F}})|{\overU}_{\tilde{\mathcal{T}}}^{(1)}(\mu) = 0$,
\begin{equation}\label{n3tac_contr_lmm4b_e1}
N(\al_1)=\blr{4a_{\hat{0}}+c_1(L_{\tilde{1}}^*),
{\overU}_{\tilde{\mathcal{T}}}^{(1)}(\mu)}
-{\mathcal{C}}_{\Bbb{P}{\mathcal{F}}|\partial
{\overU}_{\tilde{\mathcal{T}}}^{(1)}(\mu)}(\al_1^{\perp}),
\end{equation}
where $\al_1^{\perp}$ denotes the composition of $\al_1$ with the
projection map
onto the quotient ${\mathcal{O}}_1$ of \hbox{$\ga^* \otimes
\ev_{\hat{0}}^*T\PPP$}
by generic trivial line subbundle~$\Bbb{C}\bar{\nu}_1$.
Suppose ${\mathcal{T}} = (M_2,I;j,\under{d})$ is a bubble type such that
${\mathcal{T}} < \tilde{\mathcal{T}}$ and
$${\overU}_{\tilde{\mathcal{T}}}^{(1)}(\mu)\cap
{\mathcal{U}}_{\tilde{\mathcal{T}},{\mathcal{T}}} \neq\eset.$$
The section~$\al_1$ extends over
$\Bbb{P}{\mathcal{F}}|\partial{\overU}_{\tilde{\mathcal{T}}}^{(1)}(\mu)$.
Moreover, by Lemmas~\ref{n3p3_str_lmm3} and~\ref{str_lmm},
this extension does not vanish unless $d_{\tilde{1}} = 0$.
Thus, in computing the number
${\mathcal{C}}_{\Bbb{P}{\mathcal{F}}|\partial
{\overU}_{\tilde{\mathcal{T}}}^{(1)}(\mu)}(\al_1^{\perp})$,
we only need to consider bubble types ${\mathcal{T}}$ such that
$d_{\tilde{1}} = 0$.

(2)\qua If $d_{\tilde{2}} \neq 0$, then
$${\overU}_{\tilde{\mathcal{T}}}^{(1)}(\mu)\cap
{\mathcal{U}}_{\tilde{\mathcal{T}},{\mathcal{T}}}
\subset {\mathcal{U}}_{{\mathcal{T}}|\tilde{\mathcal{T}}}^{(1)}(\mu)$$
and ${\mathcal{T}} = \tilde{\mathcal{T}}(\hat{1})$ or  ${\mathcal{T}}
= \tilde{\mathcal{T}}(l)$
for some $l \in [N]\cap M_{\tilde{1}}{\mathcal{T}}$.
Moreover, by Proposition~\ref{str_prp},
$$\al_1\big(\phi_{\tilde{\mathcal{T}},{\mathcal{T}}}(\ups)\big)
=\big\{{\mathcal{D}}_{{\mathcal{T}},h}^{(1)} + \ve(\ups)\big\}\ups
\qquad\mbox{for all}\ \ups \in \mathcal{FT}_{\de},$$
if $h$ is the unique element of $\hat{I}^+$.
Thus, by the same argument as in (1) of the proof of
Lemma~\ref{n3tac_contr_lmm3}, we conclude that
$${\mathcal{C}}_{\Bbb{P}{\mathcal{F}}|
{\mathcal{U}}_{{\mathcal{T}}|\tilde{\mathcal{T}}}^{(1)}(\mu)}(\al_1^{\perp})
=N(\al_2),
\quad\hbox{where}\quad
\al_2 \in \Ga(\Bbb{P}^1 \times
{\overU}_{{\mathcal{T}}|\tilde{\mathcal{T}}}^{(1)}(\mu);
\hbox{Hom}(\Bbb{C},\Bbb{C}/\ga)\big)$$
is a nonvanishing section.
Thus, from Proposition~\ref{zeros_prp} and
the decomposition~\eqref{cart_split}, we obtain
\begin{equation}\label{n3tac_contr_lmm4b_e3}
\sum_{d_{\tilde{2}}>0}
{\mathcal{C}}_{\Bbb{P}{\mathcal{F}}|
{\mathcal{U}}_{\tilde{\mathcal{T}},{\mathcal{T}}}}(\al_1^{\perp})
=\big|{\mathcal{U}}_{\tilde{\mathcal{T}}/\hat{1}}^{(1)}(\mu)\big|
+\sum_{l\in[N]\cap M_{\tilde{1}}{\mathcal{T}}}
\big|{\mathcal{U}}_{\tilde{\mathcal{T}}/l}^{(1)}(\mu)\big|.
\end{equation}
(3)\qua If $d_{\tilde{2}} = 0$, then
$${\overU}_{\tilde{\mathcal{T}}}^{(1)}(\mu)\cap
{\mathcal{U}}_{\tilde{\mathcal{T}},{\mathcal{T}}} \subset
{\mathcal{S}}_{{\mathcal{T}}|\tilde{\mathcal{T}}}(\mu),$$
$j_{\hat{1}} = \tilde{1}$, $j_{\hat{2}} = \tilde{2}$,
$d_{\tilde{1}} = d_{\tilde{2}} = 0$,
and $|\hat{I}^+| = 2$.
Furthermore,
$$\al_1\big(\phi_{\tilde{\mathcal{T}},{\mathcal{T}}}(\ups)\big)
=\big\{{\mathcal{D}}_{{\mathcal{T}},h_1}^{(1)} + \ve(\ups)\big\}\ups_{h_1}
\qquad\mbox{for all}\ \ups = (\ups_{h_1},\ups_{h_2}) \in
\mathcal{FT}_{\de},$$
if $\hat{I}^+ = \{h_1,h_2\}$ and $h_1 = j_{\hat{1}}$.
By an argument similar to (4) of Lemma~\ref{n3p3_contr_lmm1b},
we conclude that
$${\mathcal{C}}_{\Bbb{P}{\mathcal{F}}|{\mathcal{S}}_{{\mathcal{T}}|\tilde{\mathcal{T}}}(\mu)}
(\al_1^{\perp})=N(\al_2),
\quad\hbox{where}\quad
\al_2 \in \Ga(\Bbb{P}^1 \times
{\mathcal{S}}_{{\mathcal{T}}|\tilde{\mathcal{T}}}(\mu);
\hbox{Hom}(\Bbb{C},\Bbb{C}/\ga)\big)$$
is a nonvanishing section.
Thus, from Proposition~\ref{zeros_prp} and
the decomposition~\eqref{cart_split}, we obtain
\begin{equation}\label{n3tac_contr_lmm4b_e5}
{\mathcal{C}}_{\Bbb{P}{\mathcal{F}}|
{\mathcal{U}}_{\ti{\mathcal{T}},{\mathcal{T}}}}(\al_1^{\perp})
=\big|{\mathcal{S}}_{\ti{\mathcal{T}}/\{\hat{1},\hat{2}\}}(\mu)\big|.
\end{equation}
The claim follows by plugging equations~\eqref{n3tac_contr_lmm4b_e3}
and~\eqref{n3tac_contr_lmm4b_e5} into~\eqref{n3tac_contr_lmm4b_e1}
and using~\eqref{psi_class1} and~\eqref{psi_class2}.
\end{proof}

\begin{remark}
{\rm By the second rescaling of the linear map referred to in the proof of 
Lemma~\ref{n3tac_contr_lmm4},
the number $\sum_{\tilde{\mathcal{T}}}N(\al_1)$ of Lemma~\ref{n3tac_contr_lmm4}
does not change if we replace ${\mathcal{D}}_{\ti{\mathcal{T}},\ti{1}}^{(1)}$
by~${\mathcal{D}}_{\ti{\mathcal{T}},\hat{1}}^{(1)}$.
However, a direct computation, 
ie using Propositions~\ref{zeros_prp}, \ref{euler_prp}, and~\ref{str_prp}, 
gives a slightly different answer.
As a result, we obtain yet another enumerative relationship:
$$2\blr{\eta_{\hat{0},1},{\overV}_2^{(1,1)}(\mu)}
=\big|{\mathcal{V}}_{2,(0,1)}^{(1,0,1)}(\mu)\big|.$$}
\end{remark}

\begin{crl}
\label{n3tac_contr_crl2}
The total contribution from the boundary strata 
${\mathcal{U}}_{{\mathcal{T}}_1,{\mathcal{T}}}$ such that 
$\chi_{\mathcal{T}}(\tilde{1},\hat{1}) > 0$ to the number 
${\mathcal{C}}_{\partial\Bbb{P}(L_{\tilde{1}}\oplus L_{\tilde{1}}^*)}
({\mathcal{D}}_{\tilde{1},\hat{1}})$ is given~by
\begin{equation*}\begin{split}
\sum_{\chi_{\mathcal{T}}(\tilde{1},\hat{1})=1}   
{\mathcal{C}}_{\Bbb{P}(L_{\tilde{1}}\oplus L_{\tilde{1}}^*)|
{\mathcal{U}}_{{\mathcal{T}}_1,{\mathcal{T}}}}({\mathcal{D}}_{\tilde{1},\hat{1}})
=\big\lan c_1(L_{\tilde{1}}^*),{\overV}_{1;1}^{(1)}(\mu)\big\ran
-\big\lan 4a_{\hat{0}}+\frac{1}{2}\eta_{\hat{0},1},
{\overV}_2^{(1,1)}(\mu)\big\ran ~~&\\
+2\big|{\mathcal{V}}_2^{(1)}(\mu)\big| +2\big|{\mathcal{S}}_2(\mu)\big|
-3\big|{\mathcal{S}}_{1;1}(\mu)\big|.&
\end{split}\end{equation*}
\end{crl}

\begin{proof}
This Corollary follows immediately from 
Lemmas~\ref{n3tac_contr_lmm3} and~\ref{n3tac_contr_lmm4}.
\end{proof}

\section{Level~1 numbers}
\label{level1nums_sec}

\subsection{Evaluation of cohomology classes on 
the spaces ${\overV}_1^{(1)}(\mu)$}
\label{n3p3_cohom_subs}

In this subsection, we evaluate various tautological classes on 
the space~${\overV}_1^{(1)}(\mu)$
and compute the other level~1 numbers of Lemma~\ref{n3p2_lmm}.
We again use the computational method of Section~\ref{topology_sec},
but first we represent each cohomology class by a vector-bundle section~$s$ 
on neighborhood of ${\overU}_{{\mathcal{T}}_1}(\mu)$ in~${\overU}_{{\mathcal{T}}_1}$.
We choose this section $s$ so that it is smooth and transversal
to the zero set on all the strata of ${\overU}_{{\mathcal{T}}_1}(\mu)$,
as well as on a finite number of natural submanifolds of the strata.
We will impose additional restrictions on each given section to simplify our computations.

\begin{lmm}
\label{n3p3_cohom_lmm2}
With assumptions as in (2) of Lemma~\ref{n3p2_lmm},
\begin{equation*}\begin{split}
\lan a_{\hat{0}},{\overV}_1^{(1)}(\mu)\ran=
\blr{(2d - 6)a_{\hat{0}}^3 - 4a_{\hat{0}}^2\eta_{\hat{0},1}
 - a_{\hat{0}}\eta_{\hat{0},1}^2,{\overV}_1(\mu)}
+\blr{a_{\hat{0}}^2,{\overV}_1(\mu + H^1)} \quad&\\
+\blr{a_{\hat{0}},{\overV}_2(\mu)}.&
\end{split}\end{equation*}
\end{lmm}

\begin{proof}
(1)\qua In this case, we choose a generic hyperplane $H^2$ in $\PPP$,
instead of a section of $\ev_{\hat{0}}^*{\mathcal{O}}(1_{\PPP})$.
Let $\tilde{\mu}$ be the $\tilde{M} = [N]\sqcup\{\hat{0}\}$--tuple
of constraints in $\PPP$ given~by
$$\tilde{\mu}_l=\mu_l\quad\mbox{for all}\  l \in [N];\qquad
\tilde{\mu}_{\hat{0}}=H^2.$$
By Proposition~\ref{str_prp}, ${\mathcal{U}}_{{\mathcal{T}}_1}(\tilde{\mu})$
is a pseudovariety in ${\overU}_{\tilde{\mathcal{T}}_1}$, and thus
$$\ev_{\tilde{1}} \times \ev_{\hat{1}}\co 
{\mathcal{U}}_{\tilde{\mathcal{T}}_1}(\tilde{\mu})\lra\PPP \times \PPP$$
is a $6$--pseudocycle and determines the homology-intersection number 
\begin{equation}\label{n3p3_cohom_lmm2e1}\begin{split}
\LlRr{a_{\hat{0}},{\overV}_1^{(1)}(\mu)}
&\equiv \LlRr{{\mathcal{V}}_1^{(1)}(\tilde{\mu})}
\equiv\LlRr{
\{\ev_{\tilde{1}} \times \ev_{\hat{1}}\}^{-1}
(\De_{\PPP\times\PPP}),{\overU}_{{\mathcal{T}}_1}(\tilde{\mu})}\\
&=\sum_{r+s=3}
\LlRr{\{\ev_{\tilde{1}} \times \ev_{\hat{1}}\}^{-1}
(H^r \times H^s),{\overU}_{{\mathcal{T}}_1}(\tilde{\mu})},\\
&=2d\big\lan a_{\hat{0}}^3,{\overV}_1(\mu)\big\ran
+\big\lan a_{\hat{0}}^2,{\overV}_1(\mu + H^1)\big\ran.
\end{split}\end{equation}
By the same argument as in Subsection~\ref{n3p3_sum_subs},
\begin{equation}\label{n3p3_cohom_lmm2e2}
\big\lan a_{\hat{0}},{\overV}_1^{(1)}(\mu)\big\ran=
\big|{\mathcal{V}}_1^{(1)}(\tilde{\mu})\big|
=\LlRr{a_{\hat{0}},{\overV}_1^{(1)}(\mu)}
-{\mathcal{C}}_{\partial{\overU}_{{\mathcal{T}}_1}(\tilde{\mu})}
\big(\ev_{\tilde{1}} \times \ev_{\hat{1}},\De_{\PPP\times\PPP}\big), 
\end{equation} 
where ${\mathcal{C}}_{\partial{\overU}_{{\mathcal{T}}_1}(\tilde{\mu})}
\big(\ev_{\tilde{1}} \times \ev_{\hat{1}},\De_{\PPP\times\PPP}\big)$
is the contribution of 
$\partial{\overU}_{{\mathcal{T}}_1}(\tilde{\mu})$
to $\llrr{{\mathcal{V}}_1^{(1)}(\tilde{\mu})}$
to be computed as in Subsection~\ref{n3p3_corr_subs}.

(2)\qua If ${\mathcal{T}} = (M_1,I;j,\under{d}) < {\mathcal{T}}_1$ 
is a bubble type such that $\chi_{\mathcal{T}}(\tilde{1},\hat{1}) > 0$,
the map $\ev_{\tilde{1}} \times \ev_{\hat{1}}$
is transversal to $\De_{\PPP\times\PPP}$ 
on~${\mathcal{U}}_{{\mathcal{T}}|{\mathcal{T}}_1}(\tilde{\mu})$ by Lemma~\ref{str_lmm}.
Thus, the image of ${\mathcal{U}}_{{\mathcal{T}}|{\mathcal{T}}_1}(\tilde{\mu})$
is disjoint from $\De_{\PPP\times\PPP}$, and 
${\mathcal{U}}_{{\mathcal{T}}|{\mathcal{T}}_1}$ does not contribute to
${\mathcal{C}}_{\partial{\overU}_{{\mathcal{T}}_1}(\tilde{\mu})}
\big(\ev_{\tilde{1}} \times \ev_{\hat{1}},\De_{\PPP\times\PPP}\big)$.
Thus, from now on, we assume that $\chi_{\mathcal{T}}(\tilde{1},\hat{1}) = 0$.
Note that ${\mathcal{U}}_{{\mathcal{T}}|{\mathcal{T}}_1}(\tilde{\mu}) = \eset$
unless $|\chi_{\tilde{1}}({\mathcal{T}})| \in \{1,2\}$.

(3)\qua With appropriate identifications, ${\mathcal{U}}_{{\mathcal{T}}|{\mathcal{T}}_1}(\tilde{\mu})$
is the zero set of the section
$\ev_{{\mathcal{T}}_1,\tilde{M}}$ of the bundle 
$\ev_{{\mathcal{T}}_1,\tilde{M}}^*{\mathcal{N}}\De_{{\mathcal{T}}_1}(\tilde{\mu})$ 
over an open neighborhood of 
${\mathcal{U}}_{{\mathcal{T}}|{\mathcal{T}}_1}(\tilde{\mu})$
in~${\mathcal{U}}_{{\mathcal{T}}_1,{\mathcal{T}}}$.
By Lemma~\ref{str_lmm}, this section is transversal to the zero set.
By Proposition~\ref{str_prp}, there exists 
a $C^1$--negligible map 
$$\ve_-\co \mathcal{FT}_{\de} - Y(\mathcal{FT};\hat{I}^+)
\lra\ev_{{\mathcal{T}}_1,\tilde{M}}^*{\mathcal{N}}\De_{{\mathcal{T}}_1}(\tilde{\mu})$$ 
such that
$$\ev_{{\mathcal{T}}_1,\tilde{M}}
\big(\phi_{{\mathcal{T}}_1,{\mathcal{T}}}(b;\ups)\big)=
\ev_{{\mathcal{T}},\tilde{M}}(b)+\ve_-(b;\ups)$$
for all $(b;\ups) \in \mathcal{FT}_{\de} - Y(\mathcal{FT};\hat{I}^+)$.
On the other hand, by Lemma~\ref{n3p3_str_lmm2},
\begin{gather*}
\big\{\ev_{\tilde{1}} \times \ev_{\hat{1}}\big\}
\phi_{{\mathcal{T}}_1,{\mathcal{T}}}(\ups) 
=   \sum_{h\in\chi_{\tilde{1}}({\mathcal{T}})}   
\big(y_{h;\hat{1}} - x_{\hat{1};h}\big)^{-1} \otimes 
\big\{{\mathcal{D}}_{{\mathcal{T}},h}^{(1)} + \ve_h(\ups)\big\}\rho_h(\ups),\\
\hbox{where}\qquad
\rho_h(\ups)=   
\prod_{i\in(i_{\mathcal{T}}(h,\hat{2}),h]}          \ups_i,
\end{gather*}
for all $\ups \in \mathcal{FT}_{\de}$.
Since the linear map,
$${\mathcal{F}}\equiv\bigoplus_{h\in\chi_{\tilde{1}}({\mathcal{T}})}L_h{\mathcal{T}},
\qquad \ups\lra\sum_{h\in\chi_{\tilde{1}}({\mathcal{T}})}
{\mathcal{D}}_{{\mathcal{T}},h}^{(1)}\ups_h,$$
is injective over ${\mathcal{U}}_{{\mathcal{T}}|{\mathcal{T}}_1}(\tilde{\mu})$
by Lemma~\ref{str_lmm},
${\mathcal{U}}_{{\mathcal{T}}|{\mathcal{T}}_1}$ is 
$\big(\ev_{\tilde{1}} \times \ev_{\hat{1}},\De_{\PPP\times\PPP}\big)$--hollow
unless $\chi_{\tilde{1}}({\mathcal{T}}) = \hat{I}^+$.
If $\chi_{\tilde{1}}({\mathcal{T}}) = \hat{I}^+$, by 
Proposition~\ref{euler_prp}, decomposition~\eqref{cart_split},
and a rescaling of the linear map,
\begin{gather*}
{\mathcal{C}}_{{\mathcal{U}}_{{\mathcal{T}}|\tilde{\mathcal{T}}_1}}
\big(\ev_{\tilde{1}} \times \ev_{\hat{1}};\De_{\PPP\times\PPP}\big)
=N(\al_1),\qquad\hbox{where}\\
\al_1 \in\Ga\big(\ov{\frak M}_{\{\tilde{1}\}\sqcup\chi_{\tilde{1}}({\mathcal{T}})\sqcup
                                   M_{\tilde{1}}{\mathcal{T}}} \times 
{\mathcal{U}}_{{\overT}}(\tilde{\mu});
\pi_2^*\hbox{Hom}({\mathcal{F}},\ev_{\hat{0}}^*T\PPP)\big),~~
\al_1(\ups)= \!\!\!\sum_{h\in\chi_{\tilde{1}}({\mathcal{T}})}\!\!\!
{\mathcal{D}}_{{\overT},h}^{(1)}\ups_h.
\end{gather*}
Since the linear map~$\al_1$ comes entirely from the second factor,
$N(\al_1) = 0$ unless 
$|\chi_{\tilde{1}}({\mathcal{T}})| + |M_{\tilde{1}}{\mathcal{T}}| = 2$.

(4)\qua Thus, we only need to consider the case 
$|\chi_{\tilde{1}}({\mathcal{T}})| = 1$ and 
$M_{\tilde{1}}{\mathcal{T}} = \{\hat{1}\}$ and to compute the number~$N(\al_1)$,
where
$$\al_1 = {\mathcal{D}}_{{\mathcal{T}}_0,\tilde{1}}^{(1)}
\in\Ga\big({\mathcal{U}}_{{\mathcal{T}}_0}(\tilde{\mu});
\hbox{Hom}(L_{\tilde{1}},\ev_{\hat{0}}^*T\PPP)\big).$$
Since $\al_1$ does not vanish on 
${\mathcal{U}}_{{\mathcal{T}}_0}(\tilde{\mu})$ by Lemma~\ref{str_lmm},
by Propositions~\ref{zeros_prp} and~\ref{euler_prp},
\begin{equation}\label{n3p3_cohom_lmm2e5}
N(\al_1)=\big\lan 6a_{\hat{0}}^2 + 
4a_{\hat{0}}c_1(L_{\tilde{1}}^*) + c_1^2(L_{\tilde{1}}^*),
{\overU}_{{\mathcal{T}}_0}(\tilde{\mu})\big\ran
-{\mathcal{C}}_{\partial{\overU}_{\tilde{\mathcal{T}}_0}(\tilde{\mu})}
(\al_1^{\perp}),
\end{equation}
where $\al_1^{\perp}$ denotes the composition of $\al_1$
with the projection $\pi_{\bar{\nu}_1}^{\perp}$
onto the quotient ${\mathcal{O}}_1$ of $\ev_{\hat{0}}^*T\PPP$
by a generic trivial line subbundle~$\Bbb{C}\bar{\nu}_1$.

(5)\qua If ${\mathcal{T}} = (M_0,I;j,\under{d}) < {\mathcal{T}}_0$ 
is a bubble type such that $d_{\tilde{1}} > 0$,
the section~$\al_1^{\perp}$ does not vanish over 
${\mathcal{U}}_{{\mathcal{T}}|{\mathcal{T}}_0}(\tilde{\mu})$ by Lemma~\ref{str_lmm},
if $\bar{\nu}_1$ is generic.
Thus, ${\mathcal{U}}_{{\mathcal{T}}|{\mathcal{T}}_0}$ does not contribute to
${\mathcal{C}}_{\partial{\overU}_{{\mathcal{T}}_0}(\tilde{\mu})}
(\al_1^{\perp})$.
If $d_{\tilde{1}} = 0$, by (3b) of Proposition~\ref{str_prp},
\begin{gather*}
{\mathcal{D}}_{{\mathcal{T}}_0,\tilde{1}}^{(1)}
\phi_{{\mathcal{T}}_0,{\mathcal{T}}}(\ups)
=\sum_{h\in\chi_{\tilde{1}}({\mathcal{T}})}
\big\{{\mathcal{D}}_{{\mathcal{T}},h}^{(1)} + 
\ve_{{\mathcal{T}},\hat{1};h}^{(1;1)}(\ups)\big\}
\rho_{{\mathcal{T}},\tilde{1};h}^{(1;1)}(\ups),\\   
\hbox{where}\qquad
\rho_{{\mathcal{T}},\tilde{1};h}^{(1;1)}(\ups)=
\prod_{i\in(\tilde{1},h]}  \ups_i,
\end{gather*}
for all $\ups \in \mathcal{FT}$ sufficiently small.
Thus, as before, we conclude that 
${\mathcal{U}}_{{\mathcal{T}}|{\mathcal{T}}_0}(\tilde{\mu})$
is~$\al_1^{\perp}$--hollow unless $\hat{I}^+ = \chi_{\tilde{1}}({\mathcal{T}})$.
If $\hat{I}^+ = \chi_{\tilde{1}}({\mathcal{T}})$, either
${\mathcal{T}} = {\mathcal{T}}_0(l)$ for some $l \in [N]$
or $|\hat{I}^+| = |\chi_{\tilde{1}}({\mathcal{T}})| = 2$
and $|M_{\tilde{1}}{\mathcal{T}}| = 0$.
Thus, by Proposition~\ref{euler_prp} and
decomposition~\eqref{cart_split},
\begin{gather*}
{\mathcal{C}}_{{\mathcal{U}}_{{\mathcal{T}}|{\mathcal{T}}_0}}(\al_1^{\perp})
=N(\al_2),\qquad\hbox{where}\qquad
\al_2 \in 
\Ga\big({\mathcal{U}}_{{\overT}}(\tilde{\mu});
\hbox{Hom}({\mathcal{F}}_2,{\mathcal{O}}_1)\big),\\\
{\mathcal{F}}_1=\bigoplus_{h\in\chi_{\tilde{1}}({\mathcal{T}})}L_h,\quad
{\mathcal{O}}_1=\ev_{\hat{0}}^*T\PPP/\Bbb{C}\bar{\nu}_1,\quad
\al_2(\ups)=\pi_{\bar{\nu}_1}^{\perp}\sum_{h\in\chi_{\tilde{1}}({\mathcal{T}})}
{\mathcal{D}}_{{\overT},h}^{(1)}\ups_h.
\end{gather*}
In either case, $\al_2$ does not vanish on
${\overU}_{{\overT}}(\tilde{\mu})$, and
thus, by Propositions~\ref{zeros_prp} and~\ref{euler_prp},
\begin{equation}\label{n3p3_cohom_lmm2e7}
{\mathcal{C}}_{\partial{\overU}_{{\mathcal{T}}_0}(\tilde{\mu})}
(\al_1^{\perp})=
\sum_{l\in[N]}\big\lan 4a_{\hat{0}} + c_1(L_{\tilde{1}}^*),
{\overU}_{{\mathcal{T}}_0/l}(\tilde{\mu})\big\ran
+\big|{\mathcal{V}}_2(\tilde{\mu})\big|.
\end{equation}
The lemma follows from equations
\eqref{n3p3_cohom_lmm2e1}--\eqref{n3p3_cohom_lmm2e7}
and by using~\eqref{psi_class1} and~\eqref{psi_class2}.
\end{proof}

\begin{lmm}
\label{n3p3_cohom_lmm3}
With assumptions as in (3) of Lemma~\ref{n3p2_lmm},
$$\lan a_{\hat{0}}^2,{\overV}_1^{(1)}(\mu)\ran=
2\blr{a_{\hat{0}}^3,{\overV}_1(\mu + H^1)}
-\blr{4a_{\hat{0}}^3\eta_{\hat{0},1} + 
a_{\hat{0}}^2\eta_{\hat{0},1}^2,{\overV}_1(\mu)}
+\blr{a_{\hat{0}}^2,{\overV}_2(\mu)}.$$
\end{lmm}

\begin{proof}
The proof is nearly identical to that of Lemma~\ref{n3p3_cohom_lmm2}.
\end{proof}

\begin{lmm}
\label{n3p3_cohom_lmm4}
With assumptions as in (3) of Lemma~\ref{n3p2_lmm},
\begin{equation*}\begin{split}
\blr{a_{\hat{0}}\eta_{\hat{0},1},{\overV}_1^{(1)}(\mu)}
=\blr{a_{\hat{0}}\eta_{\hat{0},1},{\overV}_1(\mu + H^0)}
+\blr{a_{\hat{0}}^2\eta_{\hat{0},1},{\overV}_1(\mu + H^1)} \qquad&\\
+d\blr{a_{\hat{0}}^3\eta_{\hat{0},1},{\overV}_1(\mu)} 
-\blr{4a_{\hat{0}}^2 + a_{\hat{0}}\eta_{\hat{0},1},{\overV}_2(\mu)}.&
\end{split}\end{equation*}
\end{lmm}

\begin{proof}
(1)\qua Let $\tilde{M} = M_0\sqcup\{\hat{0}\}$ and let 
$\tilde{\mu}$ be the $\tilde{M}$--tuple of constraints given by
$\tilde{\mu}_l = \mu_l$ for all $l \in [N]$ and $\tilde{\mu}_{\hat{0}} = H^2$.
If $s$ is a section of the bundle ${\mathcal{L}}_{\tilde{1}}^*$ over
a neighborhood of ${\overU}_{{\mathcal{T}}_1}(\tilde{\mu})$ in
${\overU}_{{\mathcal{T}}_1}$ such that $s$ is transversal 
to the zero set on all smooth strata of ${\mathcal{U}}_{{\mathcal{T}}_1}(\mu)$,
the map
$$\ev_{\tilde{1}} \times \ev_{\hat{1}}\co 
s^{-1}(0)\cap {\mathcal{U}}_{{\mathcal{T}}_1}(\tilde{\mu}) \lra 
\PPP \times \PPP$$
is a 6-pseudocycle.
In particular, we have a well-defined intersection number,
\begin{equation}\label{n3p3_cohom_lmm4e1}\begin{split}
\LlRr{a_{\hat{0}}c_1({\mathcal{L}}_{\tilde{1}}^*),{\overV}_1^{(1)}(\mu)}
&\equiv \LlRr{s^{-1}(0)\cap{\mathcal{V}}_1^{(1)}(\tilde{\mu})}\\
&\equiv\LlRr{
\{\ev_{\tilde{1}} \times \ev_{\hat{1}}\}^{-1}
(\De_{\PPP\times\PPP}),
s^{-1}(0)\cap{\overU}_{{\mathcal{T}}_1}(\tilde{\mu})}\\
&=\big\lan a_{\hat{0}}c_1({\mathcal{L}}_{\tilde{1}}^*),
{\overV}_1(\mu + H^0)\big\ran
+\big\lan a_{\hat{0}}^2c_1({\mathcal{L}}_{\tilde{1}}^*),
{\overV}_1(\mu + H^1)\big\ran\\
& \qquad\qquad\qquad\qquad\qquad\qquad~~
+d\big\lan a_{\hat{0}}^3c_1({\mathcal{L}}_{\tilde{1}}^*),
{\overV}_1(\mu)\big\ran.
\end{split}\end{equation}
As before,
\begin{equation}\label{n3p3_cohom_lmm4e2}\begin{split}
&\blr{a_{\hat{0}}c_1({\mathcal{L}}_{\tilde{1}}^*),{\overV}_1^{(1)}(\mu)}
=\, ^{\pm}\big|{s^{-1}(0)\cap\cal V}_1^{(1)}(\tilde{\mu})\big|\\
&\qquad\qquad =\LlRr{a_{\hat{0}}c_1({\mathcal{L}}_{\tilde{1}}^*),{\overV}_1^{(1)}(\mu)}
-{\mathcal{C}}_{\partial{\overU}_{{\mathcal{T}}_1}(\tilde{\mu})}
\big(\ev_{\tilde{1}} \times \ev_{\hat{1}},\De_{\PPP\times\PPP}\big), 
\end{split}\end{equation} 
where ${\mathcal{C}}_{\partial{\overU}_{{\mathcal{T}}_1}(\tilde{\mu})}
\big(\ev_{\tilde{1}} \times \ev_{\hat{1}},\De_{\PPP\times\PPP}\big)$
is the contribution of 
$s^{-1}(0)\cap\partial{\overU}_{{\mathcal{T}}_1}(\tilde{\mu})$
to $\llrr{s^{-1}(0)\cap{\mathcal{V}}_1^{(1)}(\tilde{\mu})}$.

(2)\qua If ${\mathcal{T}} = (M_1,I;j,\under{d}) < {\mathcal{T}}_1$ 
is a bubble type such that $\chi_{\mathcal{T}}(\tilde{1},\hat{1}) > 0$,
the map $\ev_{\tilde{1}} \times \ev_{\hat{1}}$
is transversal to $\De_{\PPP\times\PPP}$ 
on~${\mathcal{U}}_{{\mathcal{T}}|{\mathcal{T}}_1}(\tilde{\mu})$ by Lemma~\ref{str_lmm}.
Thus, if $s$ is chosen to be transversal to the zero set
on the set 
$\{\ev_{\tilde{1}} \times \ev_{\hat{1}}\}^{-1}(\De_{\PPP\times\PPP})
\cap{\mathcal{U}}_{{\mathcal{T}}|{\mathcal{T}}_1}(\tilde{\mu})$,
$$s^{-1}(0)\cap\{\ev_{\tilde{1}}
 \times \ev_{\hat{1}}\}^{-1}(\De_{\PPP\times\PPP})
\cap{\mathcal{U}}_{{\mathcal{T}}|{\mathcal{T}}_1}(\tilde{\mu})
=\eset,$$
and ${\mathcal{U}}_{{\mathcal{T}}|{\mathcal{T}}_1}$ does not contribute
to ${\mathcal{C}}_{\partial{\overU}_{{\mathcal{T}}_1}(\tilde{\mu})}
\big(\ev_{\tilde{1}} \times \ev_{\hat{1}},\De_{\PPP\times\PPP}\big)$.
If  $\chi_{\mathcal{T}}(\tilde{1},\hat{1}) = 0$, it can be assumed that
$s$ is transversal to the zero set on the submanifold 
${\mathcal{S}}_{{\mathcal{T}}|{\mathcal{T}}_1}(\mu)$, ie that
$s^{-1}(0)\cap{\mathcal{S}}_{{\mathcal{T}}|{\mathcal{T}}_1}(\mu) = \eset$.
Then, as in the proofs of Lemma~\ref{n3p3_cohom_lmm2} and~\ref{n3p3_cohom_lmm3}
we can conclude that 
$s^{-1}(0)\cap{\mathcal{U}}_{{\mathcal{T}}|{\mathcal{T}}_1}(\mu)$ is
$\big(\ev_{\tilde{1}} \times \ev_{\hat{1}},\De_{\PPP\times\PPP}\big)$--hollow
unless $\chi_{\tilde{1}}({\mathcal{T}}) = \hat{I}^+$.

(3)\qua We can also assume that section $\phi_{{\mathcal{T}}_1,{\mathcal{T}}}^*s$ 
is constant along the fibers of the bundle~$\mathcal{FT}$
over an open subset ${\mathcal{K}}_{\mathcal{T}}$ of ${\mathcal{U}}_{{\mathcal{T}}|{\mathcal{T}}_1}$
that contains all of the finitely many zeros of the map affine map
$${\mathcal{F}}\equiv\bigoplus_{h\in\chi_{\tilde{1}}({\mathcal{T}})}L_h
\lra\ev_{\hat{0}}^*T\PPP,\qquad
(b,\ups)\lra \nu(b)+
\sum_{h\in\chi_{\tilde{1}}({\mathcal{T}})}{\mathcal{D}}_{{\mathcal{T}},h}^{(1)}\ups_h,$$
over $s^{-1}(0)\cap{\mathcal{U}}_{{\mathcal{T}}|{\mathcal{T}}_1}(\mu)$
for a generic section 
$\nu \in \Ga({\overU}_{{\mathcal{T}}|{\mathcal{T}}_1}(\mu);\ev_{\hat{0}}^*T\PPP)$.
Then, as before,
\begin{gather*}
{\mathcal{C}}_{{\mathcal{U}}_{{\mathcal{T}}|{\mathcal{T}}_1}}
\big(\ev_{\tilde{1}} \times \ev_{\hat{1}};\De_{\PPP\times\PPP}\big)
=N(\al_1),\qquad\hbox{where}\\
\al_1\in\Ga\big(s^{-1}(0)\cap(\ov{\frak M}_{\hat{1},\chi_{\tilde{1}}({\mathcal{T}}),
                                   M_{\tilde{1}}{\mathcal{T}}} \times 
{\mathcal{U}}_{{\overT}}(\tilde{\mu}));
\hbox{Hom}({\mathcal{F}},\ev_{\hat{0}}^*T\PPP)\big),\\
\al_1(\ups)=\sum_{h\in\chi_{\tilde{1}}({\mathcal{T}})}    
{\mathcal{D}}_{{\mathcal{T}},h}^{(1)}\ups_h.
\end{gather*}
(4)\qua If
$c_1({\mathcal{L}}_{\tilde{1}}^*)|{\overU}_{{\mathcal{T}}|{\mathcal{T}}_1}(\mu)
 = 0$, we can choose $s$ so that 
$$s^{-1}(0)\cap{\overU}_{{\mathcal{T}}|{\mathcal{T}}_1}(\mu)=\eset.$$
Thus, we only need to compute contributions $N(\al_1)$
from the strata ${\mathcal{U}}_{{\mathcal{T}}|{\mathcal{T}}_1}(\mu)$
to which $c_1({\mathcal{L}}_{\tilde{1}}^*)$ restricts non-trivially.
By dimension-counting, $|\chi_{\tilde{1}}({\mathcal{T}})| \in \{1,2\}$
if ${\mathcal{U}}_{{\mathcal{T}}|{\mathcal{T}}_1}(\mu) \neq \eset$.
If $|\chi_{\tilde{1}}({\mathcal{T}})| = 1$, either 
$M_{\tilde{1}}{\mathcal{T}} = \{\hat{1}\}$ or
$M_{\tilde{1}}{\mathcal{T}} = \{\hat{1},l\}$ for some $l \in [N]$.
In either case, $c_1({\mathcal{L}}_{\tilde{1}}^*)$ restricts trivially
to ${\mathcal{U}}_{{\mathcal{T}}|{\mathcal{T}}_1}(\mu)$.
On the other hand, if $|\chi_{\tilde{1}}({\mathcal{T}})| = 2$
and $M_{\tilde{1}}{\mathcal{T}} = \{\hat{1},l\}$ for some $l \in [N]$,
$N(\al_1) = 0$, because the second factor in 
the decomposition~\eqref{cart_split} is a finite set of points,
while the map $\al_1$ comes entirely from the second factor.
In the remaining case, ie $|\chi_{\tilde{1}}({\mathcal{T}})| = 2$
and $M_{\tilde{1}}{\mathcal{T}} = \{\hat{1}\}$,
$c_1({\mathcal{L}}_{\tilde{1}}^*)$ is the pullback of the poincare dual of a point
by the projection map~$\pi_1$.
Thus, in this case,
\begin{gather*}
{\mathcal{C}}_{{\mathcal{U}}_{{\mathcal{T}}|{\mathcal{T}}_1}}
\big(\ev_{\tilde{1}} \times \ev_{\hat{1}};\De_{\PPP\times\PPP}\big)
=N(\al_1),\qquad\hbox{where}\\
\al_1\in\Ga\big({\mathcal{U}}_{{\overT}}(\tilde{\mu});
\hbox{Hom}({\mathcal{F}},\ev_{\hat{0}}^*T\PPP)\big),\quad
\al_1(\ups)=\sum_{h\in\chi_{\tilde{1}}({\mathcal{T}})}
{\mathcal{D}}_{{\overT},h}^{(1)}\ups_h.
\end{gather*}
Since $\al_1$ has full rank on all of 
${\mathcal{U}}_{{\mathcal{T}}|{\mathcal{T}}_1}(\tilde{\mu})$,
by Propositions~\ref{zeros_prp} and~\ref{euler_prp},
\begin{equation}\label{n3p3_cohom_lmm4e5}
N(\al_1)=\big\lan 4a_{\hat{0}} + \big(
c_1(L_{h_1}^*) + c_1(L_{h_2}^*)\big),
{\overU}_{{\overT}}(\tilde{\mu})\big\ran
-{\mathcal{C}}_{\Bbb{P}{\mathcal{F}}|\partial
{\overU}_{{\overT}}(\tilde{\mu})}(\tilde{\al}_1^{\perp}),
\end{equation}
where $\tilde{\al}_1^{\perp}$ denotes the composition
of the linear map 
$$\tilde{\al}_1 \in \Ga\big(\Bbb{P}{\mathcal{F}};
\hbox{Hom}(\ga_{\mathcal{F}},\pi_{\Bbb{P}{\mathcal{F}}}^*\ev_{\hat{0}}^*T\PPP)\big)$$
with the projection onto the quotient of 
$\pi_{\Bbb{P}{\mathcal{F}}}^*\ev_{\hat{0}}^*T\PPP$ by
a generic trivial line subbundle $\Bbb{C}\bar{\nu}_1$.
If 
$${\mathcal{T}}' \equiv (M_0,I';j',\under{d}')<{\overT}$$
is a bubble type such that $\tilde{\al}_1$ vanishes somewhere on 
$\Bbb{P}{\mathcal{F}}|{\mathcal{U}}_{{\mathcal{T}}'|{\overT}}(\tilde{\mu})$,
${\mathcal{T}}' = {\overT}(l)$ for some $l \in [N]$ and 
$$\tilde{\al}_1^{-1}(0)\cap 
\Bbb{P}{\mathcal{F}}|
{\mathcal{U}}_{{\mathcal{T}}'|{\overT}}(\tilde{\mu})
=\big\{(b,[\ups_{h_1},\ups_{h_2}]) :
b \in {\mathcal{U}}_{{\mathcal{T}}'|{\overT}}(\tilde{\mu}),
~\ups_{j_l'} = 0\big\},$$
as can be seen from Lemma~\ref{str_lmm}.
{}From Proposition~\ref{str_prp}, we then conclude that
$${\mathcal{C}}_{\Bbb{P}{\mathcal{F}}|
{\mathcal{U}}_{{\mathcal{T}}'|{\overT}}(\tilde{\mu})}
(\tilde{\al}_1^{\perp})=
\big|{\mathcal{U}}_{{\mathcal{T}}'|{\overT}}(\tilde{\mu})\big|.$$
Thus, summing equation~\eqref{n3p3_cohom_lmm4e5} over all
bubble type~${\mathcal{T}}$ and using~\eqref{psi_class1}
and~\eqref{psi_class2}, we~obtain
\begin{equation}\label{n3p3_cohom_lmm4e7}
{\mathcal{C}}_{\partial{\overU}_{{\mathcal{T}}_1}(\tilde{\mu})}
\big(\ev_{\tilde{1}} \times \ev_{\hat{1}},\De_{\PPP\times\PPP}\big)
=\big\lan 4a_{\hat{0}} + \big(
c_1({\mathcal{L}}_{\tilde{1}}^*) + c_1({\mathcal{L}}_{\tilde{2}}^*)\big),
{\overV}_2(\tilde{\mu})\big\ran.
\end{equation}
The claim follows from equations~\eqref{n3p3_cohom_lmm4e1},
\eqref{n3p3_cohom_lmm4e2}, and~\eqref{n3p3_cohom_lmm4e7}.
\end{proof}

\begin{lmm}
\label{n3p3_cohom_lmm5}
With assumptions as in (3) of Lemma~\ref{n3p2_lmm},
\begin{equation*}\begin{split}
\blr{\eta_{\hat{0},1}^2,{\overV}_1^{(1)}(\mu)}
=\blr{\eta_{\hat{0},1}^2,{\overV}_1(\mu + H^0)}
+\blr{a_{\hat{0}}\eta_{\hat{0},1}^2,{\overV}_1(\mu + H^1)} \qquad\quad&\\
+\blr{4a_{\hat{0}}^3\eta_{\hat{0},1}
+ d\cdot a_{\hat{0}}^2\eta_{\hat{0},1}^2,{\overV}_1(\mu)}
- \big|{\mathcal{V}}_3(\mu)\big|.&
\end{split}\end{equation*}
\end{lmm}

\begin{proof}
(1)\qua We proceed as in the proof of Lemma~\ref{n3p3_cohom_lmm4}.
Let $s$ be a section of 
\hbox{${\mathcal{L}}_{\tilde{1}}^* \oplus {\mathcal{L}}_{\tilde{1}}^*$}
with good properties.
Then, we have a well-defined homology intersection number
\begin{equation}\label{n3p3_cohom_lmm5e1}\begin{split}
&\LlRr{c_1^2({\mathcal{L}}_{\tilde{1}}^*),{\overV}_1^{(1)}(\mu)}
\equiv \LlRr{s^{-1}(0)\cap{\mathcal{V}}_1^{(1)}(\mu)} \\
&\qquad\qquad \equiv \LlRr{
\{\ev_{\tilde{1}} \times \ev_{\hat{1}}\}^{-1}
(\De_{\PPP\times\PPP}),
s^{-1}(0)\cap{\overU}_{{\mathcal{T}}_1}(\mu)}\\
&\qquad\qquad = \sum_{q+r=3}\big\lan 
c_1^2({\mathcal{L}}_{\tilde{1}}^*)a_{\hat{0}}^q a_{\hat{1}}^r,
{\overU}_{{\mathcal{T}}_1}(\mu)\big\ran\\
&\qquad\qquad =\blr{c_1^2({\mathcal{L}}_{\tilde{1}}^*),{\overV}_1(\mu + \{H^0\})} 
+\blr{a_{\hat{0}}c_1^2({\mathcal{L}}_{\tilde{1}}^*),{\overV}_1(\mu + \{H^1\})}\\
&\qquad\qquad\qquad
+d\cdot\blr{a_{\tilde{1}}^2c_1^2({\mathcal{L}}_{\tilde{1}}^*),{\overV}_1(\mu)}
+4\blr{a_{\hat{0}}^3c_1({\mathcal{L}}_{\tilde{1}}^*),{\overV}_1(\mu + \{H^0\})}.
\end{split}\end{equation}
A little care is required to obtain the last equality above.
For example, note~that
$$\blr{c_1^2({\mathcal{L}}_{\tilde{1}}^*)a_{\hat{1}}^3,{\overU}_{{\mathcal{T}}_1}(\mu)}
=\blr{c_1^2({\mathcal{L}}_{\tilde{1}}^*),{\overV}_1(\mu + \{H^0\})}
+\blr{a_{\hat{0}}^3c_1({\mathcal{L}}_{\tilde{1}}^*),{\overV}_1(\mu + \{H^0\})},$$
with our definitions; see~\eqref{psi_class1}. 
As before, 
\begin{equation}\label{n3p3_cohom_lmm5e2}\begin{split}
\blr{c_1^2({\mathcal{L}}_{\tilde{1}}^*),{\overV}_1^{(1)}(\mu)}
&=\, ^{\pm}\big|{s^{-1}(0)\cap\cal V}_1^{(1)}(\mu)\big|\\
&=\LlRr{c_1^2({\mathcal{L}}_{\tilde{1}}^*),{\overV}_1^{(1)}(\mu)}
-{\mathcal{C}}_{\partial{\overU}_{{\mathcal{T}}_1}(\mu)}
\big(\ev_{\tilde{1}} \times \ev_{\hat{1}},\De_{\PPP\times\PPP}\big), 
\end{split}\end{equation} 
where ${\mathcal{C}}_{\partial{\overU}_{{\mathcal{T}}_1}(\mu)}
\big(\ev_{\tilde{1}} \times \ev_{\hat{1}},\De_{\PPP\times\PPP}\big)$
is the contribution of 
$s^{-1}(0)\cap\partial{\overU}_{{\mathcal{T}}_1}(\mu)$
to $\llrr{s^{-1}(0)\cap{\mathcal{V}}_1^{(1)}(\mu)}$.

(2)\qua If ${\mathcal{T}} = (M_1,I;j,\under{d}) < {\mathcal{T}}_1$ 
is a bubble type such that $\chi_{\mathcal{T}}(\tilde{1},\hat{1}) > 0$,
as in the proof of Lemma~\ref{n3p3_cohom_lmm4},
the space ${\mathcal{U}}_{{\mathcal{T}}|{\mathcal{T}}_1}$ does not contribute to
$${\mathcal{C}}_{\partial{\overU}_{{\mathcal{T}}_1}(\mu)}
\big(\ev_{\tilde{1}} \times \ev_{\hat{1}},\De_{\PPP\times\PPP}\big).$$
If $\chi_{\mathcal{T}}(\tilde{1},\hat{1}) = 0$, but
$\hat{I}^+ \neq \chi_{\tilde{1}}({\mathcal{T}})$,
${\mathcal{U}}_{{\mathcal{T}}|{\mathcal{T}}_1}(\mu)$ is
$\big(\ev_{\tilde{1}} \times \ev_{\hat{1}},\De_{\PPP\times\PPP}\big)$--hollow
and again does not contribute to 
${\mathcal{C}}_{\partial{\overU}_{{\mathcal{T}}_1}(\mu)}
\big(\ev_{\tilde{1}} \times \ev_{\hat{1}},\De_{\PPP\times\PPP}\big)$.
If $\hat{I}^+ = \chi_{\tilde{1}}({\mathcal{T}})$, by dimension-counting
\begin{gather*}
|\chi_{\tilde{1}}({\mathcal{T}})| \in \{1,2\},\quad
M_{\tilde{1}}{\mathcal{T}}=\{\hat{1}\}~~\hbox{or}~~
M_{\tilde{1}}{\mathcal{T}}=\{\hat{1},l\}~\hbox{for some}~l \in [N],\\
\hbox{OR}\qquad
|\chi_{\tilde{1}}({\mathcal{T}})| = 3,\quad
M_{\tilde{1}}{\mathcal{T}}=\{\hat{1}\}.
\end{gather*}
In all cases, but the last, $c_1^2({\mathcal{L}}_{\tilde{1}}^*)$ restricts trivially
to ${\overU}_{{\mathcal{T}}|{\mathcal{T}}_1}(\mu)$.
If $|\chi_{\tilde{1}}({\mathcal{T}})| = 3$ and
$M_{\tilde{1}}{\mathcal{T}} = \{\hat{1}\}$, under 
the decomposition~\eqref{cart_split},
$c_1^2({\mathcal{L}}_{\tilde{1}}^*)$ is the pullback of the poincare dual 
of a point
by the projection map onto the second factor.
Thus, similarly to the proof of Lemma~\ref{n3p3_cohom_lmm4},
\begin{gather*}
{\mathcal{C}}_{{\mathcal{U}}_{{\mathcal{T}}|{\mathcal{T}}_1}}
\big(\ev_{\tilde{1}} \times \ev_{\hat{1}};\De_{\PPP\times\PPP}\big)
=N(\al_1),\qquad\hbox{where}\\
\al_1\in\Ga\big({\mathcal{U}}_{{\overT}}(\mu);
\hbox{Hom}({\mathcal{F}},\ev_{\hat{0}}^*T\PPP)\big),\quad
{\mathcal{F}} = \bigoplus_{h\in\chi_{\tilde{1}}({\mathcal{T}})}   L_h,\quad
\al_1(\ups)=\sum_{h\in\chi_{\tilde{1}}({\mathcal{T}})}    
{\mathcal{D}}_{{\overT},h}^{(1)}\ups_h.
\end{gather*}
Since $\al_1$ is an isomorphism on every fiber of ${\mathcal{F}}$
over the finite set~${\mathcal{U}}_{{\mathcal{T}}/M_{\tilde{1}}{\mathcal{T}}}(\mu)$,
$N(\al_1) = \big|{\mathcal{U}}_{{\overT}}(\mu)\big|$.
Thus, over all
bubble types~${\mathcal{T}}$ as above, we~obtain
\begin{equation}\label{n3p3_cohom_lmm5e7}
{\mathcal{C}}_{\partial{\overU}_{{\mathcal{T}}_1}(\mu)}
\big(\ev_{\tilde{1}} \times \ev_{\hat{1}},\De_{\PPP\times\PPP}\big)
=\big|{\mathcal{V}}_3(\mu)\big|.
\end{equation}
The claim follows from equations~\eqref{n3p3_cohom_lmm5e1},
\eqref{n3p3_cohom_lmm5e2}, and~\eqref{n3p3_cohom_lmm5e7}.
\end{proof}

\subsection{Other level~1 numbers}
\label{other_level1_subs}

In this subsection, we compute the level~1 numbers
of Lemmas~\ref{n3circ3_lmm} and~\ref{n3circ2_lmm}
and thus conclude the computation of
the enumerative numbers of Theorems~\ref{n3p3_thm} and~\ref{n3tac_thm}.

\begin{lmm}
\label{n3circ3_l}
With assumptions as in Lemma~\ref{n3circ3_lmm},
\begin{equation*}\begin{split}
\big|{\mathcal{V}}_2^{(1)}(\mu)\big|=
\big|{\mathcal{V}}_2(\mu + H^0)\big|+
\big\lan a_{\hat{0}},{\overV}_2(\mu + H^1)\big\ran+
3\big|{\mathcal{V}}_3(\mu)\big| \qquad\qquad&\\
-\big\lan (12 - d)a_{\hat{0}}^2 + 4a_{\hat{0}}\eta_{\hat{0},1}
 + 2\eta_{\hat{0},2} - \eta_{\hat{0},1}^2,{\overV}_2(\mu)\big\ran.&
\end{split}\end{equation*}
\end{lmm}

\begin{proof}
(1)\qua  By definition and the usual argument,
\begin{equation}\label{n3circ3_l_e1}\begin{split}
\big|{\mathcal{V}}_2^{(1)}(\mu)\big|
&=\sum_{\tilde{\mathcal{T}}}\Big(
\LlRr{\{\ev_{\tilde{1}} \times \ev_{\hat{1}}\}^{-1}(\De_{\PPP\times\PPP}),
{\overU}_{\tilde{\mathcal{T}}}(\mu)}\\
&\qquad\qquad\qquad\qquad
-{\mathcal{C}}_{\partial{\overU}_{\tilde{\mathcal{T}}}(\mu)}
\big(\ev_{\tilde{1}} \times \ev_{\hat{1}};\De_{\PPP\times\PPP}\big)\Big)\\
&=\big|{\mathcal{V}}_2(\mu + H^0)\big|+
\big\lan a_{\hat{0}},{\overV}_2(\mu + H^1)\big\ran+
d\big\lan a_{\hat{0}}^2,{\overV}_2(\mu)\big\ran\\
&\qquad\qquad\qquad\qquad -\sum_{\tilde{\mathcal{T}}}
{\mathcal{C}}_{\partial{\overU}_{\tilde{\mathcal{T}}}(\mu)}
\big(\ev_{\tilde{1}} \times \ev_{\hat{1}};\De_{\PPP\times\PPP}\big),
\end{split}\end{equation}
where the union is taken over all bubble types 
$\tilde{\mathcal{T}} = (M_1,I_2;\tilde{j},\tilde{d})$
such that 
$$\tilde{j}_{\hat{1}}=\tilde{1}, \qquad
\tilde{d}_{\tilde{1}},\tilde{d}_{\tilde{2}} >0, 
\quad\hbox{and}\quad 
\tilde{d}_{\tilde{1}} + \tilde{d}_{\tilde{2}}=0.$$
Let ${\mathcal{T}} = (M_1,I;j,\under{d})$ be a bubble type such that
${\mathcal{T}} < \tilde{\mathcal{T}}$ and 
${\mathcal{U}}_{{\mathcal{T}}|\tilde{\mathcal{T}}}(\mu) \neq \eset$.
If $\chi_{\mathcal{T}}({\tilde{1},\hat{1}}) > 0$,
as in the proof of Lemma~\ref{n3p3_cohom_lmm2},
${\mathcal{U}}_{\tilde{\mathcal{T}},{\mathcal{T}}}$ does not contribute~to 
$${\mathcal{C}}_{\partial{\overU}_{\tilde{\mathcal{T}}}(\mu)}
\big(\ev_{\tilde{1}} \times \ev_{\hat{1}};\De_{\PPP\times\PPP}\big).$$
Thus, we assume, $\chi_{\mathcal{T}}({\tilde{1},\hat{1}}) = 0$.

(2)\qua  By Lemma~\ref{n3p3_str_lmm2},
\begin{gather*}
\big\{\ev_{\tilde{1}} \times \ev_{\hat{1}}\big\}
\phi_{\tilde{\mathcal{T}},{\mathcal{T}}}(\ups)
=   \sum_{h\in\chi_{\tilde{1}}({\mathcal{T}})}   
\big(y_{h;\hat{1}} - x_{\hat{1};h}\big)^{-1} \otimes 
\big\{{\mathcal{D}}_{{\mathcal{T}},h}^{(1)} + \ve_h(\ups)\big\}\rho_h(\ups),\\
\hbox{where}\qquad
\rho_h(\ups)=  
\prod_{i\in(i_{\mathcal{T}}(h,\hat{1}),h]}          \ups_i,
\end{gather*}
for all $\ups \in \mathcal{FT}_{\de}$.
Since the linear map,
$${\mathcal{F}}\equiv\bigoplus_{h\in\chi_{\tilde{1}}({\mathcal{T}})}L_h{\mathcal{T}},
\qquad \ups\lra\sum_{h\in\chi_{\tilde{1}}({\mathcal{T}})}
{\mathcal{D}}_{{\mathcal{T}},h}^{(1)}\ups_h,$$
is injective over ${\mathcal{U}}_{{\mathcal{T}}|\tilde{\mathcal{T}}}(\tilde{\mu})$
by Lemma~\ref{str_lmm}, ${\mathcal{U}}_{{\mathcal{T}}|\tilde{\mathcal{T}}_1}$ is 
$\big(\ev_{\tilde{1}} \times \ev_{\hat{1}},\De_{\PPP\times\PPP}\big)$--hollow
unless $\chi_{\tilde{1}}({\mathcal{T}}) = \hat{I}^+$.
If $\chi_{\tilde{1}}({\mathcal{T}}) = \hat{I}^+$ by dimension-counting,
\begin{gather*}
|\chi_{\tilde{1}}({\mathcal{T}})| = |\hat{I}^+| = 1
\quad\hbox{and}\quad
M_{\tilde{1}}{\mathcal{T}} = \{\hat{1}\}
~~\hbox{or}~~
M_{\tilde{1}}{\mathcal{T}} = \{\hat{1},l\}
~\hbox{for some}~l \in [N]\\
\qquad\hbox{OR}\qquad
|\chi_{\tilde{1}}({\mathcal{T}})| = |\hat{I}^+| = 2
\quad\hbox{and}\quad
M_{\tilde{1}}{\mathcal{T}} = \{\hat{1}\}.
\end{gather*}
Furthermore, by Proposition~\ref{euler_prp},
a rescaling of the linear map, and 
the decomposition~\eqref{cart_split},
\begin{gather*}
{\mathcal{C}}_{{\mathcal{U}}_{{\mathcal{T}}|\tilde{\mathcal{T}}}(\mu)}
\big(\ev_{\tilde{1}} \times \ev_{\hat{1}};\De_{\PPP\times\PPP}\big)
=N(\al_1),\qquad\hbox{where}\\
\al_1\in\Ga\big(\ov{\frak M}_{\{\tilde{1}\}\sqcup
\chi_{\tilde{1}}({\mathcal{T}})\sqcup M_{\tilde{1}}({\mathcal{T}})}
 \times {\mathcal{U}}_{{\overT}}(\mu);
\hbox{Hom}({\mathcal{F}},\ev_{\hat{0}}^*T\PPP)\big),\\
{\mathcal{F}} = \bigoplus_{h\in\chi_{\tilde{1}}({\mathcal{T}})}  L_h, \qquad
\al_1(\ups)=\sum_{h\in\chi_{\tilde{1}}({\mathcal{T}})}   
{\mathcal{D}}_{{\overT},h}^{(1)}\ups_h.
\end{gather*}
Since $\al_1$ comes entirely from the second component,
$N(\al_1) = 0$ unless the first component is zero-dimensional,
ie unless $\chi_{\tilde{1}}({\mathcal{T}}) = \{h\}$ is
a single-element set and $M_{\tilde{1}}{\mathcal{T}} = \{\hat{1}\}$.
Thus, we assume that this is the case.

(3)\qua By Propositions~\ref{zeros_prp} and~\ref{euler_prp},
$$N(\al_1)=\big\lan 6a_{\hat{0}}^2 + 4a_{\hat{0}}c_1(L_h^*)
 + c_1^2(L_h^*),{\overU}_{{\overT}}(\mu)\big\ran
-{\mathcal{C}}_{\partial{\overU}_{{\overT}}(\mu)}(\al_1^{\perp}).$$
Suppose ${\mathcal{T}}' = (M_0,I';j',\under{d}')$
is a bubble type such that ${\mathcal{T}}' < {\overT}$
and ${\overU}_{{\mathcal{T}}'|{\overT}}(\mu) \neq \eset$.
Then, ${\overU}_{{\mathcal{T}}'|{\overT}}(\mu)$ does not contribute
to ${\mathcal{C}}_{\partial{\overU}_{{\overT}}(\mu)}(\al_1^{\perp})$
unless $d_h' = 0$.
If $d_h' = 0$, by Proposition~\ref{str_prp},
$$\al_1\phi_{{\overT},{\mathcal{T}}'}(\ups) 
=   \sum_{h'\in\chi_h({\mathcal{T}}')}   
\big\{{\mathcal{D}}_{{\mathcal{T}}',h'}^{(1)} + \ve_{h'}(\ups)\big\}\rho_{h'}(\ups),
~~~\hbox{where}~~~
\rho_{h'}(\ups) =   
\prod_{i\in(h,h']}  \ups_i,$$
for all $\ups \in \mathcal{FT}_{\de}'$.
Thus, ${\overU}_{{\mathcal{T}}'|{\overT}}(\mu)$ is
$\al_1^{\perp}$--hollow unless
$\chi_h({\mathcal{T}}) = {\hat{I}'}{^+}$.
In such a case, either 
${\mathcal{T}}' = {\overT}(l)$ for some 
$l \in M_h{\overT}$ or
$|\chi_h({\mathcal{T}}')| = |{\hat{I}'}{^+}| = 2$
and $M_h{\mathcal{T}}' = \eset$.
In either case,
$${\mathcal{C}}_{{\mathcal{U}}_{{\mathcal{T}}'|{\overT}}(\mu)}\big(\al_1^{\perp}\big)
=N(\al_1),\qquad\hbox{where}\quad
\al_2\in\Ga\big({\overU}_{{\mathcal{T}}'}(\mu);
\hbox{Hom}(\mathcal{FT}',\ev_{\hat{0}}^*T\PPP/\Bbb{C}\bar{\nu}_1)\big)$$
is a nonvanishing section. 
Thus, using Proposition~\ref{zeros_prp} along with identities
\eqref{psi_class1} and~\eqref{psi_class1}, we conclude~that
\begin{equation}\label{n3circ3_l_e3}
N(\al_1)=\blr{6a_{\hat{0}}^2 + 4a_{\hat{0}}c_1({\mathcal{L}}_h^*)
 + c_1^2({\mathcal{L}}_h^*),{\overU}_{{\overT}}(\mu)}
-\sum_{|\chi_h({\mathcal{T}}')|=2}
\big|{\mathcal{U}}_{{\mathcal{T}}'|{\overT}}(\mu)\big|.
\end{equation}
Summing equation~\eqref{n3circ3_l_e3} over all bubble types $\tilde{\mathcal{T}}$,
we obtain
\begin{equation}\label{n3circ3_l_e5}\begin{split}
&\sum_{\tilde{\mathcal{T}}}
{\mathcal{C}}_{\partial{\overU}_{\tilde{\mathcal{T}}}(\mu)}
\big(\ev_{\tilde{1}} \times \ev_{\hat{1}};\De_{\PPP\times\PPP}\big)\\
& \qquad\qquad\qquad\qquad
= \blr{12a_{\hat{0}}^2 + 4a_{\hat{0}}\eta_{\hat{0},1}
 + 2\eta_{\hat{0},2} - \eta_{\hat{0},1}^2,{\overV}_2(\mu)}
-3\big|{\mathcal{V}}_3(\mu)\big|.
\end{split}\end{equation}
The claim follows from~\eqref{n3circ3_l_e1} and~\eqref{n3circ3_l_e5}.
\end{proof}

\begin{lmm}
\label{n3circ2_l1}
With assumptions as in Lemma~\ref{n3circ2_lmm},
$$\big\lan a_{\hat{0}},{\overV}_2^{(1,1)}(\mu)\big\ran=
2\big\lan a_{\hat{0}},{\overV}_2(\mu + \{H^1 : H^2\})\big\ran
-\big\lan 8a_{\hat{0}}^2 + 2a_{\hat{0}}\eta_{\hat{0},1},
{\overV}_2(\mu)\big\ran.$$
\end{lmm}

\begin{proof}
(1)\qua Let $\tilde{M} = M_0\sqcup\{\hat{0}\}$
and let $\tilde{\mu}$ be the $\tilde{M}$--tuple of constraints given by
$\tilde{\mu}_l = \mu_l$ if $l \in M_0$
and $\tilde{\mu}_{\hat{0}} = H^2$,
where $H^2$ is a generic hyperplane in~$\PPP$.
By definition and the same argument as before,
\begin{equation}\label{n3circ3_l1e1}\begin{split}
&\big\lan a_{\hat{0}},{\overV}_2^{(1,1)}(\mu)\big\ran
=\sum_{\tilde{\mathcal{T}}}\Big(
\LlRr{\{\ev_{\hat{1}} \times \ev_{\hat{2}}\}^{-1}(\De_{\PPP\times\PPP}),
{\overU}_{\tilde{\mathcal{T}}}(\tilde{\mu})}\\
&\qquad\qquad\qquad\qquad\qquad\qquad
-{\mathcal{C}}_{\partial{\overU}_{\tilde{\mathcal{T}}}(\tilde{\mu})}
\big(\ev_{\hat{1}} \times \ev_{\hat{2}};\De_{\PPP\times\PPP}\big)\Big)\\
&\qquad\qquad
=2\big\lan a_{\hat{0}},{\overV}_2(\mu + \{H^1 : H^2\})\big\ran
-\sum_{\tilde{\mathcal{T}}}
{\mathcal{C}}_{\partial{\overU}_{\tilde{\mathcal{T}}}(\tilde{\mu})}
\big(\ev_{\hat{1}} \times \ev_{\hat{2}};\De_{\PPP\times\PPP}\big),
\end{split}\end{equation}
where the union is taken over all bubble types 
$\tilde{\mathcal{T}} = (M_2,I_2;\tilde{j},\tilde{d})$
such that 
$$\tilde{j}_{\hat{1}}=\tilde{1}, \qquad
\tilde{j}_{\hat{2}}=\tilde{2}, \qquad
\tilde{d}_{\tilde{1}},\tilde{d}_{\tilde{2}}>0,  \quad\hbox{and}\quad 
\tilde{d}_{\tilde{1}} + \tilde{d}_{\tilde{2}}=0.$$
Let ${\mathcal{T}} = (M_2,I;j,\under{d})$ be a bubble type such that
${\mathcal{T}} < \tilde{\mathcal{T}}$ and 
${\mathcal{U}}_{{\mathcal{T}}|\tilde{\mathcal{T}}}(\mu) \neq \eset$.
If $d_h \neq 0$ for some $h \in \hat{I}$ such that
$h \le j_{\hat{1}}$ or $h \le j_{\hat{2}}$
as in the proof of Lemma~\ref{n3p3_cohom_lmm2},
${\mathcal{U}}_{\tilde{\mathcal{T}},{\mathcal{T}}}$ does not contribute
to ${\mathcal{C}}_{\partial{\overU}_{\tilde{\mathcal{T}}}(\tilde{\mu})}
\big(\ev_{\hat{1}} \times \ev_{\hat{2}};\De_{\PPP\times\PPP}\big)$;
see also Lemma~\ref{n3p3_str_lmm3}.
Thus, we assume that $d_h = 0$ for all $h \in \hat{I}$ such that
$h \le j_{\hat{1}}$ or $h \le j_{\hat{2}}$.

(2)\qua By Lemma~\ref{n3p3_str_lmm3},
\begin{equation*}\begin{split}
&\big\{\ev_{\hat{1}} \times \ev_{\hat{2}}\big\}
\phi_{\tilde{\mathcal{T}},{\mathcal{T}}}(\ups)\\
&\qquad =\sum_{i=1,2}  (-1)^i\sum_{k=1}^{\i} 
\sum_{{}~h\in\chi_{\tilde{i}}({\mathcal{T}})}  
    \big(y_{h;\hat{i}}(\ups) - x_{\hat{i};h}(\ups)\big)^{-k}
\Big\{{\mathcal{D}}_{{\mathcal{T}},h}^{(k)} +
 \ve_{{\mathcal{T}},\hat{i};h}^{(k)}(\ups)\Big\}
\rho_{{\mathcal{T}},\hat{i};h}^{(0;k)}(\ups)
\end{split}\end{equation*}
for all $\ups \in \mathcal{FT}_{\de} - Y(\mathcal{FT};\hat{I}^+)$.
Since the map
$${\mathcal{F}} \equiv  
\bigoplus_{h\in\chi_{\tilde{1}}({\mathcal{T}})\cup\chi_{\tilde{2}}({\mathcal{T}})}
        L_h\lra\ev_{\hat{0}}^*T\PPP,\qquad
\ups\lra   
\sum_{h\in\chi_{\tilde{1}}({\mathcal{T}})\cup\chi_{\tilde{2}}({\mathcal{T}})}
        {\mathcal{D}}_{{\mathcal{T}},h}^{(1)}\ups_h,$$
is injective over ${\mathcal{U}}_{{\mathcal{T}}|\tilde{\mathcal{T}}}(\mu)$,
it follows that 
 ${\mathcal{U}}_{{\mathcal{T}}|\tilde{\mathcal{T}}}(\mu)$ is
$\big(\ev_{\hat{1}} \times \ev_{\hat{2}},\De_{\PPP\times\PPP}\big)$--hollow
unless 
$$\chi_{\tilde{1}}({\mathcal{T}})\cup\chi_{\tilde{2}}({\mathcal{T}})=\hat{I}^+.$$
In such a case, by Proposition~\ref{euler_prp},
a rescaling of the linear map, and the decomposition~\eqref{cart_split},
\begin{gather*}
{\mathcal{C}}_{{\mathcal{U}}_{{\mathcal{T}}|\tilde{\mathcal{T}}}(\tilde{\mu})}
\big(\ev_{\hat{1}} \times \ev_{\hat{2}};\De_{\PPP\times\PPP}\big)
=N(\al_1),\qquad\hbox{where}\\
\al_1 \in \Ga\big(
\ov{\frak M}_{\{\tilde{1}\}\sqcup\chi_{\tilde{1}}({\mathcal{T}})\sqcup M_{\tilde{1}}({\mathcal{T}})}
 \times 
\ov{\frak M}_{\{\tilde{2}\}\sqcup\chi_{\tilde{2}}({\mathcal{T}})\sqcup M_{\tilde{2}}({\mathcal{T}})}
 \times  {\mathcal{U}}_{{\overT}}(\tilde{\mu});
\hbox{Hom}({\mathcal{F}},\ev_{\hat{0}}^*T\PPP)\big),\\
\al(\ups)=   
\sum_{h\in\chi_{\tilde{1}}({\mathcal{T}})\cup\chi_{\tilde{2}}({\mathcal{T}})}
        {\mathcal{D}}_{{\overT},h}^{(1)}\ups_h.
\end{gather*}
Since $\al_1$ comes entirely from the third component,
$N(\al_1) = 0$ unless the first two components are zero-dimensional,
ie unless 
$$|\chi_{\tilde{1}}({\mathcal{T}})|=|\chi_{\tilde{2}}({\mathcal{T}})|=1
\quad\hbox{and}\quad
|M_{\tilde{1}}({\mathcal{T}})|=|M_{\tilde{2}}({\mathcal{T}})|=1.$$
It follows that
\begin{gather*}
\sum_{\tilde{\mathcal{T}}}
{\mathcal{C}}_{\partial{\overU}_{\tilde{\mathcal{T}}}(\tilde{\mu})}
\big(\ev_{\hat{1}} \times \ev_{\hat{2}};\De_{\PPP\times\PPP}\big)=
2N(\al_1),\qquad\hbox{where}\\
\al_1\in\Ga\big({\overV}_2(\tilde{\mu});
\hbox{Hom}(L_{\tilde{1}} \oplus L_{\tilde{2}},
\ev_{\hat{0}}^*T\PPP)\big),\quad
\al_1\big|_{{\mathcal{U}}_{\tilde{\mathcal{T}}}(\mu)}(\ups)=
{\mathcal{D}}_{\tilde{\mathcal{T}},\tilde{1}}^{(1)}\ups_1+
{\mathcal{D}}_{\tilde{\mathcal{T}},\tilde{2}}^{(1)}\ups_2,
\end{gather*}
if $\tilde{\mathcal{T}} = (M_0,I_2;\tilde{j},\tilde{d})$
bubble type such that 
$\tilde{d}_{\tilde{1}},\tilde{d}_{\tilde{2}} > 0$ and 
$\tilde{d}_{\tilde{1}} + \tilde{d}_{\tilde{2}} = 0$.
Using Propositions~\ref{zeros_prp} and~\ref{euler_prp}, 
we conclude that
\begin{equation}\label{n3circ3_l1e5}
\sum_{\tilde{\mathcal{T}}}
{\mathcal{C}}_{\partial{\overU}_{\tilde{\mathcal{T}}}(\tilde{\mu})}
\big(\ev_{\hat{1}} \times \ev_{\hat{2}};\De_{\PPP\times\PPP}\big)=
\big\lan 8a_{\hat{0}}^2 + 2a_{\hat{0}}\eta_{\hat{0},1},
{\overV}_2(\mu)\big\ran.
\end{equation}
This number is in fact computed in the proof of Lemma~5.13 in~\cite{Z1}.
The claim follows from equations~\eqref{n3circ3_l1e1} and~\eqref{n3circ3_l1e5}.
\end{proof}

\begin{lmm}
\label{n3circ2_l2}
With assumptions as in Lemma~\ref{n3circ2_lmm},
\begin{multline*}
\blr{\eta_{\hat{0},1},{\overV}_2^{(1,1)}(\mu)} 
= \blr{\eta_{\hat{0},1},{\overV}_2(\mu + \{H^1  : H^2\})}
+\big|{\mathcal{V}}_2(\mu + H^0)\big|\\
  +2\blr{a_{\hat{0}},{\overV}_2(\mu + H^1)}
  +2d\blr{a_{\hat{0}}^2,{\overV}_2(\mu)}
  -6\big|{\mathcal{V}}_3(\mu)\big|.
\end{multline*}
\end{lmm}

\begin{proof}
The proof is a mixture of the proof of
Lemma~\ref{n3circ2_l1} with the proof of  Lemma~\ref{n3p3_cohom_lmm4}.
\end{proof}

\section{Other examples}
\label{other_sec}

\subsection{Rational triple-pointed curves in $\PP$}
\label{n2p3_subs}

In this subsection, we prove Proposition~\ref{n2p3_prp},
ie the $\PP$ analogue of Theorem~\ref{n3p3_thm}.
The method is the same as in Section~\ref{n3p3_sec}, but the computation is 
significantly simpler, since there are many fewer boundary strata to consider.
Note that the formula of Proposition~\ref{n2p3_prp}
agrees with \cite[Lemma~3.2]{KQR} and \cite[Subsection~3.2]{V}.

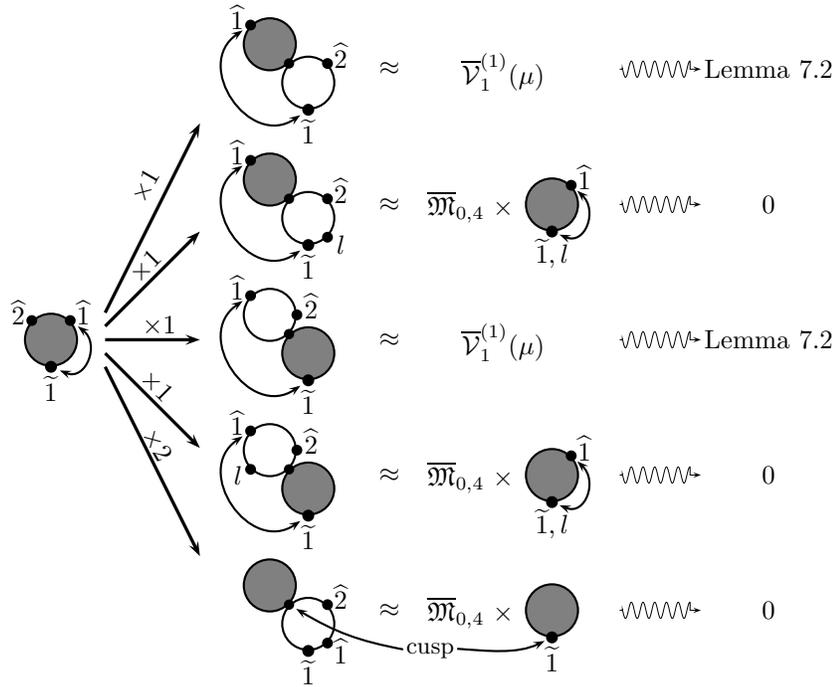
\begin{figure}[ht!]\small
\begin{pspicture}(-1,-5.7)(10,4.4)
\psset{unit=.36cm}
\pscircle[fillstyle=solid,fillcolor=gray](1,-2){1}
\pscircle*(1,-3){.22}\rput(1,-3.8){$\tilde{1}$}
\pscircle*(.29,-1.29){.2}\pscircle*(1.71,-1.29){.2}
\rput(2.2,-1){$\hat{1}$}\rput(-.2,-1){$\hat{2}$}
\pnode(1,-3){A0}\pnode(1.71,-1.29){B0}
\ncarc[nodesep=.35,arcangleA=-90,arcangleB=-75,ncurv=1.2]{<->}{A0}{B0}
% 2nd column starts here
\psline[linewidth=.12]{->}(3,-1)(6.5,6)
\rput{60}(4.5,3.7){$\times1$}
\psline[linewidth=.12]{->}(3,-1.5)(6.5,2)
\rput{45}(4.5,.7){$\times1$}
\psline[linewidth=.12]{->}(3,-2)(6.5,-2)
\rput(5,-1.5){$\times1$}
\psline[linewidth=.12]{->}(3,-2.5)(6.5,-6)
\rput{-45}(4.9,-3.7){$\times1$}
\psline[linewidth=.12]{->}(3,-3)(6.5,-10)
\rput{-60}(4.9,-5.9){$\times2$}
% 3rd column starts here
\pscircle(10.5,7.5){1}
\pscircle[fillstyle=solid,fillcolor=gray](9.08,8.92){1}
\pscircle*(10.5,6.5){.22}\rput(10.5,5.7){$\tilde{1}$}
\pscircle*(9.79,8.21){.19}\pscircle*(11.21,8.21){.2}
\rput(11.7,8.5){$\hat{2}$}
\pscircle*(8.37,9.63){.2}\rput(7.88,9.92){$\hat{1}$}
\pnode(10.5,6.5){A1}\pnode(8.37,9.63){B1}
\ncarc[nodesep=.35,arcangleA=90,arcangleB=90,ncurv=1.2]{<->}{A1}{B1}
% 2nd row starts here
\pscircle(10.5,2.5){1}
\pscircle[fillstyle=solid,fillcolor=gray](9.08,3.92){1}
\pscircle*(10.5,1.5){.22}\rput(10.5,.7){$\tilde{1}$}
\pscircle*(9.79,3.21){.19}
\pscircle*(11.21,3.21){.2}\rput(11.7,3.5){$\hat{2}$}
\pscircle*(8.37,4.63){.2}\rput(7.88,4.92){$\hat{1}$}
\pscircle*(11.21,1.79){.2}\rput(11.7,1.5){$l$}
\pnode(10.5,1.5){A2}\pnode(8.37,4.63){B2}
\ncarc[nodesep=.35,arcangleA=90,arcangleB=90,ncurv=1.2]{<->}{A2}{B2}
% 3rd row starts here
\pscircle[fillstyle=solid,fillcolor=gray](10.5,-2.5){1}
\pscircle(9.08,-1.08){1}
\pscircle*(10.5,-3.5){.22}\rput(10.5,-4.3){$\tilde{1}$}
\pscircle*(9.79,-1.79){.19}
\pscircle*(10.08,-1.08){.2}\rput(10.6,-.75){$\hat{2}$}
\pscircle*(8.37,-.37){.2}\rput(7.88,-.08){$\hat{1}$}
\pnode(10.5,-3.5){A3}\pnode(8.37,-.37){B3}
\ncarc[nodesep=.35,arcangleA=90,arcangleB=90,ncurv=1.2]{<->}{A3}{B3}
% 4th row starts here
\pscircle[fillstyle=solid,fillcolor=gray](10.5,-7.5){1}
\pscircle(9.08,-6.08){1}
\pscircle*(10.5,-8.5){.22}\rput(10.5,-9.3){$\tilde{1}$}
\pscircle*(9.79,-6.79){.19}
\pscircle*(10.08,-6.08){.2}\rput(10.6,-5.75){$\hat{2}$}
\pscircle*(8.37,-5.37){.2}\rput(7.88,-5.08){$\hat{1}$}
\pscircle*(8.37,-6.79){.2}\rput(7.9,-7){$l$}
\pnode(10.5,-8.5){A4}\pnode(8.37,-5.37){B4}
\ncarc[nodesep=.35,arcangleA=90,arcangleB=90,ncurv=1.2]{<->}{A4}{B4}
% 5th row starts here
\pscircle(10.5,-12.5){1}
\pscircle[fillstyle=solid,fillcolor=gray](9.08,-11.08){1}
\pscircle*(10.5,-13.5){.22}\rput(10.5,-14.3){$\tilde{1}$}
\pscircle*(9.79,-11.79){.19}\pscircle*(11.21,-11.79){.2}
\rput(11.7,-11.5){$\hat{2}$}
\pscircle*(11.21,-13.21){.2}\rput(11.7,-13.5){$\hat{1}$}
% 4th column starts here
\rput(13.5,8){$\approx$}\rput(13.5,3){$\approx$}\rput(13.5,-2){$\approx$}
\rput(13.5,-7){$\approx$}\rput(13.5,-12){$\approx$}
% 5th column starts here
\rput(17.7,7.9){${\overV}_1^{(1)}(\mu)$}
% 2nd row starts here
\rput(16.5,2.9){$\ov{\frak M}_{0,4}~\times$}
\pscircle[fillstyle=solid,fillcolor=gray](19.5,3){1}
\pscircle*(19.5,2){.22}\rput(19.5,1.2){$\tilde{1},l$}
\pscircle*(20.21,3.71){.2}\rput(20.7,4){$\hat{1}$}
\pnode(19.5,2){A7}\pnode(20.21,3.71){B7}
\ncarc[nodesep=.3,arcangleA=-90,arcangleB=-70,ncurv=1.2]{<->}{A7}{B7}
% 3rd row starts here
\rput(17.7,-2.1){${\overV}_1^{(1)}(\mu)$}
% 4th row starts here
\rput(16.5,-7.1){$\ov{\frak M}_{0,4}~\times$}
\pscircle[fillstyle=solid,fillcolor=gray](19.5,-7){1}
\pscircle*(19.5,-8){.22}\rput(19.5,-8.8){$\tilde{1},l$}
\pscircle*(20.21,-6.29){.2}\rput(20.7,-6){$\hat{1}$}
\pnode(19.5,-8){A9}\pnode(20.21,-6.29){B9}
\ncarc[nodesep=.3,arcangleA=-90,arcangleB=-70,ncurv=1.2]{<->}{A9}{B9}
% 5th row starts here
\rput(16.5,-12.1){$\ov{\frak M}_{0,4}~\times$}
\pscircle[fillstyle=solid,fillcolor=gray](19.5,-12){1}
\pscircle*(19.5,-13){.22}\rput(19.5,-13.8){$\tilde{1}$}
\pnode(14.3,-13.5){A2a}\pnode(15.7,-13.4){A2b}
\pnode(9.79,-11.79){B2a}\pnode(19.5,-13){B2b}
\rput(15,-13.5){\footnotesize cusp}
\ncarc[nodesep=.3,arcangleA=10,arcangleB=10,ncurv=.8]{->}{A2a}{B2a}
\ncarc[nodesep=.3,arcangleA=-10,arcangleB=-20,ncurv=.7]{->}{A2b}{B2b}
% 6th column starts here
\pscoil[linewidth=.03,coilarmA=0.1,coilarmB=0.4,coilaspect=0,coilheight=.72,
coilwidth=.6]{->}(22,8)(25,8)
\pscoil[linewidth=.03,coilarmA=0.1,coilarmB=0.4,coilaspect=0,coilheight=.72,
coilwidth=.6]{->}(22,3)(25,3)
\pscoil[linewidth=.03,coilarmA=0.1,coilarmB=0.4,coilaspect=0,coilheight=.72,
coilwidth=.6]{->}(22,-2)(25,-2)
\pscoil[linewidth=.03,coilarmA=0.1,coilarmB=0.4,coilaspect=0,coilheight=.72,
coilwidth=.6]{->}(22,-7)(25,-7)
\pscoil[linewidth=.03,coilarmA=0.1,coilarmB=0.4,coilaspect=0,coilheight=.72,
coilwidth=.6]{->}(22,-12)(25,-12)
% 7th column starts here
\rput(27.5,8){\small Lemma~\ref{n2p3contr_lmm1}}
\rput(27.5,3){$0$}
\rput(27.5,-2){\small Lemma~\ref{n2p3contr_lmm1}}
\rput(27.5,-7){$0$}
\rput(27.5,-12){$0$}
\end{pspicture}
\caption{An outline of the proof of Proposition~\ref{n2p3_prp}}
\label{n2p3_fig}
\end{figure}

Figure~\ref{n2p3_fig} outlines the computation of the boundary
contribution to the hom\-ol\-ogy-intersection 
number~$\llrr{{\mathcal{V}}_1^{(2)}(\mu)}$.
It shows all non-hollow boundary strata and the multiplicity
with which the number $N(\al)$ of zeros of an affine map
over a closure of each stratum enters  into $\llrr{{\mathcal{V}}_1^{(2)}(\mu)}$.
In three of the cases, the number~$N(\al)$ is easily seen to be zero.
Lemma~\ref{n2p3contr_lmm1} computes the number~$N(\al)$ 
in the remaining two cases.

If $d$ is a positive integer, let $n_d$ denote 
the number of degree--$d$ rational curves that pass
through \hbox{$3d - 1$} points in general position in~$\PP$.
Following~\cite{V}, we put
\begin{gather*}
A_d\equiv n_d=\big\lan a_{\hat{0}}^2,{\overV}_1(\mu)\big\ran,\\
B_d\equiv -\frac{n_d}{d}+
\frac{1}{2d}\sum_{d_1+d_2=d}   
\binom{3d - 2}{3d_1 - 1}d_1^2d_2^2n_{d_1}n_{d_2}
=\blr{a_{\hat{0}}\eta_{\hat{0},1},{\overV}_1(\mu)},\\
-C_d=\De_d\equiv 
\frac{1}{2}\sum_{d_1+d_2=d}   
\binom{3d - 2}{3d_1 - 1}d_1d_2n_{d_1}n_{d_2}
=\big|{\mathcal{V}}_2(\mu)\big|
=-\blr{\eta_{\hat{0},1}^2,{\overV}_1(\mu)},
\end{gather*}
where $\mu$ is a tuple of $3d{-}2$ points in $\PP$.
The computation of the above intersection numbers,
with essentially the same notation as in this paper,
can be found in Subsection~5.7 of~\cite{Z1}.

\begin{prp}
\label{n2p3_prp}
If $d$ is a positive integer, the number of rational one-com\-pon\-ent
degree--$d$ curves
that have a triple point and pass through a tuple $\mu$ of $3d{-}2$~points
in general position in~$\PP$ is $\frac{1}{6}|{\mathcal{V}}_1^{(2)}(\mu)|$,
where
$$\big|{\mathcal{V}}_1^{(2)}(\mu)\big|=
3(d^2 - 6d + 10)A_d-3(d - 6)B_d+6C_d.$$
\end{prp}

\begin{proof}
We use the same notation as in Section~\ref{n3p3_sec},
except now all the stable maps under consideration have values
in $\PP$, instead of~$\PPP$.
Similarly to Subsection~\ref{n3p3_sum_subs}, we have
\begin{equation}\label{n2p3_prp_e1}\begin{split}
\big|{\mathcal{V}}_1^{(2)}(\mu)\big|
&=\LlRr{ {\mathcal{V}}_1^{(2)}(\mu)}
-{\mathcal{C}}_{\partial{\overU}_{{\mathcal{T}}_2}(\mu)}
\big(\ev_{\tilde{1}} \times \ev_{\hat{2}};\De_{\PP\times\PP}\big)\\
&=2n_d^{(1)}+d
\big\lan a_{\hat{0}},{\overV}_1^{(1)}(\mu)\big\ran 
-{\mathcal{C}}_{\partial{\overU}_{{\mathcal{T}}_2}(\mu)}
\big(\ev_{\tilde{1}} \times \ev_{\hat{2}};\De_{\PP\times\PP}\big),
\end{split}\end{equation}
where $n_d^{(1)}$ denotes the number of 
degree--$d$ rational curves that pass
$3d - 1$ points in general position in~$\PP$
counted with a choice of a node, ie
\begin{equation}\label{n2p3_prp_e2}
n_d^{(1)}=\binom{d - 1}{2}n_d.
\end{equation} 
The number $\big\lan a_{\hat{0}},{\overV}_1^{(1)}(\mu)\big\ran$
is computed in Lemma~\ref{n2p2_lmm}.
In order to compute the boundary contribution
\hbox{${\mathcal{C}}_{\partial{\overU}_{{\mathcal{T}}_2}(\mu)}
\big(\ev_{\tilde{1}} \times \ev_{\hat{2}};\De_{\PP\times\PP}\big)$},
by Lemma~\ref{n3p3_exp_lmm1} it is sufficient to consider
only bubble types ${\mathcal{T}} = (M_2,I;j,\under{d})$
such that ${\mathcal{T}} < {\mathcal{T}}_2$ and
either $\chi_{\mathcal{T}}(\tilde{1},\hat{2}) = 0$ or
$\chi_{\mathcal{T}}(\hat{1},\hat{2}) = 0$.
Thus, the number
$${\mathcal{C}}_{\partial{\overU}_{{\mathcal{T}}_2}(\mu)}
\big(\ev_{\tilde{1}} \times \ev_{\hat{2}};\De_{\PP\times\PP}\big)$$
is computed by Lemmas~\ref{n2p3contr_lmm1} and~\ref{n2p3contr_lmm2}.
Finally, the numbers
\hbox{$\lan\eta_{\hat{0},1},{\overV}_1^{(1)}(\mu)\ran$} and
$|{\mathcal{S}}_1(\mu)|$ are given by Lemmas~\ref{n2p2_lmm}  and~\ref{n2cusps_lmm}.
\end{proof}

\begin{lmm}
\label{n2p3contr_lmm1}
The total contribution to the number 
$${\mathcal{C}}_{\partial{\overU}_{{\mathcal{T}}_2}(\mu)}
\big(\ev_{\tilde{1}} \times \ev_{\hat{2}};\De_{\PP\times\PP}\big)$$
from the boundary strata ${\mathcal{U}}_{{\mathcal{T}}_2,{\mathcal{T}}}$,
where ${\mathcal{T}} = (M_2,I;j,\under{d})$ is a bubble type
such that ${\mathcal{T}} < {\mathcal{T}}_2$ and
either $\chi_{\mathcal{T}}(\tilde{1},\hat{2}) = 0$ or
$\chi_{\mathcal{T}}(\hat{1},\hat{2}) = 0$, but not both, is given~by
$$    
\sum_{\chi_{\mathcal{T}}(\tilde{1},\hat{2})+\chi_{\mathcal{T}}(\hat{1},\hat{2})>0}
{\mathcal{C}}_{{\mathcal{U}}_{{\mathcal{T}}_2,{\mathcal{T}}}}
\big(\ev_{\tilde{1}} \times \ev_{\hat{2}};\De_{\PP\times\PP}\big)
=\blr{6a_{\hat{0}} + 2\eta_{\hat{0},1},{\overV}_1^{(1)}(\mu)}
-4\big|{\mathcal{S}}_1(\mu)\big|.$$
\end{lmm}

\begin{proof}
(1)\qua By symmetry, it is sufficient to consider the case
$\chi_{\mathcal{T}}(\tilde{1},\hat{2}) = 0$ and
$\chi_{\mathcal{T}}(\hat{1},\hat{2}) > 0$ and then double the answer.
By Lemma~\ref{n3p3_str_lmm1},
$${\overU}_{{\mathcal{T}}_2}^{(1)}(\mu)\cap{\mathcal{U}}_{{\mathcal{T}}_2,{\mathcal{T}}}
\subset {\mathcal{U}}_{{\mathcal{T}}|{\mathcal{T}}_2}^{(1)}(\mu).$$
By Lemma~\ref{n3p3_str_lmm2}, 
\begin{gather*}
\big\{\ev_{\tilde{1}} \times \ev_{\hat{2}}\big\}
\phi_{{\mathcal{T}}_2,{\mathcal{T}}}(\mu)
\sum_{h\in\chi_{\tilde{1}}({\mathcal{T}})}  
\big(y_{h;\hat{2}} - x_{h;\hat{2}}\big)^{-1}
\big\{{\mathcal{D}}_{{\mathcal{T}},h}^{(1)} + \ve_h(\ups)\big\}\rho_h(\ups),\\
\hbox{where}\qquad
\rho_h(\ups)=  \prod_{i\in(i_{\mathcal{T}}(h,\hat{2}),h]}       \ups_i,
\end{gather*}
for all 
$\ups \in \mathcal{FT}_{\de^*} - Y(\mathcal{FT};\hat{I}^+)$
and some $C^0$--negligible maps
$$\ve_h\co \mathcal{FT}_{\de^*} - Y(\mathcal{FT};\hat{I}^+) \lra 
\hbox{Hom}(L_h,\ev_{\hat{0}}^*T\PP).$$
The linear map
$$\al\co {\mathcal{F}} \equiv \bigoplus_{h\in\chi_{\tilde{1}}({\mathcal{T}})}   L_h
\lra\ev_{\hat{0}}^*T\PP$$
is injective over ${\mathcal{U}}_{{\mathcal{T}}|{\mathcal{T}}_2}^{(1)}(\mu)$
by Lemma~\ref{str_lmm}.
Thus, if $\hat{I}^+ \neq \chi_{\tilde{1}}({\mathcal{T}})$, then
${\mathcal{U}}_{{\mathcal{T}}|{\mathcal{T}}_2}^{(1)}(\mu)$ is
$(\ev_{\tilde{1}} \times \ev_{\hat{2}},\De_{\PP\times\PP})$--hollow, 
and
$${\mathcal{C}}_{{\mathcal{U}}_{{\mathcal{T}}_2,{\mathcal{T}}}}
\big(\ev_{\tilde{1}} \times \ev_{\hat{2}};\De_{\PP\times\PP}\big)
=0$$
by Proposition~\ref{euler_prp}, or Lemma~\ref{euler_lmm2},
and Lemma~\ref{n3p3_str_lmm1}.

(2)\qua On the other hand, if $\hat{I}^+ = \chi_{\tilde{1}}({\mathcal{T}})$, 
by the above and Lemma~\ref{n3p3_exp_lmm2}, 
${\mathcal{U}}_{{\mathcal{T}}|{\mathcal{T}}_2}^{(1)}(\mu)$ 
is $(\ev_{\tilde{1}} \times \ev_{\hat{2}},\De_{\PP\times\PP})$--regular,
and by Proposition~\ref{euler_prp}, a rescaling of the linear~map,
and the splitting~\eqref{cart_split},
\begin{gather*}
{\mathcal{C}}_{{\mathcal{U}}_{{\mathcal{T}}_2,{\mathcal{T}}}}
\big(\ev_{\tilde{1}} \times \ev_{\hat{2}};\De_{\PP\times\PP}\big)
=N(\al),\qquad\hbox{where}\\
\al \in \Ga\big(
\ov{\frak M}_{\{\tilde{1}\}\sqcup\chi_{\tilde{1}}({\mathcal{T}})\sqcup M_{\tilde{1}}({\mathcal{T}})}
 \times {\mathcal{U}}_{{\overT}}^{(1)}(\mu);
\hbox{Hom}({\mathcal{F}},\ev_{\hat{0}}^*T\PP)\big),\quad
\al(\ups)=\!\!\sum_{h\in\chi_{\tilde{1}}({\mathcal{T}})}\!\!
{\mathcal{D}}_{{\mathcal{T}},h}^{(1)}\ups_h.
\end{gather*}
Since the linear map $\al$ comes entirely from the second component
$N(\al) = 0$ unless
the first component is zero-dimensional,
ie $|\chi_{\tilde{1}}({\mathcal{T}})| = 1$ and
$M_{\tilde{1}}{\mathcal{T}} = \{\hat{2}\}$.
Thus, we conclude that
\begin{gather*}
\sum_{\chi_{\mathcal{T}}(\tilde{1},\hat{2})+\chi_{\mathcal{T}}(\hat{1},\hat{2})>1}
{\mathcal{C}}_{{\mathcal{U}}_{{\mathcal{T}}_2,{\mathcal{T}}}}
\big(\ev_{\tilde{1}} \times \ev_{\hat{2}};\De_{\PP\times\PP}\big)
=2N(\al_1),\\
\hbox{where}\qquad
\al_1 = {\mathcal{D}}_{{\mathcal{T}}_1,\tilde{1}}^{(1)}
 \in \Ga\big({\overV}_1^{(1)}(\mu);
\hbox{Hom}(L_{\tilde{1}},\ev_{\hat{0}}^*T\PP)\big).
\end{gather*}
The number $N(\al_1)$ is computed in Lemma~\ref{n2p3contr_lmm1b}.
\end{proof}

\begin{lmm}
\label{n2p3contr_lmm1b}
If $\al_1 = {\mathcal{D}}_{{\mathcal{T}}_1,\tilde{1}}^{(1)}
 \in \Ga\big({\overV}_1^{(1)}(\mu);
\hbox{Hom}(L_{\tilde{1}},\ev_{\hat{0}}^*T\PP)\big)$,
$$N(\al_1)=
\blr{3a_{\hat{0}} + c_1({\mathcal{L}}_{\tilde{1}}^*),{\overV}_1^{(1)}(\mu)}
-2\big|{\overS}_1(\mu)\big|.$$
\end{lmm}

\begin{proof}
(1)\qua Since  
$\al_1$ does not vanish on ${\mathcal{V}}_1^{(1)}(\mu)$ by Lemma~\ref{str_lmm}, 
by Propositions~\ref{zeros_prp} and~\ref{euler_prp},
\begin{equation}\label{n2p3_contr_lmm1b_e1}
N(\al_1)=\big\lan 3a_{\hat{0}}+
c_1(L_{\tilde{1}}^*),{\overV}_1^{(1)}(\mu)\big\ran
-{\mathcal{C}}_{\partial{\overV}_1^{(1)}(\mu)}(\al_1^{\perp}),
\end{equation}
where $\al_1^{\perp}$ denotes the composition of $\al_1$
with the projection $\pi_{\bar{\nu}_1}^{\perp}$
onto the quotient ${\mathcal{O}}_1$ of $\ev_{\hat{0}}^*T\PP$
by a generic trivial line subbundle~$\Bbb{C}\bar{\nu}_1$.
Suppose 
$${\mathcal{T}}=(M_1,I;j,\under{d})$$ 
is a bubble type such that ${\mathcal{T}} < {\mathcal{T}}_1$.

(2)\qua  If  $\chi_{\mathcal{T}}(\tilde{1},\hat{1}) > 0$, by
Lemma~\ref{n3p3_str_lmm1}
$${\overV}_1^{(1)}(\mu)\cap{\mathcal{U}}_{{\mathcal{T}}_1,{\mathcal{T}}}
\subset{\mathcal{U}}_{{\mathcal{T}}|{\mathcal{T}}_1}^{(1)}(\mu).$$
If in addition $d_{\tilde{1}} > 0$,
$\al_1$ does not vanish on ${\mathcal{U}}_{{\mathcal{T}}|{\mathcal{T}}_1}^{(1)}(\mu)$,
and thus ${\mathcal{U}}_{{\mathcal{T}}_1,{\mathcal{T}}}$ does not contribute
to~${\mathcal{C}}_{\partial{\overV}_1^{(1)}(\mu)}(\al_1^{\perp})$.
On the other hand, if $d_{\tilde{1}} = 0$ and
${\mathcal{U}}_{{\mathcal{T}}|{\mathcal{T}}_1}^{(1)}(\mu) \neq \eset$,
by dimension-counting via Lemma~\ref{str_lmm},
$\chi_{\tilde{1}}({\mathcal{T}}) = \{h\}$ is a single element-set
and ${\mathcal{T}} = {\mathcal{T}}_1(l)$ for some $l \in [N]$.
By Proposition~\ref{str_prp},
$${\mathcal{D}}_{{\mathcal{T}}_1,\tilde{1}}^{(1)}
\big(\phi_{{\mathcal{T}}_1,{\mathcal{T}}}(\ups)\big)
=\big\{{\mathcal{D}}_{{\mathcal{T}},h}^{(1)} + \ve(\ups)\big\}\ups,
\qquad\mbox{for all}\ 
\ups \in \mathcal{FT}_{\de^*} - Y(\mathcal{FT};\{h\}).$$
Since the section ${\mathcal{D}}_{{\mathcal{T}},h}^{(1)}$
does not vanish over ${\mathcal{U}}_{{\mathcal{T}}|{\mathcal{T}}_1}^{(1)}(\mu)$
by Lemma~\ref{str_lmm}, by Proposition~\ref{euler_prp},
$${\mathcal{C}}_{{\mathcal{U}}_{{\mathcal{T}}_1,{\mathcal{T}}}}(\al_1^{\perp})
=\big|{\mathcal{U}}_{{\overT}}(\mu)\big|.$$
Summing up for all bubble types ${\mathcal{T}} = {\mathcal{T}}_1(l)$,
we conclude~that
\begin{equation}\label{n2p3_contr_lmm1b_e3}
\sum_{\chi_{\mathcal{T}}(\tilde{1},\hat{1})>0}    
{\mathcal{C}}_{{\mathcal{U}}_{{\mathcal{T}}_1,{\mathcal{T}}}}(\al_1^{\perp})
=\big|{\mathcal{V}}_{1;1}^{(1)}(\mu)\big|.
\end{equation}
(3)\qua If ${\mathcal{T}} < {\mathcal{T}}_1$ and
$\chi_{\mathcal{T}}(\tilde{1},\hat{1}) = 0$,
by Lemma~\ref{n3p3_str_lmm2}
$${\overV}_1^{(1)}(\mu)\cap{\mathcal{U}}_{{\mathcal{T}}_1,{\mathcal{T}}}
\subset{\mathcal{S}}_{{\mathcal{T}}|{\mathcal{T}}_1}(\mu).$$
If ${\mathcal{S}}_{{\mathcal{T}}|{\mathcal{T}}_1}(\mu) \neq \eset$,
by dimension-counting via Lemma~\ref{str_lmm},
$\chi_{\tilde{1}}({\mathcal{T}}) = \{h\}$ is a single element-set
and ${\mathcal{T}} = {\mathcal{T}}_1(\hat{1})$.
 Subtracting $\big(y_{\hat{1}} - x_h\big)$
times the expansion of 
$\{\ev_{\tilde{1}} \times \ev_{\hat{1}}\} \circ 
\phi_{{\mathcal{T}}_1,{\mathcal{T}}}$ in Lemma~\ref{n3p3_str_lmm2} from
the expansion of ${\mathcal{D}}_{{\mathcal{T}}_1,\tilde{1}}^{(1)}$
in (3b) of Proposition~\ref{str_prp}, we obtain
\begin{gather*}
{\mathcal{D}}_{{\mathcal{T}}_1,\tilde{1}}^{(1)}\phi_{{\mathcal{T}}_1,{\mathcal{T}}}(\ups)=
-\big(y_{\hat{1}} - x_h\big)\otimes
\big\{{\mathcal{D}}_{{\mathcal{T}},h}^{(2)} + \ve(\ups)\big\}
\ups \otimes \ups\\
\mbox{for all}\  \ups \in \mathcal{FT}_{\de}  ~\mbox{such that}~
\phi_{{\mathcal{T}}_1,{\mathcal{T}}}(\ups) \in {\mathcal{U}}_{{\mathcal{T}}_1}^{(1)}(\mu)
\end{gather*}
and for some  $C^0$--negligible map
$$\ve\co \mathcal{FT}_{\de} - Y(\mathcal{FT};\hat{I}^+)
\lra L_h^{*\otimes2} \otimes\ev_{\hat{0}}^*T\PP.$$
By Lemma~\ref{str_lmm}, ${\mathcal{D}}_{{\mathcal{T}},h}^{(2)}$ does not vanish
on the finite set ${\mathcal{S}}_{{\mathcal{T}}|{\mathcal{T}}_1}(\mu)$.
Thus, by Proposition~\ref{euler_prp},
\begin{equation}\label{n2p3_contr_lmm1b_e5}
\sum_{\chi_{\mathcal{T}}(\tilde{1},\hat{1})=0}    
{\mathcal{C}}_{{\mathcal{U}}_{{\mathcal{T}}_1,{\mathcal{T}}}}(\al_1^{\perp})
=2\big|{\mathcal{S}}_{\mathcal{T}}(\mu)\big|
=2\big|{\mathcal{S}}_1(\mu)\big|. 
\end{equation}
The claim follows from~\eqref{n2p3_contr_lmm1b_e1}--\eqref{n2p3_contr_lmm1b_e5}
along with~\eqref{psi_class1} and~\eqref{psi_class2}.
\end{proof}

\begin{lmm}
\label{n2p3contr_lmm2}
If ${\mathcal{T}} = (M_2,I;j,\under{d}) < {\mathcal{T}}_2$ is a bubble type such that
$\chi_{\mathcal{T}}(\tilde{1},\hat{2}) = \chi_{\mathcal{T}}(\hat{1},\hat{2}) = 0$,
$${\mathcal{C}}_{{\mathcal{U}}_{{\mathcal{T}}_2,{\mathcal{T}}}}
\big(\ev_{\tilde{1}} \times \ev_{\hat{2}};\De_{\PP\times\PP}\big)=0.$$
\end{lmm}

\proof
Since $\chi_{\mathcal{T}}(\tilde{1},\hat{1}) = 0$, by Lemma~\ref{n3p3_str_lmm2},
${\overU}_{{\mathcal{T}}_2}^{(1)}(\mu)\cap{\mathcal{U}}_{{\mathcal{T}}_2,{\mathcal{T}}}
 \subset {\mathcal{S}}_{{\mathcal{T}}|{\mathcal{T}}_2}(\mu)$.
If ${\mathcal{S}}_{{\mathcal{T}}|{\mathcal{T}}_2}(\mu) \neq \eset$,
\hbox{$\chi_{\tilde{1}}({\mathcal{T}}) = \{h\}$} is a single-element set.
Subtracting the expansion of 
$\{\ev_{\tilde{1}} \times \ev_{\hat{1}}\} \circ 
\phi_{{\mathcal{T}}_2,{\mathcal{T}}}$ of Lemma~\ref{n3p3_str_lmm2}
times 
$$\big(y_{h;\hat{1}}(\ups) - x_{\hat{1};h}(\ups)\big)
                   \big(y_{h;\hat{2}}(\ups) - x_{\hat{2};h}(\ups)\big)$$
from the corresponding expansion
$\{\ev_{\tilde{1}} \times \ev_{\hat{2}}\} \circ 
\phi_{{\mathcal{T}}_2,{\mathcal{T}}}$, we obtain
$$\{\ev_{\tilde{1}} \times \ev_{\hat{2}}\} \circ  
\phi_{{\mathcal{T}}_2,{\mathcal{T}}}(\ups) 
=\{\al + \ve(\ups)\}\rho(\ups) 
~\mbox{for all}\ \ups \in \mathcal{FT}_{\de}
~\mbox{with}~\phi_{\tilde{\mathcal{T}}_2,{\mathcal{T}}}(\ups) \in 
{\mathcal{U}}_{{\mathcal{T}}_1}^{(1)}(\mu),$$
where $\rho$ is a a monomials map on $\mathcal{FT}$
with values in a line bundle $\tilde{\mathcal{F}}{\mathcal{T}}$ and
\hbox{$\al\co \tilde{\mathcal{F}}{\mathcal{T}} \lra \ev_{\hat{0}}^*T\PP$}
is a linear map.
Explicitly, if $h_1 = i_{\mathcal{T}}(h,\hat{1})$ and
$h_2 = i_{\mathcal{T}}(h,\hat{2})$,
\begin{gather*}
\rho(\ups)=  \prod_{i\in(h_1,h]}    \ups_i\otimes
  \prod_{i\in(h_2,h]}    \ups_i,\\
\al(\ups)={\mathcal{D}}_{{\mathcal{T}},h}^{(2)}\otimes
\begin{cases}
(y_{h;\hat{2}} - x_{\hat{2};h})^{-2} \otimes 
(y_{h;\hat{1}} - x_{\hat{1};h})^{-1} \otimes 
(y_{h;\hat{2}} - y_{h;\hat{1}}),&
\hbox{if}~h_1 = h_2;\\
(y_{h;\hat{2}} - x_{\hat{2};h})^{-2},& \hbox{if}~h_1 < h_2;\\
(y_{h;\hat{2}} - x_{\hat{2};h})^{-1} \otimes 
(y_{h;\hat{1}} - x_{\hat{1};h})^{-1},& \hbox{if}~h_1 > h_2.
\end{cases}
\end{gather*}
Thus, ${\mathcal{S}}_{{\mathcal{T}}|{\mathcal{T}}_2}(\mu)$ is
$(\ev_{\tilde{1}} \times \ev_{\hat{2}},\De_{\PPP\times\PPP})$--hollow
unless $\hat{I}^+ = \chi_{\tilde{1}}({\mathcal{T}})$.
On the other hand, if 
$\hat{I}^+ = \chi_{\tilde{1}}({\mathcal{T}})$,
by Proposition~\ref{euler_prp}, the decomposition~\eqref{cart_split},
and a rescaling of the linear~map,
\begin{gather*}
{\mathcal{C}}_{{\mathcal{U}}_{{\mathcal{T}}_2,{\mathcal{T}}}}
\big(\ev_{\tilde{1}} \times \ev_{\hat{2}};\De_{\PPP\times\PPP}\big)
=2N(\al_1),\qquad\hbox{where}\\
\al_1 = \pi_2^*
{\mathcal{D}}_{{\overT}}^{(2)}
\in\Ga\big(\ov{\frak M}_{\{\tilde{1},h\}\sqcup M_{\tilde{1}}{\mathcal{T}}} \times 
{\mathcal{S}}_{{\overT}}(\mu);
\pi_2^*\hbox{Hom}(L_{\tilde{1}}^*,\ev_{\hat{0}}^*T\PPP)\big).
\end{gather*}
Since $|M_{\tilde{1}}{\mathcal{T}}| \ge 2$,
the first factor is positive-dimensional, while
the linear map~$\al_1$ comes entirely from the second factor.
Thus,
$${\mathcal{C}}_{{\mathcal{U}}_{{\mathcal{T}}_2,{\mathcal{T}}}}
\big(\ev_{\tilde{1}} \times \ev_{\hat{2}};\De_{\PPP\times\PPP}\big)=
N(\al_1)=0.\eqno{\qed}$$

\begin{lmm}
\label{n2p2_lmm}
If $d$ is a positive integer 
and $\mu$ is a tuple of $3d{-}2$~points  
in general position in~$\PP$, 
the number of rational one-component degree--$d$ curves that 
pass through the constraints~$\mu$
and have a node on a generic line  
is $\frac{1}{2}\lan a_{\hat{0}},{\overV}_1^{(1)}(\mu)\ran$, where
$$\lan a_{\hat{0}},{\overV}_1^{(1)}(\mu)\ran=
\big\lan(2d - 3)a_{\hat{0}}^2 - 
a_{\hat{0}}\eta_{\hat{0},1},{\overV}_1(\mu)\big\ran.$$
Furthermore,
$$\big\lan \eta_{\hat{0},1},{\overV}_1^{(1)}(\mu)\big\ran
=\big\lan a_{\hat{0}}^2 + d\cdot
a_{\hat{0}}\eta_{\hat{0},1},{\overV}_1(\mu)\big\ran
-\big|{\mathcal{V}}_2(\mu)\big|.$$
\end{lmm}

\begin{proof}
(1)\qua In order to prove the first identity,
we take $\tilde{\mu}$ to be the 
$\tilde{M} \equiv [N]\cup\{\hat{0}\}$--tuple of constraints defined by
$\tilde{\mu}_l = \mu_l$ and $\tilde{\mu}_{\hat{0}} = H^1$,
where $H^1$ is a generic hyperplane.
Similarly to the proof of Lemma~\ref{n3p3_cohom_lmm2},
\begin{equation}\label{n2p2_lmm_e1}\begin{split}
\big\lan a_{\hat{0}},{\overV}_1^{(1)}(\mu)\big\ran
&=\LlRr{a_{\hat{0}},{\overV}_1^{(1)}(\mu)}
-{\mathcal{C}}_{\partial{\overU}_{{\mathcal{T}}_1}(\tilde{\mu})}
\big(\ev_{\tilde{1}} \times \ev_{\hat{1}},\De_{\PP\times\PP}\big)\\
&=2d\big\lan a_{\hat{0}}^2,{\overV}_1(\mu)\big\ran
-{\mathcal{C}}_{\partial{\overU}_{{\mathcal{T}}_1}(\tilde{\mu})}
\big(\ev_{\tilde{1}} \times \ev_{\hat{1}},\De_{\PP\times\PP}\big),
\end{split}\end{equation}
where ${\mathcal{C}}_{\partial{\overU}_{{\mathcal{T}}_1}(\tilde{\mu})}
\big(\ev_{\tilde{1}} \times \ev_{\hat{1}},\De_{\PP\times\PP}\big)$
is the contribution of 
$\partial{\overU}_{{\mathcal{T}}_1}(\tilde{\mu})$
to $\llrr{{\mathcal{V}}_1^{(1)}(\tilde{\mu})}$.
If \hbox{${\mathcal{T}} = (M_1,I;j,\under{d}) < {\mathcal{T}}_1$} 
is a bubble type such that $\chi_{\mathcal{T}}(\tilde{1},\hat{1}) > 0$,
the map $\ev_{\tilde{1}} \times \ev_{\hat{1}}$
is transversal to $\De_{\PP\times\PP}$ 
on~${\mathcal{U}}_{{\mathcal{T}}|{\mathcal{T}}_1}(\tilde{\mu})$ by Lemma~\ref{str_lmm}
and thus the boundary stratum 
${\mathcal{U}}_{{\mathcal{T}}_1,{\mathcal{T}}}(\tilde{\mu})$
does not contribute to the number
$${\mathcal{C}}_{\partial{\overU}_{{\mathcal{T}}_1}(\tilde{\mu})}
\big(\ev_{\tilde{1}} \times \ev_{\hat{1}},\De_{\PP\times\PP}\big).$$
If  $\chi_{\mathcal{T}}(\tilde{1},\hat{1}) = 0$
and ${\mathcal{U}}_{{\mathcal{T}}|{\mathcal{T}}_1}(\tilde{\mu}) \neq \eset$,
$\chi_{\tilde{1}}({\mathcal{T}}) = \{h\}$ and
$M_{\tilde{1}}{\mathcal{T}} = \{\hat{1}\}$
are single-element sets.
By Lemma~\ref{n3p3_str_lmm2},
$$\big\{\ev_{\tilde{1}} \times \ev_{\hat{1}}\big\}
\phi_{{\mathcal{T}}_1,{\mathcal{T}}}(\ups) 
=\big(y_{\hat{1}} - x_h\big)^{-1} \otimes 
\big\{{\mathcal{D}}_{{\mathcal{T}},h}^{(1)} + \ve(\ups)\big\}\ups
\qquad\mbox{for all}\ \ups \in \mathcal{FT}_{\de}.$$
Since the section ${\mathcal{D}}_{{\mathcal{T}},h}^{(1)}$
does not vanish on ${\mathcal{U}}_{{\mathcal{T}}|{\mathcal{T}}_1}(\mu)$
by Lemma~\ref{str_lmm},  
${\mathcal{U}}_{{\mathcal{T}}|{\mathcal{T}}_1}$ is 
\hbox{$\big(\ev_{\tilde{1}} \times \ev_{\hat{1}},
\De_{\PP\times\PP}\big)$--hollow}
unless $\hat{I}^+ = \{h\}$.
If $\hat{I}^+ = \{h\}$, by 
Proposition~\ref{euler_prp}, decomposition~\eqref{cart_split},
and a rescaling of the linear map,
\begin{gather*}
{\mathcal{C}}_{{\mathcal{U}}_{{\mathcal{T}}|\tilde{\mathcal{T}}_1}(\tilde{\mu})}
\big(\ev_{\tilde{1}} \times \ev_{\hat{1}};\De_{\PP\times\PP}\big)
=N(\al_1),\qquad\hbox{where}\\
\al_1 = {\mathcal{D}}_{{\mathcal{T}}_0,\hat{1}}^{(1)}\in
\Ga\big({\overV}_1(\tilde{\mu});
 \hbox{Hom}(L_{\tilde{1}},\ev_{\hat{0}}^*T\PP)\big).
\end{gather*}
By Propositions~\ref{zeros_prp} and~\ref{euler_prp},
$$N(\al_1)=\big\lan 3a_{\hat{0}} + c_1(L_{\tilde{1}}^*),
{\overV}_1(\tilde{\mu})\big\ran
-{\mathcal{C}}_{\partial{\overV}_1(\mu)}(\al_1^{\perp}).$$
If ${\mathcal{T}} \equiv (M_0,I;j,\under{d}) < {\mathcal{T}}_0$
is a bubble type such that $\al_1$ vanishes somewhere
on ${\mathcal{U}}_{{\mathcal{T}}|{\mathcal{T}}_0}(\tilde{\mu})$,
\hbox{${\mathcal{T}} = {\mathcal{T}}_1(l)$} for some $l \in [N]$.
{}From Proposition~\ref{str_prp}, we then obtain
$${\mathcal{C}}_{\partial{\overV}_1(\tilde{\mu})}(\al_1^{\perp})=
\big|{\mathcal{V}}_{1;1}(\tilde{\mu})\big|.$$
Putting everything together and using identities
\eqref{psi_class1} and~\eqref{psi_class2}, we conclude that
\begin{equation}\label{n2p2_lmm_e3}
{\mathcal{C}}_{\partial{\overU}_{{\mathcal{T}}_1}(\tilde{\mu})}
\big(\ev_{\tilde{1}} \times \ev_{\hat{1}},\De_{\PP\times\PP}\big)
=\big\lan 3a_{\hat{0}}^2 + a_{\hat{0}}c_1({\mathcal{L}}_{\tilde{1}}^*),
{\overV}_1(\mu)\big\ran.
\end{equation}
The first claim of the lemma follows from \eqref{n2p2_lmm_e1}
and~\eqref{n2p2_lmm_e3}.

(2)\qua Let $s$ be a section of ${\mathcal{L}}_{\tilde{1}}^*$
with good properties, ie as in the proof of Lemma~\ref{n3p3_cohom_lmm4}.
Then, 
\begin{equation}\label{n2p2_lmm_e5}\begin{split}
&\big\lan c_1({\mathcal{L}}_{\tilde{1}}^*),
{\overV}_1^{(1)}(\mu)\big\ran\\
&\quad =\LlRr{c_1({\mathcal{L}}_{\tilde{1}}^*),{\overV}_1^{(1)}(\mu)}
-{\mathcal{C}}_{s^{-1}(0)\cap\partial{\overU}_{{\mathcal{T}}_1}(\mu)}
\big(\ev_{\tilde{1}} \times \ev_{\hat{1}},\De_{\PP\times\PP}\big)\\
&\quad =\big\lan c_1({\mathcal{L}}_{\tilde{1}}^*),{\overV}_1(\mu + H^0)\big\ran
+\big\lan a_{\hat{0}}c_1({\mathcal{L}}_{\tilde{1}}^*),
{\overV}_1(\mu + H^1)\big\ran
+3\big\lan a_{\hat{0}}^2,{\overV}_1(\mu)\big\ran\\
&\quad\qquad\qquad\qquad\qquad\qquad\qquad\quad
-{\mathcal{C}}_{s^{-1}(0)\cap\partial{\overU}_{{\mathcal{T}}_1}(\mu)}
\big(\ev_{\tilde{1}} \times \ev_{\hat{1}},\De_{\PP\times\PP}\big)\\
&\quad=
\big\lan a_{\hat{0}}^2 + d\cdot a_{\hat{0}}c_1({\mathcal{L}}_{\tilde{1}}^*),
{\overV}_1(\mu + H^1)\big\ran-
{\mathcal{C}}_{s^{-1}(0)\cap\partial{\overU}_{{\mathcal{T}}_1}(\mu)}
\big(\ev_{\tilde{1}} \times \ev_{\hat{1}},\De_{\PP\times\PP}\big).
\end{split}\end{equation}
In the last equality we used Lemma~5.17 of~\cite{Z1}, 
which is essentially Lemma~2.2.2 of~\cite{P}.
If 
$${\mathcal{T}} \equiv (M_1,I;j,\under{d})<{\mathcal{T}}_1$$
is a bubble type such that $\chi_{\mathcal{T}}(\tilde{1},\hat{1}) > 0$,
the space ${\mathcal{U}}_{{\mathcal{T}}|{\mathcal{T}}_1}$
does not contribute to
${\mathcal{C}}_{\partial{\overU}_{{\mathcal{T}}_1}(\mu)}
\big(\ev_{\tilde{1}} \times \ev_{\hat{1}},\De_{\PP\times\PP}\big)$
unless $\chi_{\mathcal{T}}(\tilde{1},\hat{1}) = 0$
and $\eta_{\hat{0},1}|{\overU}_{{\mathcal{T}}|{\mathcal{T}}_1}(\mu) \neq 0$.
On the other hand,
$$\chi_{\mathcal{T}}(\tilde{1},\hat{1}) = 0, ~~
\eta_{\hat{0},1}|{\overU}_{{\mathcal{T}}|{\mathcal{T}}_1}(\mu) \neq 0
\quad\!\!\Lra\!\!\quad
\hat{I}^+ = \chi_{\tilde{1}}({\mathcal{T}}) = \{h_1,h_2\}, ~~
M_{\tilde{1}}({\mathcal{T}}) = \{\hat{1}\}.$$
By Lemma~\ref{n3p3_str_lmm2},
$$\big\{\ev_{\tilde{1}} \times \ev_{\hat{1}}\big\}
\phi_{{\mathcal{T}}_1,{\mathcal{T}}}(\ups) 
=\sum_{h\in\chi_{\tilde{1}}({\mathcal{T}})}
\big(y_{\hat{1}} - x_h\big)^{-1} \otimes 
\big\{{\mathcal{D}}_{{\mathcal{T}},h}^{(1)} + \ve_h(\ups)\big\}\ups_h,$$
for all
$\ups \in \mathcal{FT}_{\de}$.
Thus,
$${\mathcal{C}}_{{\mathcal{U}}_{{\mathcal{T}}|{\mathcal{T}}_1}(\mu)}
\big(\ev_{\tilde{1}} \times \ev_{\hat{1}};\De_{\PP\times\PP}\big)
=\, ^{\pm} \big|s^{-1}(0)\cap {\mathcal{U}}_{{\mathcal{T}}|{\mathcal{T}}_1}(\mu)\big|
=\big|{\mathcal{U}}_{{\overT}}(\mu)\big|.$$
We conclude that
\begin{equation}\label{n2p2_lmm_e7}
{\mathcal{C}}_{s^{-1}(0)\cap\partial{\overU}_{{\mathcal{T}}_1}(\mu)}
\big(\ev_{\tilde{1}} \times \ev_{\hat{1}},\De_{\PP\times\PP}\big)
=\big|{\mathcal{V}}_2(\mu)\big|.
\end{equation}
The second claim of the lemma follows from \eqref{n2p2_lmm_e5}
and~\eqref{n2p2_lmm_e7}.
\end{proof}

\begin{lmm}
\label{n2cusps_lmm}
If $d \ge 1$, the number of rational degree--$d$ cuspidal curves 
passing through a tuple~$\mu$ of $3d{-}2$ 
points in general position in~$\PP$ is given by
$$\big|{\mathcal{S}}_1(\mu)\big|=
\big\lan 3a_{\hat{0}}^2 + 3a_{\hat{0}}\eta_{\hat{0},1}
 + \eta_{\hat{0},1}^2,{\overV}_1(\mu)\big\ran-
\big|{\mathcal{V}}_2(\mu)\big|.$$
\end{lmm}

\begin{proof}
This is the $n = 2$ case of Theorem~\ref{cusps_thm};
see \cite[Lemma~5.4]{Z1} for a direct proof.
The same formula can also be found in
\cite[Subsection~4.5]{P} and \cite[Subsection~3.2]{V}.
\end{proof}

\subsection{Rational tacnodal curves in $\PP$}
\label{n2tac_subs}

In this subsection, we prove Proposition~\ref{n2tac_prp},
the $\PP$--analogue of Theorem~\ref{n3tac_thm}.
The formula we obtain agrees with previously known results;
see equation~(1.2) in~\cite{DH} and Subsection 3.2 in~\cite{V}.

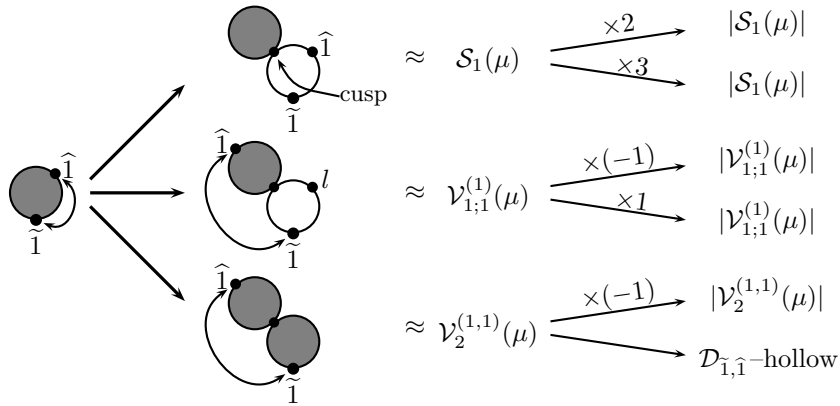
\begin{figure}[ht!]\small
\begin{pspicture}(-1,-4)(10,2.7)
\psset{unit=.36cm}
\pscircle[fillstyle=solid,fillcolor=gray](1,-2){1}
\pscircle*(1,-3){.22}\rput(1,-3.8){$\tilde{1}$}
\pscircle*(1.71,-1.29){.2}\rput(2.2,-1){$\hat{1}$}
\pnode(1,-3){A0}\pnode(1.71,-1.29){B0}
\ncarc[nodesep=.3,arcangleA=-90,arcangleB=-75,ncurv=1.2]{<->}{A0}{B0}
% 2nd column starts here
\psline[linewidth=.12]{->}(3,-1.5)(6.5,2)
\psline[linewidth=.12]{->}(3,-2)(6.5,-2)
\psline[linewidth=.12]{->}(3,-2.5)(6.5,-6)
% 3rd column starts here
\pscircle(10.5,2.5){1}
\pscircle[fillstyle=solid,fillcolor=gray](9.08,3.92){1}
\pscircle*(10.5,1.5){.22}\rput(10.5,0.7){$\tilde{1}$}
\pscircle*(9.79,3.21){.19}
\pscircle*(11.21,3.21){.2}\rput(11.7,3.5){$\hat{1}$}
\pnode(12.5,1.5){A1}\rput(13.1,1.5){\footnotesize cusp}\pnode(9.79,3.21){B1}
\ncarc[nodesep=.3,arcangleA=20,arcangleB=30,ncurv=1]{->}{A1}{B1}
% 2nd row starts here
\pscircle(10.5,-2.5){1}
\pscircle[fillstyle=solid,fillcolor=gray](9.08,-1.08){1}
\pscircle*(10.5,-3.5){.22}\rput(10.5,-4.3){$\tilde{1}$}
\pscircle*(9.79,-1.79){.19}
\pscircle*(11.21,-1.79){.2}\rput(11.7,-1.5){$l$}
\pscircle*(8.37,-.37){.2}\rput(7.88,-.08){$\hat{1}$}
\pnode(10.5,-3.5){A2}\pnode(8.37,-.37){B2}
\ncarc[nodesep=.35,arcangleA=90,arcangleB=90,ncurv=1.2]{<->}{A2}{B2}
% 3rd row starts here
\pscircle[fillstyle=solid,fillcolor=gray](10.5,-7.5){1}
\pscircle[fillstyle=solid,fillcolor=gray](9.08,-6.08){1}
\pscircle*(10.5,-8.5){.22}\rput(10.5,-9.3){$\tilde{1}$}
\pscircle*(9.79,-6.79){.19}
\pscircle*(8.37,-5.37){.2}\rput(7.88,-5.08){$\hat{1}$}
\pnode(10.5,-8.5){A3}\pnode(8.37,-5.37){B3}
\ncarc[nodesep=.35,arcangleA=90,arcangleB=90,ncurv=1.2]{<->}{A3}{B3}
% 4th column starts here
\rput(15,3){$\approx$}\rput(15,-2){$\approx$}\rput(15,-7){$\approx$}
% 5th column starts here
\rput(17.7,2.9){${\mathcal{S}}_1(\mu)$}
\rput(17.7,-2.1){${\mathcal{V}}_{1;1}^{(1)}(\mu)$}
\rput(17.7,-7.1){${\mathcal{V}}_2^{(1,1)}(\mu)$}
% 6th column starts here
\psline{->}(20,3.25)(25,4)\rput{8}(22.5,4.2){$\times2$}
\psline{->}(20,2.75)(25,2)\rput{-8}(23,2.7){$\times3$}
\psline{->}(20,-1.75)(25,-1)\rput{8}(22.5,-.8){$\times(-1)$}
\psline{->}(20,-2.25)(25,-3)\rput{-8}(23,-2.3){$\times1$}
\psline{->}(20,-6.75)(25,-6)\rput{8}(22.5,-5.8){$\times(-1)$}
\psline{->}(20,-7.25)(25,-8)
% 7th column starts here
\rput(28,4.25){$|{\mathcal{S}}_1(\mu)|$}
\rput(28,2){$|{\mathcal{S}}_1(\mu)|$}
\rput(28,-.75){$|{\mathcal{V}}_{1;1}^{(1)}(\mu)|$}
\rput(28,-3){$|{\mathcal{V}}_{1;1}^{(1)}(\mu)|$}
\rput(28,-5.75){$|{\mathcal{V}}_2^{(1,1)}(\mu)|$}
\rput(28,-8.3){${\mathcal{D}}_{\tilde{1},\hat{1}}$--\small{hollow}}  
\end{pspicture}
\caption{An outline of the proof of Proposition~\ref{n2tac_prp}}
\label{n2tac_fig}
\end{figure}

Figure~\ref{n2tac_fig} shows the three types of boundary strata 
${\overV}_1^{(1)}(\mu)\cap{\mathcal{U}}_{{\mathcal{T}}_1,{\mathcal{T}}}$
such that
$$\Bbb{P}(L_{\tilde{1}} \oplus L_{\tilde{1}}^*)|
{\overV}_1^{(1)}(\mu)\cap{\mathcal{U}}_{{\mathcal{T}}_1,{\mathcal{T}}}$$
is not contained in a finite union of 
${\mathcal{D}}_{\tilde{1},\hat{1}}$--hollow sets.
For such boundary strata, 
$$\Bbb{P}(L_{\tilde{1}} \oplus L_{\tilde{1}}^*)|
{\overV}_1^{(1)}(\mu)\cap{\mathcal{U}}_{{\mathcal{T}}_1,{\mathcal{T}}}$$
is a union of one 
${\mathcal{D}}_{\tilde{1},\hat{1}}$--regular or hollow subset
and one ${\mathcal{D}}_{\tilde{1},\hat{1}}$--regular subset:
a section over the base 
${\overV}_1^{(1)}(\mu)\cap{\mathcal{U}}_{{\mathcal{T}}_1,{\mathcal{T}}}$
and its complement.
The second-to-last column of Figure~\ref{n2tac_fig}
shows the multiplicity with which each number $N(\al)$
of zeros of an affine map over a closure 
of the larger and the smaller subset, if it is regular, 
enters into the euler class
of the bundle 
$\ga_{L_{\tilde{1}}\oplus L_{\tilde{1}}^*}^* \otimes \ev_{\hat{0}}^*T\PP$
as computed via the section~${\mathcal{D}}_{\tilde{1},\hat{1}}$.
The last column gives the number $N(\al)$ for each regular subset 
of the boundary strata.
Contributions from the boundary strata as in the first row of 
Figure~\ref{n2tac_fig} are computed in Lemma~\ref{n2tac_contr_lmm1}.
Lemma~\ref{n2tac_contr_lmm2} deals with 
the boundary strata as in the last two rows of Figure~\ref{n2tac_fig}.

\begin{prp}
\label{n2tac_prp}
If $d$ is a nonnegative integer, the number of rational one-component degree--$d$ curves 
that have a tacnodal point and pass through a tuple $\mu$ of $3d{-}2$~points  
in general position in~$\PP$ is $\frac{1}{2}|{\mathcal{S}}_1^{(1)}(\mu)|$, where
$$\big|{\mathcal{S}}_1^{(1)}(\mu)\big| = 2(3d - 11)A_d+2(d - 9)B_d-8C_d.$$
\end{prp}

\begin{proof}
Similarly to Subsection~\ref{n3tac_sum_subs},
\begin{equation}\label{n2tac_prp_e1}
\big|{\mathcal{S}}_1^{(1)}(\mu)\big|=
\blr{3a_{\hat{0}},{\overV}_1^{(1)}(\mu)}
- {\mathcal{C}}_{\partial\Bbb{P}(L_{\tilde{1}}\oplus L_{\tilde{1}}^*)}
({\mathcal{D}}_{\tilde{1},\hat{1}}),
\end{equation}
where ${\mathcal{C}}_{\partial\Bbb{P}(L_{\tilde{1}}\oplus L_{\tilde{1}}^*)}
({\mathcal{D}}_{\tilde{1},\hat{1}})$ is the contribution 
from the boundary strata of the space
$\Bbb{P}(L_{\tilde{1}} \oplus L_{\tilde{1}}^*)$.
This contribution is computed in Lemmas~\ref{n2tac_contr_lmm1} 
and~\ref{n2tac_contr_lmm2}.
The numbers 
\hbox{$\lan a_{\hat{0}},{\overV}_1^{(1)}(\mu)\ran$}
and $|{\mathcal{S}}_1(\mu)|$ are given by 
Lemmas~\ref{n2p2_lmm} and~\ref{n2cusps_lmm}.
Finally, the number of two-component rational curves that
pass through $3d{-}2$ points in general position in~$\PP$,
counted a with choice of an ordered pair of distinct nodes
at which the two components intersect and with a choice of a branch
at one of these nodes, is easily seen to be
\begin{equation}\label{n2tac_prp_e3}\begin{split}
\big|{\mathcal{V}}_2^{(1,1)}(\mu)\big|
&=\sum_{d_1+d_2=d}\binom{3d - 2}{3d_1 - 1}
d_1d_2\big(d_1d_2 - 1\big)n_{d_1}n_{d_2}\\
&=2A_d+2d B_d+2C_d
\end{split}\end{equation}
which proves the claim.
\end{proof}

\begin{lmm}
\label{n2tac_contr_lmm1} 
The total contribution to 
${\mathcal{C}}_{\partial\Bbb{P}(L_{\tilde{1}}\oplus L_{\tilde{1}}^*)}
({\mathcal{D}}_{\tilde{1},\hat{1}})$
from the boundary strata 
\hbox{$\Bbb{P}(L_{\tilde{1}} \oplus L_{\tilde{1}}^*)
|{\mathcal{U}}_{{\mathcal{T}}_1,{\mathcal{T}}}$},
where ${\mathcal{T}} = (M_1,I;j,\under{d})$ is a bubble type
such that ${\mathcal{T}} < {\mathcal{T}}_1$ and
$\chi_{\mathcal{T}}(\tilde{1},\hat{1}) = 0$, is given~by
$$\sum_{\chi_{\mathcal{T}}(\tilde{1},\hat{1})=0}    
{\mathcal{C}}_{\partial\Bbb{P}(L_{\tilde{1}}\oplus L_{\tilde{1}}^*)|
{\mathcal{U}}_{{\mathcal{T}}_2,{\mathcal{T}}}} ({\mathcal{D}}_{\tilde{1},\hat{1}})
=5\big|{\mathcal{S}}_1(\mu)\big|.$$
\end{lmm}

\begin{proof}
(1)\qua By Lemma~\ref{n3p3_str_lmm2},
$${\overU}_{{\mathcal{T}}_1}^{(1)}(\mu)\cap{\mathcal{U}}_{{\mathcal{T}}_1,{\mathcal{T}}}
\subset{\mathcal{S}}_{{\mathcal{T}}|{\mathcal{T}}_1}(\mu).$$
If ${\mathcal{S}}_{{\mathcal{T}}|{\mathcal{T}}_1}(\mu) \neq \eset$,
$\hat{I}^+ = \chi_{\tilde{1}}({\mathcal{T}}) = \{h\}$ 
and $M_{\tilde{1}}{\mathcal{T}} = \{\hat{1}\}$
are a single-element sets.
Subtracting the expansion of 
$\{\ev_{\tilde{1}} \times \ev_{\hat{1}}\} \circ \phi_{{\mathcal{T}}_1,{\mathcal{T}}}$
of Lemma~\ref{n3p3_str_lmm2} multiplied by
$(y_{\hat{1}} - x_h)$  and by $-(y_{\hat{1}} - x_h)^{-1}$
from the expansions of 
${\mathcal{D}}_{{\mathcal{T}}_1,\tilde{1}}^{(1)} \circ \phi_{{\mathcal{T}}_1,{\mathcal{T}}}$
and ${\mathcal{D}}_{{\mathcal{T}}_1,\hat{1}}^{(1)} \circ \phi_{{\mathcal{T}}_1,{\mathcal{T}}}$,
respectively, given by Proposition~\ref{str_prp},
we~obtain
\begin{equation*}\begin{split}
{\mathcal{D}}_{\tilde{1},\hat{1}}\phi_{{\mathcal{T}}_1,{\mathcal{T}}}
([\ups_{\tilde{1}},\ups_{\hat{1}}];\ups)=-
\big\{\big( (y_{\hat{1}} - x_h)^{-1} \otimes \ups_{\tilde{1}} 
+(y_{\hat{1}} - x_h)^{-3} \otimes \ups_{\hat{1}}\big)
\otimes{\mathcal{D}}_{{\mathcal{T}},h}^{(2)} \quad&\\
+\ve(\ups)\}\ups \otimes \ups&
\end{split}\end{equation*}
for all
$([\ups_{\tilde{1}},\ups_{\hat{1}}];\ups) \in \mathcal{FT}_{\de}$
such that
$\phi_{{\mathcal{T}}_1,{\mathcal{T}}}(\ups) \in {\mathcal{U}}_{{\mathcal{T}}_1}^{(1)}$.
Let 
$${\mathcal{Z}}_{\mathcal{T}}=\big\{
\big[\ups_{\tilde{1}},\ups_{\hat{1}}\big]
 \in \Bbb{P}(L_{\tilde{1}} \oplus L_{\tilde{1}}^*)
|{\mathcal{U}}_{{\mathcal{T}}|{\mathcal{T}}_1} :
(y_{\hat{1}} - x_h)^{-1}\ups_{\tilde{1}}
+(y_{\hat{1}} - x_h)^{-3}\ups_{\hat{1}} = 0\big\}.$$
Since the section~${\mathcal{D}}_{{\mathcal{T}},h}^{(2)}$ does not vanish
over ${\mathcal{S}}_{{\mathcal{T}}|{\mathcal{T}}_1}(\mu)$,
by Proposition~\ref{euler_prp}, the decomposition~\eqref{cart_split}, and 
a rescaling of the linear~map,
\begin{gather*}
{\mathcal{C}}_{\Bbb{P}(L_{\tilde{1}}\oplus L_{\tilde{1}}^*)|
{\mathcal{S}}_{{\mathcal{T}}|{\mathcal{T}}_1}(\mu)-{\mathcal{Z}}_{\mathcal{T}}}
({\mathcal{D}}_{\tilde{1},\hat{1}})
=2N(\al_1),\qquad\hbox{where}\\
\al_1\in\Ga\big(\Bbb{P}^1 \times {\mathcal{S}}_1(\mu);
\hbox{Hom}(\ga^*,\ga^* \otimes \ev_{\hat{0}}^*T\PP)\big)~~
\end{gather*}
is a nonvanishing section.
Thus, by Proposition~\ref{zeros_prp},
\begin{equation}\label{n2tac_contr_lmm1e3}
\sum_{\chi_{\mathcal{T}}(\tilde{1},\hat{1})=0}    
{\mathcal{C}}_{\Bbb{P}(L_{\tilde{1}}\oplus L_{\tilde{1}}^*)|
{\mathcal{S}}_{{\mathcal{T}}|{\mathcal{T}}_1}(\mu)-{\mathcal{Z}}_{\mathcal{T}}}
({\mathcal{D}}_{\tilde{1},\hat{1}})=
2\big\lan 2\la - \la,\Bbb{P}^1\big\ran \big|{\mathcal{S}}_1(\mu)\big|=
2\big|{\mathcal{S}}_1(\mu)\big|.
\end{equation}
(3)\qua  In order to compute the contribution from the
space~${\mathcal{Z}}_{\mathcal{T}}$,
we model a neighborhood of ${\mathcal{Z}}_{\mathcal{T}}$ in
\hbox{$\Bbb{P}(L_{\tilde{1}} \oplus L_{\tilde{1}}^*)$} by the~map
$$L_{\tilde{1}}^*\otimes L_{\tilde{1}}^*\lra
\Bbb{P}(L_{\tilde{1}} \oplus L_{\tilde{1}}^*),\qquad
\big([\ups_{\tilde{1}},\ups_{\hat{1}}],u\big)\lra 
\big[\ups_{\tilde{1}},\ups_{\hat{1}} + u(\ups_{\tilde{1}})\big].$$
Near ${\mathcal{Z}}_{\mathcal{T}}$,
\begin{equation*}\begin{split}
{\mathcal{D}}_{\tilde{1},\hat{1}}\phi_{{\mathcal{T}}_1,{\mathcal{T}}}
([\ups_{\tilde{1}},\ups_{\hat{1}}];u,\ups)=&
-(y_{\hat{1}} - x_h)^{-3} \otimes 
\big({\mathcal{D}}_{{\mathcal{T}},h}^{(2)}+\ve_2(u,\ups)\big)
\ups \otimes \ups \otimes u\\
&-(y_{\hat{1}} - x_h)^{-2} \otimes \ups_{\tilde{1}} \otimes 
\big({\mathcal{D}}_{{\mathcal{T}},h}^{(3)}+\ve_3(u,\ups)\big)
\ups \otimes \ups \otimes \ups
\end{split}\end{equation*}
for all
$([\ups_{\tilde{1}},\ups_{\hat{1}}];\ups) \in \mathcal{FT}_{\de}$
such that
$\phi_{{\mathcal{T}}_1,{\mathcal{T}}}(\ups) \in {\mathcal{U}}_{{\mathcal{T}}_1}^{(1)}$.
By Lemma~\ref{str_lmm}, the images of
${\mathcal{D}}_{{\mathcal{T}},h}^{(2)}$ and ${\mathcal{D}}_{{\mathcal{T}},h}^{(3)}$
are distinct over ${\mathcal{S}}_{{\mathcal{T}}|{\mathcal{T}}_1}(\mu)$.
Thus, by Proposition~\ref{euler_prp},
the decomposition~\eqref{cart_split}, and 
a rescaling of the linear~map,
\begin{equation}\label{n2tac_contr_lmm1e5}
\sum_{\chi_{\mathcal{T}}(\tilde{1},\hat{1})=0}    
{\mathcal{C}}_{{\mathcal{Z}}_{\mathcal{T}}}({\mathcal{D}}_{\tilde{1},\hat{1}})
=3\big|{\mathcal{S}}_1\big|.
\end{equation}
The claim follows from equations
\eqref{n2tac_contr_lmm1e3} and \eqref{n2tac_contr_lmm1e5}.
\end{proof}

\begin{lmm}
\label{n2tac_contr_lmm2} 
The total contribution to 
${\mathcal{C}}_{\partial\Bbb{P}(L_{\tilde{1}}\oplus L_{\tilde{1}}^*)}
({\mathcal{D}}_{\tilde{1},\hat{1}})$
from the boundary strata 
$\Bbb{P}(L_{\tilde{1}}\oplus L_{\tilde{1}}^*)|{\mathcal{U}}_{{\mathcal{T}}_1,{\mathcal{T}}}$,
where ${\mathcal{T}} = (M_1,I;j,\under{d})$ is a bubble type
such that ${\mathcal{T}} < {\mathcal{T}}_1$ and
$\chi_{\mathcal{T}}(\tilde{1},\hat{1}) > 0$, is given~by
$$\sum_{\chi_{\mathcal{T}}(\tilde{1},\hat{1})>0}  
{\mathcal{C}}_{\partial\Bbb{P}(L_{\tilde{1}}\oplus L_{\tilde{1}}^*)|
{\mathcal{U}}_{{\mathcal{T}}_2,{\mathcal{T}}}} ({\mathcal{D}}_{\tilde{1},\hat{1}})
=-\big|{\mathcal{V}}_2^{(1,1)}(\mu)\big|.$$
\end{lmm}

\begin{proof}
(1)\qua By Lemma~\ref{n3p3_str_lmm1},
$${\overU}_{{\mathcal{T}}_1}^{(1)}(\mu)\cap{\mathcal{U}}_{{\mathcal{T}}_1,{\mathcal{T}}}
\subset{\mathcal{U}}_{{\mathcal{T}}|{\mathcal{T}}_1}^{(1)}(\mu).$$
If $j_{\hat{1}} = \tilde{1}$ or $j_{\hat{1}} > \tilde{1}$
and $d_{j_{\hat{1}}} = 0$, the section 
${\mathcal{D}}_{\tilde{1},\hat{1}}$ has a nonvanishing extension
over~${\mathcal{U}}_{{\mathcal{T}}|{\mathcal{T}}_1}^{(1)}(\mu)$.
Thus, we only need to consider bubble types~${\mathcal{T}}$
such that 
$$j_{\hat{1}} \equiv h >\tilde{1} \qquad\hbox{and}\qquad d_h > 0.$$
Furthermore, if ${\mathcal{U}}_{{\mathcal{T}}|{\mathcal{T}}_1}^{(1)}(\mu) \neq \eset$,
$I^+ = \{\tilde{1},h\}$.

(2)\qua If $d_{\tilde{0}} = 0$ and 
${\mathcal{U}}_{{\mathcal{T}}|{\mathcal{T}}_1}^{(1)}(\mu) \neq \eset$,
${\mathcal{T}} = {\mathcal{T}}_1(l)$ for some $l \in [N]$.
By Proposition~\ref{str_prp},
\begin{gather*}
{\mathcal{D}}_{\tilde{1},\hat{1}}\phi_{{\mathcal{T}}_1,{\mathcal{T}}}
([\ups_{\tilde{1}},\ups_{\hat{1}}];\ups)=
-(y_{\hat{1}} - x_h)^{-2}
\big\{\ups_{\hat{1}} \otimes {\mathcal{D}}_{{\mathcal{T}},h}^{(1)}
 + \ve(\ups)\big\}\ups^*\\
\quad\mbox{for all}\  \ups\in\mathcal{FT}_{\de} - Y(\mathcal{FT};\hat{I}^+).
\end{gather*}
Let ${\mathcal{Z}}_{\mathcal{T}} = \Bbb{P}L_{\tilde{1}}$.
Since the section ${\mathcal{D}}_{{\mathcal{T}},h}^{(1)}$ does not vanish
on ${\mathcal{U}}_{{\mathcal{T}}|{\mathcal{T}}_1}^{(1)}(\mu)$,
by Proposition~\ref{euler_prp} and a rescaling of the linear map,
\begin{gather*}
{\mathcal{C}}_{\Bbb{P}(L_{\tilde{1}}\oplus L_{\tilde{1}}^*)|
{\mathcal{U}}_{{\mathcal{T}}|{\mathcal{T}}_1}^{(1)}(\mu)-{\mathcal{Z}}_{\mathcal{T}}}
({\mathcal{D}}_{\tilde{1},\hat{1}})=-N(\al_1),
\qquad\hbox{where}\\
\al_1\in\Ga\big(\Bbb{P}^1 \times {\mathcal{S}}_1(\mu);
\hbox{Hom}(\ga^*,\ga^* \otimes \ev_{\hat{0}}^*T\PP)\big),
\end{gather*}
is a nonvanishing section.
Thus, by Proposition~\ref{zeros_prp},
$${\mathcal{C}}_{\Bbb{P}(L_{\tilde{1}}\oplus L_{\tilde{1}}^*)|
{\mathcal{U}}_{{\mathcal{T}}|{\mathcal{T}}_1}^{(1)}(\mu)-{\mathcal{Z}}_{\mathcal{T}}}
({\mathcal{D}}_{\tilde{1},\hat{1}})
=-\big\lan 2\la - \la,\Bbb{P}^1\big\ran
\big|{\mathcal{U}}_{{\mathcal{T}}|{\mathcal{T}}_1}^{(1)}(\mu)\big|
=-\big|{\mathcal{U}}_{{\mathcal{T}}|{\mathcal{T}}_1}^{(1)}(\mu)\big|.$$
On the other hand, with the same notation as in (2) 
of the proof of Lemma~\ref{n2tac_contr_lmm1}, near ${\mathcal{Z}}_{\mathcal{T}}$,
\begin{equation*}\begin{split}
{\mathcal{D}}_{\tilde{1},\hat{1}}\phi_{{\mathcal{T}}_1,{\mathcal{T}}}
([\ups_{\tilde{1}},\ups_{\hat{1}}];u,\ups)=
-(y_{\hat{1}} - x_h)^{-2}\big\{{\mathcal{D}}_{{\mathcal{T}},h}^{(1)}
 + \ve_{\hat{1}}(\ups)\big\}u \otimes \ups^* \qquad&\\
+\big\{{\mathcal{D}}_{{\mathcal{T}},h}^{(1)}
 + \ve_{\tilde{1}}(\ups)\big\}\ups.&
\end{split}\end{equation*}
Thus, by Proposition~\ref{euler_prp},
$${\mathcal{C}}_{{\mathcal{Z}}_{\mathcal{T}}}({\mathcal{D}}_{\tilde{1},\hat{1}})
=\big|{\mathcal{Z}}_{\mathcal{T}}\big|
=\big|{\mathcal{U}}_{{\mathcal{T}}|{\mathcal{T}}_1}^{(1)}(\mu)\big|.$$
We conclude that
\begin{equation}\label{n2tac_contr_e3}
\sum_{d_{\tilde{1}}=0}
{\mathcal{C}}_{\Bbb{P}(L_{\tilde{1}}\oplus L_{\tilde{1}}^*)|
{\mathcal{U}}_{{\mathcal{T}}_1,{\mathcal{T}}}}
({\mathcal{D}}_{\tilde{1},\hat{1}}) =0.
\end{equation}
(3)\qua Finally, suppose $d_{\tilde{1}} > 0$.
The same argument as in (2) above shows that
$${\mathcal{C}}_{\Bbb{P}(L_{\tilde{1}}\oplus L_{\tilde{1}}^*)|
{\mathcal{U}}_{{\mathcal{T}}|{\mathcal{T}}_1}^{(1)}(\mu)-{\mathcal{Z}}_{\mathcal{T}}}
({\mathcal{D}}_{\tilde{1},\hat{1}})
=-\big|{\mathcal{U}}_{{\mathcal{T}}|{\mathcal{T}}_1}^{(1)}(\mu)\big|,$$
but ${\mathcal{Z}}_{\mathcal{T}}$ is
${\mathcal{D}}_{\tilde{1},\hat{1}}$--hollow.
Thus, summing up over all bubble types ${\mathcal{T}}$
of appropriate form, we obtain
\begin{equation}\label{n2tac_contr_e5}
\sum_{d_{\tilde{1}}>0}
{\mathcal{C}}_{\Bbb{P}(L_{\tilde{1}}\oplus L_{\tilde{1}}^*)|
{\mathcal{U}}_{{\mathcal{T}}_1,{\mathcal{T}}}}
({\mathcal{D}}_{\tilde{1},\hat{1}})
=-\big|{\mathcal{V}}_2^{(1,1)}(\mu)\big|.
\end{equation}
The claim follows equations~\eqref{n2tac_contr_e3} and~\eqref{n2tac_contr_e5}.
\end{proof}

\subsection[Rational cuspidal curves in Pn]{Rational cuspidal curves in $\P$}
\label{cusps_subs}

In this subsection, we prove Theorem~\ref{cusps_thm}.
In particular, we construct a tree of chern classes,
as mentioned in the third-to-last paragraph of Subsection~\ref{method_subs}.
The sum of these chern classes, with an appropriate sign, 
is the number that appears on the right-hand
side of the equation in Theorem~\ref{cusps_thm}.
The tree is very similar to that constructed in Subsection~3.1 of~\cite{Z4};
the main difference is that here we focus on intersection numbers,
instead of zeros of polynomial maps.
Theorem~\ref{cusps_thm} follows immediately from 
Corollary~\ref{cusps_crl1} and Lemma~\ref{cusps_lmm4}.

We first introduce a little more notation.
If $d$, $N$, and $\mu$ are as in the statement of Theorem~\ref{cusps_thm},
$k \ge 1$, and $m \ge 0$,
let ${\mathcal{V}}_{k,m}(\mu)$ and ${\overV}_{k,m}(\mu)$ denote
the quotients of the disjoint unions of the spaces 
${\mathcal{U}}_{\mathcal{T}}(\tilde{\mu})$ and~${\overU}_{\mathcal{T}}(\tilde{\mu})$, 
respectively, taken over all bubble types
\begin{gather*}
{\mathcal{T}} = ([N] - M_0,I_k;j,\under{d})
\quad\mbox{such that}\quad
M_0 \subset [N],~|M_0| = m,\\
I_k = \{\hat{0}\}\cup\{\tilde{1},\ldots,\tilde{k}\},\quad
d_{\tilde{1}},\ldots,d_{\tilde{k}} > 0,  \quad
\sum d_i = d,
\end{gather*}
by the natural action of the symmetric group~$S_k$.
Here $\tilde{\mu}$ is the $([N]{-}M_0)\sqcup\{\hat{0}\}$--tuple
of constraints defined by
$$\tilde{\mu}_l=\mu_l\quad\hbox{if}~~l \in [N] - M_0;
\qquad
\tilde{\mu}_{\hat{0}}=\bigcap_{l\in M_0}\mu_l.$$
By dimension-counting, the spaces  ${\mathcal{V}}_{k,m}(\mu)$
are smooth manifolds.
We define the vector bundle \hbox{$E_{k,m} \to {\overV}_{k,m}(\mu)$}
and homomorphism \hbox{$\al_{k,m}\co E_{k,m} \to \ev_{\hat{0}}^*T\P$} 
over ${\overV}_{k,m}(\mu)$ by
$$E_{k,m}|{\overU}_{\mathcal{T}}(\mu)=
\bigoplus_{i\in I^+}L_i,\quad
\al_{k,m}\big((\ups_i)_{i\in I^+}\big)
=\sum_{i\in I^+}{\mathcal{D}}_{{\mathcal{T}},i}^{(1)}\ups_i,$$
whenever ${\mathcal{T}}$ is a bubble type as above.

We now construct the tree mentioned above. Each node is a tuple 
$$\si=(r;k,m;\phi),$$ 
where $r \ge 0$ is the distance to the root 
$\si_0 = (0;1,0;\cdot)$, $k \ge 1$, and $m \ge 0$.
The tree satisfies the following properties.
If $r > 0$ and 
$$\si^*\equiv(r - 1;k^*,m^*;\phi^*)$$
is the node from which $\si$ is directly descendant, 
we require that $k^* \le k$, $m^* \le m$, and at least
one of the inequalities is strict.
Furthermore, $\phi$ specifies a splitting of the set~$[k]$
into $k^*$ disjoint subsets and an assignment of $m - m^*$
of the elements of the set 
$$[m]\equiv\big\{(1,1),\ldots,(1,m)\big\}$$ 
to these subsets.
This description inductively constructs an infinite tree.
However, we will need to consider only the nodes 
$$\si\equiv(r;k,m;\phi) \quad\mbox{such that}\quad 2k + m\le n + 2.$$
We will write $\si \vdash \si^*$  to indicate that 
$\si$ is directly descendant from~$\si^*$.

For each node in the above tree, except for the root, 
we now define a linear map between vector bundles.
If $\si = (r;k,m;\phi)$ and $s \ge 1$, 
let 
$$\{\si_s = (s;k_s,m_s;\phi_s) :0 \le s\le r\}$$ 
be the sequence of nodes such that
$\si_r = \si$ and $\si_s \vdash \si_{s-1}$ for all $s > 0$.
Put
\begin{gather*}
{\overV}_{\si}={\overV}_{k,m}(\mu),\qquad
E_{\si} = E_{k,m}\lra{\overV}_{\si},\qquad
\al_{\si}=\al_{k,m},\\
{\mathcal{X}}_{\si}={\mathcal{Y}}_{\si} \times {\overV}_{\si},\quad
{\mathcal{X}}_{\si,s}={\mathcal{Y}}_{\si,s} \times {\overV}_{\si},
\end{gather*}
where
\begin{gather*}
{\mathcal{Y}}_{\si}={\mathcal{Y}}_{\si,r},\quad
{\mathcal{Y}}_{\si,0}=\{pt\},\quad
{\mathcal{Y}}_{\si,s}=\Bbb{P}F_{\si_s} \times {\mathcal{Y}}_{\si,s-1}
\hbox{~~if~~}s > 0,\\
\ov{\frak M}_{\si}=
\prod_{i\in\Im~ \phi}  \ov{\frak M}_{i+\phi^{-1}(i)}, \quad
F_{\si}=\bigoplus_{i\in\Im~ \phi}  \ga_{\si;i}
\lra\ov{\frak M}_{\si}.
\end{gather*}
For the purposes of the last line above,
we view $\phi$ as a map from $[k] - [k^*]$ and 
a subset of~$[m]$ to~$[k^*]$ in the notation of the previous paragraph.
Then, 
$$\ga_{\si;i}\lra\ov{\frak M}_{i+\phi^{-1}(i)}$$
is the ``tautological'' line bundle,
ie the universal tangent bundle at the marked point~$i$. 
Let
$${\mathcal{O}}_{\si}={\mathcal{O}}_{\si,r},\quad
{\mathcal{O}}_{\si,1}=\ev_{\hat{0}}^*T\P,\quad
{\mathcal{O}}_{\si,s}={\mathcal{O}}_{\si,s-1}\big/\hbox{Im}~ \bar{\nu}_{\si,s-1}
~~\hbox{if}~s > 1,$$
where
$\bar{\nu}_{\si,s}\in\Ga\big({\mathcal{X}}_{\si,s};
\hbox{Hom}(\ga_{F_{\si_s}},{\mathcal{O}}_{\si,s})\big)$
is a generic section.
Since $k_{s-1} \le k_s$, $m_{s-1} \le m_s$, 
and one of the inequalities is strict,
\begin{equation*}\begin{split}
\frac{1}{2}\dim{\mathcal{X}}_{\si,s}
\le \frac{1}{2}\dim{\mathcal{X}}_{\si}
&= \big(n + 2 - 2k - m\big)+
\sum_{s=1}^{s=r}\big(\big|\Im~ \phi_s\big| - 1\big)\\
&=n+1-k-r<\rk~ {\mathcal{O}}_{\si,1}-(r - 1).
\end{split}\end{equation*}
Thus, we see inductively that 
each bundle ${\mathcal{O}}_{\si,s}$ is well-defined
and a generic section $\bar{\nu}_{\si,s}$ of 
$\hbox{Hom}(\ga_{F_{\si,s}},{\mathcal{O}}_{\si,s})$
does not vanish.
Let 
$$\pi_{\si}\co \ev_{\hat{0}}^*T\P\lra{\mathcal{O}}_{\si}$$ 
be the projection map. We define 
\begin{gather*}
\tilde{\al}_{\si}\in\Ga\big({\mathcal{X}}_{\si};
\hbox{Hom}(\ga_{F_{\si}}^* \otimes E_{\si};
\ga_{F_{\si}}^* \otimes {\mathcal{O}}_{\si})\big)
\qquad\hbox{by}\\
\big\{\tilde{\al}_{\si}(\tau \otimes \ups)\big\}(w)
=\tau(w)\cdot\pi_{\si}\al_{\si}(\ups)
\in{\mathcal{O}}_{\si}.
\end{gather*}

\begin{lmm}
\label{cusps_lmm1}
With notation as above,
$$\big|{\mathcal{S}}_1(\mu)\big|=
\big\lan c\big(L_{\tilde{1}}^* \otimes \ev_{\hat{0}}^*T\P\big),
{\overV}_1(\mu)\big\ran
-\sum_{\si\vdash\si_0}N\big(\tilde{\al}_{\si}\big).$$
Furthermore, for every node $\si^* \neq \si_0$,
$$N\big(\tilde{\al}_{\si^*}\big)=
\big\lan c\big(\ga_{F_{\si^*}}^* \otimes {\mathcal{O}}_{\si^*}\big)
c\big(\ga_{F_{\si^*}}^* \otimes E_{\si^*}\big)^{-1},
{\mathcal{X}}_{\si^*}\big\ran
-\sum_{\si\vdash\si^*}N\big(\tilde{\al}_{\si}\big).$$
\end{lmm}

\begin{proof}
This lemma is obtained by the usual argument 
from the estimate (3b) of Proposition~\ref{str_prp}
via Propositions~\ref{zeros_prp} and~\ref{euler_prp}.
If $\si^* \neq \si_0$, the proof is the same as
the proof of Lemma~3.3 in~\cite{Z4}.
For the first identity, apply the proof of Lemma~3.3
with $\al_{\si_0}^{\perp} = {\mathcal{D}}_{{\mathcal{T}},\tilde{1}}^{(1)}$.
\end{proof}

\begin{lmm}
\label{cusps_lmm2}
For every node $\si \neq \si_0$,
$$\big\lan c\big(\ga_{F_{\si}}^* \otimes {\mathcal{O}}_{\si}\big)
c\big(\ga_{F_{\si}}^* \otimes E_{\si}\big)^{-1},
{\mathcal{X}}_{\si}\big\ran
=\big\lan c(\ev_{\hat{0}}^*T\P)c(E_{k,m})^{-1},
{\overV}_{k,m}(\tilde{\mu})\big\ran.$$
Furthermore,
$$\big\lan c_1\big(L_{\tilde{1}}^* \otimes \ev_{\hat{0}}^*T\P\big),
{\overV}_1(\mu)\big\ran
=\big\lan c(\ev_{\hat{0}}^*T\P)c(E_{1,0})^{-1},
{\overV}_{1,0}(\tilde{\mu})\big\ran.$$
\end{lmm}

\begin{proof}
For the first identity,
see the proof of \cite[Corollary~3.5]{Z4}.
The second equality is clear from the fact that
$\dim{\overV}_1(\mu) = \rk\ev_{\hat{0}}^*T\P$
\end{proof}

\begin{crl}
\label{cusps_crl1}
\begin{equation*}\begin{split}
\big|{\mathcal{S}}_1(\mu)\big|
&=\sum_{(1,0)\le(k,m)}\Bigg\{(-1)^{k+m-1}k^m(m - 1)! \\
&\qquad\qquad\qquad   \times\sum_{l=0}^{n+2-(2k+m)}  
\binom{n + 1}{l} \big\lan a^l\tilde{\eta}_{\hat{0},n+2-(2k+m)-l},
{\overV}_{k,m}(\tilde{\mu})\big\ran\Bigg\}.
\end{split}\end{equation*}
\end{crl}

\begin{proof}
This corollary follows from Lemma~\ref{cusps_lmm1}
and Lemma~\ref{cusps_lmm2} via straightforward combinatorics;
see \cite[Corollary~3.6 and Lemma~3.7]{Z4}.
\end{proof}

\begin{lmm}
\label{cusps_lmm4}
For all $k \ge 1$ and $l \ge 0$,
\begin{equation*}\begin{split}
\sum_{m\ge0}(-1)^m k^m(m - 1)!
\blr{a^l\tilde{\eta}_{\hat{0},n+2-(2k+m)-l},{\overV}_{k,m}(\tilde{\mu})}
\qquad\qquad\qquad\qquad&\\
=(k - 1)!\blr{a^l\eta_{\hat{0},n+2-2k-l},{\overV}_k(\mu)}.&
\end{split}\end{equation*}
\end{lmm}

\begin{proof}
See the proof of Corollary~3.10 in~\cite{Z4},
which uses~\eqref{psi_class1} along with
\cite[Subsection~3.2]{P}.
\end{proof}

\section{Low-degree numbers}
\label{tables_sec}

We now give some low-degree enumerative numbers
for rational curves in projective spaces.
In all five tables, the top row lists the degree~$d$ of the map.
In Tables~\ref{n2p3_table} and~\ref{n2tac_table},
the constraints are assumed to be $3d{-}2$ points in general position
in~$\P$.
In Tables~\ref{n3p3_table} and~\ref{n3tac_table},
the constraints are $p$~points and $q$~lines in~$\PPP$, as specified 
by the second row.
Similarly, in Table~\ref{n4cusps_table},
the constraints are $p$~points, $q$~lines, and $r$ two-planes
in~$\Bbb{P}^4$.

The formulas of Theorems~\ref{n3p3_thm} and~\ref{n3tac_thm}
give zeros in degrees one, two, and three.
{}From classical algebraic geometry,
one would expect these low-degree numbers, as well as
the first three degree-four numbers listed in
Tables~\ref{n3p3_table} and~\ref{n3tac_table}, to vanish.
In addition, as expected, the fourth number in 
Table~\ref{n3p3_table} (Table~\ref{n3tac_table})
is the same as the degree-four number of 
Table~\ref{n2p3_table} (Table~\ref{n2tac_table}).
Similarly, all degree-one and -two numbers $|{\mathcal{S}}_1(\mu)|$
and several degree-three and -four numbers,
as listed in Table~\ref{n4cusps_table}, are zero, as the case should~be.
Finally, observe that the third number of Table~\ref{n4cusps_table}
is the same as the long-known number of plane cubic cuspidal curves
that pass through seven general points.

\begin{table}[ht!]
\begin{footnotesize}\begin{center}
\begin{tabular}{||c|c|c|c|c|c|c|c|c||}
\hline\hline
$d$& 1& 2& 3& 4& 5& 6& 7& 8\\
\hline
$|{\cal V}_1^{(2)}(\mu)|$& 0& 0& 0& 60& 56,400& 49,177,440& 
56,784,765,120&  91,466,185,097,280\\
\hline\hline
\end{tabular}
\end{center}\end{footnotesize}
\vskip-8pt
\caption{One-component rational triple-pointed curves in $\PP$}
\label{n2p3_table}
\end{table}

\begin{table}[ht!]
\begin{footnotesize}\begin{center}
\begin{tabular}{||c|c|c|c|c|c|c|c|c||}
\hline\hline
$d$&  ${}\!\!1\!\!{}$&  ${}\!\!2\!\!{}$&  ${}\!\!3\!\!{}$& 4& 5& 6& 7& 8\\
\hline
$|{\cal S}_1^{(1)}(\mu)|$& ${}\!\!0\!\!{}$& ${}\!\!0\!\!{}$& 
${}\!\!0\!\!{}$& ${}\!\!1,296\!\!{}$& 
${}\!\!499,680\!\!{}$&  ${}\!\!271,751,040\!\!{}$& 
${}\!\!227,509,931,520\!\!{}$& ${}\!\!287,190,836,432,640\!\!{}$\\
\hline\hline
\end{tabular}
\end{center}\end{footnotesize}
\vskip-8pt
\caption{One-component rational tacnodal curves in $\PP$}
\label{n2tac_table}
\end{table}

\begin{table}[ht!]
\begin{footnotesize}\begin{center}
\begin{tabular}{||c|c|c|c|c|c|c|c|c|c|c||}
\hline\hline
$d$& 4&4&4& 4&4&4& 5&5&5& 6\\
\hline
$(p,q)$& (6,1)&(5,3)&(4,5)& (3,7)&(2,9)&(1,11)&  (8,1)&
${}\!(7,3)\!{}$& ${}\!(6,5)\!{}$& ${}\!(10,1)\!{}$\\
\hline
$\frac{1}{6}|{\cal V}_1^{(2)}(\mu)|$& 0&0&0& 60&1,280& ${}\!19,640\!{}$& 
${}\!8\!{}$& ${}\!264\!{}$& ${}\!4,360\!{}$& ${}\!4,680\!{}$\\
\hline\hline
\end{tabular}
\end{center}\end{footnotesize}
\vskip-8pt
\caption{One-component rational  triple-pointed curves in $\PPP$}
\label{n3p3_table}
\end{table}

\begin{table}[ht!]
\begin{footnotesize}\begin{center}
\begin{tabular}{||c|c|c|c|c|c|c|c|c|c|c||}
\hline\hline
$d$& 4&4&4& 4&4&4& 5&5&5& 6\\
\hline
$(p,q)$& ${}\!\!\!(6,1)\!\!\!$& ${}\!\!\!(5,3)\!\!\!$& ${}\!\!\!(4,5)\!\!\!$& 
${}\!\!(3,7)\!\!$& ${}\!\!(2,9)\!\!$& ${}\!\!(1,11)\!\!$&  
(8,1)&(7,3)&(6,5)& (10,1)\\
\hline
$\frac{1}{2}|{\cal S}_1^{(1)}(\mu)|$& 0&0&0& 
${}\!\!1,296\!\!$& ${}\!\!27,648\!\!$& ${}\!\!426,672\!\!$& 
${}\!960\!{}$& ${}\!9,792\!{}$& ${}\!111,840\!{}$& ${}\!112,320\!{}$\\
\hline\hline
\end{tabular}
\end{center}\end{footnotesize}
\vskip-8pt
\caption{One-component rational tacnodal curves in $\PPP$}
\label{n3tac_table}
\end{table}

\begin{table}[ht!]
\begin{footnotesize}\begin{center}
\begin{tabular}{||c|c|c|c|c|c|c|c|c|c|c||}
\hline\hline
$d$& 3&3& 3&3&  4&4&4&  4& 5\\
\hline
$(p,q,r)$& ${}\!(4,0,1)\!{}$& ${}\!(3,1,2)\!{}$& ${}\!(3,0,4)\!{}$& 
${}\!(2,1,5)\!{}$&  ${}\!\!(6,0,0)\!\!{}$& ${}\!\!(5,1,1)\!\!{}$& 
${}\!\!(5,0,3)\!\!{}$& ${}\!\!(4,1,4)\!\!{}$& ${}\!\!(7,1,0)\!\!{}$\\
\hline
$|{\cal S}_1(\mu)|$& 0&0& 24&240& 
${}\!0\!{}$& ${}\!0\!{}$& ${}\!0\!{}$& ${}\!\!1,680\!\!{}$& ${}\!120\!{}$\\
\hline\hline
\end{tabular}
\end{center}\end{footnotesize}
\vskip-8pt
\caption{One-component rational cuspidal curves in $\Bbb{P}^4$}
\label{n4cusps_table}
\end{table}


\begin{thebibliography}

\bibitem{BT}
\textbf{Raoul Bott}, \textbf{Loring~W Tu}, \emph{Differential forms in
  algebraic topology}, volume~82 of \emph{Graduate Texts in Mathematics},
  Springer--Verlag, New York (1982)
  \MR{0658304}

\bibitem{DH}
\textbf{Steven Diaz}, \textbf{Joe Harris}, \emph{Geometry of the {S}everi
  variety}, Trans. Amer. Math. Soc. 309 (1988) 1--34
  \MR{0957060}

\bibitem{FP}
\textbf{W Fulton}, \textbf{R Pandharipande}, \emph{Notes on stable maps and
  quantum cohomology}, from: ``Algebraic geometry---Santa Cruz 1995'', Proc.
  Sympos. Pure Math. 62, Amer. Math. Soc., Providence, RI (1997)  45--96
  \MR{1492534}

\bibitem{GH}
\textbf{Phillip Griffiths}, \textbf{Joseph Harris}, \emph{Principles of
  algebraic geometry}, Wiley Classics Library, John Wiley \& Sons Inc., New
  York (1994)
  \MR{1288523}

\bibitem{G}
\textbf{M Gromov}, \emph{Pseudoholomorphic curves in symplectic manifolds},
  Invent. Math. 82 (1985) 307--347
  \MR{0809718}

\bibitem{KQR}
\textbf{Sheldon Katz}, \textbf{Zhenbo Qin}, \textbf{Yongbin Ruan},
  \emph{Enumeration of nodal genus--{$2$} plane curves with fixed complex
  structure}, J. Algebraic Geom. 7 (1998) 569--587
  \MR{1618136}

\bibitem{KM}
\textbf{Maxim Kontsevich}, \textbf{Yu Manin}, \emph{Gromov-{W}itten classes,
  quantum cohomology, and enumerative geometry}, Comm. Math. Phys. 164 (1994)
  525--562
  \MR{1291244}

\bibitem{K}
\textbf{Maxim Kontsevich}, \emph{Enumeration of rational curves via torus
  actions}, from: ``The moduli space of curves (Texel Island, 1994)'', Progr.
  Math. 129, Birkh\"auser Boston, Boston, MA (1995)  335--368
  \MR{1363062}

\bibitem{LT}
\textbf{Jun Li}, \textbf{Gang Tian}, \emph{Virtual moduli cycles and
  {G}romov-{W}itten invariants of general symplectic manifolds}, from: ``Topics
  in symplectic $4$-manifolds (Irvine, CA, 1996)'', First Int. Press Lect.
  Ser., I, Internat. Press, Cambridge, MA (1998)  47--83
  \MR{1635695}

\bibitem{MS}
\textbf{Dusa McDuff}, \textbf{Dietmar Salamon}, \emph{{$J$}-holomorphic curves
  and quantum cohomology}, volume~6 of \emph{University Lecture Series},
  American Mathematical Society, Providence, RI (1994)
  \MR{1286255}

\bibitem{P}
\textbf{Rahul Pandharipande}, \emph{Intersections of {$\mathbf{Q}$}-divisors on
  {K}ontsevich's moduli space {$\overline M\sb {0,n}(\mathbf{P}\sp r,d)$} and
  enumerative geometry}, Trans. Amer. Math. Soc. 351 (1999) 1481--1505
  \MR{1407707}

\bibitem{R1}
\textbf{Ziv Ran}, \emph{On the quantum cohomology of the plane, old and new, and
  a {$K3$} analogue}, Collect. Math. 49 (1998) 519--526
  \MR{1677096}

\bibitem{R2}
\textbf{Ziv Ran}, \emph{Enumerative geometry of divisorial families of rational
  curves}, Ann. Sc. Norm. Super. Pisa Cl. Sci. (5) 3 (2004) 67--85
  \MR{2064968}

\bibitem{RT}
\textbf{Yongbin Ruan}, \textbf{Gang Tian}, \emph{A mathematical theory of
  quantum cohomology}, J. Differential Geom. 42 (1995) 259--367
  \MR{1366548}

\bibitem{V}
\textbf{R Vakil}, 
\emph{Enumerative Geometry of Plane Curve of Low Genus}, preprint,
\arxiv{math.AG/9803007}

\bibitem{Ze}
\textbf{H Zeuthen}, 
\emph{Almindelige Egenskaber ved Systemer af Plane Kurver},
Kongelige Danske Videnskabernes Selskabs Skrifter,
10 (1873) 285--393 (Danish)

\bibitem{Z1}
\textbf{Aleksey Zinger}, \emph{Enumeration of genus-two curves with a fixed
  complex structure in {$\Bbb P\sp 2$} and {$\Bbb P\sp 3$}}, J. Differential
  Geom. 65 (2003) 341--467
  \MR{2064428}

\bibitem{Z2}
\textbf{Aleksey Zinger}, 
\textit{Enumerative vs. Symplectic Invariants and Obstruction Bundles},
to appear in J. Symplectic Geom.

\bibitem{Z3}
\textbf{Aleksey Zinger}, \emph{Enumeration of genus-three plane curves with a
  fixed complex structure}, J. Algebraic Geom. 14 (2005) 35--81
  \MR{2092126}

\bibitem{Z4}
\textbf{Aleksey Zinger}, \emph{Enumeration of one-nodal rational curves in
  projective spaces}, Topology 43 (2004) 793--829
  \MR{2061208}

\end{thebibliography}
\end{document}